\documentclass[11pt,preprint]{elsarticle}
\usepackage[top=1in,bottom=1in,left=0.5in,right=0.5in]{geometry}
\usepackage{stmaryrd}
\usepackage{lineno}
\usepackage[colorlinks]{hyperref}
\usepackage{graphicx}
\usepackage{subfigure}
\usepackage{color}
\usepackage{amsmath}
\usepackage{amssymb,bm}
\usepackage{amsfonts}
\usepackage{amstext}
\usepackage{amsxtra}
\usepackage{mathrsfs}
\usepackage{amscd}
\usepackage{amsgen}
\usepackage{amsopn}
\usepackage{cases}
\usepackage{ntheorem}
\usepackage[utf8]{inputenc}
\usepackage{booktabs}
\usepackage{listings}
\usepackage[normalem]{ulem}
\usepackage{tikz}
\usepackage{tikz-3dplot}
\usepackage{tcolorbox}
\usepackage{verbatim}
\usetikzlibrary{calc}
\usetikzlibrary{plotmarks}
\usetikzlibrary{circuits}
\usetikzlibrary{decorations.shapes,shapes.geometric}
\usetikzlibrary{3d}
\usetikzlibrary{positioning}
\usetikzlibrary{arrows,shapes,backgrounds}
\numberwithin{equation}{section}

\newcommand{\topcaption}{%
\setlength{\abovecaptionskip}{0pt}%
\setlength{\belowcaptionskip}{10pt}%
\caption}

\journal{...}

\begin{document}

\newtheorem{theorem}{Theorem}[section]
\newtheorem{definition}{Definition}[section]
\newtheorem{example}{Example}[section]
\newtheorem{lemma}[theorem]{Lemma}
\newtheorem{corollary}[theorem]{Corollary}
\newtheorem{proposition}[theorem]{Proposition}
\newtheorem{remark}{Remark}[section]
\theorembodyfont{\normalfont}
\newtheorem*{pf}{\it Proof.}
\def\endpf{\hspace*{\fill}~$\square$\par\endtrivlist\unskip}

\begin{frontmatter}

\title{Analysis and computation of a pressure-robust method for the rotation form of the stationary incompressible Navier--Stokes equations by using high-order finite elements\tnoteref{mytitlenote}}
\tnotetext[mytitlenote]{This research was supported in part by National Science Foundation of China under Grant No.11771348.}
\author[add1]{Di Yang}\ead{ydxjtu0226@gmail.com}
\author[add1]{Yinnian He\corref{correspondingauthor}}\ead{heyn@mail.xjtu.edu.cn}
\cortext[correspondingauthor]{Corresponding author.}
\address[add1]{School of Mathematics and Statistics, Xi'an Jiaotong University, Xi'an, Shaanxi 710049, PR China }

\begin{abstract}
In this work, we develop a high-order pressure-robust method for the rotation form of the stationary incompressible Navier--Stokes equations. The original idea is to change the velocity test functions in the discretization of trilinear and right hand side terms by using an $\bm{H}(\mathrm{div})$-conforming velocity reconstruction operator. In order to match the rotation form and error analysis, a novel skew-symmetric discrete trilinear form containing the reconstruction operator is proposed, in which not only the velocity test function is changed. The corresponding well-posed discrete weak formulation stems straight from the classical inf-sup stable mixed conforming high-order finite elements, and it is proven to achieve the pressure-independent velocity errors. Optimal convergence rates of $H^1$, $L^2$-error for the velocity and $L^2$-error for the Bernoulli pressure are completely established. Adequate numerical experiments are presented to demonstrate the theoretical results and the remarkable performance of the proposed method.
\end{abstract}


\end{frontmatter}


\section{Introduction}
\label{sec:intro}

The Navier--Stokes equations are widely followed by flow issues, like flows in pipes and channels, flows around objects such as a cylinder and wings of a plane, to name just a few. The most common stationary incompressible Navier--Stokes equations consist of a momentum balance including velocity and pressure two fundamental variables, and a mass conservation which is also called the velocity divergence-free constraint, with appropriate boundary conditions. In the numerical world, it is popular with mixed finite element methods used to obtain a non-positive algebra system, such that the mixed finite elements satisfying the famous discrete inf-sup condition guarantees stability and the unique solution. In classical mixed methods, one tends to relax the divergence-free constraint to access the discrete inf-sup stability more easily. In this case, a serious cost of this relaxation is that the velocity errors depend on a pressure-dependent error contribution \cite{Linke2016, Olshanskii20041699}, that is,

$$\nu^{-1}\inf_{q_h\in Q_h}\|p-q_h\|_{L^2},$$
where $\nu$ is the kinematic viscosity, $p$ is pressure and $Q_h$ is the discrete pressure space. Hence, if the pressure error is large (i.e., strong gradient fields in the momentum balance), then the velocity errors of classical inf-sup stable mixed element methods will be possibly severe. This lack of robustness is also called poor mass conservation in \cite{Galvin2012166, Linke2016, Manica201103}, and it was adequately demonstrated by Fig. 1.1--Fig. 1.3 in \cite{John2017}.

Recently, as an efficient method of improving pressure robustness which means delivering a velocity error independent of pressure $p$, changing the velocity test functions in the discretization of right hand side terms by using an $\bm{H}(\mathrm{div})$-conforming velocity reconstruction operator was proposed in \cite{Linke2012837, Linke2014782, Linke2016289}. It is based on the inf-sup stable mixed finite elements. In particular, Linke et al. \cite{Linke2016289} constructed families of conforming and nonconforming mixed finite elements of arbitrary order suitable for velocity reconstruction operators in two and three space dimensions and successfully applied them to the Stokes equations in theory and numerical performance. Also, they provided some indications to deal with the pressure-robust problems in the incompressible Navier--Stokes equations. In the later literature, this pressure-robust method has been widely coordinated with some novel finite element methods recently proposed such as weak Galerkin (WG) finite element methods and virtual element methods (VEM), for solving the Stokes equations and Brinkman problems \cite{Liu2020, Mu20201422, Mu2020B608, Mu2021A2614, Wang2021}.

In the incompressible Navier--Stokes equations, the nonlinear convection term has various equivalent forms and especially if it is taken as the rotation form $(\mathrm{rot}\,\bm{u})\times\bm{u}$ then the Bernoulli pressure, a mixture of the velocity and the kinematic pressure, will be employed to replace the kinematic pressure. Layton et al. discussed the accuracy of the classical mixed finite elements for the rotation form of the incompressible Navier--Stokes equations in \cite{Layton20093433}, and pointed out that the resolution of the Bernoulli pressure and the strong coupling between the velocity and the Bernoulli pressure approximations especially for high Reynolds numbers possibly lead to a loss of the accuracy, which is also mentioned in \cite{Layton2010916, Olshanskii20025515}. In \cite{Linke2016}, the authors applied the pressure-robust method straight to the rotation form of the incompressible Navier--Stokes equations by using the standard discrete trilinear form with respect to the nonlinear term $(\mathrm{rot}\,\bm{u})\times\bm{u}$ with an additional nonzero term in the discrete level to preserve the skew-symmetry. They provided a priori $H^1$-error estimates for the velocity, and the advantage of the pressure-robust method was shown via several numerical examples including potential flows, irrotational flows with nonlinear convection and rigid body rotation with nonlinear convection. Later, Quiroz and Di Pietro \cite{Quiroz20202655} embedded this pressure-robust method into the hybrid high-order (HHO) method, which is a nonconforming finite element method, for the rotation form of the stationary incompressible Navier--Stokes equations. Likewise, the intended results were successfully achieved by the proposed pressure-robust HHO method. In addition, as a recent published result \cite{Liu2021375}, a modified nonconforming VEM with the pressure-robust method was proposed for the convection form (i.e., $(\bm{u}\cdot\nabla)\bm{u}$) of the stationary incompressible Navier--Stokes equations. The authors established an optimal convergence results for $H^1$, $L^2$-velocity and $L^2$-pressure for the proposed pressure-robust VEM method, and achieved the pressure-independence of velocity errors and the effectiveness of small viscosities.

It is indisputable, to our best knowledge, that the classical conforming finite element methods are more popular from the views of scheme complexity and implementation by codes. Hence, our work is based on the inf-sup stable mixed conforming finite elements with high order. We propose a novel skew-symmetric discrete trilinear form containing a velocity reconstruction operator for the nonlinear term $(\mathrm{rot}\,\bm{u})\times\bm{u}$, in which not only the velocity test function is changed to preserve the skew-symmetry. The most advantage is that it is discretized pressure-robustly straight from the continuous trilinear form without additional terms. Besides, compared to the discrete rotation form proposed in Remark 7.3 \cite{Ahmed2018}, it benefits the error analysis much. As the theory contributions, both optimal $H^1$, $L^2$-error estimates for the velocity and an optimal $L^2$-error estimate for the Bernoulli pressure are completely established. Furthermore, contrast to the numerical examples shown in \cite{Linke2016, Liu2021375, Quiroz20202655}, various practical and more complicated numerical experiments are presented to demonstrate the efficiency of the proposed high-order pressure-robust method, and by them one can also find that there is no obvious loss of the accuracy for the velocity and pressure at high Reynolds numbers.

The rest of this article is organized as follows. In section \ref{sec:setting}, we introduce some useful notations, lemmas and necessary preliminaries. In section \ref{sec:discrete}, an $\bm{H}(\mathrm{div})$-conforming velocity reconstruction operator $\mathcal R$ is constructed with respect to the high-order inf-sup stable mixed elements $\mathbf{P}_k^{\mathrm{bubble}}$-$\mathrm{P}_{k-1}^{\mathrm{dc}}$ with any integer $k\geqslant 2$. Then, a novel skew-symmetric discrete trilinear form containing the reconstruction operator $\mathcal R$ is proposed, and the well-posedness of the corresponding discrete weak formulation is also proved. In section \ref{sec:error:estimate}, a complete error analysis of the $H^1$, $L^2$-error estimates for the velocity and the $L^2$-error estimate for the Bernoulli pressure is carried out, and all of their convergence orders are optimal. The theoretical results also show that the velocity error is independent of the Bernoulli pressure indeed. In section \ref{sec:numer}, several numerical experiments are provided to demonstrate the theoretical results and the performance of the proposed high-order pressure-robust method, including a comparison with the classical mixed method. Finally, the conclusion follows in section \ref{sec:conc}.

\section{Settings of Problem}
\label{sec:setting}

\subsection{Some notations}
\label{subsec:prelim}

Throughout this paper we use the following standard function spaces. Spaces and variables for vector-valued functions are both indicated with bold letters. For a Lipschitz domain $D\subset\mathbb{R}^d\ (d=1,2,3)$, we denote by $W^{s,p}(D)$ the Sobolev space with indexes $s\geqslant 0$, $1\leqslant p\leqslant\infty$ of real-valued functions defined on $D$, endowed with the seminorm $|\cdot|_{W^{s,p}(D)}$ denoted by $|\cdot|_{s,p,D}$ and norm $\|\cdot\|_{W^{s,p}(D)}$ denoted by $\|\cdot\|_{s,p,D}$. When $p=2$, $H^s(D)$ is denoted as $W^{s,2}(D)$ and the corresponding seminorm and norm are written as $|\cdot|_{s,D}$ and $\|\cdot\|_{s,D}$, respectively. Furthermore, with $|D|$ we denote the $d$-dimensional Hausdorff measure of $D$.

Given an open polygonal domain $\Omega\subset\mathbb{R}^2$ (for ease of analysis) with the boundary $\Gamma:=\partial\Omega$, we define some necessary Hilbert spaces as follows.

\begin{align*}
\bm{X}&=\bm{H}_0^1(\Omega):=\left\{\bm{v}\in \bm{H}^1(\Omega):\,\bm{v}=\bm{0}\ \text{on}\ \Gamma\right\},\quad \bm{Y}:=\bm{L}^2(\Omega),\\
Q&=L_0^2(\Omega):=\left\{q\in L^2(\Omega):\,|\Omega|^{-1}\int_\Omega q\,\mathrm{d}\bm{x}=0\right\}.
\end{align*}
We also need the closed subspace $\bm{V}$ of $\bm{X}$ given by

\[
\bm{V}:=\left\{\bm{w}\in\bm{X}:\,\mathrm{div}\,\bm{w}=0\ \text{in}\ \Omega\right\},
\]
and its dual space $\bm{V}^\ast$ equipped with the following dual norm:

\[
\|\bm{f}\|_{V^\ast}:=\sup_{\bm{w}\in\bm{V}\setminus\{\bm{0}\}}
\frac{\int_\Omega\bm{f}\cdot\bm{w}\,\mathrm{d}\bm{x}}{\|\nabla\bm{w}\|_{0,\Omega}},
\quad\forall\,\bm{f}\in\bm{V}^\ast.
\]
We denote by $\bm{H}$ the closed subspace of $\bm{Y}$, which is defined as:

\[
\bm{H}:=\left\{\bm{v}\in\bm{Y}:\,\mathrm{div}\,\bm{v}=0\ \text{in}\ \Omega,\ \bm{v}\cdot\bm{n}=0\ \text{on}\ \Gamma\right\},
\]
where $\bm{n}$ represents the outward unit normal vectors to $\Gamma$. In addition, $(\cdot,\cdot)_D$ represents the $L^2$ inner product in arbitrary bounded domain $D\subset\mathbb{R}^2$, and $\langle\cdot,\cdot\rangle_{\partial D}$ represents the $L^2$ inner product (or duality pairing) on its boundary $\partial D$. In particular, the subscript will be dropped for convenience when $D=\Omega$ or $\partial D=\Gamma$.

For the 2D gradient operator $\nabla$, we will define the vector-valued and tensor-valued results when it is taken on any scalar function $\psi\in H^1$ and any vector-valued function $\bm{v}=(v_1,v_2)^\top\in\bm{H}^1$, respectively.

\[
\nabla \psi:=\left(\frac{\partial\psi}{\partial x_1},\frac{\partial\psi}{\partial x_2}\right),\quad
\nabla\bm{v}:=\left[\frac{\partial v_i}{\partial x_j}\right]_{i,j}.
\]
For arbitrary two vectors $\bm{w},\bm{z}\in\mathbb{R}^2$ and arbitrary two matrices $\underline{\sigma},\underline{\delta}\in\mathbb{R}^{2\times 2}$, we define $\underline{\sigma}:\underline{\delta}=\sum_{i,j=1}^2\sigma_{ij}\delta_{ij}$, and as defined in \cite{Cockburn2002319}, let $\bm{w}\otimes\bm{z}$ denote the matrix whose $ij$th component is $w_iz_j$. Then contrary to the three-dimensional case, in the two-dimensional case two rotation operators $\mathrm{rot}:\,\bm{H}^1\rightarrow\mathbb{R}$ and $\bm{\mathrm{curl}}:\,H^1\rightarrow\mathbb{R}^2$ are defined as:

\[
\mathrm{rot}\,\bm{v}:=\frac{\partial v_2}{\partial x_1}-\frac{\partial v_1}{\partial x_2},\quad\forall\,\bm{v}=(v_1,v_2)^\top\in\bm{H}^1,
\]
and

\[
\bm{\mathrm{curl}}\,\psi:=\left(\frac{\partial\psi}{\partial x_2}, -\frac{\partial\psi}{\partial x_1}\right)^\top,\quad\forall\,\psi\in H^1,
\]
respectively. We also define a vector-valued function about rotation:

\[
\omega\times\bm{v}:=(-\omega v_2,\,\omega v_1)^\top,\quad\forall\,\omega\in\mathbb{R},\ \forall\,\bm{v}=(v_1,v_2)^\top\in\bm{H}^1.
\]

Based on the above definitions, we can obtain the following identity:

\begin{equation}\label{relation:conv:rot}
  (\mathrm{rot}\,\bm{v})\times\bm{v}+\frac{1}{2}\nabla|\bm{v}|^2=(\bm{v}\cdot\nabla)\bm{v},\quad
  \forall\,\bm{v}\in\bm{H}^1(\Omega).
\end{equation}
Besides, the following lemma holds that will play an important role on latter analysis and computation.

\begin{lemma}\label{lem:rot:skew}
Let $\mathcal D$ denote a simply connected open polygonal subset of $\Omega$, then for all $\bm{w},\bm{z},\bm{v}\in\bm{H}^1(\mathcal D)$, it holds

\begin{equation}\label{relation:rot:skew}
\int_{\mathcal D}((\mathrm{rot}\,\bm{w})\times\bm{z})\cdot\bm{v}\,\mathrm{d}\bm{x}
=\int_{\mathcal D}((\nabla\bm{w})\bm{z})\cdot\bm{v}\,\mathrm{d}\bm{x}
-\int_{\mathcal D}((\nabla\bm{w})\bm{v})\cdot\bm{z}\,\mathrm{d}\bm{x}.
\end{equation}
\end{lemma}

\begin{pf}
Let $\bm{w}=(w_1,w_2)^\top$, $\bm{z}=(z_1,z_2)^\top$, and $\bm{v}=(v_1,v_2)^\top$. Since

$$
\begin{aligned}
(\mathrm{rot}\,\bm{w})\times\bm{z}&=\left(z_2\frac{\partial w_1}{\partial x_2}-z_2\frac{\partial w_2}{\partial x_1},\;z_1\frac{\partial w_2}{\partial x_1}-z_1\frac{\partial w_1}{\partial x_2}\right)^\top,\\
(\nabla\bm{w})\bm{z}-\mathrm{div}\,(\bm{w}\cdot\bm{z})+(\nabla\bm{z})^\top\bm{w}
&=\left(z_2\frac{\partial w_1}{\partial x_2}-z_2\frac{\partial w_2}{\partial x_1},\;z_1\frac{\partial w_2}{\partial x_1}-z_1\frac{\partial w_1}{\partial x_2}\right)^\top,
\end{aligned}
$$
we have the following identity

$$
((\mathrm{rot}\,\bm{w})\times\bm{z},\bm{v})_{\mathcal D}=((\nabla\bm{w})\bm{z},\bm{v})_{\mathcal D}-(\mathrm{div}\,(\bm{w}\cdot\bm{z}),\bm{v})_{\mathcal D}+((\nabla\bm{z})^\top\bm{w},\bm{v})_{\mathcal D}.
$$
Applying integration by parts and matrix transformations we arrive at

$$
\begin{aligned}
-(\mathrm{div}\,(\bm{w}\cdot\bm{z}),\bm{v})_{\mathcal D}&=(\mathrm{div}\,\bm{v},\bm{w}\cdot\bm{z})_{\mathcal D}-\langle\bm{w}\cdot\bm{z},\bm{v}\cdot\bm{n}_{\mathcal D}\rangle_{\partial\mathcal D},\\
((\nabla\bm{z})^\top\bm{w},\bm{v})_{\mathcal D}&=\int_{\mathcal D}\nabla\bm{z}:(\bm{w}\otimes\bm{v})\,\mathrm{d}\bm{x}=
-(\mathrm{div}\,(\bm{w}\otimes\bm{v}),\bm{z})_{\mathcal D}+\langle\bm{w}\cdot\bm{z},\bm{v}\cdot\bm{n}_{\mathcal D}\rangle_{\partial\mathcal D}\\
&=-(\mathrm{div}\,\bm{v},\bm{w}\cdot\bm{z})_{\mathcal D}-((\nabla\bm{w})\bm{v},\bm{z})_{\mathcal D}+\langle\bm{w}\cdot\bm{z},\bm{v}\cdot\bm{n}_{\mathcal D}\rangle_{\partial\mathcal D},
\end{aligned}
$$
where $\bm{n}_{\mathcal D}$ is denoted as the outward unit normal vectors to $\partial\mathcal D$. Hence, the identity \eqref{relation:rot:skew} then holds.

\end{pf}

\begin{remark}\label{rem:Hdiv:test}
From the above details of proving Lemma \ref{lem:rot:skew}, it follows that the regularity of $\bm{v}$ can be relaxed to belonging to the space $\bm{H}(\mathrm{div};\mathcal D)$.
\end{remark}

\subsection{Navier--Stokes equations and continuous Helmholtz projection}
\label{subsec:cont:NS}

In this work, we shall consider the rotation form of the stationary incompressible Navier--Stokes equations reading as: find a pair $(\bm{u},p)\in\bm{X}\times Q$ such that

\begin{equation}\label{eq:con:NS}
\begin{split}
-\nu\Delta\bm{u}+(\mathrm{rot}\,\bm{u})\times\bm{u}+\nabla p\,&=\,\bm{f}\quad\mathrm{in}\ \Omega,\\
\mathrm{div}\,\bm{u}\,&=\,0\,\quad\mathrm{in}\ \Omega,\\
\bm{u}\,&=\,\bm{0}\,\quad\mathrm{on}\ \Gamma,
\end{split}
\end{equation}
where $\nu$ is a positive constant representing the kinematic viscosity, and let $\bm{f}\in\bm{Y}$ be a given body force.
\begin{remark}
Note that the pressure $p$ here is called the Bernoulli pressure satisfying the equation $p=p^{\mathrm{kin}}+\frac{1}{2}|\bm{u}|^2$ with $p^{\mathrm{kin}}$ the kinematic pressure, based on \eqref{relation:conv:rot}.
\end{remark}
The bilinear forms $a(\cdot,\cdot)$, $d(\cdot,\cdot)$ and the trilinear form $b(\cdot;\cdot,\cdot)$ are defined as

\begin{align*}
a(\bm{w},\bm{v})&:=(\nabla\bm{w},\nabla\bm{v}),
\qquad\qquad\qquad\qquad\qquad\qquad\quad\quad\quad\ \
\forall\,\bm{w},\bm{v}\in\bm{X},\\
d(\bm{v},q)&:=-(q,\mathrm{div}\,\bm{v}),
\qquad\qquad\qquad\qquad\qquad\quad\quad\quad\quad\quad\,\;
\forall\,\bm{v}\in\bm{X},\ \forall\,q\in Q,\\
b(\bm{w};\bm{z},\bm{v})&:=((\mathrm{rot}\,\bm{w}\times\bm{z}),\bm{v})
=((\nabla\bm{w})\bm{z},\bm{v})-((\nabla\bm{w})\bm{v},\bm{z}),\quad\forall\,
\bm{w},\bm{z},\bm{v}\in\bm{H}^1(\Omega).
\end{align*}
Hence, the Galerkin weak formulation of \eqref{eq:con:NS} reads as: find $(\bm{u},p)\in\bm{X}\times Q$ such that $\forall\,(\bm{v},q)\in\bm{X}\times Q$,

\begin{align}
\label{cont:rot:NS:momentum}
  \nu a(\bm{u},\bm{v})+b(\bm{u};\bm{u},\bm{v})+d(\bm{v},p)&=(\bm{f},\bm{v}),\\
\label{cont:rot:NS:mass}
  d(\bm{u},q)&=0.
\end{align}

\begin{remark}\label{rem:comp:bcs}
If the boundary $\Gamma$ is dissected into two parts $\Gamma_D$ and $\Gamma_N$ with $\Gamma=\Gamma_D\cup\Gamma_N$, $\Gamma_D\cap\Gamma_N=\emptyset$, and a more complicated boundary condition is employed as follows:

\[
\bm{u}=\bm{u}_D\quad\mathrm{on}\ \Gamma_D;\quad (p^{\mathrm{kin}}\mathbb{I}-\nu\nabla\bm{u})\bm{n}=\bm{u}_N
\quad\mathrm{on}\ \Gamma_N,
\]
where $\mathbb{I}$ denotes the $2\times 2$ identity matrix and $\bm{u}_D$, $\bm{u}_N$ are given functions, then \eqref{cont:rot:NS:momentum} will be changed into

\begin{equation}
\label{cont:rot:NS:momentum:change}
\nu a(\bm{u},\bm{v})
+b(\bm{u};\bm{u},\bm{v})
+\frac{1}{2}\left\langle|\bm{u}|^2,
\bm{v}\cdot\bm{n}\right\rangle_{\Gamma_N}
+d(\bm{v},p)=(\bm{f},\bm{v}),
\end{equation}
for any $\bm{v}\in\bm{H}_{0,D}^1(\Omega):=\{\bm{w}\in\bm{H}^1(\Omega):\,\bm{w}|_{\Gamma_D}=\bm{0}\}$ due to the integration $(\frac{1}{2}\nabla|\bm{u}|^2,\bm{v})$. The discretised form of \eqref{cont:rot:NS:momentum:change} will appear in the numerical experiments of section \ref{sec:numer}, and the discretised term of $\frac{1}{2}\left\langle|\bm{u}|^2,
\bm{v}\cdot\bm{n}\right\rangle_{\Gamma_N} $ will play a decisive role to get correct numerical results.
\end{remark}

The bilinear form $a(\cdot,\cdot)$ is continuous and coercive on $\bm{X}\times\bm{X}$; the bilinear form $d(\cdot,\cdot)$ is continuous on the couple $\bm{X}\times Q$ and satisfies the inf-sup condition \cite{Girault1986}, i.e., there is a positive constant $\beta$ such that, for any $q\in Q$,

\begin{equation}\label{cont:inf:sup}
  \sup_{\forall\,\bm{w}\in\bm{X}}\frac{d(\bm{w},q)}{\|\nabla\bm{w}\|_{0,\Omega}}
  \geqslant\beta\|q\|_{0,\Omega}.
\end{equation}
It follows from Remark 3.35 in \cite{John2016} that

\begin{equation}\label{ineq:div:grad}
  \|\mathrm{div}\,\bm{w}\|_{0,\Omega}\leqslant\|\nabla\bm{w}\|_{0,\Omega},\quad\forall\,\bm{w}\in\bm{X}.
\end{equation}
According to chap. 4 in \cite{Boffi2013}, we define the norm of the bilinear form $\|d\|$ as

\[
\|d\|:=\sup_{\forall\,\bm{w}\in\bm{X},\,\forall\,q\in Q}\frac{d(\bm{w},q)}
{\|\nabla\bm{w}\|_{0,\Omega}\|q\|_{0,\Omega}},
\]
and the bound that $\|d\|\leqslant 1$ derives straight from \eqref{ineq:div:grad}. In addition, the classical Sobolev inequality reads as: if $\mathcal D$ is denoted as arbitrary simply connected open polygonal subset of $\Omega$, then it holds

\begin{equation}\label{cont:Sobolev:L4}
  \|\psi\|_{0,q,\mathcal D}\leqslant C_{\mathrm{s}}\|\psi\|_{1,\mathcal D},\quad
  \forall\,\psi\in H^1(\mathcal D),\quad 1\leqslant q\leqslant 6,
\end{equation}
where the positive constant $C_{\mathrm{s}}$ depends only on $\mathcal D$. As shown below, the trilinear form $b(\cdot;\cdot,\cdot)$ is skew-symmetric on the product space $\bm{H}^1(\Omega)\times\bm{H}^1(\Omega)\times\bm{H}^1(\Omega)$, and it follows from the H\"{o}lder inequality with exponents $(\frac{1}{2},\frac{1}{4},\frac{1}{4})$, the Sobolev inequality \eqref{cont:Sobolev:L4} in $L^4(\Omega)$ and the Poincar\'{e}--Friedrichs inequality that the continuity holds on $\bm{X}\times\bm{X}\times\bm{X}$.

\begin{align}
\label{cont:skew:sym}
b(\bm{w};\bm{z},\bm{v})
&=-b(\bm{w};\bm{v},\bm{z}),
\qquad\qquad\qquad\qquad\,\forall\,\bm{w},\bm{z},\bm{v}\in\bm{H}^1(\Omega),\\
\label{cont:bound}
|b(\bm{w};\bm{z},\bm{v})|
&\leqslant N\|\nabla\bm{w}\|_{0,\Omega}\|\nabla\bm{z}\|_{0,\Omega}
\|\nabla\bm{v}\|_{0,\Omega},
\quad\forall\,\bm{w},\bm{z},\bm{v}\in\bm{X},
\end{align}
where $N$ is a positive constant independent of $h$, $\bm{w}$, $\bm{z}$ and $\bm{v}$.

Next, without proof we introduce the following lemma about the famous Helmholtz-Hodge decomposition. One can read \cite{Ahmed2018, Gauger201989, Girault1986, John2017, Linke2016} and the references therein for more details.

\begin{lemma}\label{lem:hhdec}
Every vector field $\bm{f}\in\bm{Y}$ can be uniquely decomposed into $\nabla\phi$ with $\phi\in H^1(\Omega)$, and a divergence-free vector field $\bm{f}_0\in\bm{H}$, i.e., $\bm{f}=\bm{f}_0+\nabla\phi$.
\end{lemma}
By integration by parts, we know that for any $\bm{w}\in\bm{H}$ and any $\psi\in H^1(\Omega)$,

\[
(\bm{w},\nabla\psi)
=\langle\bm{w}\cdot\bm{n},\psi\rangle-(\psi,\mathrm{div}\,\bm{w})
=0.
\]
Hence, this decomposition is $L^2$-orthogonal. The divergence-free remainder $\bm{f}_0$ is called the Helmholtz projector, which is denoted as $\mathbb{P}(\bm{f})$. In addition, it holds that $\mathbb{P}(\nabla\psi)=0$ for all $\psi\in H^1(\Omega)$.

\begin{remark}\label{remark:relation:HPconv:rot}
Recalling the identity \eqref{relation:conv:rot}, one can find that
$\mathbb{P}((\bm{v}\cdot\nabla)\bm{v})=\mathbb{P}((\mathrm{rot}\,\bm{v})\times\bm{v})$ holds for any $\bm{v}\in\bm{X}$.
\end{remark}

Owing to the Helmholtz projector, an a priori estimate of the continuous velocity satisfying \eqref{cont:rot:NS:momentum}--\eqref{cont:rot:NS:mass} are not quite the same as the classical results provided in \cite{Girault1986, He2009, He2008, Temam1984}.

\begin{theorem}\label{thm:cont:wellposed}
  Given $\bm{f}\in\bm{Y}$, there exists at least a solution pair $(\bm{u},p)\in\bm{X}\times Q$ which satisfies \eqref{cont:rot:NS:momentum}--\eqref{cont:rot:NS:mass} and

  \begin{equation}
  \label{cont:ineq:stable}
    \|\nabla\bm{u}\|_{0,\Omega}\leqslant\nu^{-1}\|\mathbb{P}(\bm{f})\|_{V^\ast}.
  \end{equation}
  Furthermore, if the condition $\nu^{-2}N\|\mathbb{P}(\bm{f})\|_{V^\ast}=:\sigma<1$ holds, the solution pair $(\bm{u},p)$ is unique.
\end{theorem}

\begin{pf}
The existence and uniqueness of solution to \eqref{cont:rot:NS:momentum}--\eqref{cont:rot:NS:mass} are classical results, see \cite{Girault1986, Temam1984}. To avoid repeating, we just verify \eqref{cont:ineq:stable}. Here we assume that $\bm{u}\in\bm{V}$ is the continuous velocity solution.

Recalling Lemma \ref{lem:hhdec}, for the given $\bm{f}\in\bm{Y}$, we have $\bm{f}=\mathbb{P}(\bm{f})+\nabla\phi$. Then plugging it into \eqref{cont:rot:NS:momentum}, and taking $\bm{v}=\bm{u}$ in \eqref{cont:rot:NS:momentum}, $q=p-\phi$ in \eqref{cont:rot:NS:mass} yield $\nu a(\bm{u},\bm{u})+b(\bm{u};\bm{u},\bm{u})=(\mathbb{P}(\bm{f}),\bm{u})$. Based on \eqref{cont:skew:sym}, we obtain $\nu\|\nabla\bm{u}\|_{0,\Omega}^2\leqslant \|\mathbb{P}(\bm{f})\|_{V^\ast}\|\nabla\bm{u}\|_{0,\Omega}$, which directly leads to \eqref{cont:ineq:stable}.

\end{pf}

\section{Discrete Problem}
\label{sec:discrete}

Let the mesh $\mathcal T_h$ be a subdivision of $\Omega$ with disjoint triangles. We set

\[
h=\max_{T\in\mathcal T_h}h_T,\quad h_T=\mathrm{diam}(T),\quad\forall\,T\in\mathcal T_h.
\]
We assume that $\mathcal T_h$ is shape-regular, i.e., $h_T/\rho_T\leqslant\gamma$ for all triangle $T\in\mathcal T_h$ where $\rho_T$ is the radius of the largest ball inscribed in $T$, and $\gamma$ is a positive constant independent of $h$. In addition, we denote by $L$ the diameter of $\Omega$.

\subsection{High-order inf-sup stable mixed conforming finite element spaces}

Now referring to chap. 2., section 2.2 in \cite{Girault1986}, we define the finite element spaces for the discrete velocity field and pressure respectively as follows.

\[
\bm{X}_h:=\left[\mathcal L_k^1\oplus B_{k+1}\right]^2\cap\bm{X},\quad
Q_h:=\mathcal L_{k-1}^0\cap Q,
\]
where $k\geqslant 2$ is a positive integer. The spaces $\mathcal L_k^s$ and $B_{k+1}$ are respectively defined as

\begin{align*}
\mathcal L_k^s&:=\left\{v\in H^s(\Omega):\,v|_T\in P_k(T),\ \forall\,T\in\mathcal T_h\right\},\\
B_{k+1}&:=\left\{v:\,v|_T\in \widetilde{P}_{k+1}(T)\cap H_0^1(T),\ \forall\,T\in\mathcal T_h\right\},
\end{align*}
where $P_l(D)$ denotes the space of all polynomials with total degree not greater than $l$ defined on domain $D$, and $\widetilde{P}_l(D)$ denotes the space of homogeneous polynomials of degree $l$, i.e., $\widetilde{P}_l=\mathrm{span}\{x_1^ix_2^{l-i};0\leqslant i\leqslant l\}$. Hence, every $v$ in $B_{k+1}$, on arbitrary triangle $T\in\mathcal T_h$, has the form $\alpha(T)\lambda_1\lambda_2\lambda_3$ with $\alpha(T)\in \widetilde{P}_{k-2}(T)$ and $\lambda_i$s the triangular area coordinates, $i=1,2,3$. Vividly, the mixed finite element space $\bm{X}_h\times Q_h$ is denoted as $\mathbf{P}_k^{\mathrm{bubble}}$-$\mathrm{P}_{k-1}^{\mathrm{dc}}$ hereafter.

Then, the spaces $\bm{X}_h$ and $Q_h$ have the following approximation property and projection property, respectively.

\begin{enumerate}[({M}1)]
\item $\forall\,\bm{v}\in\bm{X}\cap \bm{H}^{k+1}(\Omega)$, there exist an approximation $\bm{\Pi}_h \bm{v}\in\bm{X}_h$ and a positive constants $C_\Pi$ independent of $h$ and $\bm{v}$, such that

    \begin{equation}\label{appro:vel}
      \|\bm{v}-\bm{\Pi}_h \bm{v}\|_{0,\Omega}
      +h|\bm{v}-\bm{\Pi}_h \bm{v}|_{1,\Omega}\leqslant C_\Pi h^{k+1}|\bm{v}|_{k+1,\Omega}.
    \end{equation}

\item For each $T\in\mathcal T_h$ and some nonnegative integer $m$ satisfying $0\leqslant m\leqslant k-1$, let $\pi_{h,T}^m$ be the local $L^2$-projection onto $P_m(T)$ such that

    \begin{equation}\label{proj:pressure}
      \|q-\pi_{h,T}^m q\|_{0,T}+h_T|q-\pi_{h,T}^m q|_{1,T}\leqslant C_\pi h_T^{m+1}|q|_{m+1,T},\quad \forall\,q\in H^{m+1}(T),
    \end{equation}
    where $C_\pi$ is a positive constant independent of $h$. It will sometimes be combined to the global $L^2$-projection $\pi_h:\,Q\rightarrow Q_h$ when $m=k-1$.
\end{enumerate}
By the above definition of the $L^2$-projecton $\pi_h$, we can define the discrete divergence as $\mathrm{div}_h\,\bm{w}:=\pi_h(\mathrm{div}\,\bm{w})$ for any $\bm{w}\in\bm{X}$. The corresponding space of discretely divergence-free functions is given by

\[
\bm{V}_h:=\{\bm{w}_h\in\bm{X}_h:\;\mathrm{div}_h\,\bm{w}_h=0\}.
\]

\begin{remark}\label{rem:H1Pk:appro}
For any $\bm{v}\in\bm{X}$, let $\bm{\mathrm{I}}_h^{k}\bm{v}$ be its standard Lagrange interpolation in the Lagrange finite element space $\bm{W}_h:=[\mathcal L_k^1]^2\cap\bm{X}$. From affine transformations, scaling arguments and the famous Bramble--Hilbert's Lemma, it follows that for all integers $m$ and real $s$ with $0\leqslant m\leqslant s+1$, $1\leqslant s\leqslant k$,

\begin{equation}\label{H1Pk:appro:vel}
|\bm{v}-\bm{\mathrm{I}}_h^{k}\bm{v}|_{m,\Omega}\leqslant C_I h^{s+1-m}|\bm{v}|_{s+1,\Omega},
\quad\forall\,\bm{v}\in\bm{X}\cap\bm{H}^{s+1}(\Omega),
\end{equation}
where the constant $C_I$ is positive and independent of $h$ and $\bm{v}$. By comparison between \eqref{appro:vel} and \eqref{H1Pk:appro:vel}, one can find that the approximation ``effect'' of $\bm{\mathrm{I}}_h^{k}\bm{v}$ is the same as $\bm{\Pi}_h\bm{v}$ in the norm of $\|\cdot\|_{0,\Omega}$ or $\|\cdot\|_{1,\Omega}$ but only the estimate constants when $s=k$ and $m\in\{0,1\}$.
\end{remark}

\begin{remark}\label{remark:disc:divfree:equiv}
Note that the domain $\Omega$ is assumed as a two-dimensional polygon so that there always exists a shape-regular triangulation $\mathcal T_h$ satisfying

\[
\bigcup_{T\in\mathcal T_h}\overline{T}=\overline{\Omega}.
\]
Hence, the $L^2$ inner product $(\cdot,\cdot)$ on $\Omega$ has another equivalent form $\sum_{T\in\mathcal T_h}(\cdot,\cdot)_T$. As a result, in this work the equivalent form of $\bm{V}_h$ is given by

\begin{equation}\label{equi:weak:divfree}
  \bm{V}_h=\{\bm{w}_h\in\bm{X}_h:\;-(q_h,\mathrm{div}\,\bm{w}_h)=0,\ \forall\,q_h\in Q_h\},
\end{equation}
where the divergence term $\mathrm{div}\,\bm{w}_h$ will be replaced by $\mathrm{div}_h\,\bm{w}_h$ if needs to be distinguished.
\end{remark}

Arising from \eqref{equi:weak:divfree}, we shall consider the discrete inf--sup stability. According to chap. 2., Lemma 2.6 in \cite{Girault1986} or Proposition 8.6.2. in \cite{Boffi2013}, for any integer $k\geqslant 2$, $\mathbf{P}_k^{\mathrm{bubble}}$-$\mathrm{P}_{k-1}^{\mathrm{dc}}$ is inf-sup stable, i.e., there is a positive constant $\beta_0$ independent of $h$ such that, for any $q_h\in Q_h$,

\begin{equation}\label{discrete:inf:sup}
  \sup_{\forall\bm{w}_h\in\bm{X}_h}
  \frac{-(q_h,\mathrm{div}\,\bm{w}_h)}{\|\nabla\bm{w}_h\|_{0,\Omega}}
  \geqslant\beta_0\|q_h\|_{0,\Omega}.
\end{equation}

\subsection{Pressure-robustness by velocity reconstruction}

We first define the following two discrete dual norms with respect to the discrete space $\bm{V}_h$:

\begin{align*}
\|\mathscr{L}\|_{V_h^\ast}&:=\sup_{\bm{w}_h\in\bm{V}_h\setminus\{\bm{0}\}}
\frac{\mathscr{L}(\bm{w}_h)}{\|\nabla\bm{w}_h\|_{0,\Omega}},\quad\forall\,\mathscr{L}\in\bm{V}_h^\ast,\\
\|\bm{g}\|_{V_h^\ast}&:=\sup_{\bm{w}_h\in\bm{V}_h\setminus\{\bm{0}\}}\frac{(\bm{g},\bm{w}_h)}{\|\nabla\bm{w}_h\|_{0,\Omega}},
\quad\forall\,\bm{g}\in\bm{Y},
\end{align*}
where $\bm{V}_h^\ast$ is denoted as the dual space of $\bm{V}_h$. Then, like the continuous Helmholtz projector introduced in section \ref{subsec:cont:NS}, the classical discrete Helmholtz projecton $\widetilde{\mathbb{P}}_h$ is given by: for any given $\bm{f}\in\bm{Y}$,

\begin{equation}\label{discrete:class:HP}
\widetilde{\mathbb{P}}_h(\bm{f}):=
\mathop{\arg\min}_{\bm{w}_h\in\bm{V}_h}
\|\bm{f}-\bm{w}_h\|_{0,\Omega},\quad
\text{s.t.}\quad
 (\widetilde{\mathbb{P}}_h(\bm{f}),\bm{w}_h)
  =(\bm{f},\bm{w}_h),
  \quad\forall\,\bm{w}_h\in\bm{V}_h.
\end{equation}
Hence, for some $\bm{f}\in\bm{Y}$ with the continuous Helmholtz projection $\bm{f}=\mathbb{P}(\bm{f})+\nabla\phi$, by \eqref{discrete:class:HP} we have

\[
\|\mathbb{P}(\bm{f})-\widetilde{\mathbb{P}}_h(\bm{f})\|_{V_h^\ast}
=\sup_{\bm{w}_h\in\bm{V}_h\setminus\mathbb{R}}
\frac{(\mathbb{P}(\bm{f}),\bm{w}_h)-(\widetilde{\mathbb{P}}_h(\bm{f}),\bm{w}_h)}
{\|\nabla\bm{w}_h\|_{0,\Omega}}
=\sup_{\bm{w}_h\in\bm{V}_h\setminus\mathbb{R}}
\frac{-(\nabla\phi,\bm{w}_h)}{\|\nabla\bm{w}_h\|_{0,\Omega}}
=\|\widetilde{\mathbb{P}}_h(\nabla\phi)\|_{V_h^\ast}.
\]
According to \eqref{equi:weak:divfree}, for every $\bm{w}_h\in\bm{V}_h$, it holds $(q_h,\mathrm{div}\,\bm{w}_h)=0$ for any $q_h\in Q_h$. As a result, for every $\bm{w}_h\in\bm{V}_h$ and any $q_h\in Q_h$, it holds $(\widetilde{\mathbb{P}}_h(\nabla\phi),\bm{w}_h)
=-(\phi-q_h,\mathrm{div}\,\bm{w}_h)$. Hence, it follows from \eqref{ineq:div:grad} and \eqref{proj:pressure} that if $\phi\in H^k(\Omega)$,

\begin{equation}\label{diff:cont:disc:HP}
  \|\mathbb{P}(\bm{f})-\widetilde{\mathbb{P}}_h(\bm{f})\|_{V_h^\ast}
=\|\widetilde{\mathbb{P}}_h(\nabla\phi)\|_{V_h^\ast}
\leqslant \min_{q_h\in Q_h}\|\phi-q_h\|_{0,\Omega}\leqslant C_\pi h^k|\phi|_{k,\Omega}.
\end{equation}

The above inequality \eqref{diff:cont:disc:HP} implies the lack of pressure-robustness. In short, it lies at the root of the fact that $\widetilde{\mathbb{P}}_h(\nabla\psi)\neq 0$ for any $\psi\in H^1(\Omega)\setminus Q_h$. So, in order to make the discrete Helmholtz projection zero when applied to gradients, a natural method is to employ a reconstruction operator $\mathcal R$ that maps discretely divergence-free test functions onto divergence-free test functions \cite{Ahmed2018, Linke2016289, Linke2016}, i.e., $\mathcal R\bm{V}_h\subset\bm{H}$.

To this end, as the image space of $\mathcal R$, we shall choose the $k$-th order $\bm{H}(\mathrm{div})$-conforming Brezzi-Douglas-Marini element space $\mathbf{BDM}_k$, which is denoted as $\bm{R}_h$ hereafter, with respect to shape-regular triangular meshes. We now define a velocity reconstruction operator $\mathcal R:\,\bm{X}+\bm{X}_h\rightarrow\bm{R}_h$, such that for any $\bm{w}\in\bm{X}+\bm{X}_h$, any $T\in\mathcal T_h$ and any $e\subset\partial T$,

\begin{align}
\label{BDM:property:1}
\left(\bm{w}-\mathcal R(\bm{w}),\bm{v}_h\right)_T&=0,
\quad\forall\,\bm{v}_h\in\bm{M}_k(T);\\
\label{BDM:property:2}
\left\langle(\bm{w}-\mathcal R(\bm{w}))\cdot\bm{n}_e, q_h\right\rangle_e&=0,
\quad\forall\,q_h\in P_k(e),
\end{align}
where $\bm{n}_e$ denotes the outward unit normal vectors to the edge $e$,

\[
\bm{M}_k(T):=\left\{\nabla z_h+\mathbf{curl}\,(b_T\omega_h):\,(z_h,\omega_h)\in P_{k-1}(T)\times P_{k-2}(T)\right\},
\]
and $b_T=\lambda_1\lambda_2\lambda_3\in B_3(T)$ is the bubble function on $T$. It is well known that $\mathrm{dim}(\bm{\mathrm{BDM}}_k(T))=\mathrm{dim}(\bm{P}_k(T))$ for any triangle $T\in\mathcal T_h$ in the two dimensional case. Hence, from Proposition 2.3.1. and Lemma 2.3.2. in \cite{Boffi2013}, and \eqref{BDM:property:1}--\eqref{BDM:property:2}, it follows that for any $T\in\mathcal T_h$,

\begin{equation}\label{equal:Pk:Rec}
  \mathcal R(\bm{w}_h)|_T=\bm{w}_h|_T,
  \quad\forall\,\bm{w}_h\in\bm{W}_h.
\end{equation}

\begin{remark}
In fact, if $\bm{w}_h\in\bm{W}_h$, it will satisfy that $\bm{w}_h|_T\in\bm{P}_k(T)$ and $(\bm{w}_h\cdot\bm{n}_T)|_{\partial T}\in P_k(\partial T)$ for any $T\in\mathcal T_h$.
\end{remark}
As a result, combining \eqref{equal:Pk:Rec} and the famous Bramble--Hilbert's Lemma yields the following estimate: for any $T\in\mathcal T_h$,

\begin{equation}\label{BDM:property:3}
  \|\bm{w}-\mathcal R(\bm{w})\|_{s,T}\leqslant C_{\mathcal R}h_T^{m-s}|\bm{w}|_{m,T},
  \quad\forall\,\bm{w}\in\bm{X}+\bm{X}_h,\ m\in\{1,2\},\ s\in\{0,1\},
\end{equation}
where $C_{\mathcal R}$ is a positive constant depending only on the mesh regularity parameter $\gamma$.

Next, we should verify that $\mathcal R(\bm{w}_h)\in\bm{H}$ holds for any $\bm{w}_h\in\bm{V}_h$. In fact, for any $\bm{w}_h\in\bm{V}_h$ and any $q_h\in Q_h$, it follows from integration by parts and \eqref{BDM:property:1}--\eqref{BDM:property:2} that

\begin{align*}
0&=(q_h,\mathrm{div}_h\,\bm{w}_h)=
\sum_{T\in\mathcal T_h}(q_h,\mathrm{div}\,\bm{w}_h)_T
=\sum_{T\in\mathcal T_h}
\int_T\mathrm{div}(q_h\bm{w}_h)\,\mathrm{d}\bm{x}
-\int_T\nabla q_h\cdot\bm{w}_h\,\mathrm{d}\bm{x}\\
&=\sum_{T\in\mathcal T_h}\langle q_h,\bm{w}_h\cdot\bm{n}_T\rangle_{\partial T}-(\nabla q_h,\bm{w}_h)_T
=\sum_{T\in\mathcal T_h}\langle q_h,\mathcal R(\bm{w}_h)\cdot\bm{n}_T\rangle_{\partial T}-(\nabla q_h,\mathcal R(\bm{w}_h))_T\\
&=(q_h,\mathrm{div}\,\mathcal R(\bm{w}_h)),
\end{align*}
where $\bm{n}_T$ represents the outward unit normal vectors to the element boundary $\partial T$ for each $T\in\mathcal T_h$. Due to $\mathcal R(\bm{w}_h)\in\bm{R}_h$, it holds that $\mathrm{div}\,\mathcal R(\bm{w}_h)\in\mathrm{P}_{k-1}^{\mathrm{dc}}=Q_h$. We can then take $q_h=\mathrm{div}\,\mathcal R(\bm{w}_h)$ to obtain $\mathrm{div}\,\mathcal R(\bm{w}_h)=0$ in the sense of $\bm{H}(\mathrm{div};\Omega)$. Besides, we obtain from $\mathcal R(\bm{w}_h)\in\bm{R}_h$ again that $\mathcal R(\bm{w}_h)\cdot\bm{n}_e\in P_k(e)$, and if we take $q_h=\mathcal R(\bm{w}_h)\cdot\bm{n}_e$ in \eqref{BDM:property:2}, then it holds that $R(\bm{w}_h)\cdot\bm{n}_e=0$ on all $e\subset\Gamma$. Hence, the purpose that $\mathcal R\bm{V}_h\subset\bm{H}$ then realizes.

Further with the help of the reconstruction operator $\mathcal R$, the classical discrete Helmholtz projection $\widetilde{\mathbb{P}}_h$ can be repaired as: for any given $\bm{f}\in\bm{Y}$,

\begin{equation}\label{discrete:new:HP}
\mathbb{P}_h(\bm{f}):=
\mathop{\arg\min}_{\bm{w}_h\in\bm{V}_h}
\|\bm{f}-\mathcal R(\bm{w}_h)\|_{0,\Omega},\quad
\text{s.t.}\quad
 (\bm{f}-\mathbb{P}_h(\bm{f}),\mathcal R(\bm{w}_h))=0,
  \quad\forall\,\bm{w}_h\in\bm{V}_h.
\end{equation}
The following lemma rewrites that the above method accesses pressure-robustness in terms of the dual norm $\|\cdot\|_{V_h^\ast}$.

\begin{lemma}
\label{lem:fix:pres:robust}
For some $\bm{f}\in\bm{Y}$ with the continuous Helmholtz projection $\bm{f}=\mathbb{P}(\bm{f})+\nabla\phi$, it holds

\begin{align*}
&\mathbb{P}_h(\nabla\psi)=\mathbb{P}(\nabla\psi)=0,
\quad\forall\,\psi\in H^1(\Omega);\\
&\|(\mathbb{P}_h(\bm{f})-\mathbb{P}(\bm{f}))
\circ\mathcal R\|_{V_h^\ast}=0.
\end{align*}
\end{lemma}

\begin{pf}
An integration by parts, \eqref{discrete:new:HP}, and the fact that $\mathcal R(\bm{w}_h)\in\bm{H}$ for all $\bm{w}_h\in\bm{V}_h$ yield that for any $\psi\in H^1(\Omega)$,

\[
(\mathbb{P}_h(\nabla\psi),\mathcal R(\bm{w}_h))=
(\nabla\psi,\mathcal R(\bm{w}_h))
=-(\psi,\mathrm{div}\,\mathcal R(\bm{w}_h))=0,\quad\forall\,\bm{w}_h\in\bm{V}_h.
\]
The second property can be proved by:

\[
\|(\mathbb{P}_h(\bm{f})-\mathbb{P}(\bm{f}))
\circ\mathcal R\|_{V_h^\ast}
=\sup_{\bm{w}_h\in\bm{V}_h\setminus\{\bm{0}\}}
\frac{(\bm{f}-\mathbb{P}(\bm{f}),\mathcal R(\bm{w}_h))}{\|\nabla\bm{w}_h\|_{0,\Omega}}
=\sup_{\bm{w}_h\in\bm{V}_h\setminus\{\bm{0}\}}
\frac{\mathbb{P}_h(\nabla\phi)}{\|\nabla\bm{w}_h\|_{0,\Omega}}
=0.
\]
\end{pf}

We shall end this subsection by proposing some useful results about the reconstruction operator $\mathcal R$ as follows. We denote by $\|\bm{f}\|_{\ast,h}$ the discrete dual norm $\|\mathbb{P}(\bm{f})\circ\mathcal R\|_{V_h^\ast}$ for any $\bm{f}\in\bm{Y}$, and $C_P$ the standard Poincar\'{e}--Friedrichs inequality constant hereafter.

\begin{lemma}
\label{lem:discrete:new:rh}
For any $\bm{f}\in\bm{Y}$, it holds

\begin{equation}\label{ineq:discrete:new:rh}
  \|\bm{f}\|_{\ast,h}\leqslant\|\mathbb{P}(\bm{f})\|_{V^\ast}+C_{\mathcal R}L\|\mathbb{P}(\bm{f})\|_{0,\Omega}.
\end{equation}
\end{lemma}

\begin{pf}
It follows the Cauchy--Schwarz inequality, \eqref{BDM:property:3}, the fact that $\sum_{T\in\mathcal T_h}|\bm{w}|_{1,T}^2
\leqslant\|\nabla\bm{w}\|_{0,\Omega}^2$ for any $\bm{w}\in\bm{H}^1(\Omega)$, and the fact that $h<L$, that for any $\bm{f}\in\bm{Y}$,

\begin{align*}
  \|\bm{f}\|_{\ast,h}
  &=\sup_{\bm{w}_h\in\bm{V}_h\setminus\{\bm{0}\}}
  \frac{\int_\Omega\mathbb{P}(\bm{f})\cdot\mathcal R(\bm{w}_h)\,\mathrm{d}\bm{x}}{\|\nabla\bm{w}_h\|_{0,\Omega}}\\
  &=\sup_{\bm{w}_h\in\bm{V}_h\setminus\{\bm{0}\}}
  \frac{\int_\Omega\mathbb{P}(\bm{f})\cdot\bm{w}_h\,\mathrm{d}\bm{x}}
  {\|\nabla\bm{w}_h\|_{0,\Omega}}+\sup_{\bm{w}_h\in\bm{V}_h\setminus\{\bm{0}\}}
  \frac{\int_\Omega\mathbb{P}(\bm{f})\cdot(\mathcal R(\bm{w}_h)-\bm{w}_h)\,\mathrm{d}\bm{x}}{\|\nabla\bm{w}_h\|_{0,\Omega}}\\
  &\leqslant \|\mathbb{P}(\bm{f})\|_{V^\ast}+C_{\mathcal R}L\|\mathbb{P}(\bm{f})\|_{0,\Omega}.
\end{align*}
\end{pf}

\begin{lemma}[Discrete Sobolev Inequalities in $L^4(\Omega)$]
\label{lem:disc:sobolev:L4}
For any $\bm{w}_h\in\bm{X}_h$, it holds that

\begin{align}
\label{discrete:Sobolev:ineq:1}
\|\mathcal R(\bm{w}_h)-\bm{w}_h\|_{0,4,\Omega}
&\leqslant C_{1,S}\|\nabla\bm{w}_h\|_{0,\Omega},\\
\label{discrete:Sobolev:ineq:2}
\|\mathcal R(\bm{w}_h)\|_{0,4,\Omega}
&\leqslant C_{2,S}\|\nabla\bm{w}_h\|_{0,\Omega},
\end{align}
where $C_{1,S}$ and $C_{2,S}$ are two positive constants both independent of $h$ and $\bm{w}_h$, but possibly depend on $\Omega$, $k$ and $\gamma$.
\end{lemma}

\begin{pf}
First of all, based on Lemma 5.1 in \cite{DiPietro20172159} or Proposition 3 in \cite{Quiroz20202655}, for any $T\in\mathcal T_h$ the discrete Lebesgue embedding from $L^4(T)$ onto $L^2(T)$ reads as: given an integer $l\geqslant 0$ fixed, it holds that

\begin{equation}\label{discrete:Sobolev:L4}
  \|q_h\|_{0,4,T}\leqslant C_{\mathrm{ds}}h_T^{-1/2}\|q_h\|_{0,T},\quad
  \forall\,q_h\in P_l(T),
\end{equation}
where $C_{\mathrm{ds}}$ is a positive constant depends only on $l$ and $\gamma$. Then from \eqref{discrete:Sobolev:L4} and \eqref{BDM:property:3}, it follows that for any $T\in\mathcal T_h$ and any $\bm{w}_h\in\bm{X}_h$,

\[
  \|\mathcal R(\bm{w}_h)-\bm{w}_h\|_{0,4,T}\leqslant
  C_{\mathrm{ds}}h_T^{-1/2}\|\mathcal R(\bm{w}_h)-\bm{w}_h\|_{0,T}
  \leqslant C_{\mathcal R}C_{\mathrm{ds}}h_T^{1/2}|\bm{w}_h|_{1,T}.
\]
Summing the above inequality over all $T\in\mathcal T_h$ yields

\begin{align*}
&\|\mathcal R(\bm{w}_h)-\bm{w}_h\|_{0,4,\Omega}^4
=\sum_{T\in\mathcal T_h}\|\mathcal R(\bm{w}_h)-\bm{w}_h\|_{0,4,T}^4\leqslant (C_{\mathcal R}C_{\mathrm{ds}})^4\sum_{T\in\mathcal T_h}h_T^2|\bm{w}_h|_{1,T}^4\\
&\leqslant (C_{\mathcal R}C_{\mathrm{ds}})^4L^2
\|\nabla\bm{w}_h\|_{0,\Omega}^2
\left(\sum_{T\in\mathcal T_h}|\bm{w}_h|_{1,T}^2\right)
\leqslant(C_{\mathcal R}C_{\mathrm{ds}})^4L^2
\|\nabla\bm{w}_h\|_{0,\Omega}^4,
\end{align*}
where we have used the fact that $h_T\leqslant L$ along with $|\bm{w}_h|_{1,T}\leqslant\|\nabla\bm{w}_h\|_{0,\Omega}$ for any $T\in\mathcal T_h$, and another fact that $\sum_{T\in\mathcal T_h}|\bm{w}_h|_{1,T}^2
\leqslant\|\nabla\bm{w}_h\|_{0,\Omega}^2$ to conclude. Hence, if we set $C_{1,S}=C_{\mathcal R}C_{\mathrm{ds}}L^{1/2}$, the result \eqref{discrete:Sobolev:ineq:1} then holds.

The second result \eqref{discrete:Sobolev:ineq:2} follows from using the triangle inequality, \eqref{discrete:Sobolev:ineq:1}, the Sobolev inequality \eqref{cont:Sobolev:L4} in $L^4(\Omega)$, the Poincar\'{e}--Friedrichs inequality in turn, and setting $C_{2,S}=C_{1,S}+C_{\mathrm{s}}(C_P^2+1)^{1/2}$.

\end{pf}

\begin{lemma}[Consistency Error]
\label{lem:consistency:error}
For any $\bm{w}\in\bm{H}^{k+1}(\Omega)$ and any $\bm{v}\in\bm{X}+\bm{X}_h$, it holds that

\begin{equation}\label{disc:consistency:error}
  \left|(\nabla\bm{w},\nabla\bm{v})
  +(\Delta\bm{w},\mathcal R(\bm{v}))\right|\leqslant C_E h^k|\bm{w}|_{k+1,\Omega}
  \|\nabla\bm{v}\|_{0,\Omega},
\end{equation}
where $C_E$ is a positive constant independent of $h$, $\bm{w}$ and $\bm{v}$.
\end{lemma}

\begin{pf}
By adding and subtract $(\Delta\bm{w},\bm{v})$ and using an integration by parts, we obtain

$$
  (\nabla\bm{w},\nabla\bm{v})
  +(\Delta\bm{w},\mathcal R(\bm{v}))=
  (\Delta\bm{w},\mathcal R(\bm{v})-\bm{v})
  +(\nabla\bm{w},\nabla\bm{v})
  +(\Delta\bm{w},\bm{v})
  =(\Delta\bm{w},\mathcal R(\bm{v})-\bm{v}).
$$
The property \eqref{BDM:property:1} allows to subtract the local $L^2$-projection $\pi_{h,T}^{k-2}$, and then we shall use the Cauchy--Schwarz inequality, \eqref{proj:pressure} and \eqref{BDM:property:3} to estimate the only term $(\Delta\bm{w},\mathcal R(\bm{v})-\bm{v})$ as follows.

\begin{align*}
\left|(\Delta\bm{w},\mathcal R(\bm{v})-\bm{v})\right|&
=\left|\sum_{T\in\mathcal T_h}\left(\Delta\bm{w}-\pi_{h,T}^{k-2}\Delta\bm{w},
\mathcal R(\bm{v})-\bm{v}\right)_T\right|\\
&\leqslant C_\pi C_{\mathcal R}\left(\sum_{T\in\mathcal T_h}h_T^{2k-2}|\Delta\bm{w}|_{k-1,T}^2\right)^{1/2}
\left(\sum_{T\in\mathcal T_h}h_T^{2}\|\nabla\bm{v}|_{0,T}^2\right)^{1/2}\\
&\leqslant C_E h^k|\bm{w}|_{k+1,\Omega}
  \|\nabla\bm{v}\|_{0,\Omega},
\end{align*}
where we set $C_E=C_\pi C_{\mathcal R}$, which is independent of $h$, $\bm{w}$ and $\bm{v}$ indeed.

\end{pf}

\subsection{A novel discrete weak formulation and well-posedness}
\label{subsec:discrete:problem}

By Lemma \ref{lem:rot:skew} and the reconstruction operator $\mathcal R$, a different discrete trilinear form $b_h(\cdot;\cdot,\cdot)$ is defined as

\begin{equation}\label{discrete:tril}
b_h(\bm{w};\bm{z},\bm{v}):=
((\nabla\bm{w})\mathcal R(\bm{z}),\mathcal R(\bm{v}))-((\nabla\bm{w})\mathcal R(\bm{v}),\mathcal R(\bm{z})),\quad\forall\,\bm{w},\bm{z},\bm{v}\in\bm{X}+\bm{X}_h.
\end{equation}
Obviously it is skew-symmetric, i.e.,

\begin{equation}\label{skew:sym:discrete:tril}
b_h(\bm{w};\bm{z},\bm{v})=-b_h(\bm{w};\bm{v},\bm{z}),
\quad\forall\,\bm{w},\bm{z},\bm{v}\in\bm{X}+\bm{X}_h,
\end{equation}
and clearly \eqref{discrete:tril} stems from the continuous trilinear form $b(\cdot;\cdot,\cdot)$ without additional forms. Then, a novel discrete weak formulation of \eqref{eq:con:NS} reads as: find $(\bm{u}_h,p_h)\in\bm{X}_h\times Q_h$ such that $\forall\,(\bm{v}_h,q_h)\in\bm{X}_h\times Q_h$,

\begin{align}
\label{dis:rot:NS:momentum}
  \nu a(\bm{u}_h,\bm{v}_h)+b_h(\bm{u}_h;\bm{u}_h,\bm{v}_h)+d(\bm{v}_h,p_h)
  &=(\bm{f},\mathcal R(\bm{v}_h)),\\
\label{dis:rot:NS:mass}
  d(\bm{u}_h,q_h)&=0.
\end{align}

\begin{remark}\label{rem:num:kin:pres}
If necessary, in this case the discrete kinematic pressure $p_h^{\mathrm{kin}}$ can be computed from the Bernoulli pressure $p_h$ satisfying \eqref{dis:rot:NS:momentum} a posteriori by

\begin{equation}
\label{kin:pres:posteriori}
p_h^{\mathrm{kin}}:=p_h-\frac{1}{2}\,\mathop{\arg\min}_{q_h\in Q_h}\left\|q_h-|\bm{u}_h|^2\right\|_{0,\Omega}.
\end{equation}
\end{remark}

\begin{lemma}
\label{lem:disc:tri:bound}
For any $\bm{w}\in\bm{X}+\bm{X}_h$ and any $\bm{z}_h,\bm{v}_h\in\bm{X}_h$, there exists a positive constant $\mathcal N$ independent of $h$, such that

\begin{equation}\label{disc:tri:bound}
  |b_h(\bm{w};\bm{z}_h,\bm{v}_h)|
  \leqslant\mathcal N\|\nabla\bm{w}\|_{0,\Omega}
  \|\nabla\bm{z}_h\|_{0,\Omega}
  \|\nabla\bm{v}_h\|_{0,\Omega}.
\end{equation}
\end{lemma}

\begin{pf}
The proof follows in turn from the triangle inequality, a H\"{o}lder inequality with exponents $(\frac{1}{2},\frac{1}{4},\frac{1}{4})$, and \eqref{discrete:Sobolev:ineq:2}. To avoid repeating, we omit the details. In addition, one can derive that $\mathcal N=2C_{2,S}^2$.

\end{pf}

Now, we start to prove the well-posedness of the discrete problem \eqref{dis:rot:NS:momentum}--\eqref{dis:rot:NS:mass}. In fact, by Lemma \ref{lem:fix:pres:robust} it can be reformulated as seeking $\bm{u}_h\in\bm{V}_h$ such that

\begin{equation}\label{dis:rot:NS:momentum:reform}
  \nu a(\bm{u}_h,\bm{v}_h)+b_h(\bm{u}_h;\bm{u}_h,\bm{v}_h)
  =(\mathbb{P}(\bm{f}),\mathcal R(\bm{v}_h)),\quad\forall\,\bm{v}_h\in\bm{V}_h.
\end{equation}
Let $\mathcal F:\,\bm{V}_h\rightarrow\bm{V}_h$ be a nonlinear map such that for each $\bm{w}_h\in\bm{V}_h$, $\widetilde{\bm{u}}_h:=\mathcal F(\bm{w}_h)\in\bm{V}_h$ is given as the solution of the following linear problem:

\begin{equation}\label{dis:rot:NS:momentum:linear:1}
  \nu a(\widetilde{\bm{u}}_h,\bm{v}_h)+b_h(\bm{w}_h;\widetilde{\bm{u}}_h,\bm{v}_h)
  =(\mathbb{P}(\bm{f}),\mathcal R(\bm{v}_h)),\quad\forall\,\bm{v}_h\in\bm{V}_h.
\end{equation}
We know that the map $\mathcal F$ is continuous and compact in the finite dimensional space $\bm{V}_h$. Hence, if $\lambda>0$ and $\bm{w}_h$ satisfies $\mathcal F(\bm{w}_h)=\lambda\bm{w}_h$, then from above it follows that

\begin{equation}\label{dis:rot:NS:momentum:linear:2}
  \lambda\nu a(\bm{w}_h,\bm{v}_h)+\lambda b_h(\bm{w}_h;\bm{w}_h,\bm{v}_h)
  =(\mathbb{P}(\bm{f}),\mathcal R(\bm{v}_h)),\quad\forall\,\bm{v}_h\in\bm{V}_h.
\end{equation}
By choosing $\bm{v}_h=\bm{w}_h$ in \eqref{dis:rot:NS:momentum:linear:2} and using \eqref{skew:sym:discrete:tril}, we obtain

\[
\lambda\nu\|\nabla\bm{w}_h\|_{0,\Omega}^2=(\mathbb{P}(\bm{f}),\mathcal R(\bm{w}_h)).
\]
Based on the definition of $\|\cdot\|_{\ast,h}$ in Lemma \ref{lem:discrete:new:rh}, it holds

\[
\lambda\leqslant\frac{\|\bm{f}\|_{\ast,h}}{\nu\|\nabla\bm{w}_h\|_{0,\Omega}}.
\]
Thus, the condition $\lambda<1$ holds true for any $\bm{w}_h$ being on the boundary of the ball in $\bm{V}_h$ centered at the origin with radius $\rho>\nu^{-1}\|\bm{f}\|_{\ast,h}$. In addition, the famous Leray--Schauder fixed point theorem \cite{Cioranescu2016} implies that the nonlinear map $\mathcal F$ defined by \eqref{dis:rot:NS:momentum:linear:1} has a fixed point denoted as $\bm{u}_h$ such that $\mathcal F(\bm{u}_h)=\bm{u}_h$ in any ball centered at the origin with radius $\rho>\nu^{-1}\|\bm{f}\|_{\ast,h}$. As a result, the fixed point $\bm{u}_h\in\bm{V}_h$ is also a solution of problem \eqref{dis:rot:NS:momentum:reform}, and further satisfies an a prior estimate that $\|\nabla\bm{u}_h\|_{0,\Omega}\leqslant\nu^{-1}\|\bm{f}\|_{\ast,h}$. Then the last task is to obtain the global uniqueness of the solution pair $(\bm{u}_h,p_h)\in\bm{X}_h\times Q_h$ to the original problem \eqref{dis:rot:NS:momentum}--\eqref{dis:rot:NS:mass}. To this end, we shall use the discrete inf-sup condition \eqref{discrete:inf:sup} and the Bab\u{u}ska--Brezzi's theory \cite{Babuska1973, Brezzi1974, Girault1986, Temam1984} to conclude the following theorem.

\begin{theorem}
\label{thm:disc:wellposed}
Under the condition that $\nu^{-2}\mathcal N\|\bm{f}\|_{\ast,h}=:\sigma_h<1$, problem \eqref{dis:rot:NS:momentum}--\eqref{dis:rot:NS:mass} admits a unique solution pair $(\bm{u}_h,p_h)\in\bm{X}_h\times Q_h$ such that

\begin{equation}\label{disc:vel:bound}
  \|\nabla\bm{u}_h\|_{0,\Omega}\leqslant\nu^{-1}\|\bm{f}\|_{\ast,h}.
\end{equation}
\end{theorem}

\begin{pf}
We have proved that the equivalent problem \eqref{dis:rot:NS:momentum:reform} has at least one solution $\bm{u}_h\in\bm{V}_h$ which satisfies \eqref{disc:vel:bound}. If we assume that there are two solutions $\bm{u}_h^1$, $\bm{u}_h^2\in\bm{V}_h$ to problem \eqref{dis:rot:NS:momentum:reform}, then the difference $\bm{\delta}_h:=\bm{u}_h^1-\bm{u}_h^2$ satisfies that for any $\bm{v}_h\in\bm{V}_h$,

\begin{equation}\label{disc:error:Vh}
\nu a(\bm{\delta}_h,\bm{v}_h)
+b_h(\bm{u}_h^1;\bm{\delta}_h,\bm{v}_h)
+b_h(\bm{\delta}_h;\bm{u}_h^2,\bm{v}_h)
=0.
\end{equation}
Taking $\bm{v}_h=\bm{\delta}_h$ in \eqref{disc:error:Vh}, using the skew-symmetric property of $b_h(\cdot;\cdot,\cdot)$, and using  \eqref{disc:tri:bound} in Lemma \ref{lem:disc:tri:bound} yield

\begin{equation}\label{disc:error:Vh:2}
\left(\nu-\mathcal N\|\nabla\bm{u}_h^2\|_{0,\Omega}\right)
\|\nabla\bm{\delta}_h\|_{0,\Omega}^2
\leqslant 0.
\end{equation}
Note that $\bm{u}_h^2$ is a solution to problem \eqref{dis:rot:NS:momentum:reform}, i.e., the bound \eqref{disc:vel:bound} holds for $\bm{u}_h^2$. Hence, if we set $\sigma_h:=\nu^{-2}\mathcal N\|\bm{f}\|_{\ast,h}$ with $\sigma_h<1$, then \eqref{disc:error:Vh:2} shows that the solution of problem \eqref{dis:rot:NS:momentum:reform} is unique.

For the unique solution $\bm{u}_h\in\bm{V}_h$ of problem \eqref{dis:rot:NS:momentum:reform}, the following mapping:

\[
\bm{v}_h\in\bm{X}_h\ \mapsto (\bm{f},\mathcal R(\bm{v}_h))-\nu a(\bm{u}_h,\bm{v}_h)
-b_h(\bm{u}_h;\bm{u}_h,\bm{v}_h)
\]
defines an element $\ell$ on the dual space $\bm{X}_h^\ast$, and further $\ell$ vanishes on $\bm{V}_h$. Consequently, the discrete inf-sup condition \eqref{discrete:inf:sup} implies that there exists exactly one $p_h\in Q_h$ such that

\begin{equation}\label{unique:disc:pressure}
  \ell(\bm{v}_h)=d(\bm{v}_h,p_h),\quad
  \forall\,\bm{v}_h\in\bm{X}_h.
\end{equation}
Hence, it follows from the fact that $\bm{u}_h\in\bm{X}_h$ and \eqref{unique:disc:pressure} that the original problem \eqref{dis:rot:NS:momentum}--\eqref{dis:rot:NS:mass} admits a unique solution pair $(\bm{u}_h,p_h)\in\bm{X}_h\times Q_h$.

\end{pf}

Recalling \eqref{ineq:discrete:new:rh} from Lemma \ref{lem:discrete:new:rh}, we denote by $\mathcal M$ the right hand bound $\|\mathbb{P}(\bm{f})\|_{V^\ast}+C_{\mathcal R}L\|\mathbb{P}(\bm{f})\|_{0,\Omega}$, which is exactly independent of $h$ and any physical parameter. In fact, one can derive that $\|\bm{f}\|_{\ast,h}\leqslant\mathcal M$ for any mesh $\mathcal T_h$. Hence in this case, the a priori bound \eqref{disc:vel:bound} can be rewritten into

\begin{equation}\label{disc:other:vel:bound}
  \|\nabla\bm{u}_h\|_{0,\Omega}\leqslant \nu^{-1}\mathcal M.
\end{equation}
Simultaneously, the small data assumption $\sigma_h<1$ for the solution uniqueness needs to be replaced by

\[
\nu^{-2}\mathcal N\mathcal M=:\widetilde{\sigma}<1.
\]
Moreover, this fixed constant $\widetilde{\sigma}$ will be employed to later analysis of error estimates due to its independency of $h$.

\section{A Prior Error Estimates}
\label{sec:error:estimate}

In this section, we shall discuss some convergence results of the new discrete weak formulation \eqref{dis:rot:NS:momentum}--\eqref{dis:rot:NS:mass} proposed in section \ref{subsec:discrete:problem}. Let $(\bm{u},p)$ and $(\bm{u}_h,p_h)$ be the solutions of the continuous problem \eqref{cont:rot:NS:momentum}--\eqref{cont:rot:NS:mass} and the discrete problem \eqref{dis:rot:NS:momentum}--\eqref{dis:rot:NS:mass}, respectively.

We first propose an estimate about the trilinear form, see the lemma below.

\begin{lemma}
\label{lem:error:tril}
For any $\bm{v}_h\in\bm{X}_h$, it holds that

\begin{equation}\label{discrete:tril:error}
\begin{split}
&\left|b_h(\bm{u}_h;\bm{u}_h,\bm{v}_h)
-b(\bm{u};\bm{u},\mathcal R(\bm{v}_h))\right|
\leqslant\mathcal N\|\nabla\bm{v}_h\|_{0,\Omega}^2
\|\nabla\bm{u}_h\|_{0,\Omega}\\
&\quad+\mathcal N\|\nabla\bm{v}_h\|_{0,\Omega}
\|\nabla(\bm{u}-\bm{u}_h+\bm{v}_h)\|_{0,\Omega}
\left(\|\nabla\bm{u}\|_{0,\Omega}+\|\nabla\bm{u}_h\|_{0,\Omega}\right)\\
&\quad+(C_{2,S}C_{\mathrm{s}}C_P+\mathcal N)\|\nabla\bm{v}_h\|_{0,\Omega}
\|\nabla\bm{u}\|_{0,\Omega}\|\nabla(\bm{\mathrm{I}}_h^{k}\bm{u}
-\bm{u})\|_{0,\Omega}.
\end{split}
\end{equation}
\end{lemma}

\begin{pf}
From the definitions of $b(\cdot;\cdot,\cdot)$ and $b_h(\cdot;\cdot,\cdot)$, and Remark \ref{rem:Hdiv:test}, it follows that

\begin{align}
&b_h(\bm{u}_h;\bm{u}_h,\bm{v}_h)
-b(\bm{u};\bm{u},\mathcal R(\bm{v}_h))=b_h(\bm{u}_h-\bm{u};\bm{u}_h,\bm{v}_h)
+b(\bm{u};\mathcal R(\bm{u}_h)-\bm{u},\mathcal R(\bm{v}_h))\notag\\
\label{equiv:tri:diff}
=\,&b_h(\bm{v}_h;\bm{u}_h,\bm{v}_h)
-b(\bm{u};\bm{u}-\mathcal R(\bm{u}_h-\bm{v}_h),\mathcal R(\bm{v}_h))
-b_h(\bm{u}-\bm{u}_h+\bm{v}_h;\bm{u}_h,\bm{v}_h).
\end{align}
Based on \eqref{disc:tri:bound} in Lemma \ref{lem:disc:tri:bound}, the first and third terms of \eqref{equiv:tri:diff} are bounded by

\begin{align*}
|b_h(\bm{v}_h;\bm{u}_h,\bm{v}_h)|
&\leqslant \mathcal N\|\nabla\bm{v}_h\|_{0,\Omega}^2
\|\nabla\bm{u}_h\|_{0,\Omega},\\
|b_h(\bm{u}-\bm{u}_h+\bm{v}_h;\bm{u}_h,\bm{v}_h)|
&\leqslant \mathcal N\|\nabla\bm{v}_h\|_{0,\Omega}
\|\nabla\bm{u}_h\|_{0,\Omega}\|\nabla(\bm{u}-\bm{u}_h+\bm{v}_h)\|_{0,\Omega},
\end{align*}
respectively. By introducing $\bm{\mathrm{I}}_h^{k}\bm{u}$ from Remark \ref{rem:H1Pk:appro}, and using \eqref{disc:tri:bound}, a H\"{o}lder inequality with exponents $(\frac{1}{2},\frac{1}{4},\frac{1}{4})$, \eqref{discrete:Sobolev:ineq:2}, \eqref{equal:Pk:Rec}, the Sobolev inequality \eqref{cont:Sobolev:L4} in $L^4(\Omega)$ and the Poincar\'{e}--Friedrichs inequality, the second term of \eqref{equiv:tri:diff} can be estimated by

\begin{align*}
&|b(\bm{u};\bm{u}-\mathcal R(\bm{u}_h-\bm{v}_h),\mathcal R(\bm{v}_h))|\leqslant |b_h(\bm{u};\bm{\mathrm{I}}_h^{k}\bm{u}-\bm{u}_h+\bm{v}_h,\bm{v}_h)|
+|b(\bm{u};\bm{u}-\mathcal R(\bm{\mathrm{I}}_h^{k}\bm{u}),\mathcal R(\bm{v}_h))|\\
&\leqslant\|\nabla\bm{v}_h\|_{0,\Omega}
\|\nabla\bm{u}\|_{0,\Omega}\left\{\mathcal N\|\nabla(\bm{u}-\bm{u}_h+\bm{v}_h)\|_{0,\Omega}
+(C_{2,S}C_{\mathrm{s}}C_P+\mathcal N)\|\nabla(\bm{\mathrm{I}}_h^{k}\bm{u}
-\bm{u})\|_{0,\Omega}\right\}.
\end{align*}
Hence, the result \eqref{discrete:tril:error} then follows by adding the above results together.

\end{pf}

With the help of Lemma \ref{lem:error:tril}, we can straight obtain the a priori $H^1$-error estimate of the velocity.

\begin{theorem}\label{thm:H1error:vel}
  Under the uniqueness conditions that $\sigma<1$ and $\widetilde{\sigma}<1$, if $\bm{u}\in\bm{X}\cap\bm{H}^{k+1}(\Omega)$ and $p\in Q\cap H^1(\Omega)$, then it holds

  \begin{equation}\label{H1error:vel}
    \|\nabla(\bm{u}-\bm{u}_h)\|_{0,\Omega}\leqslant \kappa_1h^k|\bm{u}|_{k+1,\Omega},
  \end{equation}
  where the positive constant $\kappa_1$ is independent of $h$, $\bm{u}$, $\bm{u}_h$, $p$ and $p_h$.
\end{theorem}

\begin{pf}
Since $\bm{u}$ and $\bm{u}_h$ are respectively the solutions of the continuous problem \eqref{cont:rot:NS:momentum}--\eqref{cont:rot:NS:mass} and the discrete problem \eqref{dis:rot:NS:momentum}--\eqref{dis:rot:NS:mass}, it holds for any $\bm{v}_h\in\bm{V}_h$ and any $\bm{w}_h\in\bm{u}_h+\bm{V}_h$ that

\begin{align}
&\nu a(\bm{u}_h-\bm{w}_h,\bm{v}_h)=\nu a(\bm{u}-\bm{w}_h,\bm{v}_h)+\nu a(\bm{u}_h-\bm{u},\bm{v}_h)\notag\\
=\,&\nu a(\bm{u}-\bm{w}_h,\bm{v}_h)
-\left(\nu a(\bm{u},\bm{v}_h)
+b_h(\bm{u}_h;\bm{u}_h,\bm{v}_h)
-(\bm{f},\mathcal R(\bm{v}_h)\right)\notag\\
\label{equvi:disc:blinear:form}
=\,&\nu a(\bm{u}-\bm{w}_h,\bm{v}_h)
-\left[\nu a(\bm{u},\bm{v}_h)
+b(\bm{u};\bm{u},\mathcal R(\bm{v}_h))
-(\bm{f},\mathcal R(\bm{v}_h)\right]
-\left[b_h(\bm{u}_h;\bm{u}_h,\bm{v}_h)
-b(\bm{u};\bm{u},\mathcal R(\bm{v}_h))\right].
\end{align}
For the second term of \eqref{equvi:disc:blinear:form}, it follows from Lemma \ref{lem:fix:pres:robust} that

\begin{equation}
\label{equvi:disc:blinear:form:2ndterm}
\nu a(\bm{u},\bm{v}_h)
+b(\bm{u};\bm{u},\mathcal R(\bm{v}_h))
-(\bm{f},\mathcal R(\bm{v}_h))
=\nu\left[(\nabla\bm{u},\nabla\bm{v}_h)
+(\Delta\bm{u},\mathcal R(\bm{v}_h))\right].
\end{equation}
Then, plugging \eqref{equvi:disc:blinear:form:2ndterm} into \eqref{equvi:disc:blinear:form}, taking $\bm{v}_h=\bm{u}_h-\bm{w}_h$ in \eqref{equvi:disc:blinear:form}, and using
\eqref{disc:consistency:error}, \eqref{discrete:tril:error} and \eqref{H1Pk:appro:vel} yield

\begin{align*}
(\nu-\mathcal N\|\nabla\bm{u}_h\|_{0,\Omega})\|\nabla(\bm{u}_h-\bm{w}_h)\|_{0,\Omega}
&\leqslant
(\nu+\mathcal N\|\nabla\bm{u}_h\|_{0,\Omega}+\mathcal N\|\nabla\bm{u}\|_{0,\Omega})
\|\nabla(\bm{u}-\bm{w}_h)\|_{0,\Omega}\\
&\quad+\left(\nu C_E+C_I(C_{2,S}C_{\mathrm{s}}C_P+\mathcal N)\|\nabla\bm{u}\|_{0,\Omega}\right)
h^k|\bm{u}|_{k+1,\Omega}.
\end{align*}
Following the triangle inequality, \eqref{cont:ineq:stable}, \eqref{disc:other:vel:bound}, Proposition 5.1.3. in \cite{Boffi2013} or section 3 in \cite{Layton200201}, and \eqref{appro:vel}, the $H^1$-error estimate \eqref{H1error:vel} then holds, in which the positive constant $\kappa_1$ is taken as

\[
\kappa_1=\frac{2C_{\Pi}(3N+\mathcal N)+\beta_0(C_E+C_IC_{2,S}C_{\mathrm{s}}C_P+C_I\mathcal N)}{\beta_0 N(1-\widetilde{\sigma})}.
\]
\end{pf}

Next, we shall prove the $L^2$-error estimate of the pressure. Since $(\bm{u}_h,p_h)\in\bm{X}_h\times Q_h$ is the solution pair of problem \eqref{dis:rot:NS:momentum}--\eqref{dis:rot:NS:mass}, we obtain that for any $\bm{v}_h\in\bm{X}_h$,

\[
-(p_h,\mathrm{div}\,\bm{v}_h)=d(\bm{v}_h,p_h)
=(\bm{f},\mathcal R(\bm{v}_h))-\nu a(\bm{u}_h,\bm{v}_h)-b(\bm{u}_h;\bm{u}_h,\bm{v}_h).
\]
Hence, if $\bm{u}$ and $p$ are smooth and let $\bm{f}=-\nu\Delta\bm{u}+(\mathrm{rot}\,\bm{u})\times\bm{u}+\nabla p$, then it holds that for any $\bm{v}_h\in\bm{X}_h$,

\begin{equation}\label{error:eq:pres}
-(p_h,\mathrm{div}\,\bm{v}_h)=(\nabla p,\mathcal R(\bm{v}_h))
-\nu\left[(\nabla\bm{u},\nabla\bm{v}_h)
+(\Delta\bm{u},\mathcal R(\bm{v}_h))\right]
+\left[b(\bm{u};\bm{u},\mathcal R(\bm{v}_h))-b_h(\bm{u}_h;\bm{u}_h,\bm{v}_h)\right]
\end{equation}

\begin{theorem}\label{thm:L2error:pres}
Under the uniqueness conditions that $\sigma<1$ and $\widetilde{\sigma}<1$, if $\bm{u}\in\bm{X}\cap\bm{H}^{k+1}(\Omega)$ and $p\in Q\cap H^k(\Omega)$, then it holds

\begin{equation}\label{L2error:pres}
  \|p-p_h\|_{0,\Omega}\leqslant h^k\left(\kappa_2|\bm{u}|_{k+1,\Omega}+C_\pi|p|_{k,\Omega}\right),
\end{equation}
where the positive constant $\kappa_2$ is independent of $h$, $\bm{u}$, $\bm{u}_h$, $p$ and $p_h$.
\end{theorem}

\begin{pf}
Let us start with \eqref{error:eq:pres}, and then we have for any $\bm{v}_h\in\bm{X}_h$,

\begin{equation}\label{error:eq:proj:pres}
\begin{split}
(\pi_h p-p_h,\mathrm{div}\,\bm{v}_h)
&=\left[(\pi_h p,\mathrm{div}\,\bm{v}_h)+(\nabla p,\mathcal R(\bm{v}_h))\right]
-\nu\left[(\nabla\bm{u},\nabla\bm{v}_h)+(\Delta\bm{u},\mathcal R(\bm{v}_h))\right]\\
&\quad+\left[b(\bm{u};\bm{u},\mathcal R(\bm{v}_h))-b_h(\bm{u}_h;\bm{u}_h,\bm{v}_h)\right].
\end{split}
\end{equation}
For the first term of \eqref{error:eq:proj:pres}, due to integration by parts and the fact that $\mathrm{div}\,\mathcal R(\bm{v}_h)\in Q_h$, it holds

\[
(\pi_h p,\mathrm{div}\,\bm{v}_h)+(\nabla p,\mathcal R(\bm{v}_h))
=(\pi_h p,\mathrm{div}\,\bm{v}_h)-(\pi_h p,\mathrm{div}\,\mathcal R(\bm{v}_h))=0.
\]
For the second term of \eqref{error:eq:proj:pres}, it follows from \eqref{disc:consistency:error} that

\[
-\nu\left[(\nabla\bm{u},\nabla\bm{v}_h)+(\Delta\bm{u},\mathcal R(\bm{v}_h))\right]
\leqslant C_E \nu h^k|\bm{u}|_{k+1,\Omega}\|\nabla\bm{v}_h\|_{0,\Omega}.
\]
For the last term of \eqref{error:eq:proj:pres}, its estimate result will differ from \eqref{discrete:tril:error} and some changes are made on it such that

\[
b(\bm{u};\bm{u},\mathcal R(\bm{v}_h))-b_h(\bm{u}_h;\bm{u}_h,\bm{v}_h)
=b(\bm{u};\bm{u}-\mathcal R(\bm{\mathrm{I}}_h^k\bm{u}),\mathcal R(\bm{v}_h))
+b_h(\bm{u};\bm{\mathrm{I}}_h^k\bm{u}-\bm{u}_h,\bm{v}_h)
-b_h(\bm{u}_h-\bm{u};\bm{u}_h,\bm{v}_h).
\]
Then following the similar technique of proving the estimates of the second term of \eqref{equiv:tri:diff} in Lemma \ref{lem:error:tril} and using \eqref{H1error:vel}, we obtain

$$b(\bm{u};\bm{u},\mathcal R(\bm{v}_h))-b_h(\bm{u}_h;\bm{u}_h,\bm{v}_h)
\leqslant \left\{(\kappa_1\mathcal N+C_I\mathcal N+C_IC_{2,S}C_{\mathrm{s}}C_P)\|\nabla\bm{u}\|_{0,\Omega}
+\kappa_1\mathcal N\|\nabla\bm{u}_h\|_{0,\Omega}\right\}h^k|\bm{u}|_{k+1,\Omega}\|\nabla\bm{v}_h\|_{0,\Omega}.$$
Combining the above results, and using the discrete inf-sup condition \eqref{discrete:inf:sup}, \eqref{cont:ineq:stable} and \eqref{disc:other:vel:bound} yield

\begin{equation}\label{L2error:proj:pres}
\|\pi_h p-p_h\|_{0,\Omega}\leqslant \frac{1}{\beta_0}\sup_{\bm{v}_h\in\bm{X}_h}
\frac{(\pi_h p-p_h,\mathrm{div}\,\bm{v}_h)}{\|\nabla\bm{v}_h\|_{0,\Omega}}
\leqslant\kappa_2 h^k|\bm{u}|_{k+1,\Omega},
\end{equation}
where the constant $\kappa_2$ is given by

\[
\kappa_2=\frac{\nu}{\beta_0 N}\left((\kappa_1+C_I)\mathcal N+(\kappa_1+C_E)N+C_IC_{2,S}C_{\mathrm{s}}C_P\right).
\]
Hence, the result \eqref{L2error:pres} is straight derived from the triangle inequality, \eqref{proj:pressure} and \eqref{L2error:proj:pres}.

\end{pf}

Let us end this section by proving the $L^2$-error estimate of the velocity, which is a little tricky task. With the additional rotation term $(\mathrm{rot}\,\bm{u})\times\bm{u}$, the proof procedure will differ much from which shown in \cite{Linke2016289} for solving the Stokes equations. Also, due to different trilinear forms and conforming properties, it will differ slightly from the proof procedure of $L^2$-error estimate of the velocity shown in \cite{Liu2021375} for solving the stationary Navier--Stokes equations with $\mathbf{BDM}$-like reconstruction by a nonconforming VEM.

To this end, we first consider the dual problem of \eqref{eq:con:NS}, that is, for any given $\bm{g}\in\bm{Y}$, to find the unique solution pair $(\bm{\phi},\psi)\in\bm{U}\times Z:=[\bm{H}^2(\Omega)\cap\bm{X}]\times[H^1(\Omega)\cap Q]$ such that

\begin{equation}\label{dual:NS}
\begin{split}
-\nu\Delta\bm{\phi}-(\bm{u}\cdot\nabla)\bm{\phi}+(\nabla\bm{u})^\top\bm{\phi}+\nabla \psi\,&=\,\bm{g}\quad\;\mathrm{in}\ \Omega,\\
\mathrm{div}\,\bm{\phi}\,&=\,0\ \quad\mathrm{in}\ \Omega,\\
\bm{\phi}\,&=\,\bm{0}\;\quad\mathrm{on}\ \Gamma,
\end{split}
\end{equation}
where $\bm{u}\in\bm{V}$ is the unique solution of \eqref{eq:con:NS}, i.e., the uniqueness condition $\sigma<1$ holds. Then, it is derived from \eqref{cont:bound}, \eqref{cont:ineq:stable} and the Poincar\'{e}-Friedrichs inequality that

\begin{equation}\label{H1bound:phig}
  \|\nabla\bm{\phi}\|_{0,\Omega}\leqslant\frac{C_P}{\nu(1-\sigma)}\|\bm{g}\|_{0,\Omega}.
\end{equation}

According to \cite{He2009, He2008, Temam1984}, if the domain $\Omega$ is convex, the unique solution $(\bm{v},q)\in\bm{U}\times Z$ of the following Stokes equations

\[
-\Delta\bm{v}+\nabla q=\bm{g},\quad \mathrm{div}\,\bm{v}=0\quad\text{in}\ \Omega,\quad\bm{v}|_\Gamma=\bm{0}
\]
for any given $\bm{g}\in\bm{Y}$ then exists and satisfies

\begin{equation}\label{bound:auxiliary}
\|\bm{v}\|_{2,\Omega}+\|q\|_{1,\Omega}\leqslant C_{\mathrm{st}}\|\bm{g}\|_{0,\Omega},
\end{equation}
where the positive constant $C_{\mathrm{st}}$ depends only on $\Omega$. In addition, based on (2.7) in \cite{He2008}, it holds that for any $\bm{w}\in\bm{U}$,

\begin{equation}\label{Linfty:H1:H2}
  \|\bm{w}\|_{0,\infty,\Omega}\leqslant
  C_\infty\|\bm{w}\|_{2,\Omega}^{1/2}
  \|\nabla\bm{w}\|_{0,\Omega}^{1/2},
\end{equation}
where $C_\infty$ is a positive constant depending only on $\Omega$. Hence for the $H^2$-regularity estimate of $\bm{\phi}$, from \eqref{bound:auxiliary}, the H\"{o}lder inequality, the Littlewood's inequality, \eqref{cont:Sobolev:L4}, the Poincar\'{e}-Friedrichs inequality, \eqref{Linfty:H1:H2}, \eqref{cont:ineq:stable}, \eqref{H1bound:phig} and the uniqueness condition that $\sigma<1$, it follows that

\begin{align*}
\nu\|\bm{\phi}\|_{2,\Omega}+\|\psi\|_{1,\Omega}&\leqslant C_{\mathrm{st}}\|\bm{g}\|_{0,\Omega}
+C_{\mathrm{st}}\left(C_{\mathrm{s}}^{3/2}(1+C_P^2)^{1/2}+C_\infty\right)
\|\nabla\bm{u}\|_{0,\Omega}\|\nabla\bm{\phi}\|_{0,\Omega}^{1/2}\|\bm{\phi}\|_{2,\Omega}^{1/2}\\
&\leqslant \left(C_{\mathrm{st}}
+\frac{C_{\mathrm{st}}^2C_P\left(C_{\mathrm{s}}^3(1+C_P^2)+C_\infty^2\right)}{N^2(1-\sigma)}\right)
\|\bm{g}\|_{0,\Omega}+\frac{\nu}{2}\|\bm{\phi}\|_{2,\Omega},
\end{align*}
which yields that

\begin{equation}\label{H2bound:phi}
  \nu\|\bm{\phi}\|_{2,\Omega}+\|\psi\|_{1,\Omega}\leqslant \kappa_3\|\bm{g}\|_{0,\Omega},
\end{equation}
where the positive constant $\kappa_3$, independent of $\bm{u}$, $\bm{\phi}$, $p$ and $\psi$, is taken as $\kappa_3=2C_{\mathrm{st}}+\frac{2C_{\mathrm{st}}^2C_P\left(C_{\mathrm{s}}^3(1+C_P^2)+C_\infty^2\right)}{N^2(1-\sigma)}$. Furthermore, we imitate the above process to obtain the $H^2$-regularity estimate of the velocity $\bm{u}$ that if $(\bm{u},p)\in\bm{U}\times Z$, then

\begin{equation}\label{H2bound:vel}
  \nu\|\bm{u}\|_{2,\Omega}+\|p\|_{1,\Omega}\leqslant \kappa_4\|\bm{f}\|_{0,\Omega},
\end{equation}
where $\kappa_4=2C_{\mathrm{st}}+\frac{C_{\mathrm{st}}^2C_{\mathrm{s}}^2(1+C_P^2)C_P}{N^2}$.

Then, we shall propose the following theorem about the $L^2$-error estimate of the velocity, which is not shown in \cite{Linke2016, Quiroz20202655}.

\begin{theorem}\label{thm:L2error:vel}
Let $\bm{u}\in\bm{X}\cap\bm{H}^{k+1}(\Omega)$ and $p\in Q\cap H^k(\Omega)$. If the assumptions $\sigma<1$, $\widetilde{\sigma}<1$ and $\frac{2\kappa_3^{1/2}C_P^{1/2}C_\infty\mathcal M}{\nu^2(1-\sigma)^{1/2}}<1$ hold, then it holds

\begin{equation}\label{L2error:vel}
 \|\bm{u}-\bm{u}_h\|_{0,\Omega}\leqslant h^{k+1}
 \left(
 \kappa_5|\bm{u}|_{k+1,\Omega}+\kappa_6\|\bm{u}\|_{k,\Omega}\|\bm{u}\|_{k+1,\Omega}+\kappa_7\|\bm{f}\|_{k-1,\Omega}
 \right),
\end{equation}
where $\kappa_5$, $\kappa_6$ and $\kappa_7$ are three positive constants, all independent of $h$, $\bm{u}$, $p$, $\bm{u}_h$, $p_h$, $\bm{\phi}$ and $\psi$.
\end{theorem}

\begin{pf}
For any given $\bm{g}\in\bm{Y}$, by choosing appropriate interpolation $(\bm{\phi}_h,\psi_h)\in[\bm{V}_h\cap\bm{W}_h]\times Q_h$ for the solution pair $(\bm{\phi},\psi)$ and taking $\bm{v}=\bm{\phi}$ and $\bm{v}_h=\bm{\phi}_h$ in \eqref{cont:rot:NS:momentum} and \eqref{dis:rot:NS:momentum}, respectively, one can see that

\begin{align}
(\bm{g},\bm{u}-\bm{u}_h)&=\left(\bm{g},(\bm{u}-\bm{u}_h)-\mathcal R(\bm{u}-\bm{u}_h)\right)\notag\\
&\quad\ +\left\{
\left(\bm{g},\mathcal R(\bm{u}-\bm{u}_h)\right)
-\nu a(\bm{\phi},\bm{u}-\bm{u}_h)
+\left(((\bm{u}\cdot\nabla)\bm{\phi},\bm{u}-\bm{u}_h\right)
-\left((\nabla\bm{u})^\top\bm{\phi},\bm{u}-\bm{u}_h\right)
\right\}\notag\\
&\quad\ +\Big\{
\nu a(\bm{u}-\bm{u}_h,\bm{\phi})
+b(\bm{u};\bm{u}-\bm{u}_h,\bm{\phi})
+b(\bm{u}-\bm{u}_h;\bm{u},\bm{\phi})\notag\\
&\qquad\quad+\nu a(\bm{u}_h,\bm{\phi}_h)
+b_h(\bm{u}_h;\bm{u}_h,\bm{\phi}_h)
-(\bm{f},\mathcal R(\bm{\phi}_h))
-\nu a(\bm{u},\bm{\phi}_h)-b(\bm{u};\bm{u},\mathcal R(\bm{\phi}_h))
\Big\}\notag\\
&\quad\ +\Big\{
\nu a(\bm{u},\bm{\phi}_h)+b(\bm{u};\bm{u},\bm{\phi}_h)
-\nu a(\bm{u},\bm{\phi})-b(\bm{u};\bm{u},\bm{\phi})
+(\bm{f},\bm{\phi})
\Big\},\notag\\
&:=\left(\bm{g},(\bm{u}-\bm{u}_h)-\mathcal R(\bm{u}-\bm{u}_h)\right)
+\mathcal I_1+\mathcal I_2+\mathcal I_3,\notag
\end{align}
where we have used the fact $\bm{\phi}\in\bm{V}$ such that the identity $((\bm{\phi}\cdot\nabla)\bm{u},\bm{v})+((\bm{\phi}\cdot\nabla)\bm{v},\bm{u})=0$ holds for any $\bm{v}\in\bm{H}^1(\Omega)$, the fact $\bm{u}\in\bm{V}$ such that the identity $-((\bm{u}\cdot\nabla)\bm{\phi},\bm{v})=((\bm{u}\cdot\nabla)\bm{v},\bm{\phi})$ holds for any $\bm{v}\in\bm{H}^1(\Omega)$, the identity $((\nabla\bm{u})^\top\bm{\phi},\bm{v})=((\bm{v}\cdot\nabla)\bm{u},\bm{\phi})$ holds for any $\bm{v}\in\bm{H}^1(\Omega)$, and \eqref{equal:Pk:Rec}.

For $\mathcal I_1$, plugging the first equation of \eqref{dual:NS} into it and following the similar technique of proving \eqref{disc:consistency:error} in Lemma \ref{lem:consistency:error}, we obtain

\begin{equation}\label{estimate:I1}
  |\mathcal I_1|\leqslant C_{\mathcal I_1}
  h^{k+1}|\bm{u}|_{k+1,\Omega}\|\bm{g}\|_{0,\Omega},
\end{equation}
where the constant $C_{\mathcal I_1}$ is given by

\[
C_{\mathcal I_1}=\kappa_1\kappa_3\left( C_E
+C_EC_{\mathrm{s}}^2(1+C_P^2)^{1/2}\frac{\|\mathbb{P}(\bm{f})\|_{V^\ast}}{\nu}
+C_EC_\infty C_P^{1/2}\frac{\|\mathbb{P}(\bm{f})\|_{V^\ast}}{\nu^2\kappa_3^{1/2}(1-\sigma)^{1/2}}\right),
\]
and we have used \eqref{disc:consistency:error}, \eqref{BDM:property:1}, \eqref{proj:pressure}, \eqref{BDM:property:3}, the H\"{o}lder inequality, the Littlewood's inequality, the Poincar\'{e}-Friedrichs inequality, \eqref{Linfty:H1:H2}, \eqref{cont:ineq:stable}, \eqref{H2bound:phi} and \eqref{H1error:vel}.

For $\mathcal I_2$ and $\mathcal I_3$, we rearrange them into

\begin{align*}
\mathcal I_2+\mathcal I_3
&=\nu a(\bm{u}-\bm{u}_h,\bm{\phi}-\bm{\phi}_h)
+\left[
b_h(\bm{u}_h;\bm{u}_h,\bm{\phi}_h-\bm{\phi})
-b(\bm{u};\bm{u},\mathcal R(\bm{\phi}_h-\bm{\phi}))
\right]\\
&\quad\ +\big[
b(\bm{u}-\bm{u}_h;\bm{u}-\bm{u}_h,\bm{\phi})
+b(\bm{u}_h-\bm{u};\mathcal R(\bm{u}_h),\mathcal R(\bm{\phi})-\bm{\phi})
+b(\bm{u};\mathcal R(\bm{u}_h)-\bm{u}_h,\mathcal R(\bm{\phi})-\bm{\phi})\\
&\qquad\ \
+b(\bm{u};\bm{u}_h-\bm{u},\mathcal R(\bm{\phi})-\bm{\phi})
+b(\bm{u}_h;\mathcal R(\bm{u}_h)-\bm{u}_h,\bm{\phi})
\big]\\
&\quad\ +(\bm{f},\bm{\phi}-\mathcal R(\bm{\phi}))
+\left[
\nu a(\bm{u},\bm{\phi}_h-\bm{\phi})
+b(\bm{u};\bm{u},\bm{\phi}_h-\bm{\phi})
-(\bm{f},\mathcal R(\bm{\phi}_h-\bm{\phi}))
\right]\\
&:=\mathcal I_{23}^1+\mathcal I_{23}^2+\mathcal I_{23}^3+\mathcal I_{23}^4+\mathcal I_{23}^5.
\end{align*}
Hence, for $\mathcal I_{23}^1$, from \eqref{H1Pk:appro:vel} and \eqref{H2bound:phi} it follows that

\begin{equation}\label{estimate:I23:1}
  |\mathcal I_{23}^1|\leqslant C_{\mathcal I_{23}^1}h^{k+1}|\bm{u}|_{k+1,\Omega}\|\bm{g}\|_{0,\Omega}.
\end{equation}
where $C_{\mathcal I_{23}^1}=\kappa_1\kappa_3 C_I$. Like the conclusions of Lemma \ref{lem:disc:sobolev:L4}, here we need a similar estimate of $\|\mathcal R(\bm{w})-\bm{w}\|_{0,4,\Omega}$ for any $\bm{w}\in\bm{X}$ or $\bm{w}\in\bm{X}\cap\bm{H}^2(\Omega)$. In fact, by using the Sobolev inequality \eqref{cont:Sobolev:L4} in $L^4(T)$ for any $T\in\mathcal T_h$ and \eqref{BDM:property:3}, one can obtain

\begin{align}
\label{rec:cont:sobolev:L4}
  \|\mathcal R(\bm{w})-\bm{w}\|_{0,4,\Omega}&\leqslant C_{3,S}h^{m-1}|\bm{w}|_{m,\Omega},
  \quad m\in\{1,2\},\\
\label{rec:cont:sobolev:L4:2}
  \|\mathcal R(\bm{w})\|_{0,4,\Omega}&\leqslant C_{4,S}\|\nabla\bm{w}\|_{0,\Omega},\qquad\ \forall\,
  \bm{w}\in\bm{X},
\end{align}
where $C_{3,S}=C_{\mathrm{s}}C_{\mathcal R}$ and $C_{4,S}=C_{\mathrm{s}}C_{\mathcal R}+C_{\mathrm{s}}(1+C_P^2)^{1/2}$, both independent of $h$ and $\bm{w}$. For $\mathcal I_{23}^2$, we imitate the proof process in Lemma \ref{thm:L2error:pres} to obtain

\begin{align}
|\mathcal I_{23}^2|
&\leqslant |b(\bm{u};\bm{u}-\bm{\mathrm{I}}_h^k\bm{u},\mathcal R(\bm{\phi}_h-\bm{\phi}))|
+|b_h(\bm{u};\bm{\mathrm{I}}_h^k\bm{u}-\bm{u},\bm{\phi}_h-\bm{\phi})|\notag\\
&\quad\ +|b_h(\bm{u};\bm{u}-\bm{u}_h,\bm{\phi}_h-\bm{\phi})|
+|b_h(\bm{u}_h-\bm{u};\bm{u}_h,\bm{\phi}_h-\bm{\phi})|\notag\\
\label{estimate:I23:2}
&\leqslant C_{\mathcal I_{23}^2}h^{k+1}|\bm{u}|_{k+1,\Omega}\|\bm{g}\|_{0,\Omega},
\end{align}
where the constant $C_{\mathcal I_{23}^2}$ is set as

\[
C_{\mathcal I_{23}^2}=\frac{\kappa_3C_{4,S}C_I}{\nu^2}\left(
C_IC_{\mathrm{s}}(1+C_P^2)^{1/2}\|\mathbb{P}(\bm{f})\|_{V^\ast}
+2C_IC_{4,S}\|\mathbb{P}(\bm{f})\|_{V^\ast}
+2\kappa_3C_{4,S}\|\mathbb{P}(\bm{f})\|_{V^\ast}
+2\kappa_3C_{2,S}\mathcal M
\right).
\]
For $\mathcal I_{23}^3$, by using repeatedly the H\"{o}lder inequality with exponents $(\frac{1}{2},\frac{1}{4},\frac{1}{4})$ for all terms but the last one where another H\"{o}lder inequality

$$
|b(\bm{u}_h;\mathcal R(\bm{u}_h)-\bm{u}_h,\bm{\phi})|\leqslant \|\nabla\bm{u}_h\|_{0,\Omega}\|\mathcal R(\bm{u}_h)-\bm{u}_h\|_{0,\Omega}\|\bm{\phi}\|_{0,\infty,\Omega}
$$
is employed with \eqref{H1bound:phig}, \eqref{Linfty:H1:H2} and \eqref{H2bound:phi}, and having in mind the identity that $\mathcal R(\bm{u}_h)-\bm{u}_h=\mathcal R(\bm{u}_h-\bm{u}+\bm{u}-\bm{\mathrm{I}}_h^k\bm{u})+\bm{\mathrm{I}}_h^k\bm{u}-\bm{u}+\bm{u}-\bm{u}_h$ due to \eqref{equal:Pk:Rec}, we obtain

\begin{equation}\label{estimate:I23:3}
  |\mathcal I_{23}^3|\leqslant C_{\mathcal I_{23}^3}h^{k+1}|\bm{u}|_{k+1,\Omega}\|\bm{g}\|_{0,\Omega}
  +\frac{2\kappa_3^{1/2}C_P^{1/2}C_\infty\mathcal M}{\nu^2(1-\sigma)^{1/2}}\|\bm{u}-\bm{u}_h\|_{0,\Omega}\|\bm{g}\|_{0,\Omega},
\end{equation}
where the constant $C_{\mathcal I_{23}^3}$ is given by

\[
C_{\mathcal I_{23}^3}=\frac{\kappa_3}{\nu^2}\left(
\frac{\kappa_4 C_P\|\bm{f}\|_{0,\Omega}}{1-\sigma}+\kappa_1C_{2,S}C_{3,S}\mathcal M
+(C_{4,S}+C_{\mathrm{s}})(\kappa_1+C_I)\mathcal M+\kappa_1C_{\mathrm{s}}C_{3,S}\|\mathbb{P}(\bm{f})\|_{V^\ast}
+\frac{2C_P^{1/2}C_\infty C_I\mathcal M}{\kappa_3^{1/2}(1-\sigma)^{1/2}}
\right).
\]
For $\mathcal I_{23}^4$, because $\bm{u}\in\bm{X}\cap\bm{H}^{k+1}(\Omega)$, $p\in Q\cap H^k(\Omega)$ and $(\bm{u},p)$ satisfies \eqref{cont:rot:NS:momentum}, it holds $\bm{f}\in\bm{H}^{k-1}(\Omega)$. Hence, imitating the proof process in Lemma \ref{thm:L2error:pres} we have

\begin{equation}\label{estimate:I23:4}
|\mathcal I_{23}^4|\leqslant C_{\mathcal I_{23}^4}h^{k+1}\|\bm{f}\|_{k-1,\Omega}\|\bm{g}\|_{0,\Omega},
\end{equation}
where the constant $C_{\mathcal I_{23}^4}$ is given by $C_{\mathcal I_{23}^4}=\nu^{-1}\kappa_3 C_E$. The same technique can be used for $\mathcal I_{23}^5$ to obtain

\begin{align}
|\mathcal I_{23}^5|&\leqslant \nu|(\Delta\bm{u},\mathcal R(\bm{\phi}_h-\bm{\phi}))-(\nabla\bm{u},\nabla(\bm{\phi}_h-\bm{\phi}))|
+|((\mathrm{rot}\,\bm{u})\times\bm{u},(\bm{\phi}_h-\bm{\phi})-\mathcal R(\bm{\phi}_h-\bm{\phi}))|\notag\\
\label{estimate:I23:5}
&\leqslant C_{\mathcal I_{23}^{5,1}}h^{k+1}|\bm{u}|_{k+1,\Omega}\|\bm{g}\|_{0,\Omega}+C_{\mathcal I_{23}^{5,2}}h^{k+1}\|\bm{u}\|_{k,\Omega}\|\bm{u}\|_{k+1,\Omega}\|\bm{g}\|_{0,\Omega},
\end{align}
where the two constants $C_{\mathcal I_{23}^{5,1}}$ and $C_{\mathcal I_{23}^{5,2}}$ are given by $C_{\mathcal I_{23}^{5,1}}=\kappa_3 C_E$ and $C_{\mathcal I_{23}^{5,2}}=\nu^{-1}\kappa_3 C_EC_{\mathrm{s}}^2$.

Now, by combining \eqref{estimate:I1}, \eqref{estimate:I23:1}, \eqref{estimate:I23:2}--\eqref{estimate:I23:5}, the fact that

\[
\|\bm{u}-\bm{u}_h\|_{0,\Omega}=\sup_{\bm{g}\in\bm{Y}}\frac{(\bm{g},\bm{u}-\bm{u}_h)}{\|\bm{g}\|_{0,\Omega}},
\]
and the following estimate

\[
|(\bm{g},(\bm{u}-\bm{u}_h)-\mathcal R(\bm{u}-\bm{u}_h))|\leqslant C_{\mathcal I_0}h^{k+1}|\bm{u}|_{k+1,\Omega}\|\bm{g}\|_{0,\Omega}
\]
with $C_{\mathcal I_0}=\kappa_1 C_{\mathcal R}$ due to \eqref{BDM:property:3} and \eqref{H1error:vel}, we can deduce the result \eqref{L2error:vel} with

\[\kappa_5=C_{\mathcal I_0}+C_{\mathcal I_1}+C_{\mathcal I_{23}^1}+C_{\mathcal I_{23}^2}+C_{\mathcal I_{23}^3}+C_{\mathcal I_{23}^{5,1}},\quad\kappa_6=C_{\mathcal I_{23}^{5,2}},\quad\kappa_7=C_{\mathcal I_{23}^4},
\]
and an additional assumption that $\frac{2\kappa_3^{1/2}C_P^{1/2}C_\infty\mathcal M}{\nu^2(1-\sigma)^{1/2}}<1$ stemming from deduction for the term $\mathcal I_{23}^3$ in \eqref{estimate:I23:3}.

\end{pf}

\section{Numerical Experiments}
\label{sec:numer}

In this section, we present various numerical experiments to verify the convergence rates and demonstrate the performance of the proposed high-order pressure-robust method. In order to linearize the nonlinear discrete problem \eqref{dis:rot:NS:momentum}--\eqref{dis:rot:NS:mass}, the standard Newton's method is employed. Hence, we obtain an iterative algorithm given by: for $n=0,1,\cdots$, to find $(\bm{u}_h^{n+1},p_h^{n+1})\in\bm{X}_h\times Q_h$ such that for all $(\bm{v}_h,q_h)\in\bm{X}_h\times Q_h$,

\begin{equation}\label{linear:Newton}
\nu a(\bm{u}_h^{n+1},\bm{v}_h)+b_h(\bm{u}_h^n;\bm{u}_h^{n+1},\bm{v}_h)
+b_h(\bm{u}_h^{n+1};\bm{u}_h^n,\bm{v}_h)+d(\bm{v}_h,p_h^{n+1})+d(\bm{u}_h^{n+1},q_h)
=(\bm{f},\mathcal R(\bm{v}_h))+b_h(\bm{u}_h^n;\bm{u}_h^n,\bm{v}_h),
\end{equation}
with the initial data $(\bm{u}_h^0,p_h^0)$ taken to satisfy the corresponding Stokes equations:

\[
\nu a(\bm{u}_h^0,\bm{v}_h)+d(\bm{v}_h,p_h^0)+d(\bm{u}_h^0,q_h)=(\bm{f},\mathcal R(\bm{v}_h)).
\]
We proceed the iteration \eqref{linear:Newton} until the stopping criterion that

\[
\left( \|\bm{u}_h^{n+1}-\bm{u}_h^n\|_{0,\Omega}^2+\|p_h^{n+1}-p_h^n\|_{0,\Omega}^2 \right)^{1/2}<10^{-10}
\]
or $n>20$. For comparison, we also introduce the classical mixed conforming method and the corresponding iterative algorithm reads as: for $n=0,1,\cdots$, to find $(\bm{u}_h^{n+1},p_h^{n+1})\in\bm{X}_h\times Q_h$ such that for all $(\bm{v}_h,q_h)\in\bm{X}_h\times Q_h$,

\begin{equation}\label{linear:Newton:classical}
\nu a(\bm{u}_h^{n+1},\bm{v}_h)+b(\bm{u}_h^n;\bm{u}_h^{n+1},\bm{v}_h)
+b(\bm{u}_h^{n+1};\bm{u}_h^n,\bm{v}_h)+d(\bm{v}_h,p_h^{n+1})+d(\bm{u}_h^{n+1},q_h)
=(\bm{f},\bm{v}_h)+b(\bm{u}_h^n;\bm{u}_h^n,\bm{v}_h),
\end{equation}
with the initial data $(\bm{u}_h^0,p_h^0)$ taken to satisfy the corresponding Stokes equations:

\[
\nu a(\bm{u}_h^0,\bm{v}_h)+d(\bm{v}_h,p_h^0)+d(\bm{u}_h^0,q_h)=(\bm{f},\bm{v}_h).
\]
Likewise, the same stopping criterion is taken for the classical algorithm. To benefit the statements, for $k\geqslant 2$, hereafter we denote the proposed method with the iteration \eqref{linear:Newton} as $\mathbf{P}_k^{\mathrm{bubble}}$-$\mathrm{P}_{k-1}^{\mathrm{dc}}$-$\mathbf{BDM}_k$ while the classical method with the iteration \eqref{linear:Newton:classical} is denoted as $\mathbf{P}_k^{\mathrm{bubble}}$-$\mathrm{P}_{k-1}^{\mathrm{dc}}$. In particular, all of the numerical experiments are performed by using the NGSolve software \cite{Schoberl2014}.

\subsection{Example 1: Kovasznay flow}

\begin{table}[htbp]
    \centering
    \topcaption{\emph{Example 1}. Errors and convergence rates for $\mathbf{P}_k^{\mathrm{bubble}}$-$\mathrm{P}_{k-1}^{\mathrm{dc}}$-$\mathbf{BDM}_k$ with $\nu=0.1$ and $k\in\{2,3,4\}$}
    \begin{tabular*}{\hsize}{@{}@{\extracolsep{\fill}}ccccccc@{}}
        \toprule
        $h$    & $\|\bm{u}-\bm{u}_h\|_{0,\Omega}$ & Rate & $\|\nabla(\bm{u}-\bm{u}_h)\|_{0,\Omega}$ & Rate & $\frac{\|p-p_h\|_{0,\Omega}}{\|p\|_{0,\Omega}}$ & Rate\\
        \hline
        \vspace{3pt}
        $k=2$  &&&&&&\\
        $1/8$  & 3.32e$-$02    &  --     &  1.20e$+$00  & --     &  3.27e$-$02  &  --   \\
        $1/16$ & 3.70e$-$03    &  3.16   &  3.18e$-$01  & 1.91   &  8.53e$-$03  &  1.93 \\
        $1/32$ & 3.55e$-$04    &  3.38   &  8.15e$-$02  & 1.96   &  2.17e$-$03  &  1.98 \\
        $1/64$ & 3.97e$-$05    &  3.16   &  2.05e$-$02  & 1.99   &  5.45e$-$04  &  2.00 \\
        \hline
        \vspace{3pt}
        $k=3$  &&&&&&\\
        $1/8$  & 7.40e$-$04    &  --     &  6.08e$-$02  & --     &  2.84e$-$03  &  --   \\
        $1/16$ & 3.69e$-$05    &  4.33   &  7.28e$-$03  & 3.06   &  3.76e$-$04  &  2.92 \\
        $1/32$ & 2.04e$-$06    &  4.18   &  8.87e$-$04  & 3.04   &  4.77e$-$05  &  2.98 \\
        $1/64$ & 1.21e$-$07    &  4.07   &  1.09e$-$04  & 3.02   &  5.98e$-$06  &  2.99 \\
        \hline
        \vspace{3pt}
        $k=4$  &&&&&&\\
        $1/8$  & 4.06e$-$05    &  --     &  2.25e$-$03  & --     &  8.85e$-$05  &  --   \\
        $1/16$ & 1.12e$-$06    &  5.18   &  1.35e$-$04  & 4.06   &  5.87e$-$06  &  3.92 \\
        $1/32$ & 3.37e$-$08    &  5.06   &  8.20e$-$06  & 4.04   &  3.72e$-$07  &  3.98 \\
        $1/64$ & 1.03e$-$09    &  5.02   &  5.04e$-$07  & 4.02   &  2.33e$-$08  &  3.99 \\
        \bottomrule
    \end{tabular*}
    \label{tab:ex1:new}
\end{table}

\begin{table}[!htbp]
    \centering
    \topcaption{\emph{Example 1}. Errors and convergence rates for $\mathbf{P}_k^{\mathrm{bubble}}$-$\mathrm{P}_{k-1}^{\mathrm{dc}}$ with $\nu=0.1$ and $k\in\{2,3,4\}$}
    \begin{tabular*}{\hsize}{@{}@{\extracolsep{\fill}}ccccccc@{}}
        \toprule
        $h$    & $\|\bm{u}-\bm{u}_h\|_{0,\Omega}$ & Rate & $\|\nabla(\bm{u}-\bm{u}_h)\|_{0,\Omega}$ & Rate & $\frac{\|p-p_h\|_{0,\Omega}}{\|p\|_{0,\Omega}}$ & Rate\\
        \hline
        \vspace{3pt}
        $k=2$  &&&&&&\\
        $1/8$  & 3.73e$+$00    &  --     &  2.06e$+$02  & --     &  3.64e$-$02  &  --   \\
        $1/16$ & 2.92e$-$01    &  3.67   &  3.77e$+$01  & 2.45   &  9.66e$-$03  &  1.91 \\
        $1/32$ & 2.17e$-$02    &  3.75   &  5.88e$+$00  & 2.68   &  2.46e$-$03  &  1.97 \\
        $1/64$ & 1.45e$-$03    &  3.90   &  8.33e$-$01  & 2.82   &  6.19e$-$04  &  1.99 \\
        \hline
        \vspace{3pt}
        $k=3$  &&&&&&\\
        $1/8$  & 4.95e$-$01    &  --     &  5.28e$+$01  & --     &  3.17e$-$03  &  --   \\
        $1/16$ & 3.43e$-$02    &  3.85   &  7.28e$+$00  & 2.86   &  4.17e$-$04  &  2.93 \\
        $1/32$ & 2.32e$-$03    &  3.88   &  9.53e$-$01  & 2.93   &  5.26e$-$05  &  2.99 \\
        $1/64$ & 1.52e$-$04    &  3.93   &  1.22e$-$01  & 2.97   &  6.59e$-$06  &  3.00 \\
        \hline
        \vspace{3pt}
        $k=4$  &&&&&&\\
        $1/8$  & 2.84e$-$02    &  --     &  4.64e$+$00  & --     &  1.47e$-$04  &  --   \\
        $1/16$ & 9.46e$-$04    &  4.91   &  3.07e$-$01  & 3.92   &  9.71e$-$06  &  3.92 \\
        $1/32$ & 3.04e$-$05    &  4.96   &  1.95e$-$02  & 3.98   &  6.14e$-$07  &  3.98 \\
        $1/64$ & 9.65e$-$07    &  4.98   &  1.23e$-$03  & 3.99   &  3.85e$-$08  &  4.00 \\
        \bottomrule
    \end{tabular*}
    \label{tab:ex1:old}
\end{table}

The first example we use to verify the convergence rates for the proposed method was introduced by Kovasznay \cite{Kovasznay194858}, and the original pressure solution will be modified here. On the square $\Omega=(-0.5,1.5)\times (0,2)$, the right-hand body force $\bm{f}(\bm{x})$ is chosen so that the exact solution $(\bm{u},p)$ of problem \eqref{eq:con:NS} is given by

\[
u_1(\bm{x})=1-e^{\lambda x_1}\cos(2\pi x_2),\quad
u_2(\bm{x})=\frac{\lambda}{2\pi} e^{\lambda x_1} \sin(2\pi x_2),\quad
p(\bm{x})= -500 e^{2\lambda x_1}+\frac{1}{2}(u_1^2+u_2^2),
\]
where $\bm{u}(\bm{x})=(u_1(\bm{x}),u_2(\bm{x}))^\top$, $\lambda:=\frac{\mathrm{Re}}{2}-\sqrt{\frac{\mathrm{Re}^2}{4}+4\pi^2}$ with the global Reynolds number $\mathrm{Re}=\frac{1}{2\nu}$ and the vector-valued function $((u_1(\bm{x}),u_2(\bm{x}))^\top$ is imposed as Dirichlet boundary condition on $\Gamma$. We take $\nu=0.1$ and consider computations with polynomial degrees $k\in\{2,3,4\}$ performed on a sequence of uniformly $h$-refined triangular meshes with the coarsest mesh size $1/8$.

The $H^1$, $L^2$-errors of the velocity, the $L^2$-error of the pressure and their convergence rates for the proposed method are collected in Table \ref{tab:ex1:new}. It is observed that for each polynomial degree $k$, all of these errors convergence as the optimal orders as shown in Theorems \ref{thm:H1error:vel}, \ref{thm:L2error:pres} and \ref{thm:L2error:vel}. Further in order to show the advantage of the proposed method, we collect the corresponding errors for the classical mixed method in Table \ref{tab:ex1:old}. Clearly, the classical $\mathbf{P}_k^{\mathrm{bubble}}$-$\mathrm{P}_{k-1}^{\mathrm{dc}}$ are optimally convergent but their velocity errors are quite terrible compared to $\mathbf{P}_k^{\mathrm{bubble}}$-$\mathrm{P}_{k-1}^{\mathrm{dc}}$-$\mathbf{BDM}_k$ for each fixed $k$. In addition, we can also find that even the $L^2$-errors of the pressure are slightly larger. In fact, if noticing the proof of Theorem \ref{thm:L2error:pres}, one can observe that the error estimate of $\|\pi_h p-p_h\|_{0,\Omega}$ is also pressure-independent while it depends on the pressure in the classical method.

\subsection{Example 2: No flow test}

\begin{table}[!htbp]
    \centering
    \topcaption{\emph{Example 2}. Errors and convergence rates for $\mathbf{P}_k^{\mathrm{bubble}}$-$\mathrm{P}_{k-1}^{\mathrm{dc}}$ with $\nu=0.01$ and $k\in\{2,3,4\}$}
    \begin{tabular*}{\hsize}{@{}@{\extracolsep{\fill}}ccccccc@{}}
        \toprule
        $h$    & $\|\bm{u}-\bm{u}_h\|_{0,\Omega}$ & Rate & $\|\nabla(\bm{u}-\bm{u}_h)\|_{0,\Omega}$ & Rate & $\frac{\|p-p_h\|_{0,\Omega}}{\|p\|_{0,\Omega}}$ & Rate\\
        \hline
        \vspace{3pt}
        $k=2$  &&&&&&\\
        $1/8$  & 5.21e$-$05    &  --     &  3.43e$-$03  & --     &  2.49e$-$02  &  --   \\
        $1/16$ & 3.67e$-$06    &  3.83   &  4.74e$-$04  & 2.86   &  6.32e$-$03  &  1.98 \\
        $1/32$ & 2.41e$-$07    &  3.93   &  6.19e$-$05  & 2.94   &  1.59e$-$03  &  1.99 \\
        $1/64$ & 1.53e$-$08    &  3.97   &  7.89e$-$06  & 2.97   &  3.97e$-$04  &  2.00 \\
        \hline
        \vspace{3pt}
        $k=3$  &&&&&&\\
        $1/8$  & 1.01e$-$05    &  --     &  9.17e$-$04  & --     &  1.42e$-$03  &  --   \\
        $1/16$ & 6.41e$-$07    &  3.98   &  1.17e$-$04  & 2.97   &  1.77e$-$04  &  3.00 \\
        $1/32$ & 4.07e$-$08    &  3.98   &  1.49e$-$05  & 2.98   &  2.21e$-$05  &  3.00 \\
        $1/64$ & 2.57e$-$09    &  3.98   &  1.88e$-$06  & 2.99   &  2.75e$-$06  &  3.00 \\
        \hline
        \vspace{3pt}
        $k=4$  &&&&&&\\
        $1/8$  & 3.87e$-$07    &  --     &  4.53e$-$05  & --     &  2.96e$-$05  &  --   \\
        $1/16$ & 1.19e$-$09    &  5.03   &  2.87e$-$06  & 3.98   &  1.81e$-$06  &  4.03 \\
        $1/32$ & 3.67e$-$10    &  5.01   &  1.81e$-$07  & 3.98   &  1.12e$-$07  &  4.02 \\
        $1/64$ & 1.14e$-$11    &  5.01   &  1.14e$-$08  & 3.99   &  6.92e$-$09  &  4.01 \\
        \bottomrule
    \end{tabular*}
    \label{tab:ex2:old}
\end{table}

In the above example, the velocity errors by $\mathbf{P}_k^{\mathrm{bubble}}$-$\mathrm{P}_{k-1}^{\mathrm{dc}}$-$\mathbf{BDM}_k$ are satisfactory by comparison with the classical method but the pressure-independency is inapparent. Hence, in this example the velocity exact solution is taken zero so that the velocity errors should be up to the machine precision by $\mathbf{P}_k^{\mathrm{bubble}}$-$\mathrm{P}_{k-1}^{\mathrm{dc}}$-$\mathbf{BDM}_k$ for any given $k\geqslant 2$. On the square $\Omega=(0,1)\times (0,1)$, the right-hand body force $\bm{f}(\bm{x})$ is chosen so that the exact solution $(\bm{u},p)$, with the homogeneous Dirichlet boundary condition, is given by

\[
u_1(\bm{x})=0,\quad
u_2(\bm{x})=0,\quad
p(\bm{x})= 2x_1^2(1-x_1)x_2(1-x_2).
\]
Note that $\bm{f}(\bm{x})$ is an irrotational function now. We take $\nu=0.01$ and consider computations with polynomial degrees $k\in\{2,3,4\}$ performed on a sequence of uniformly $h$-refined triangular meshes with the coarsest mesh size $1/8$.

The approximation results for zero velocity by $\mathbf{P}_k^{\mathrm{bubble}}$-$\mathrm{P}_{k-1}^{\mathrm{dc}}$ and $\mathbf{P}_k^{\mathrm{bubble}}$-$\mathrm{P}_{k-1}^{\mathrm{dc}}$-$\mathbf{BDM}_k$ are collected in Tables \ref{tab:ex2:old} and \ref{tab:ex2:new}, respectively. The classical $\mathbf{P}_k^{\mathrm{bubble}}$-$\mathrm{P}_{k-1}^{\mathrm{dc}}$ are optimally convergent but their velocity errors are not pressure-independent. The performance of $\mathbf{P}_k^{\mathrm{bubble}}$-$\mathrm{P}_{k-1}^{\mathrm{dc}}$-$\mathbf{BDM}_k$ satisfies the expectation, and likewise, the $L^2$-errors of the pressure are slightly smaller than $\mathbf{P}_k^{\mathrm{bubble}}$-$\mathrm{P}_{k-1}^{\mathrm{dc}}$.

\begin{table}[!htbp]
    \centering
    \topcaption{\emph{Example 2}. Errors and convergence rates for $\mathbf{P}_k^{\mathrm{bubble}}$-$\mathrm{P}_{k-1}^{\mathrm{dc}}$-$\mathbf{BDM}_k$ with $\nu=0.01$ and $k\in\{2,3,4\}$}
    \begin{tabular*}{\hsize}{@{}@{\extracolsep{\fill}}ccccccc@{}}
        \toprule
        $h$    & $\|\bm{u}-\bm{u}_h\|_{0,\Omega}$ & Rate & $\|\nabla(\bm{u}-\bm{u}_h)\|_{0,\Omega}$ & Rate & $\frac{\|p-p_h\|_{0,\Omega}}{\|p\|_{0,\Omega}}$ & Rate\\
        \hline
        \vspace{3pt}
        $k=2$  &&&&&&\\
        $1/8$  & 2.42e$-$17    &  --   &  1.48e$-$15  & --   &  2.21e$-$02  &  --   \\
        $1/16$ & 1.57e$-$17    &  --   &  1.42e$-$15  & --   &  5.56e$-$03  &  1.99 \\
        $1/32$ & 2.08e$-$17    &  --   &  1.42e$-$15  & --   &  1.39e$-$03  &  2.00 \\
        $1/64$ & 3.05e$-$17    &  --   &  1.43e$-$15  & --   &  3.49e$-$04  &  2.00 \\
        \hline
        \vspace{3pt}
        $k=3$  &&&&&&\\
        $1/8$  & 2.03e$-$16    &  --   &  2.06e$-$14  & --   &  1.26e$-$03  &  --   \\
        $1/16$ & 1.01e$-$16    &  --   &  2.06e$-$14  & --   &  1.58e$-$04  &  2.99 \\
        $1/32$ & 5.05e$-$17    &  --   &  2.06e$-$14  & --   &  1.98e$-$05  &  3.00 \\
        $1/64$ & 5.52e$-$17    &  --   &  2.06e$-$14  & --   &  2.48e$-$06  &  3.00 \\
        \hline
        \vspace{3pt}
        $k=4$  &&&&&&\\
        $1/8$  & 5.01e$-$17    &  --   &  1.03e$-$14  & --   &  1.33e$-$05  &  --   \\
        $1/16$ & 4.95e$-$17    &  --   &  1.03e$-$14  & --   &  8.27e$-$07  &  4.00 \\
        $1/32$ & 2.71e$-$17    &  --   &  1.03e$-$14  & --   &  5.16e$-$08  &  4.00 \\
        $1/64$ & 1.57e$-$17    &  --   &  1.03e$-$14  & --   &  3.23e$-$09  &  4.00 \\
        \bottomrule
    \end{tabular*}
    \label{tab:ex2:new}
\end{table}

\begin{table}[!htbp]
    \centering
    \topcaption{\emph{Example 2}. Errors and convergence rates for $\mathbf{P}_2^{\mathrm{bubble}}$-$\mathrm{P}_1^{\mathrm{dc}}$ and $\mathbf{P}_2^{\mathrm{bubble}}$-$\mathrm{P}_1^{\mathrm{dc}}$-$\mathbf{BDM}_2$ with $\nu=0.01$ and a non-polynomial function for the body force}
    \begin{tabular*}{\hsize}{@{}@{\extracolsep{\fill}}ccccccc@{}}
        \toprule
        $h$    & $\|\bm{u}-\bm{u}_h\|_{0,\Omega}$ & Rate & $\|\nabla(\bm{u}-\bm{u}_h)\|_{0,\Omega}$ & Rate & $\frac{\|p-p_h\|_{0,\Omega}}{\|p\|_{0,\Omega}}$ & Rate\\
        \hline
        \vspace{3pt}
        $\mathbf{P}_2^{\mathrm{bubble}}$-$\mathrm{P}_1^{\mathrm{dc}}$  &&&&&&\\
        $1/32$  & 1.41e$-$05    &  --     &  4.50e$-$03  & --     &  1.08e$-$03  &  --   \\
        $1/64$  & 8.90e$-$07    &  3.99   &  5.73e$-$04  & 2.97   &  2.71e$-$04  &  2.00 \\
        $1/128$ & 5.58e$-$08    &  4.00   &  7.21e$-$05  & 2.99   &  6.78e$-$05  &  2.00 \\
        $1/256$ & 3.49e$-$09    &  4.00   &  9.05e$-$06  & 2.99   &  1.70e$-$05  &  2.00 \\
        \hline
        \vspace{3pt}
        $\mathbf{P}_2^{\mathrm{bubble}}$-$\mathrm{P}_1^{\mathrm{dc}}$-$\mathbf{BDM}_2$  &&&&&&\\
        $1/32$  & 2.76e$-$12    &  --     &  5.83e$-$10  & --     &  9.52e$-$04  &  --   \\
        $1/64$  & 1.09e$-$14    &  7.98   &  4.67e$-$12  & 6.96   &  2.38e$-$04  &  2.00 \\
        $1/128$ & 5.77e$-$16    &  4.24   &  1.00e$-$13  & 5.54   &  5.95e$-$05  &  2.00 \\
        $1/256$ & 5.88e$-$16    &  --     &  9.29e$-$14  & --     &  1.49e$-$05  &  2.00 \\
        \bottomrule
    \end{tabular*}
    \label{tab:ex2:nonpoly}
\end{table}

Furthermore, we shall check the performance of a more complicated irrotational function for the right-hand body force. To this end, the body force $\bm{f}(\bm{x})$ is chosen so that the exact solution of the pressure is given by a non-polynomial function that

\[
p(\bm{x})=\exp(x_1+x_2)+\sin(2\pi x_1)\cos(2\pi x_2).
\]
Here we consider computations with polynomial degree $k=2$ performed on a sequence of uniformly $h$-refined triangular meshes with the coarsest mesh size adjusted to 1/32. The corresponding approximation results by $\mathbf{P}_2^{\mathrm{bubble}}$-$\mathrm{P}_1^{\mathrm{dc}}$ and $\mathbf{P}_2^{\mathrm{bubble}}$-$\mathrm{P}_1^{\mathrm{dc}}$-$\mathbf{BDM}_2$ are collected in Table \ref{tab:ex2:nonpoly}. In this case, while the classical method performs as usual, the velocity errors for the proposed method are far from the machine precision unlike the results shown in Table \ref{tab:ex2:new} under the same mesh size but fortunately they are close to the machine precision in the finest mesh. This phenomenon is caused by the conspicuous form change of the right-hand body force indeed. A more complicated form of $\bm{f}$ may lead to larger values of the high-order semi-norms $|\bm{f}|_{m,\Omega}$ for $m\geqslant 1$, and the improvement is almost impossible to achieve because $\bm{f}$ is preset.

\subsection{Example 3: Robustness for irrotational body forces}

This example, arising from Benchmark 3.3 in \cite{Linke2016654}, is computed to show the robustness of the proposed method for large irrotational body forces. On the square $\Omega=(0,1)\times (0,1)$, the exact solution $(\bm{u},p)$ is given by

\[
u_1(\bm{x})=-x_2,\quad
u_2(\bm{x})=x_1,\quad
p(\bm{x})= \lambda\,x_1^6+x_1^2+x_2^2,\quad\lambda\in\left\{10^2,10^6\right\},
\]
where the vector-valued function $(-x_2,x_1)^\top$ is imposed as Dirichlet boundary condition on $\Gamma$. In this case, we can figure out the right-hand body force $\bm{f}(\bm{x})=(6\lambda x_1^5,0)^\top$ which is exactly irrotational. We take $\nu=1$ and consider computations with polynomial degrees $k\in\{2,3,4\}$ performed on a sequence of uniformly $h$-refined triangular meshes with the coarsest mesh size $1/8$.

The approximation results by $\mathbf{P}_k^{\mathrm{bubble}}$-$\mathrm{P}_{k-1}^{\mathrm{dc}}$ and $\mathbf{P}_k^{\mathrm{bubble}}$-$\mathrm{P}_{k-1}^{\mathrm{dc}}$-$\mathbf{BDM}_k$ with $\lambda=10^2$ are collected in Tables \ref{tab:ex3:old:1e2} and \ref{tab:ex3:new:1e2}, respectively. While the approximation results by $\mathbf{P}_k^{\mathrm{bubble}}$-$\mathrm{P}_{k-1}^{\mathrm{dc}}$ and $\mathbf{P}_k^{\mathrm{bubble}}$-$\mathrm{P}_{k-1}^{\mathrm{dc}}$-$\mathbf{BDM}_k$ with $\lambda=10^6$ are collected in Tables \ref{tab:ex3:old:1e6} and \ref{tab:ex3:new:1e6}, respectively. Unsurprisingly, we witness the significant differences between the two methods again. Whenever $\lambda=10^2$ or $\lambda=10^6$, the velocity field is exactly reproduced by $\mathbf{P}_k^{\mathrm{bubble}}$-$\mathrm{P}_{k-1}^{\mathrm{dc}}$-$\mathbf{BDM}_k$ for each fixed polynomial degree $k$, whereas the classical $\mathbf{P}_k^{\mathrm{bubble}}$-$\mathrm{P}_{k-1}^{\mathrm{dc}}$ are out of reach; especially as $\lambda=10^6$, the velocity $H^1$-errors by $\mathbf{P}_k^{\mathrm{bubble}}$-$\mathrm{P}_{k-1}^{\mathrm{dc}}$ are disastrous despite the optimal convergence rates. Besides, it is also observed that the velocity approximation by $\mathbf{P}_k^{\mathrm{bubble}}$-$\mathrm{P}_{k-1}^{\mathrm{dc}}$-$\mathbf{BDM}_k$ is independent of the value of $\lambda$, which implies the robustness of the velocity errors for large irrotational body forces.

\begin{table}[htbp]
    \centering
    \topcaption{\emph{Example 3}. Errors and convergence rates for $\mathbf{P}_k^{\mathrm{bubble}}$-$\mathrm{P}_{k-1}^{\mathrm{dc}}$ with $\nu=1$, $\lambda=10^2$ and $k\in\{2,3,4\}$}
    \begin{tabular*}{\hsize}{@{}@{\extracolsep{\fill}}ccccccc@{}}
        \toprule
        $h$    & $\|\bm{u}-\bm{u}_h\|_{0,\Omega}$ & Rate & $\|\nabla(\bm{u}-\bm{u}_h)\|_{0,\Omega}$ & Rate & $\frac{\|p-p_h\|_{0,\Omega}}{\|p\|_{0,\Omega}}$ & Rate\\
        \hline
        \vspace{3pt}
        $k=2$  &&&&&&\\
        $1/8$  & 2.87e$-$04    &  --     &  1.81e$-$02  & --     &  2.08e$-$02  &  --   \\
        $1/16$ & 2.15e$-$05    &  3.74   &  2.86e$-$03  & 2.67   &  5.30e$-$03  &  1.97 \\
        $1/32$ & 1.47e$-$06    &  3.87   &  4.04e$-$04  & 2.82   &  1.33e$-$03  &  1.99 \\
        $1/64$ & 9.60e$-$08    &  3.94   &  5.39e$-$05  & 2.91   &  3.34e$-$04  &  2.00 \\
        \hline
        \vspace{3pt}
        $k=3$  &&&&&&\\
        $1/8$  & 3.30e$-$05    &  --     &  3.47e$-$03  & --     &  9.65e$-$04  &  --   \\
        $1/16$ & 2.25e$-$06    &  3.88   &  4.58e$-$04  & 2.92   &  1.21e$-$04  &  3.00 \\
        $1/32$ & 1.47e$-$07    &  3.93   &  5.87e$-$05  & 2.96   &  1.51e$-$05  &  3.00 \\
        $1/64$ & 9.42e$-$09    &  3.97   &  7.44e$-$06  & 2.98   &  1.88e$-$06  &  3.00 \\
        \hline
        \vspace{3pt}
        $k=4$  &&&&&&\\
        $1/8$  & 7.87e$-$07    &  --     &  1.22e$-$04  & --     &  1.90e$-$05  &  --   \\
        $1/16$ & 2.49e$-$08    &  4.98   &  7.79e$-$06  & 3.97   &  1.18e$-$06  &  4.01 \\
        $1/32$ & 7.81e$-$10    &  4.99   &  4.91e$-$07  & 3.99   &  7.35e$-$08  &  4.01 \\
        $1/64$ & 2.45e$-$11    &  5.00   &  3.08e$-$08  & 3.99   &  4.59e$-$09  &  4.00 \\
        \bottomrule
    \end{tabular*}
    \label{tab:ex3:old:1e2}
\end{table}

\begin{table}[htbp]
    \centering
    \topcaption{\emph{Example 3}. Errors and convergence rates for $\mathbf{P}_k^{\mathrm{bubble}}$-$\mathrm{P}_{k-1}^{\mathrm{dc}}$-$\mathbf{BDM}_k$ with $\nu=1$, $\lambda=10^2$ and $k\in\{2,3,4\}$}
    \begin{tabular*}{\hsize}{@{}@{\extracolsep{\fill}}ccccccc@{}}
        \toprule
        $h$    & $\|\bm{u}-\bm{u}_h\|_{0,\Omega}$ & Rate & $\|\nabla(\bm{u}-\bm{u}_h)\|_{0,\Omega}$ & Rate & $\frac{\|p-p_h\|_{0,\Omega}}{\|p\|_{0,\Omega}}$ & Rate\\
        \hline
        \vspace{3pt}
        $k=2$  &&&&&&\\
        $1/8$  & 2.70e$-$16    &  --   &  1.40e$-$14  & --   &  1.86e$-$02  &  --   \\
        $1/16$ & 4.22e$-$16    &  --   &  1.50e$-$14  & --   &  4.70e$-$03  &  1.98 \\
        $1/32$ & 9.48e$-$16    &  --   &  1.95e$-$14  & --   &  1.18e$-$03  &  2.00 \\
        $1/64$ & 3.62e$-$15    &  --   &  3.71e$-$14  & --   &  2.95e$-$04  &  2.00 \\
        \hline
        \vspace{3pt}
        $k=3$  &&&&&&\\
        $1/8$  & 1.26e$-$15    &  --   &  1.20e$-$13  & --   &  8.72e$-$04  &  --   \\
        $1/16$ & 7.52e$-$16    &  --   &  1.41e$-$13  & --   &  1.10e$-$04  &  2.99 \\
        $1/32$ & 1.44e$-$15    &  --   &  1.53e$-$13  & --   &  1.38e$-$05  &  3.00 \\
        $1/64$ & 4.24e$-$15    &  --   &  1.63e$-$13  & --   &  1.72e$-$06  &  3.00 \\
        \hline
        \vspace{3pt}
        $k=4$  &&&&&&\\
        $1/8$  & 7.45e$-$16    &  --   &  7.30e$-$14  & --   &  1.14e$-$05  &  --   \\
        $1/16$ & 5.29e$-$16    &  --   &  8.16e$-$14  & --   &  7.16e$-$07  &  3.99 \\
        $1/32$ & 1.29e$-$15    &  --   &  8.61e$-$14  & --   &  4.48e$-$08  &  4.00 \\
        $1/64$ & 4.19e$-$15    &  --   &  9.85e$-$14  & --   &  2.80e$-$09  &  4.00 \\
        \bottomrule
    \end{tabular*}
    \label{tab:ex3:new:1e2}
\end{table}

\begin{table}[htbp]
    \centering
    \topcaption{\emph{Example 3}. Errors and convergence rates for $\mathbf{P}_k^{\mathrm{bubble}}$-$\mathrm{P}_{k-1}^{\mathrm{dc}}$ with $\nu=1$, $\lambda=10^6$ and $k\in\{2,3,4\}$}
    \begin{tabular*}{\hsize}{@{}@{\extracolsep{\fill}}ccccccc@{}}
        \toprule
        $h$    & $\|\bm{u}-\bm{u}_h\|_{0,\Omega}$ & Rate & $\|\nabla(\bm{u}-\bm{u}_h)\|_{0,\Omega}$ & Rate & $\frac{\|p-p_h\|_{0,\Omega}}{\|p\|_{0,\Omega}}$ & Rate\\
        \hline
        \vspace{3pt}
        $k=2$  &&&&&&\\
        $1/8$  & 2.87e$+$00    &  --     &  1.81e$+$02  & --     &  2.10e$-$02  &  --   \\
        $1/16$ & 2.15e$-$01    &  3.74   &  2.86e$+$01  & 2.67   &  5.35e$-$03  &  1.97 \\
        $1/32$ & 1.47e$-$02    &  3.87   &  4.04e$+$00  & 2.82   &  1.35e$-$03  &  1.99 \\
        $1/64$ & 9.60e$-$04    &  3.94   &  5.39e$-$01  & 2.91   &  3.37e$-$04  &  2.00 \\
        \hline
        \vspace{3pt}
        $k=3$  &&&&&&\\
        $1/8$  & 3.30e$-$01    &  --     &  3.47e$+$01  & --     &  9.78e$-$04  &  --   \\
        $1/16$ & 2.25e$-$02    &  3.88   &  4.58e$+$00  & 2.92   &  1.22e$-$04  &  3.00 \\
        $1/32$ & 1.47e$-$03    &  3.93   &  5.87e$-$01  & 2.96   &  1.53e$-$05  &  3.00 \\
        $1/64$ & 9.42e$-$05    &  3.97   &  7.44e$-$02  & 2.98   &  1.91e$-$06  &  3.00 \\
        \hline
        \vspace{3pt}
        $k=4$  &&&&&&\\
        $1/8$  & 7.87e$-$03    &  --     &  1.22e$+$00  & --     &  1.92e$-$05  &  --   \\
        $1/16$ & 2.49e$-$04    &  4.98   &  7.79e$-$02  & 3.97   &  1.19e$-$06  &  4.01 \\
        $1/32$ & 7.81e$-$06    &  4.99   &  4.91e$-$03  & 3.99   &  7.44e$-$08  &  4.01 \\
        $1/64$ & 2.45e$-$07    &  5.00   &  3.08e$-$04  & 3.99   &  4.64e$-$09  &  4.00 \\
        \bottomrule
    \end{tabular*}
    \label{tab:ex3:old:1e6}
\end{table}

\begin{table}[htbp]
    \centering
    \topcaption{\emph{Example 3}. Errors and convergence rates for $\mathbf{P}_k^{\mathrm{bubble}}$-$\mathrm{P}_{k-1}^{\mathrm{dc}}$-$\mathbf{BDM}_k$ with $\nu=1$, $\lambda=10^6$ and $k\in\{2,3,4\}$}
    \begin{tabular*}{\hsize}{@{}@{\extracolsep{\fill}}ccccccc@{}}
        \toprule
        $h$    & $\|\bm{u}-\bm{u}_h\|_{0,\Omega}$ & Rate & $\|\nabla(\bm{u}-\bm{u}_h)\|_{0,\Omega}$ & Rate & $\frac{\|p-p_h\|_{0,\Omega}}{\|p\|_{0,\Omega}}$ & Rate\\
        \hline
        \vspace{3pt}
        $k=2$  &&&&&&\\
        $1/8$  & 2.24e$-$12    &  --   &  1.39e$-$10  & --   &  1.87e$-$02  &  --   \\
        $1/16$ & 7.58e$-$13    &  --   &  1.19e$-$10  & --   &  4.72e$-$03  &  1.98 \\
        $1/32$ & 5.62e$-$13    &  --   &  1.10e$-$10  & --   &  1.18e$-$03  &  2.00 \\
        $1/64$ & 6.20e$-$13    &  --   &  1.10e$-$10  & --   &  2.96e$-$04  &  2.00 \\
        \hline
        \vspace{3pt}
        $k=3$  &&&&&&\\
        $1/8$  & 1.18e$-$11    &  --   &  1.16e$-$09  & --   &  8.77e$-$04  &  --   \\
        $1/16$ & 6.58e$-$12    &  --   &  1.36e$-$09  & --   &  1.11e$-$04  &  2.99 \\
        $1/32$ & 3.62e$-$12    &  --   &  1.48e$-$09  & --   &  1.39e$-$05  &  3.00 \\
        $1/64$ & 3.75e$-$12    &  --   &  1.54e$-$09  & --   &  1.73e$-$06  &  3.00 \\
        \hline
        \vspace{3pt}
        $k=4$  &&&&&&\\
        $1/8$  & 6.53e$-$12    &  --   &  6.98e$-$10  & --   &  1.15e$-$05  &  --   \\
        $1/16$ & 3.57e$-$12    &  --   &  7.52e$-$10  & --   &  7.20e$-$07  &  3.99 \\
        $1/32$ & 2.23e$-$12    &  --   &  7.79e$-$10  & --   &  4.50e$-$08  &  4.00 \\
        $1/64$ & 1.73e$-$12    &  --   &  7.64e$-$10  & --   &  2.82e$-$09  &  4.00 \\
        \bottomrule
    \end{tabular*}
    \label{tab:ex3:new:1e6}
\end{table}

\subsection{Example 4: Lid-driven cavity flow at high Reynolds numbers}

In this example, we consider the famous lid-driven cavity flow problem at high Reynolds numbers. On the square $\Omega=(0,1)\times (0,1)$, we choose the right-hand body force $\bm{f}(\bm{x})=\bm{0}$, and set the nonhomogeneous Dirichlet boundary condition as

\begin{equation*}
\bm{u}(\bm{x})=
  \begin{cases}
  (1,0)^\top, & \text{on}\ (0,1)\times\{x_2=1\},\\
  (0,0)^\top, & \text{otherwise}.
  \end{cases}
\end{equation*}

\begin{figure}[htbp]
\centering
\subfigure{
\includegraphics[width=5cm]{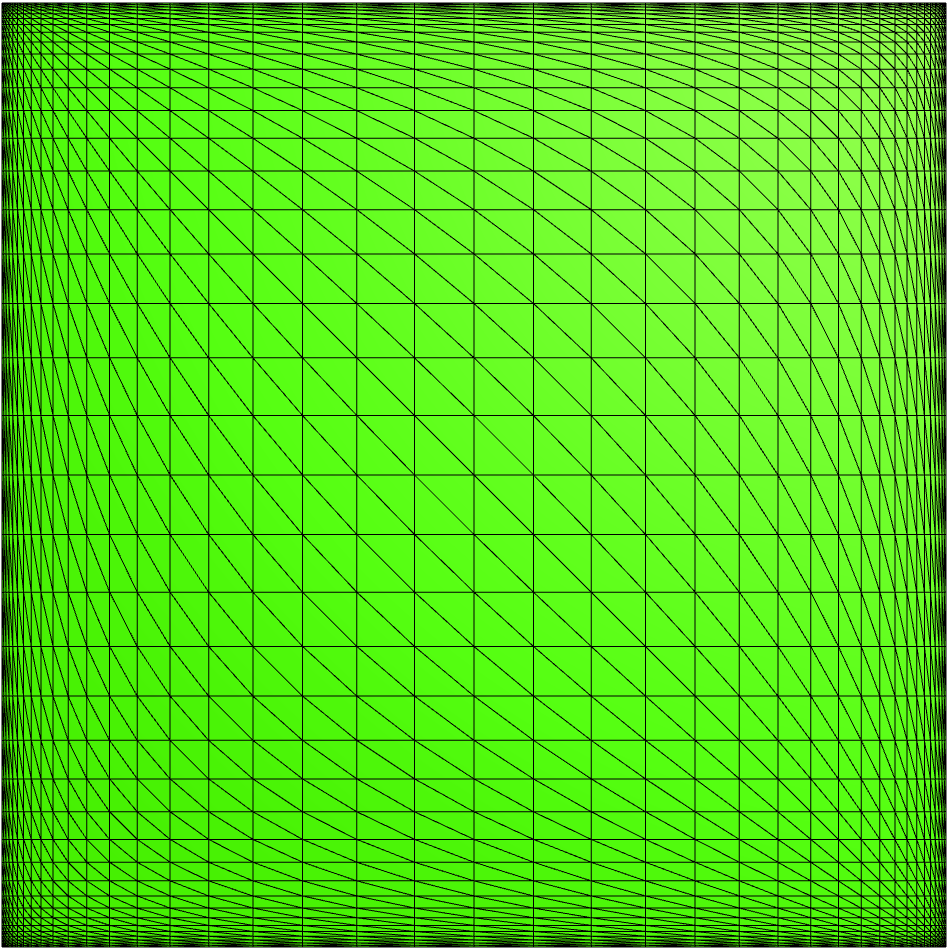}
}
\qquad\qquad
\subfigure{
\includegraphics[width=5cm]{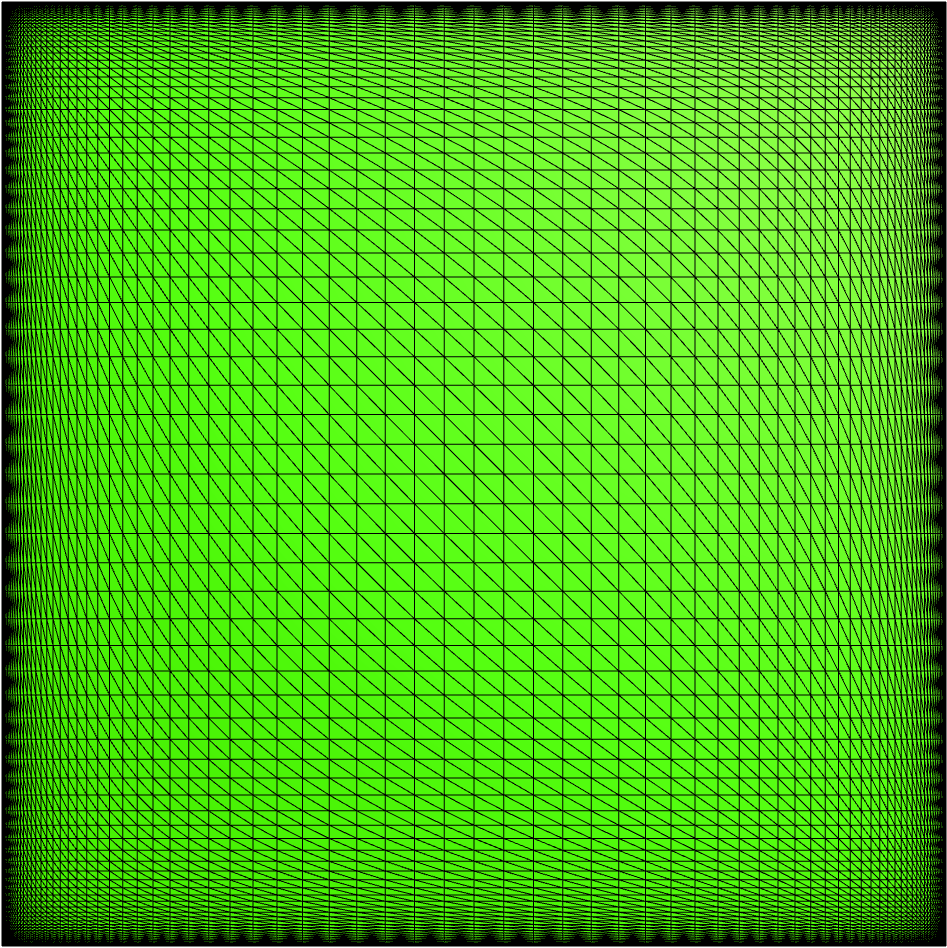}
}
\caption{\emph{Example 4.} The coarse mesh plot with $M=40$, $\gamma=2.5$ (left) and the refined mesh plot with $M=80$, $\gamma=2.5$ (right)}\label{fig:ex4:mesh40:80}
\end{figure}

We shall resolve sub-grid scale vortices with considerable computations. To this end, we first denote by $(\xi_i,\eta_j)$ the Cartesian coordinates under the uniform grid over $\Omega$ with the same number of cells $M$ in the $x$ and $y$ directions, i.e., $\xi_i=(i-1)/M$, $\eta_j=(j-1)/M$ for $i,j=1,\cdots,M+1$. Then the mesh grid point distance in the $x$ and $y$ directions is set to be a hyperbolic tangent profile. Hence, as shown below, a stretching function with an adjustable parameter $\gamma$ is implemented to get the final computational coordinates $(x_i,y_j)$.

\[
x_i=0.5+\frac{\tanh(2\gamma(\xi_i-0.5))}{2\tanh(\gamma)},\quad
y_j=0.5+\frac{\tanh(2\gamma(\eta_j-0.5))}{2\tanh(\gamma)}.
\]
Here two meshes are employed, which are a coarse mesh with $M=40$, $\gamma=2.5$ and a refined mesh with $M=80$, $\gamma=2.5$, respectively, see Figure \ref{fig:ex4:mesh40:80}.

Before achieving the following presented results, we tried to use the Newton iteration \eqref{linear:Newton} or \eqref{linear:Newton:classical} with the Stokes initial data but both did not converge at higher Reynolds numbers $\frac{1}{\nu}=:\mathrm{Re}\geqslant 1000$. Hence, to solve this $\mathrm{Re}$-restraint problem for the standard Newton iteration method and avoid all factitious terms in the numerical schemes, we decide to follow the so-called continuation method, which was successfully applied in \cite{Cadou2001825, Farrell2021075, Farrell2019A3073, Layton200201}. The specific approach is that the problem is first solved for $\mathrm{Re}=100$, then $\mathrm{Re}=400$, $1000$, $1800$, $2500$, $3200$, $5000$, and then in steps of $2500$ until $\mathrm{Re}=15000$ under the coarse mesh and until $\mathrm{Re}=20000$ under the refined mesh, with the solution for the previous value of $\mathrm{Re}$ used as initial guess for the next; the Stokes equations are solved to provide the initial guess used at $\mathrm{Re}=100$. A necessary remark is that this demanding but effective method is also used in the rest of our numerical experiments.

\begin{figure}[htbp]
\centering
\subfigure{
\includegraphics[width=5.5cm]{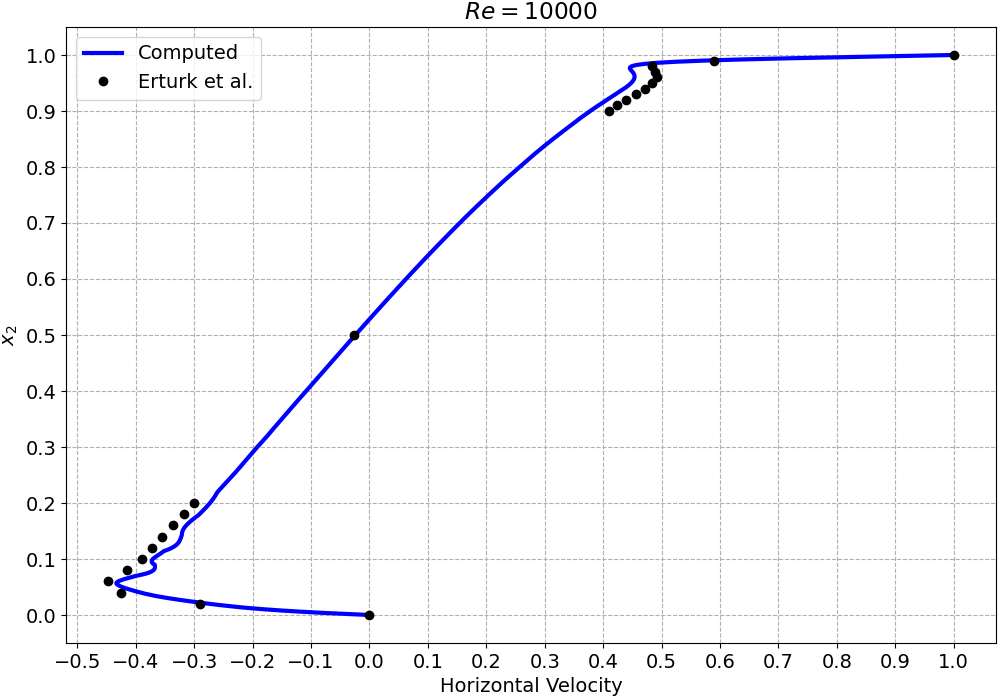}
}
\quad
\subfigure{
\includegraphics[width=5.5cm]{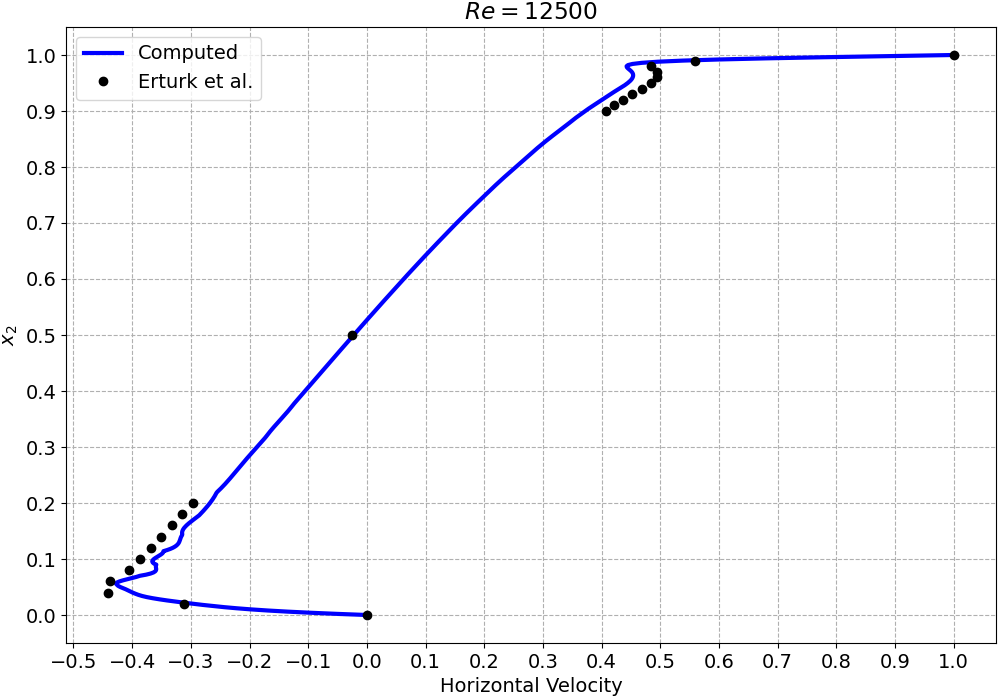}
}
\quad
\subfigure{
\includegraphics[width=5.5cm]{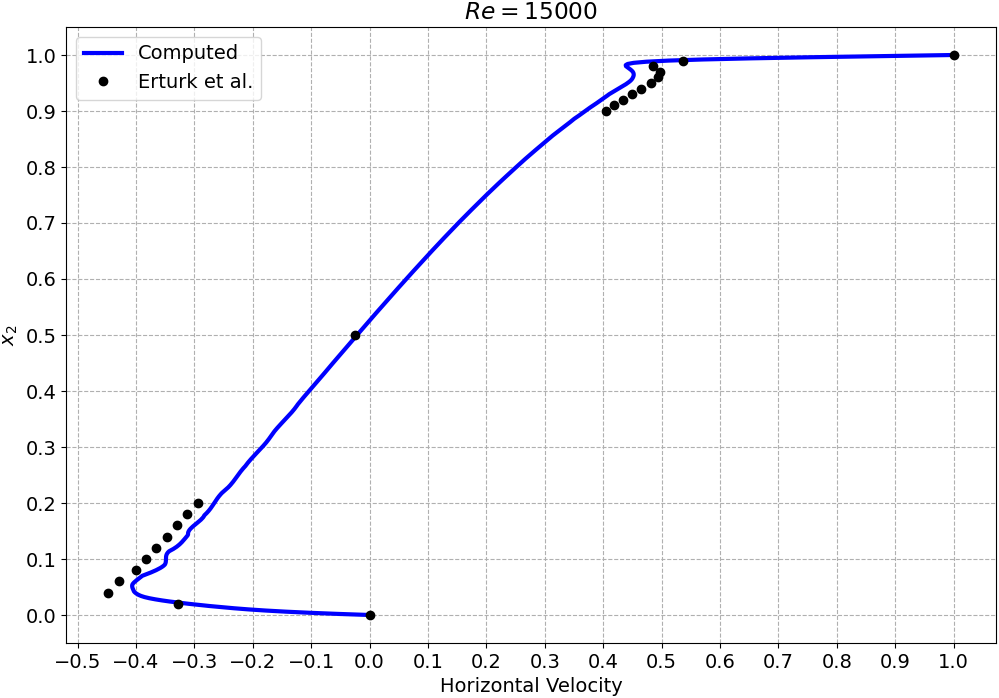}
}

\subfigure{
\includegraphics[width=5.5cm]{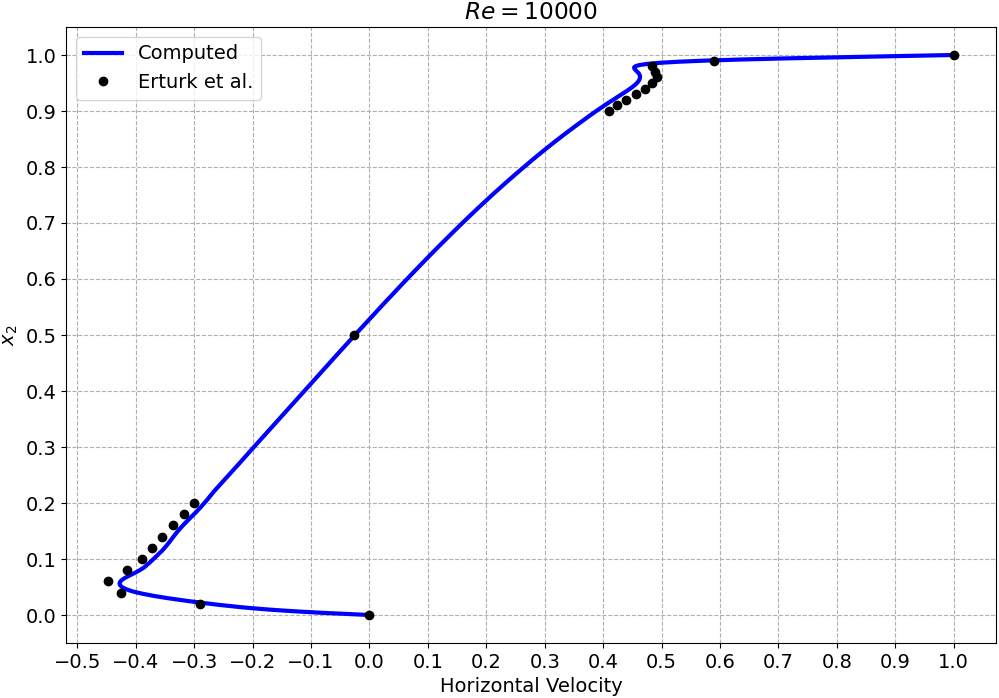}
}
\quad
\subfigure{
\includegraphics[width=5.5cm]{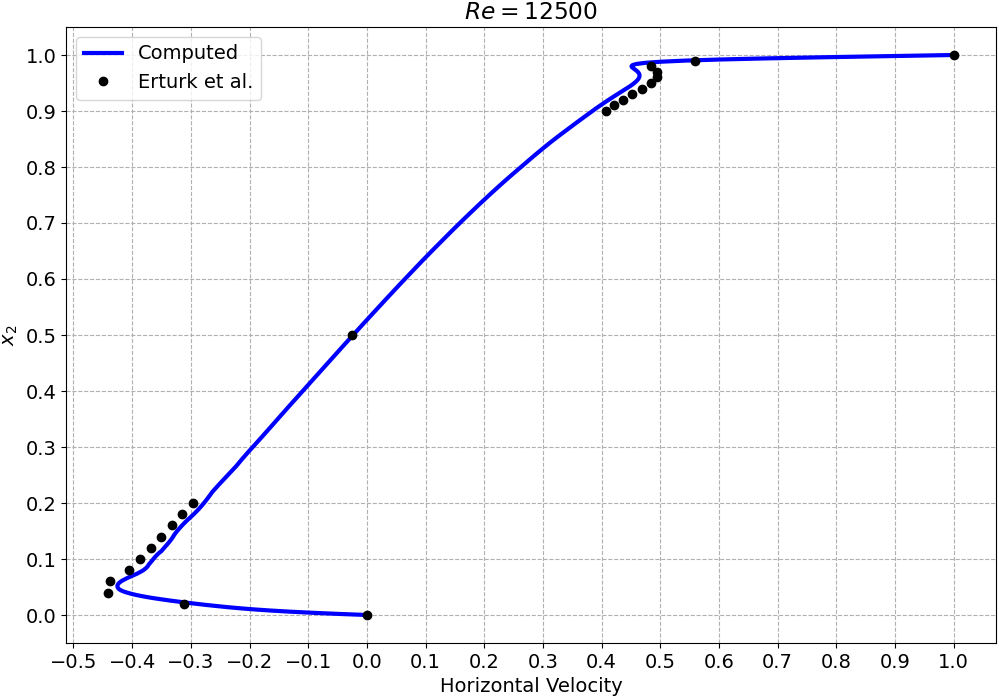}
}
\quad
\subfigure{
\includegraphics[width=5.5cm]{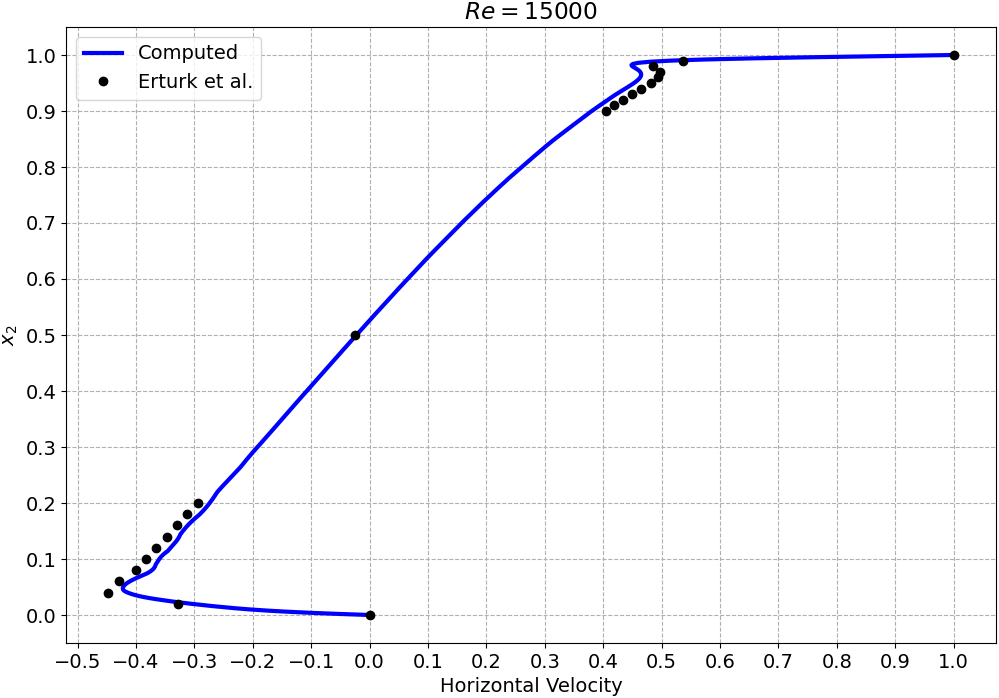}
}
\centering
\caption{\emph{Example 4.} Under the coarse mesh, computed horizontal component $u_1$ of the velocity along the vertical centreline $x_1=0.5$ compared with the reference results by a high-order accurate finite difference scheme under a fine mesh $601\times 601$ from Erturk et al. \cite{Erturk2005747}, for $\mathrm{Re}=10000$, $12500$ and $15000$ from left to right: $\bm{\mathrm{P}}_3^{\mathrm{bubble}}$-$\mathrm{P}_2^{\mathrm{dc}}$ (upper row) and $\bm{\mathrm{P}}_3^{\mathrm{bubble}}$-$\mathrm{P}_2^{\mathrm{dc}}$-$\bm{\mathrm{BDM}}_3$ (lower row)}\label{fig:ex4:mesh40:u1:0.5}
\end{figure}

\begin{figure}[htbp]
\centering
\subfigure{
\includegraphics[width=5.5cm]{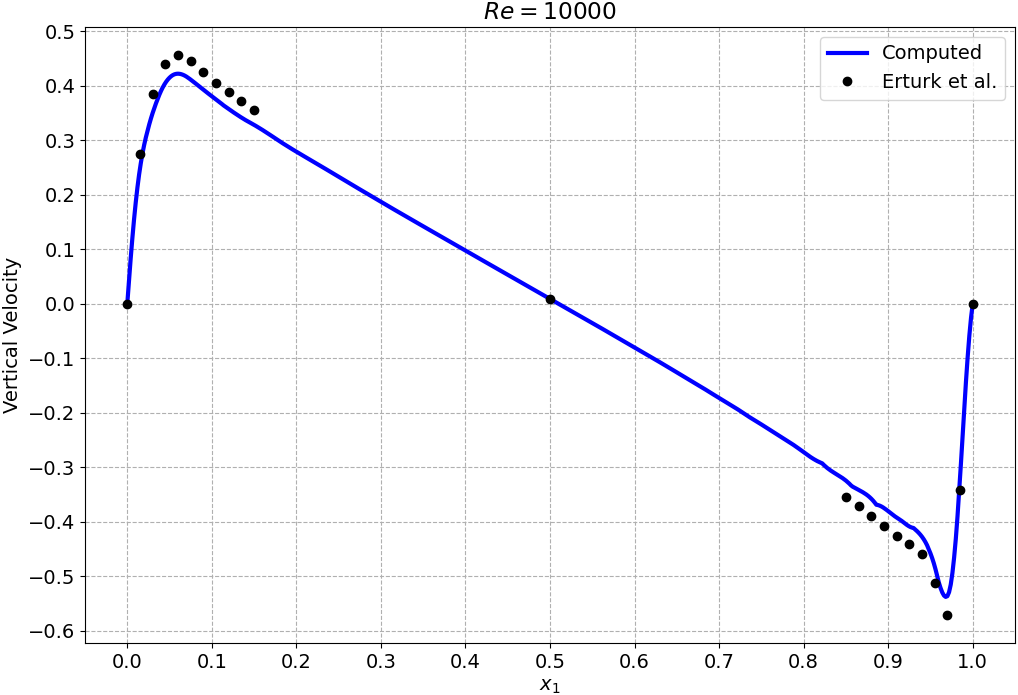}
}
\quad
\subfigure{
\includegraphics[width=5.5cm]{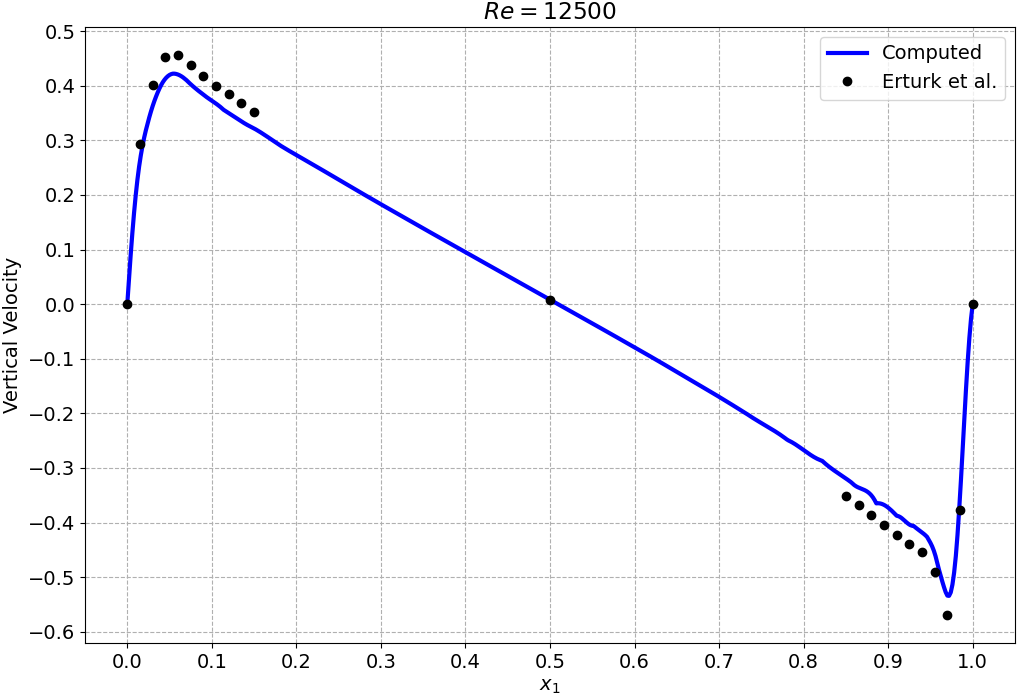}
}
\quad
\subfigure{
\includegraphics[width=5.5cm]{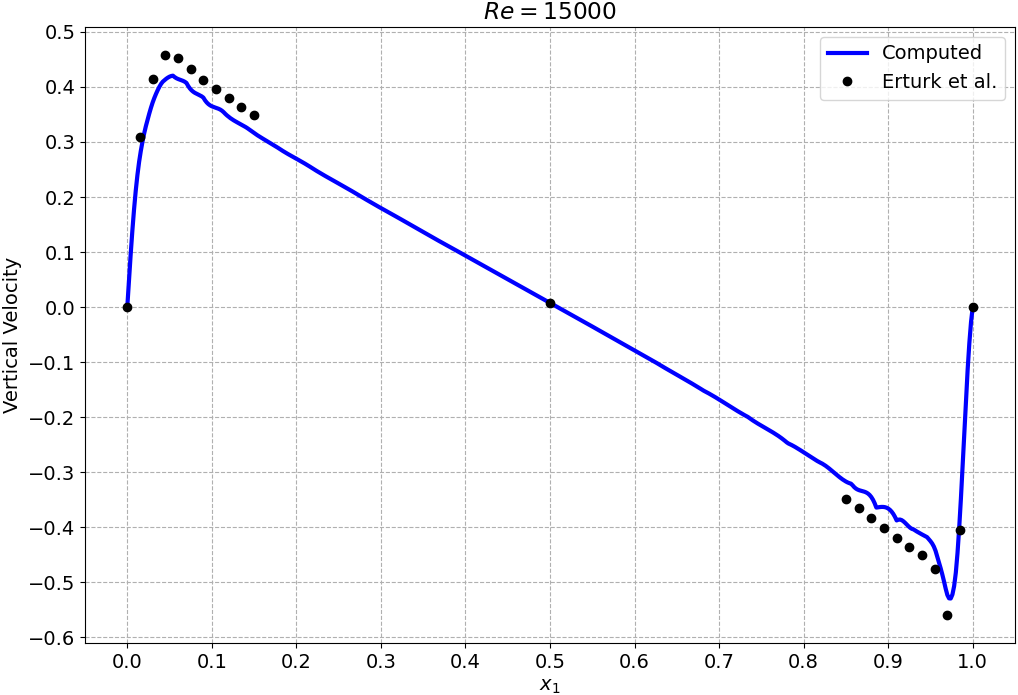}
}

\subfigure{
\includegraphics[width=5.5cm]{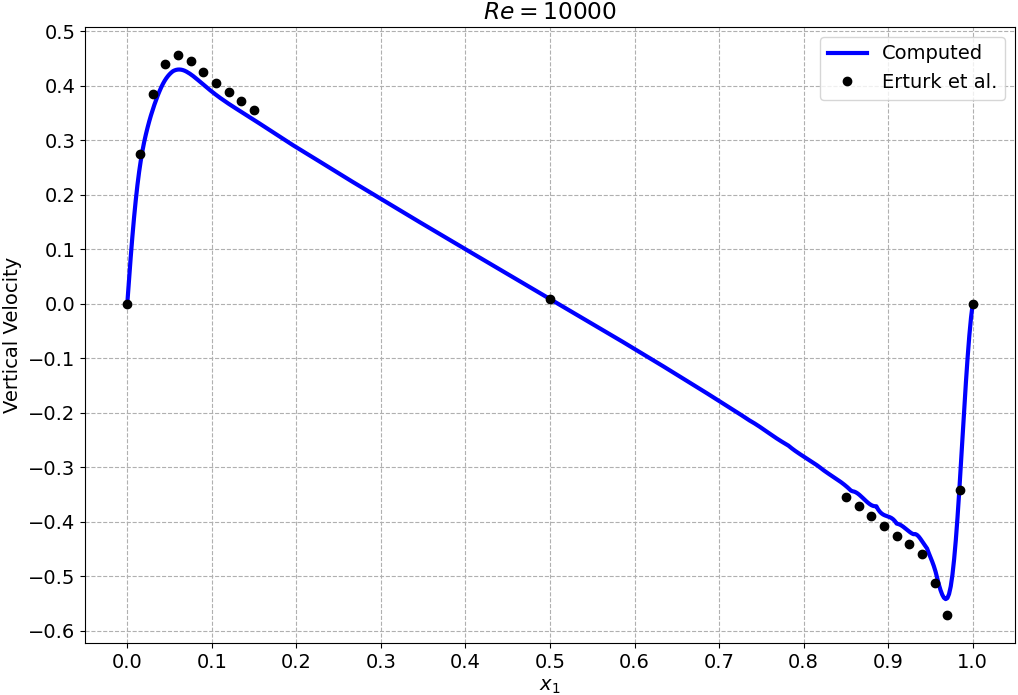}
}
\quad
\subfigure{
\includegraphics[width=5.5cm]{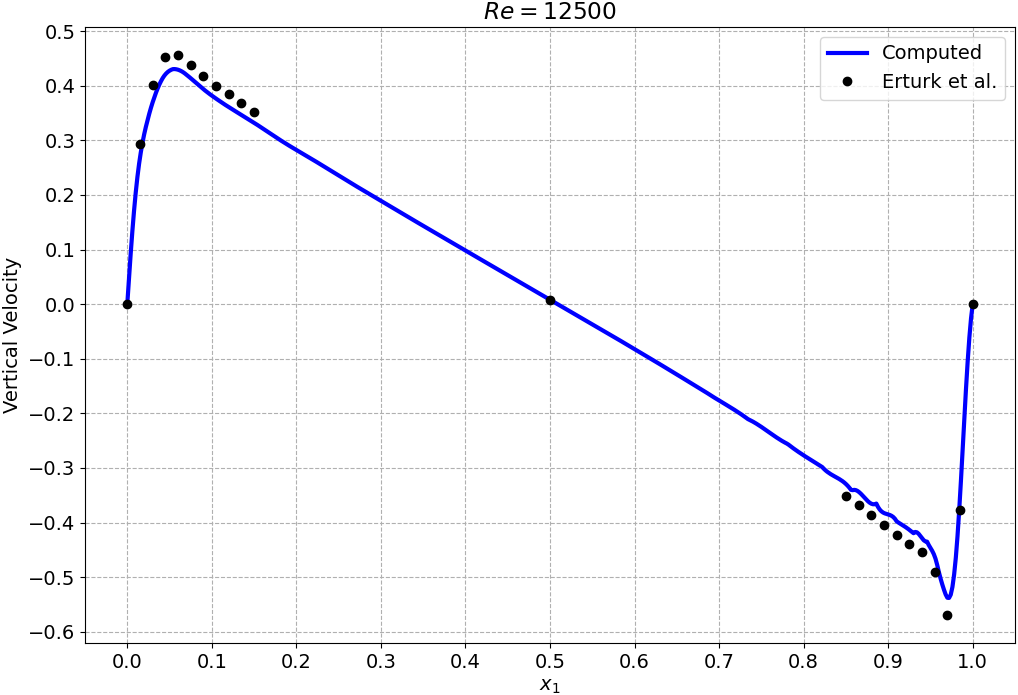}
}
\quad
\subfigure{
\includegraphics[width=5.5cm]{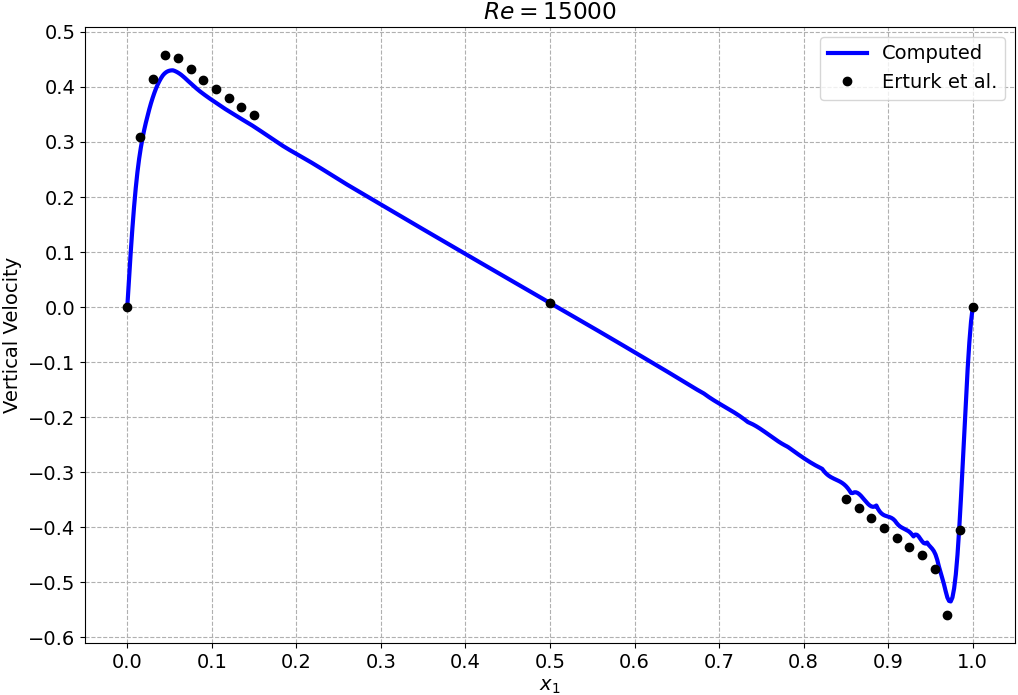}
}
\centering
\caption{\emph{Example 4.} Under the coarse mesh, computed vertical component $u_2$ of the velocity along the horizontal centreline $x_2=0.5$ compared with the reference results by a high-order accurate finite difference scheme under a fine mesh $601\times 601$ from Erturk et al. \cite{Erturk2005747}, for $\mathrm{Re}=10000$, $12500$ and $15000$ from left to right: $\bm{\mathrm{P}}_3^{\mathrm{bubble}}$-$\mathrm{P}_2^{\mathrm{dc}}$ (upper row) and $\bm{\mathrm{P}}_3^{\mathrm{bubble}}$-$\mathrm{P}_2^{\mathrm{dc}}$-$\bm{\mathrm{BDM}}_3$ (lower row)}\label{fig:ex4:mesh40:u2:0.5}
\end{figure}

\begin{figure}[htbp]
\centering
\subfigure[The reference result \cite{Erturk2005747}]{
\includegraphics[width=5.5cm]{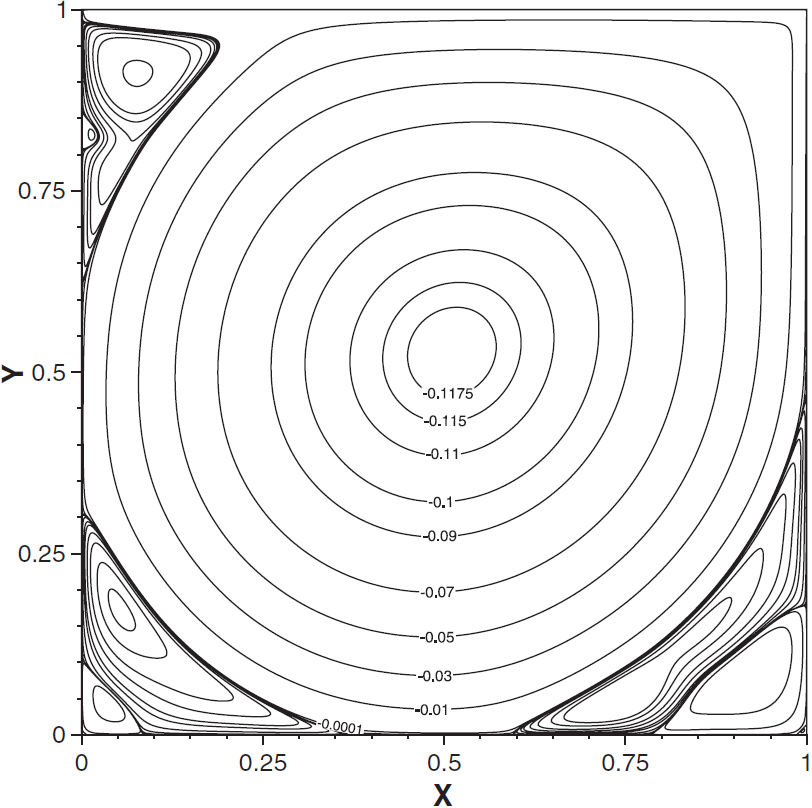}
}
\quad
\subfigure[$\bm{\mathrm{P}}_3^{\mathrm{bubble}}$-$\mathrm{P}_2^{\mathrm{dc}}$-$\bm{\mathrm{BDM}}_3$]{
\includegraphics[width=5.5cm]{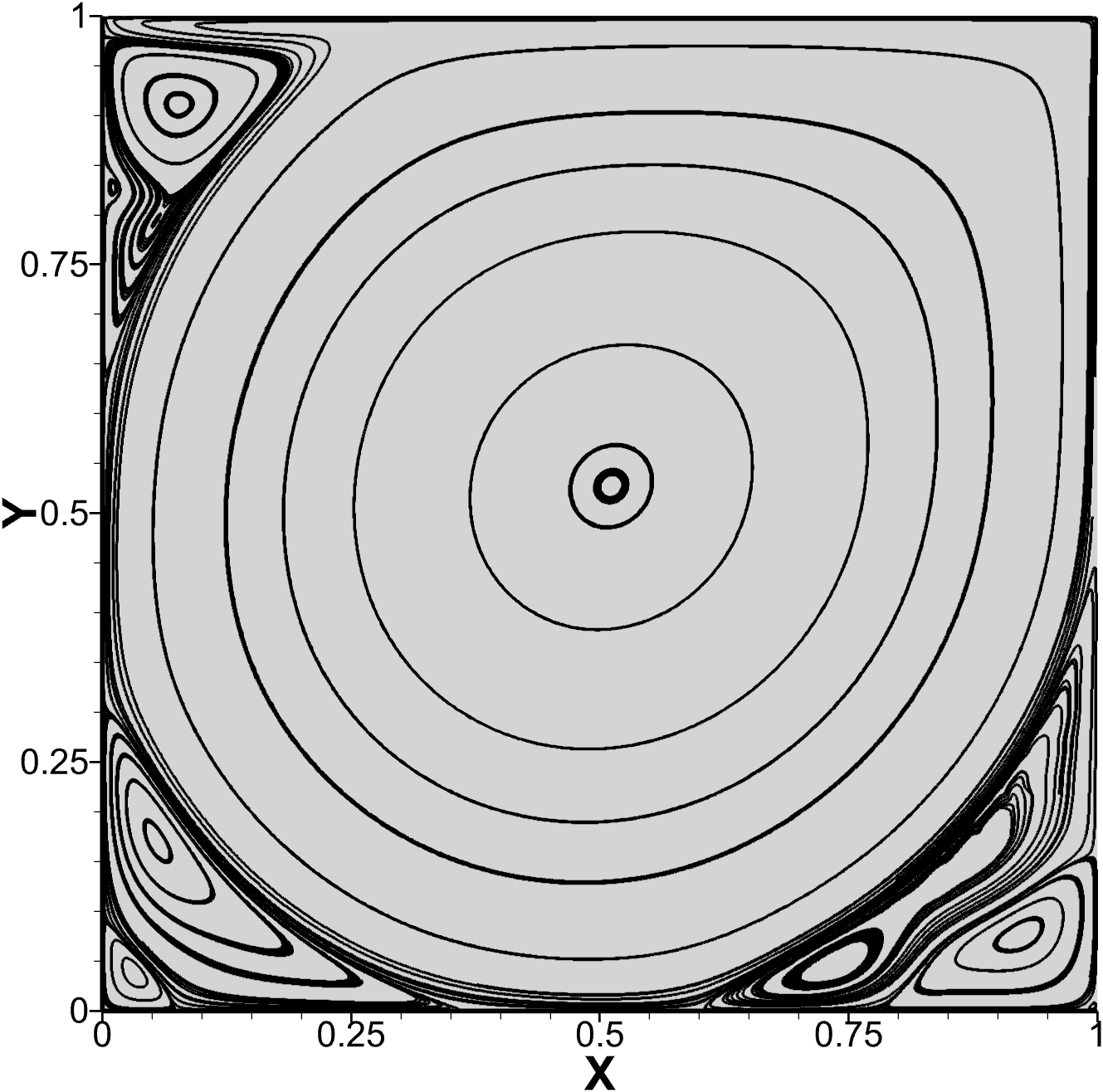}
}
\quad
\subfigure[$\bm{\mathrm{P}}_3^{\mathrm{bubble}}$-$\mathrm{P}_2^{\mathrm{dc}}$]{
\includegraphics[width=5.5cm]{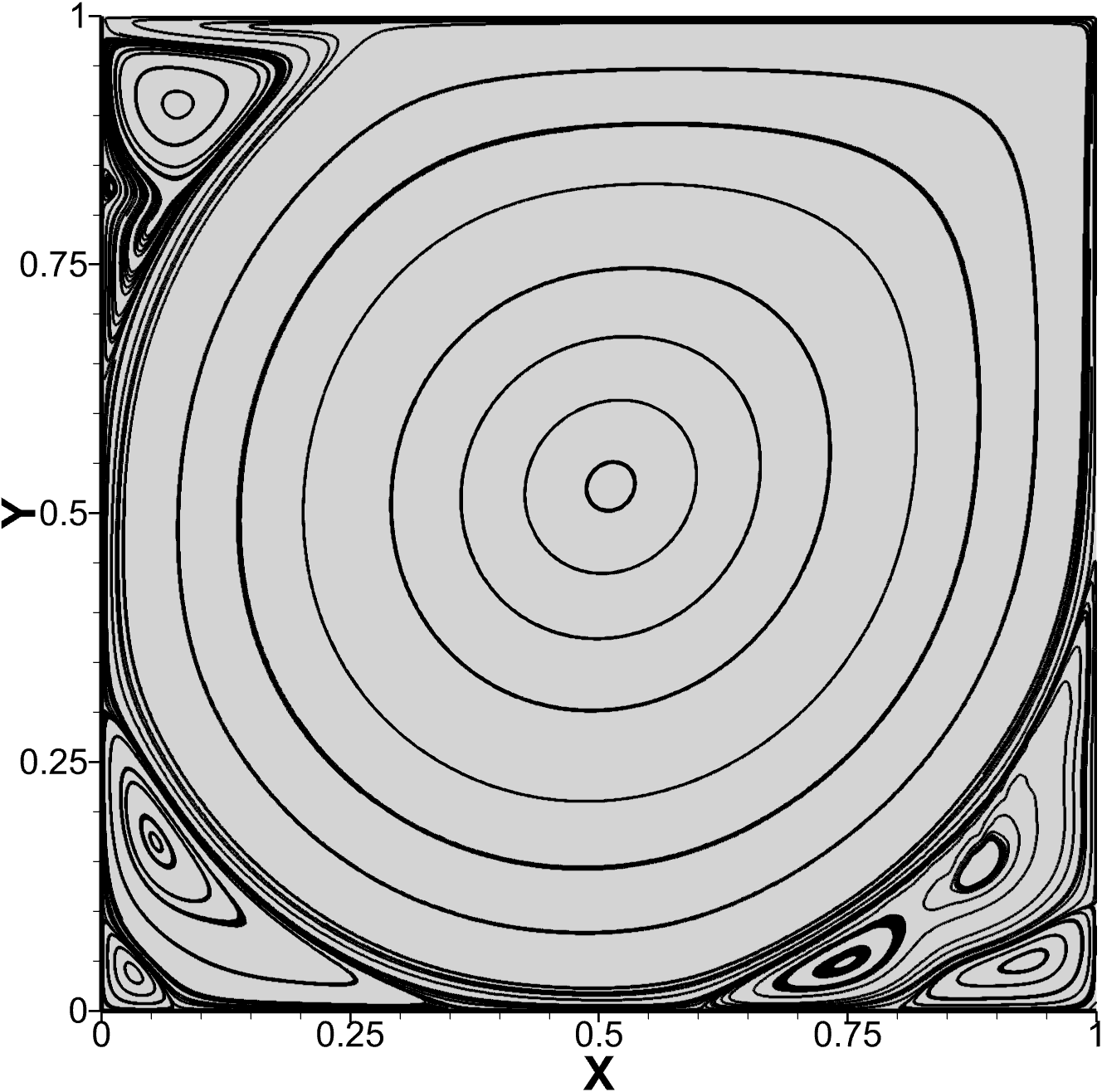}
}
\centering
\caption{\emph{Example 4.} Under the coarse mesh, streamline contours for the primary and secondary vortices for $\mathrm{Re}=15000$: (a) a high-order accurate finite difference scheme under a fine mesh $601\times 601$ from Erturk et al. \cite{Erturk2005747}, (b) $\bm{\mathrm{P}}_3^{\mathrm{bubble}}$-$\mathrm{P}_2^{\mathrm{dc}}$-$\bm{\mathrm{BDM}}_3$ and (c) $\bm{\mathrm{P}}_3^{\mathrm{bubble}}$-$\mathrm{P}_2^{\mathrm{dc}}$}\label{fig:ex4:mesh40:sl:15000}
\end{figure}

Next, the performance of the two methods will be compared under the coarse mesh. In this example, the evaluation system consists of three items, which are the horizontal component $u_1$ of the velocity along the vertical centreline $x_1=0.5$, vertical component $u_2$ of the velocity along the horizontal centreline $x_2=0.5$ and streamline contours for vortices, respectively. As results, Figures \ref{fig:ex4:mesh40:u1:0.5}, \ref{fig:ex4:mesh40:u2:0.5} and \ref{fig:ex4:mesh40:sl:15000} are presented to show the performance differences between the two high-order methods $\bm{\mathrm{P}}_3^{\mathrm{bubble}}$-$\mathrm{P}_2^{\mathrm{dc}}$ and $\bm{\mathrm{P}}_3^{\mathrm{bubble}}$-$\mathrm{P}_2^{\mathrm{dc}}$-$\bm{\mathrm{BDM}}_3$ for $\mathrm{Re}\geqslant 10^4$ under the coarse mesh. Here we choose the numerical solutions computed by a high-order accurate finite difference scheme under a fine mesh $601\times 601$ provided in \cite{Erturk2005747} as the reference results. It follows from Figures \ref{fig:ex4:mesh40:u1:0.5}, \ref{fig:ex4:mesh40:u2:0.5} and \ref{fig:ex4:mesh40:sl:15000} that $\bm{\mathrm{P}}_3^{\mathrm{bubble}}$-$\mathrm{P}_2^{\mathrm{dc}}$-$\bm{\mathrm{BDM}}_3$ does distinctly better than $\bm{\mathrm{P}}_3^{\mathrm{bubble}}$-$\mathrm{P}_2^{\mathrm{dc}}$ in the aspects of approximations, weakening pseudo numerical oscillations and capturing correct vortex locations. It also means that the proposed method is not only pressure-robust but also an efficient stabilization method without additional factitious stabilization terms in some practical problems.

\begin{figure}[htbp]
\centering
\subfigure{
\includegraphics[width=5.5cm]{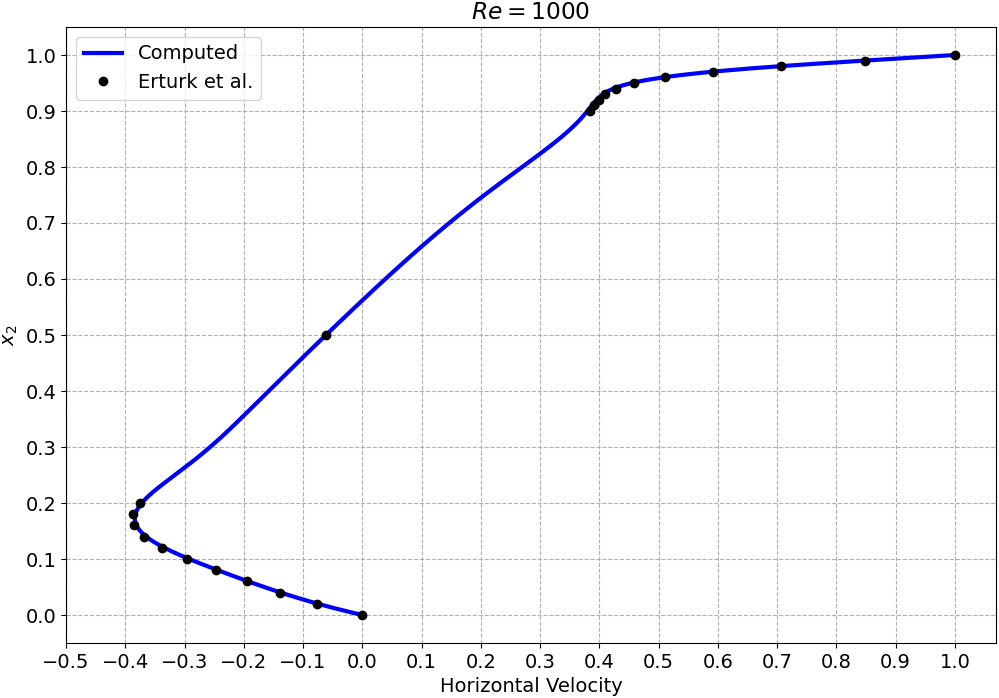}
}
\quad
\subfigure{
\includegraphics[width=5.5cm]{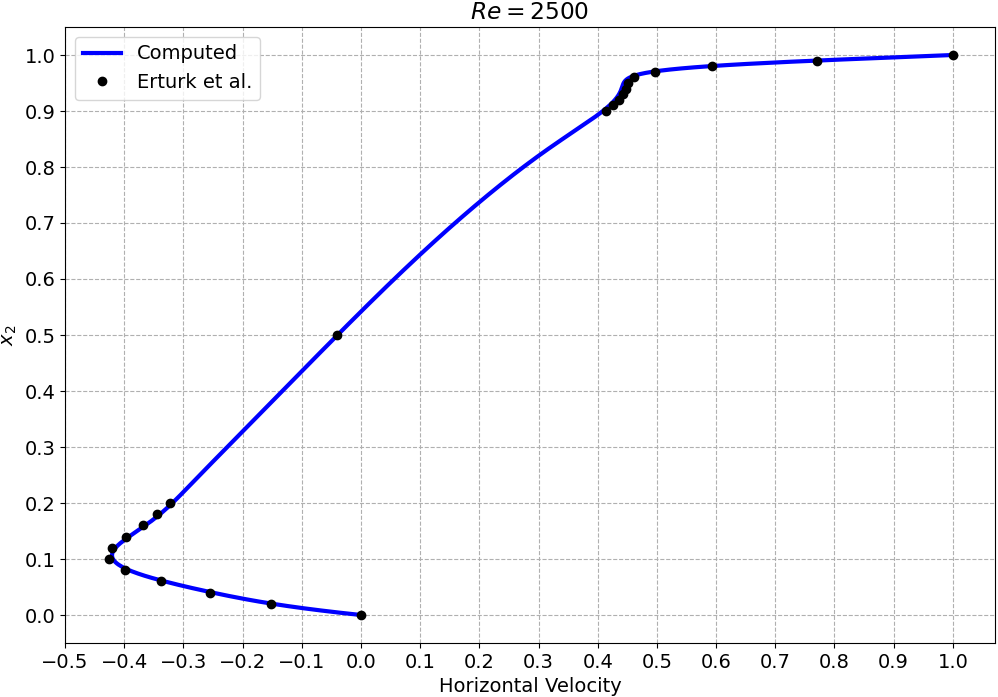}
}
\quad
\subfigure{
\includegraphics[width=5.5cm]{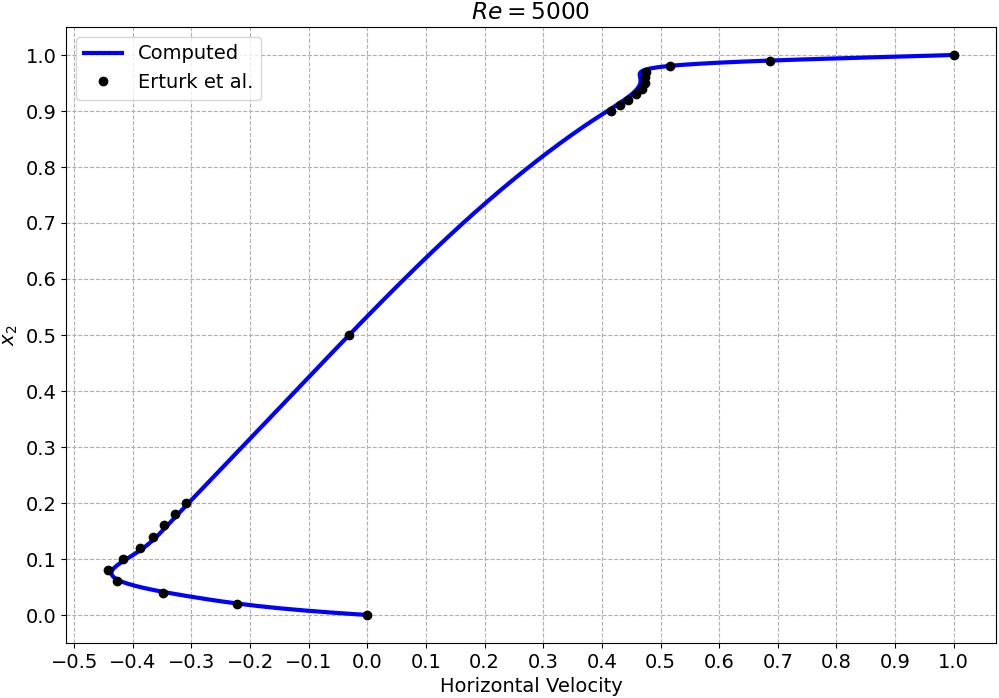}
}

\subfigure{
\includegraphics[width=5.5cm]{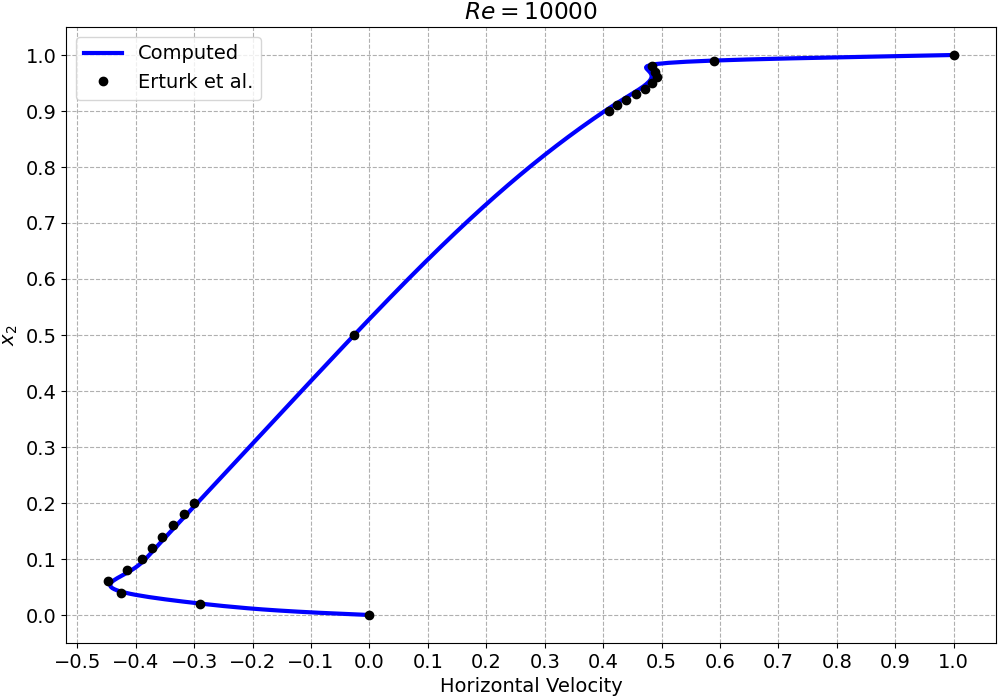}
}
\quad
\subfigure{
\includegraphics[width=5.5cm]{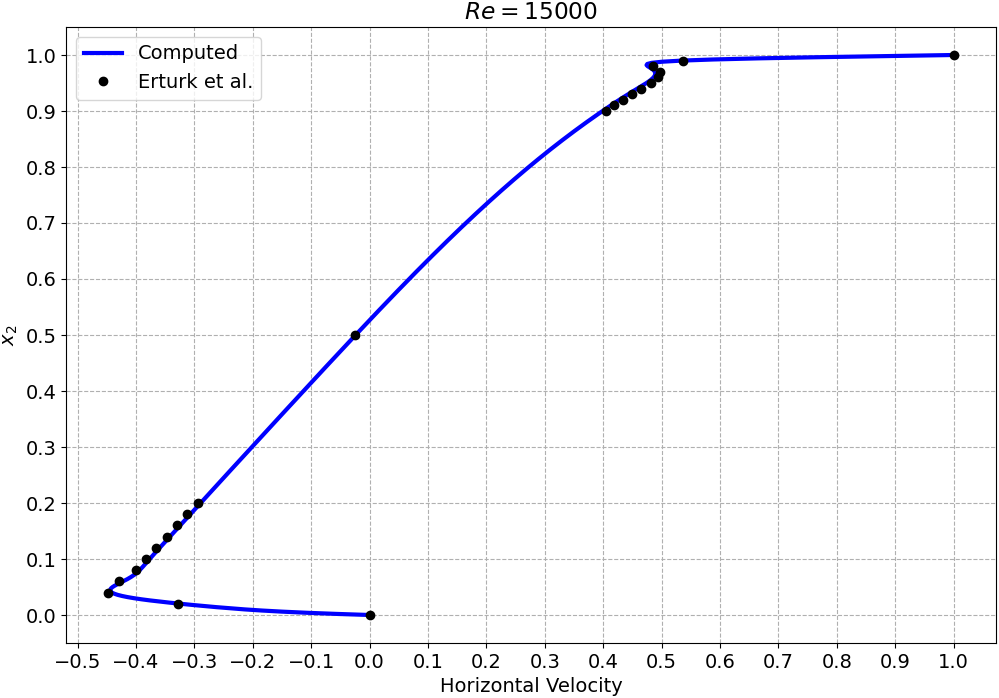}
}
\quad
\subfigure{
\includegraphics[width=5.5cm]{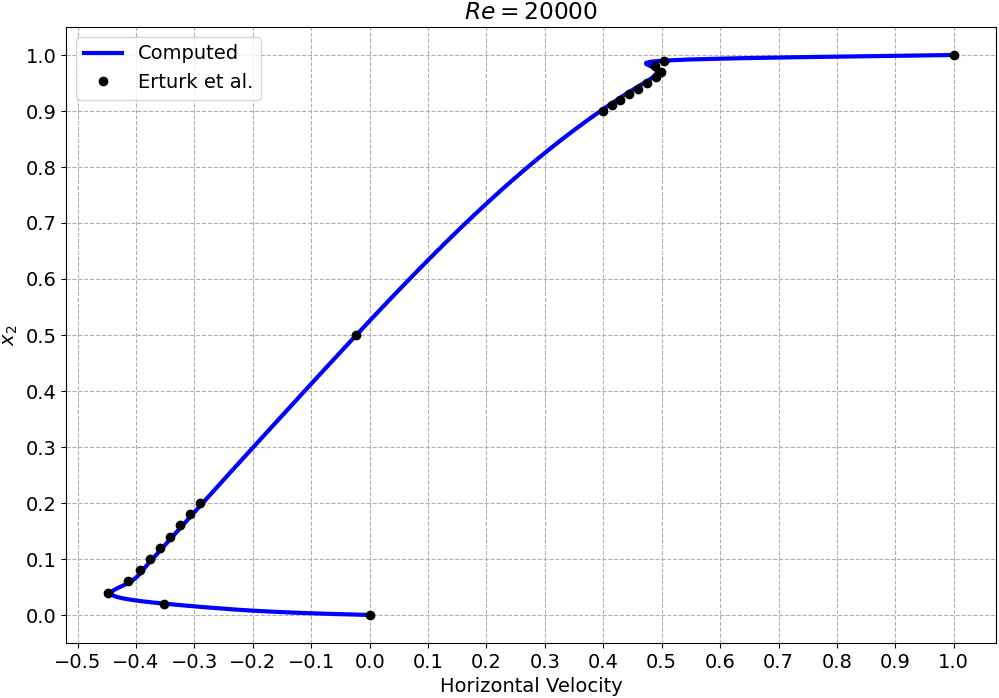}
}
\centering
\caption{\emph{Example 4.} Under the refined mesh, computed horizontal component $u_1$ of the velocity along the vertical centreline $x_1=0.5$ by $\bm{\mathrm{P}}_3^{\mathrm{bubble}}$-$\mathrm{P}_2^{\mathrm{dc}}$-$\bm{\mathrm{BDM}}_3$, compared with the reference results by a high-order accurate finite difference scheme under a fine mesh $601\times 601$ from Erturk et al. \cite{Erturk2005747}, for $\mathrm{Re}=1000$, $2500$, $5000$, $10000$, $15000$ and $20000$ from top-left to bottom-right}\label{fig:ex4:mesh80:u1:0.5}
\end{figure}

\begin{figure}[htbp]
\centering
\subfigure{
\includegraphics[width=5.5cm]{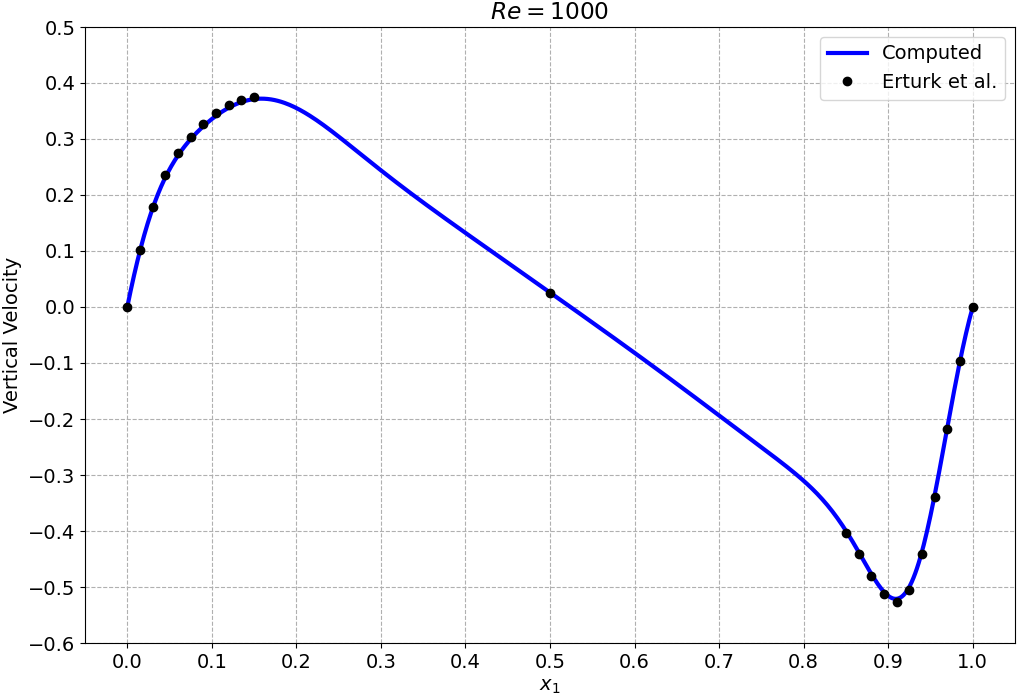}
}
\quad
\subfigure{
\includegraphics[width=5.5cm]{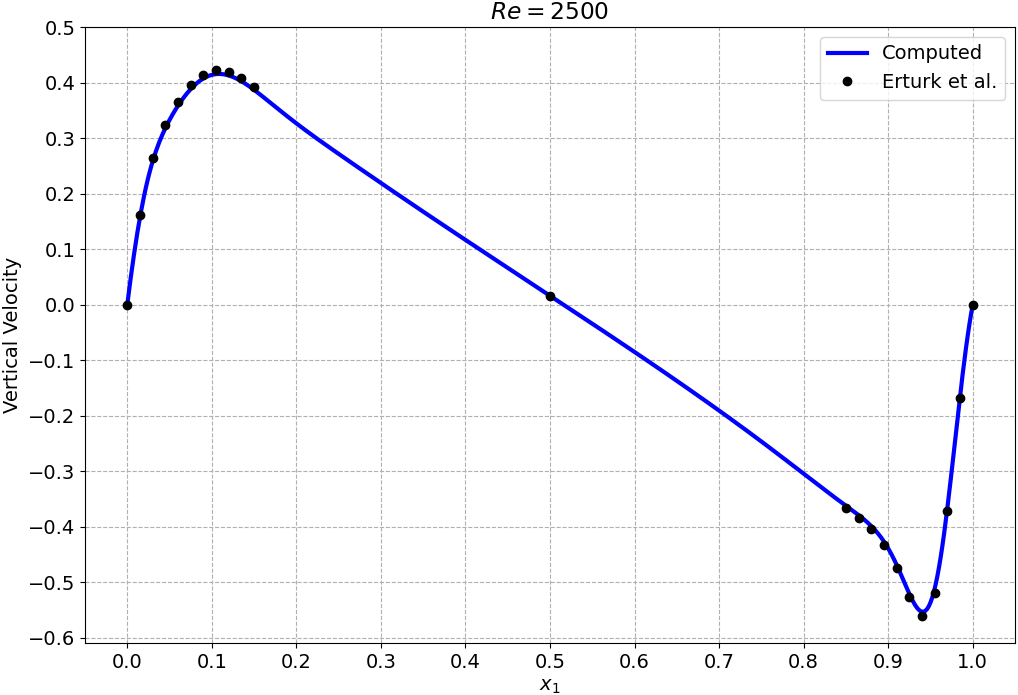}
}
\quad
\subfigure{
\includegraphics[width=5.5cm]{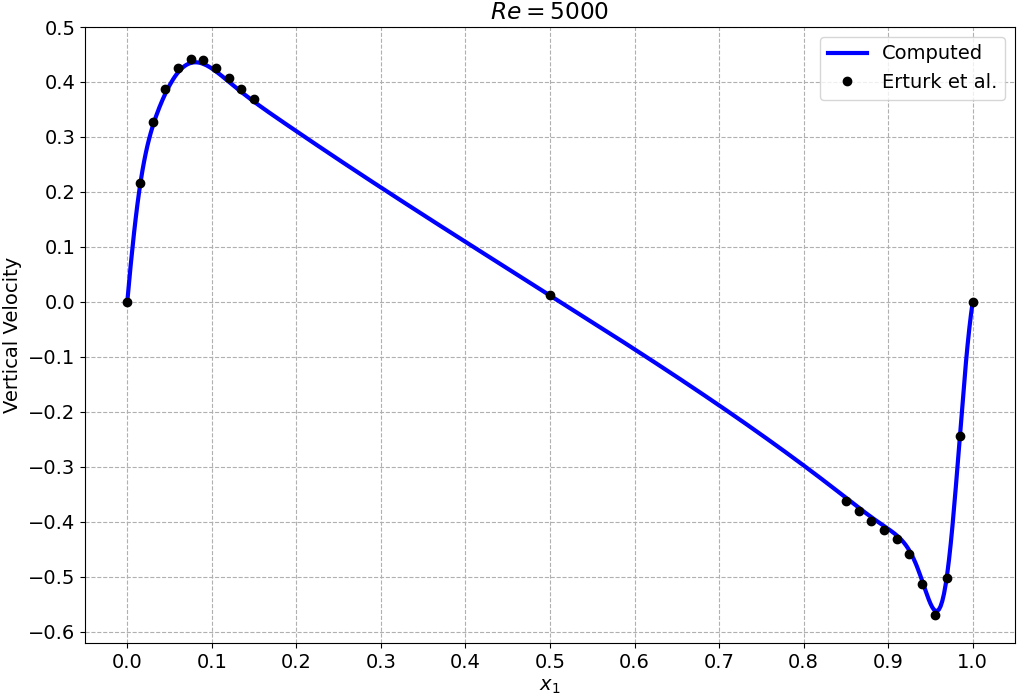}
}

\subfigure{
\includegraphics[width=5.5cm]{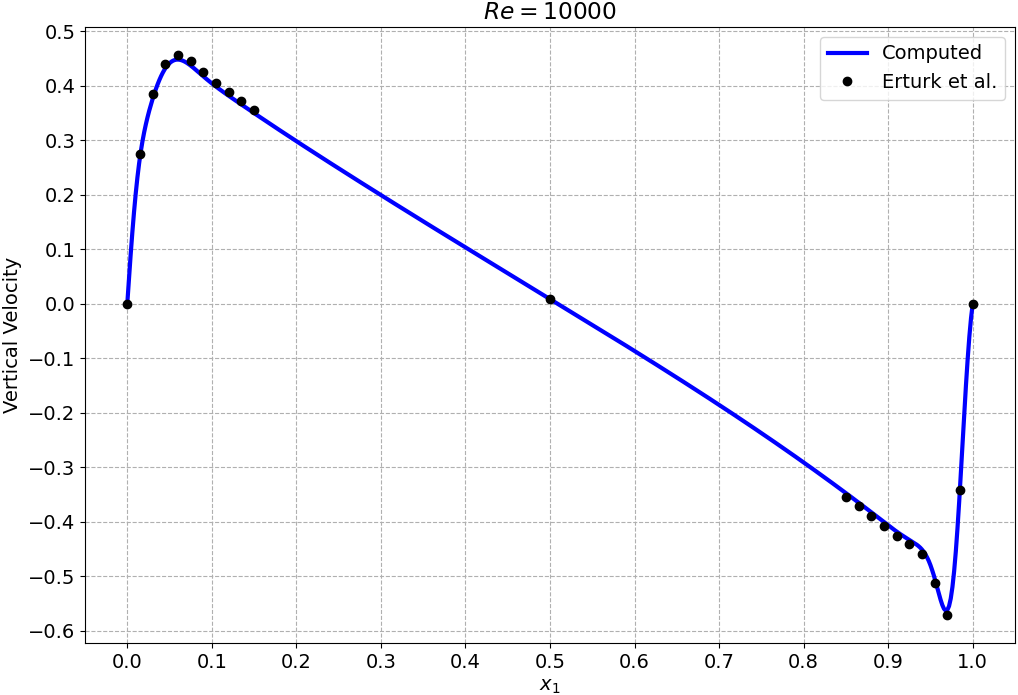}
}
\quad
\subfigure{
\includegraphics[width=5.5cm]{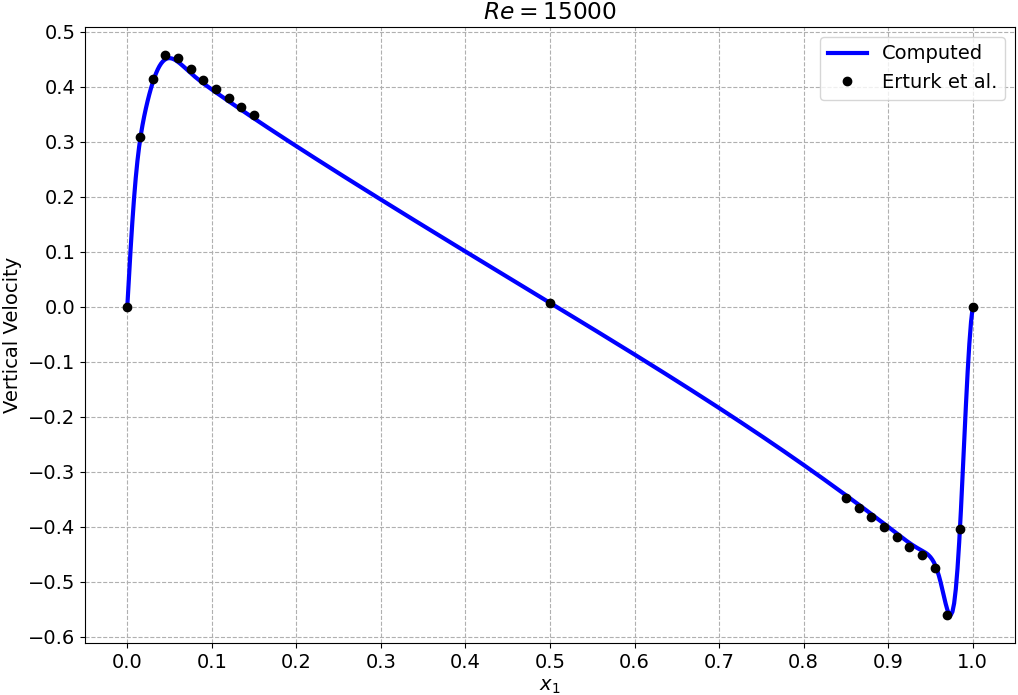}
}
\quad
\subfigure{
\includegraphics[width=5.5cm]{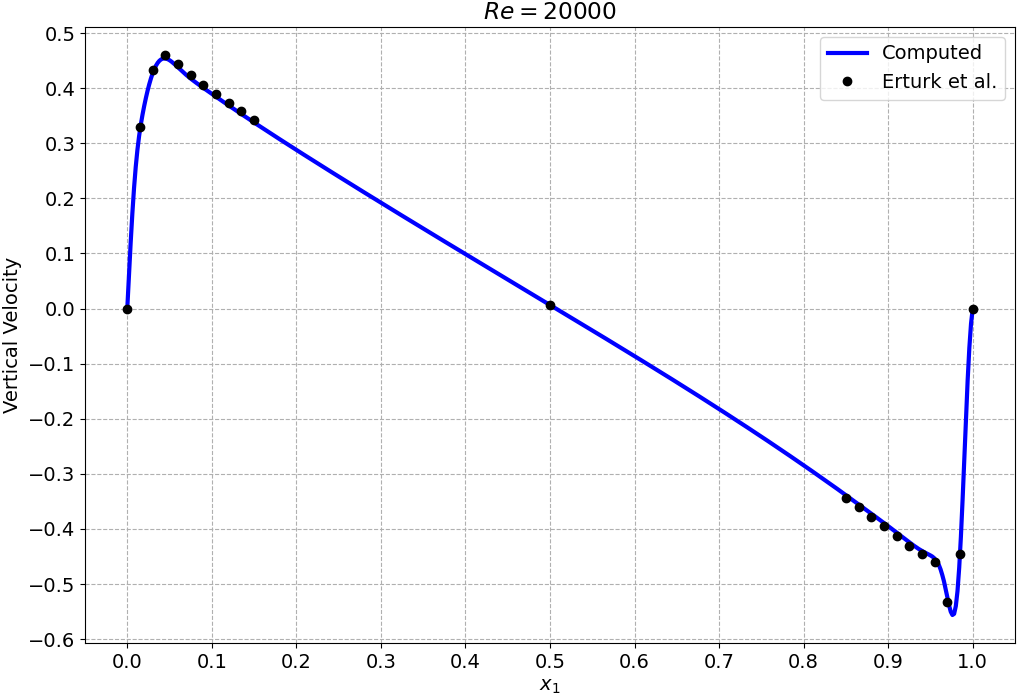}
}
\centering
\caption{\emph{Example 4.} Under the refined mesh, computed vertical component $u_2$ of the velocity along the horizontal centreline $x_2=0.5$ by $\bm{\mathrm{P}}_3^{\mathrm{bubble}}$-$\mathrm{P}_2^{\mathrm{dc}}$-$\bm{\mathrm{BDM}}_3$, compared with the reference results by a high-order accurate finite difference scheme under a fine mesh $601\times 601$ from Erturk et al. \cite{Erturk2005747}, for $\mathrm{Re}=1000$, $2500$, $5000$, $10000$, $15000$ and $20000$ from top-left to bottom-right}\label{fig:ex4:mesh80:u2:0.5}
\end{figure}

\begin{figure}[htbp]
\centering
\subfigure{
\includegraphics[width=5cm]{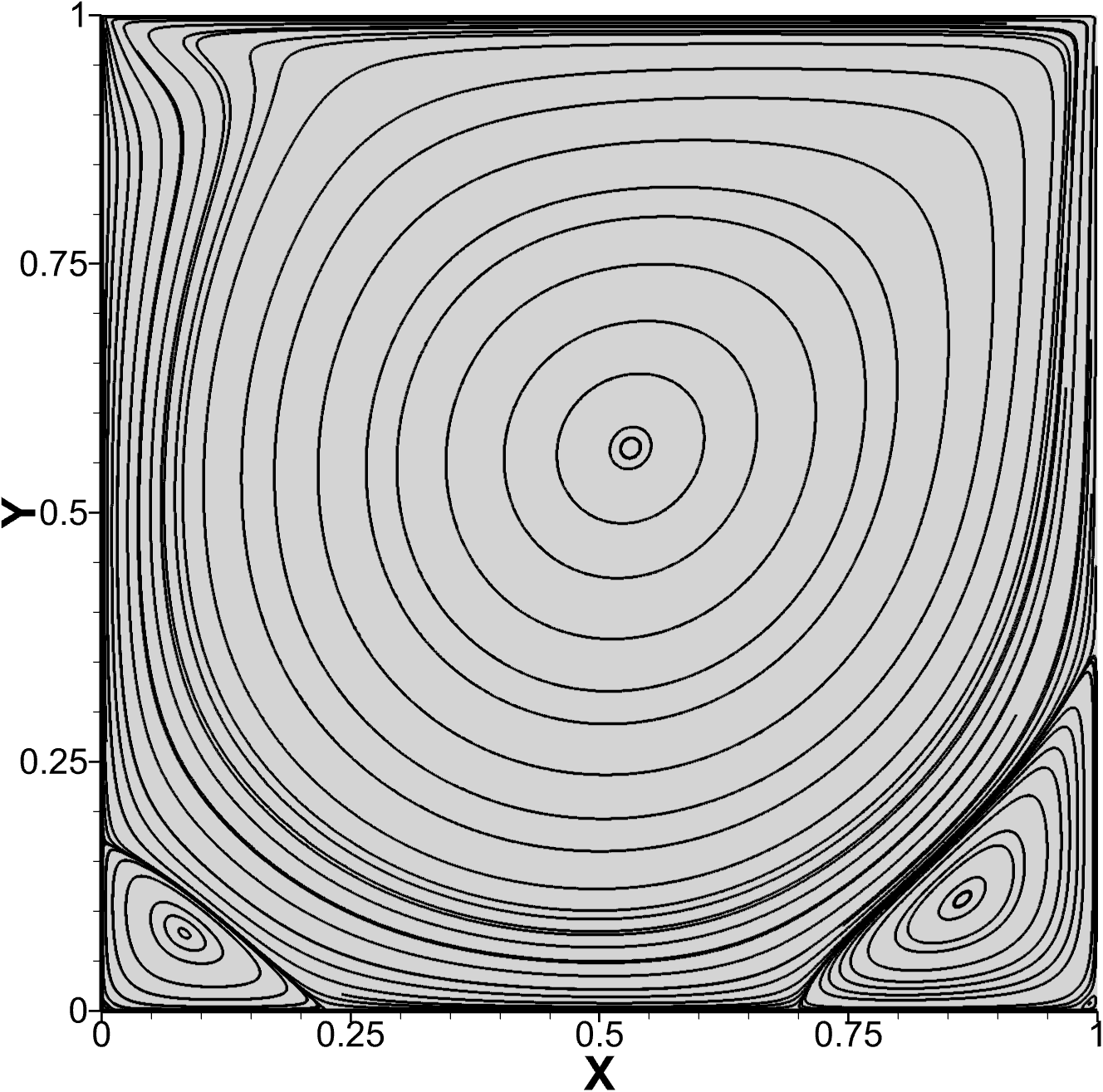}
}
\quad
\subfigure{
\includegraphics[width=5cm]{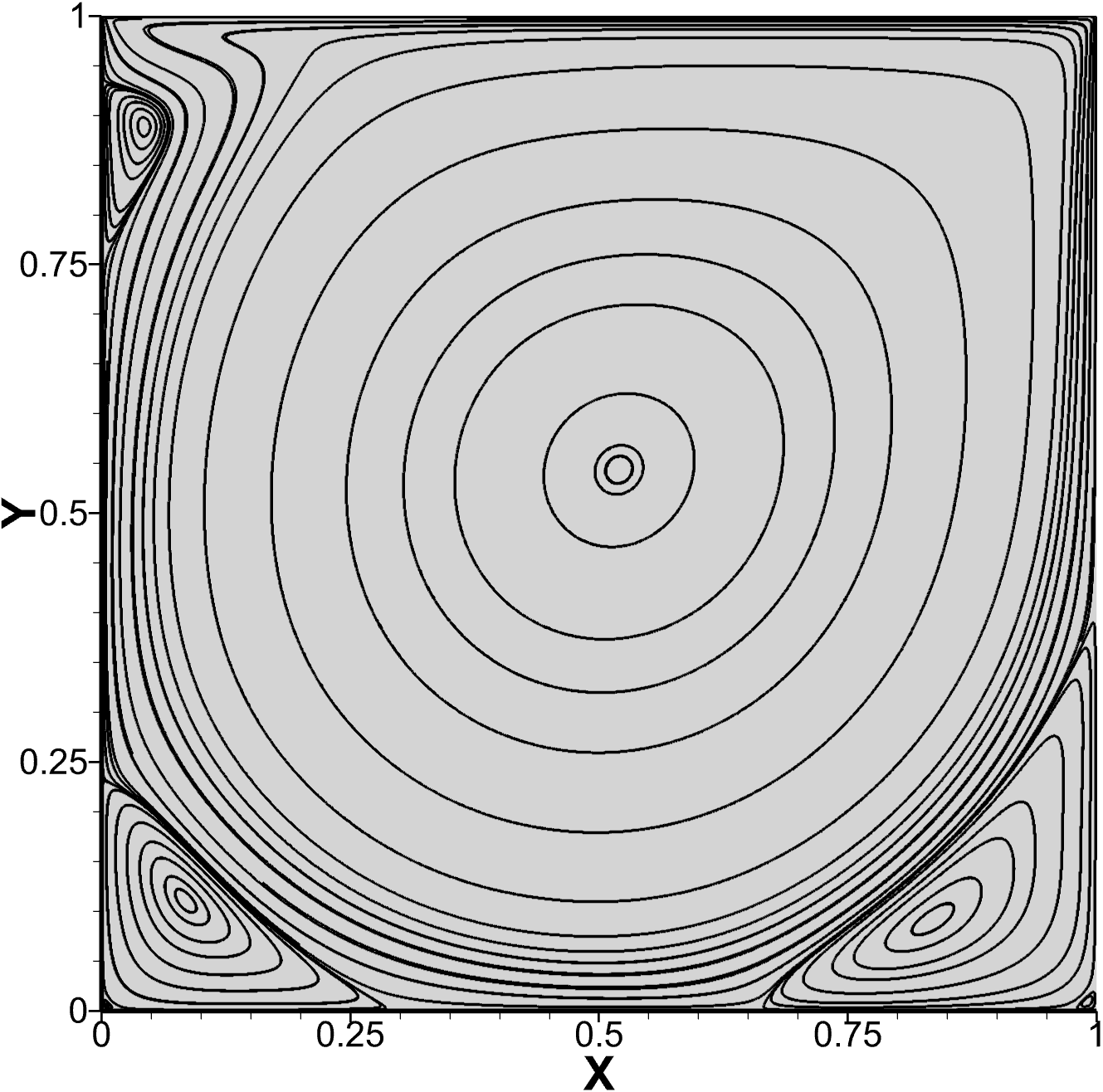}
}
\quad
\subfigure{
\includegraphics[width=5cm]{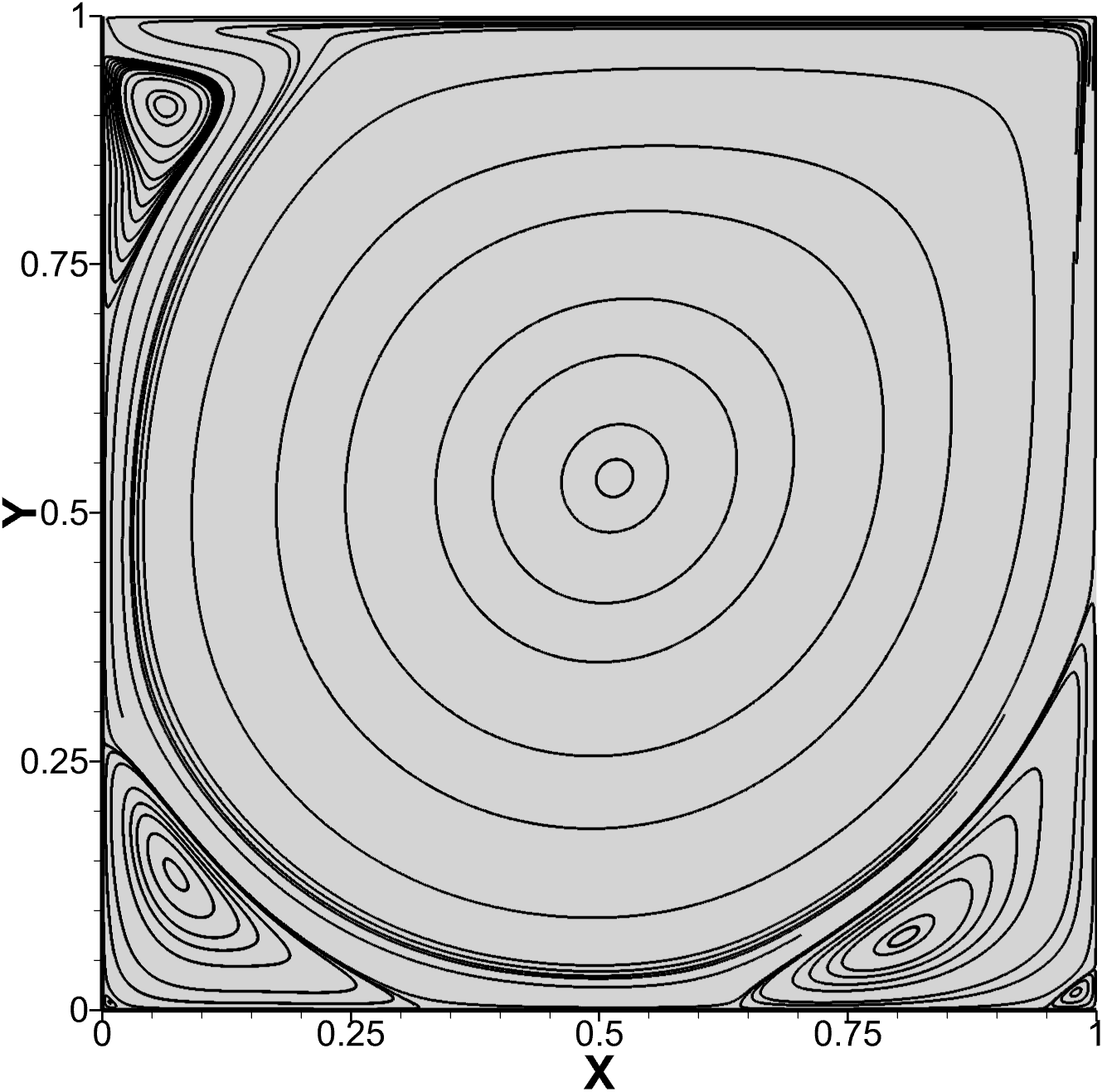}
}

\subfigure{
\includegraphics[width=5cm]{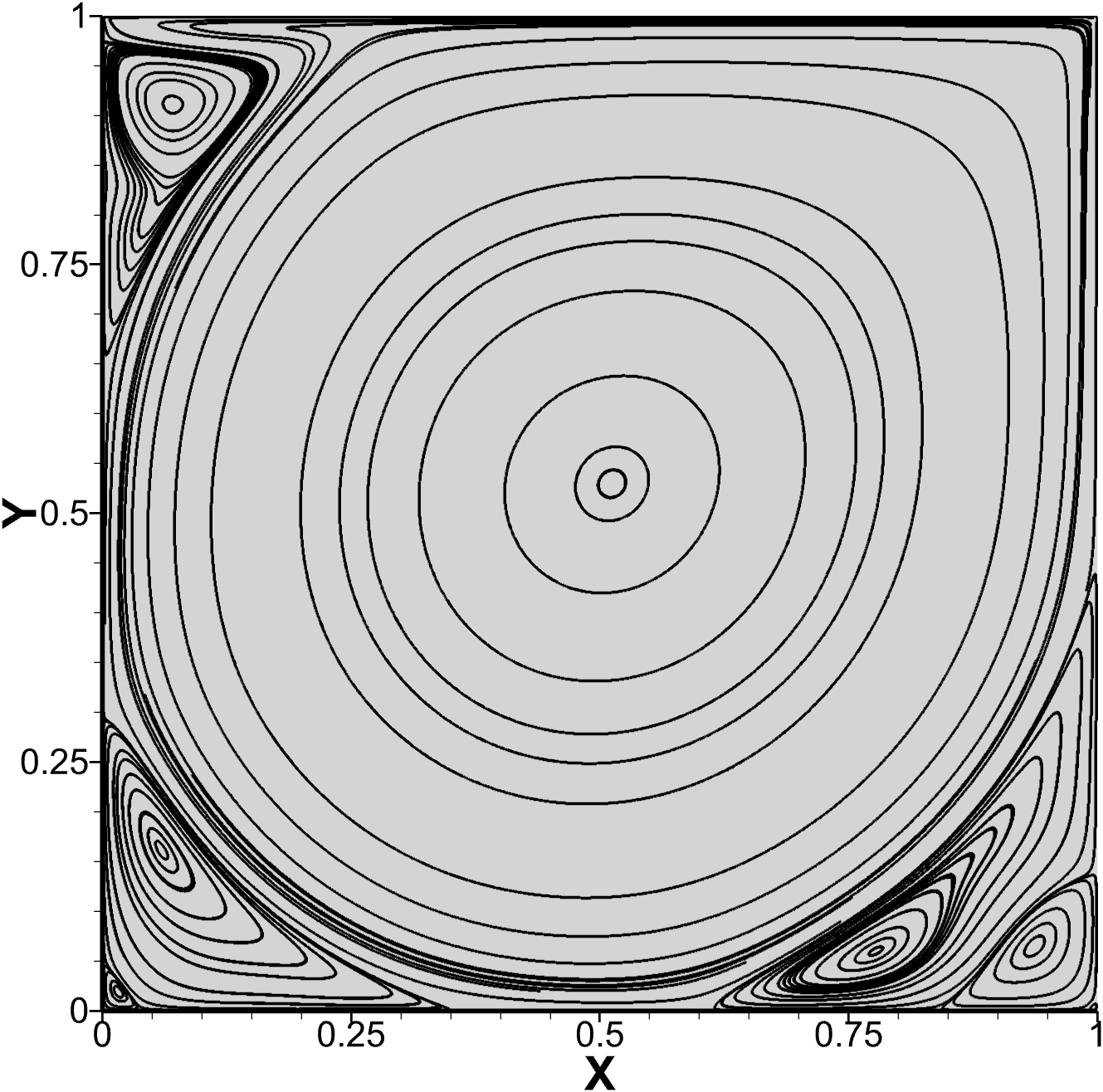}
}
\quad
\subfigure{
\includegraphics[width=5cm]{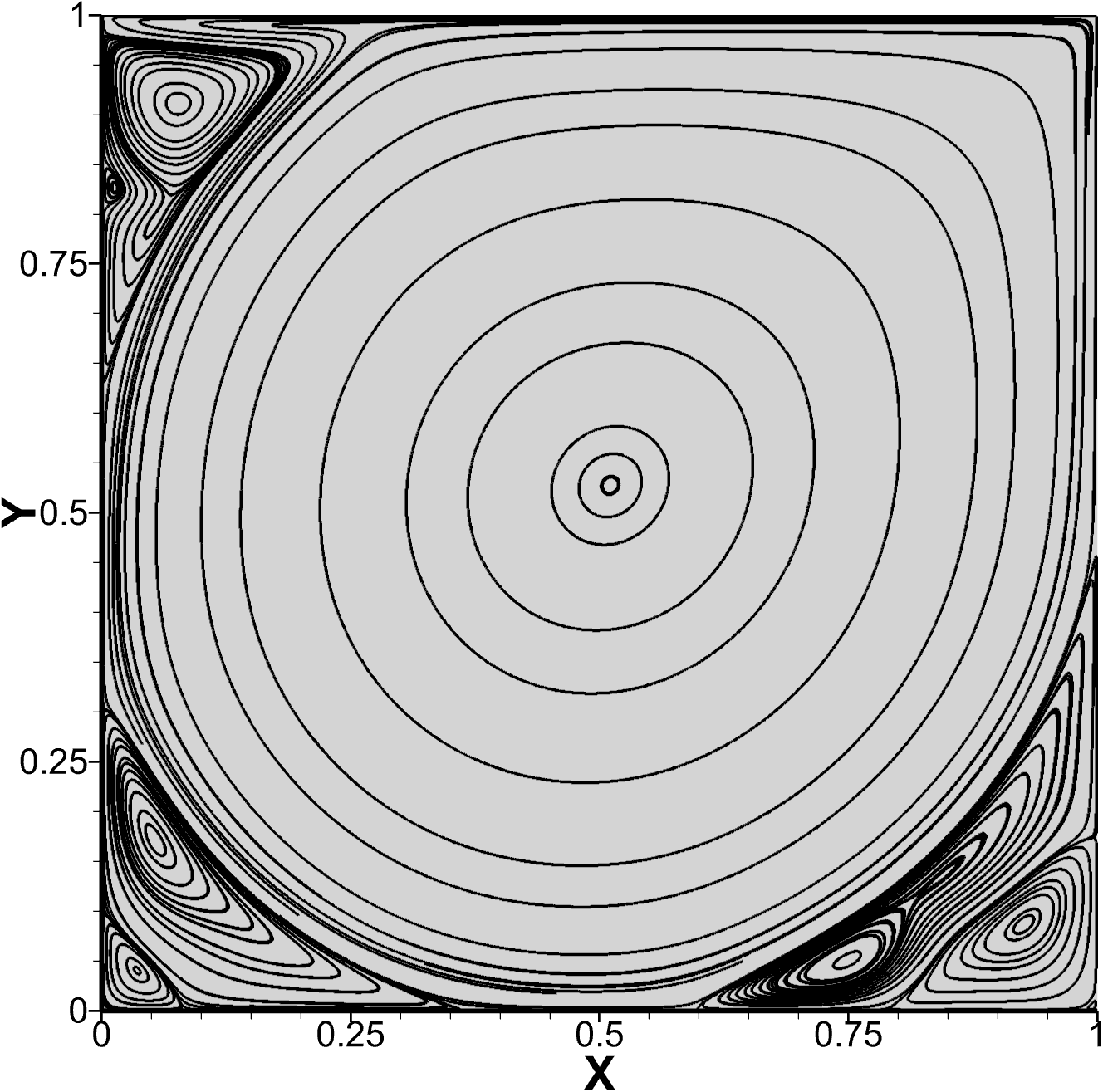}
}
\quad
\subfigure{
\includegraphics[width=5cm]{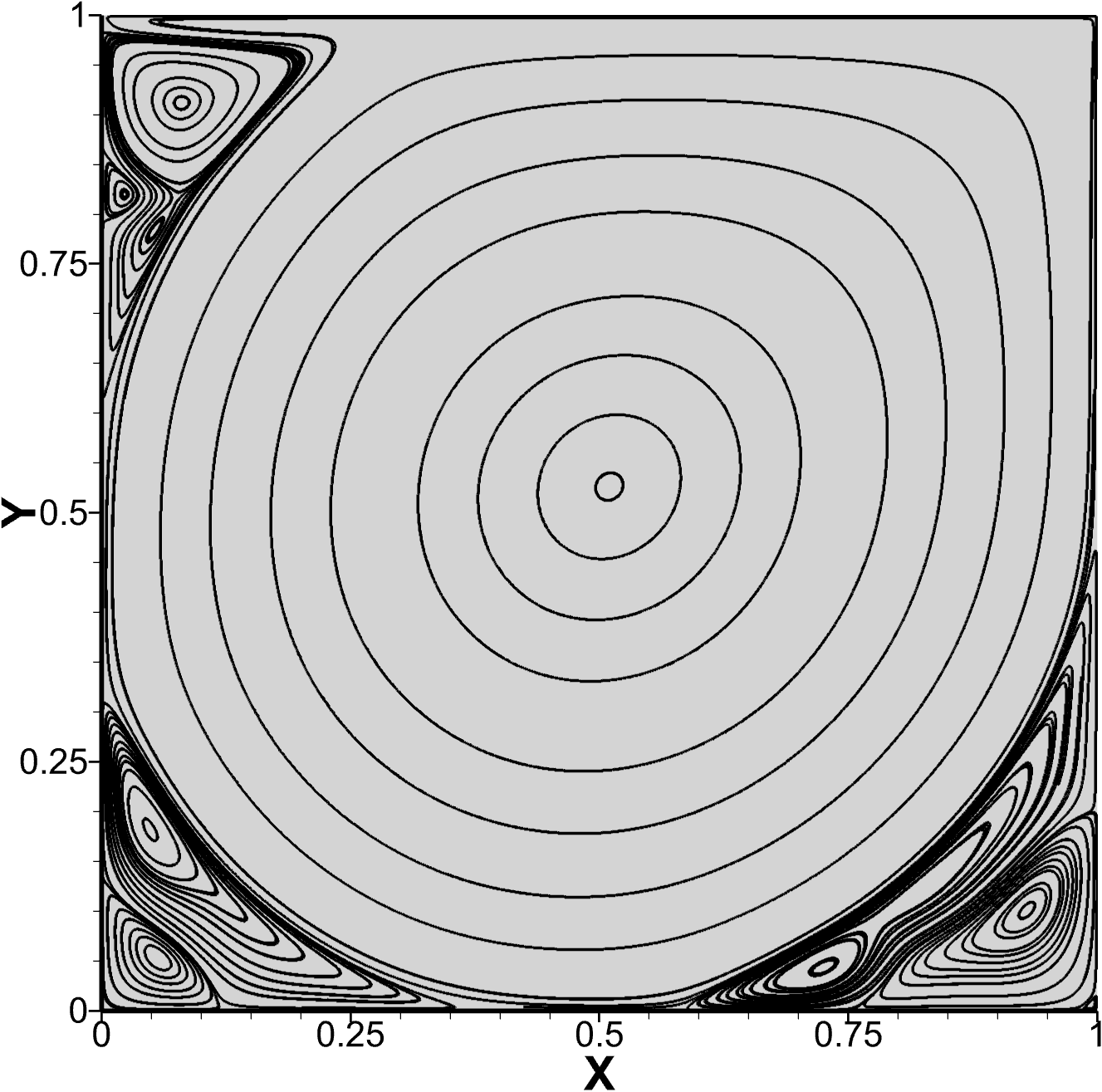}
}
\centering
\caption{\emph{Example 4.} Under the refined mesh, streamline contours for primary and secondary vortices by $\bm{\mathrm{P}}_3^{\mathrm{bubble}}$-$\mathrm{P}_2^{\mathrm{dc}}$-$\bm{\mathrm{BDM}}_3$ for $\mathrm{Re}=1000$, $2500$, $5000$, $10000$, $15000$ and $20000$ from top-left to bottom-right}\label{fig:ex4:mesh80:global:sl}
\end{figure}

Hence, we make use of $\bm{\mathrm{P}}_3^{\mathrm{bubble}}$-$\mathrm{P}_2^{\mathrm{dc}}$-$\bm{\mathrm{BDM}}_3$ to achieve more precise results under the refined mesh. Figures \ref{fig:ex4:mesh80:u1:0.5} and \ref{fig:ex4:mesh80:u2:0.5} present the $u_1$-velocity profiles along a vertical line $x_1=0.5$ and the $u_2$-velocity profiles along a horizontal line $x_2=0.5$ passing through the cavity's geometric centre respectively. These profiles are in good agreement with that of Erturk et al. \cite{Erturk2005747} shown by dark spots in Figures \ref{fig:ex4:mesh80:u1:0.5} and \ref{fig:ex4:mesh80:u2:0.5}. Figure \ref{fig:ex4:mesh80:global:sl} exhibits streamline contours for the primary and secondary vortices at various Reynolds numbers. Likewise, the formation of the counter-rotating secondary vortices, which appear as the Reynolds number increases, almost matches that shown in \cite{Erturk2005747}. Hughes et al. \cite{Hughes2000135} pointed out that the number of unknowns of a discretized problem is a proper indicator for the efficiency of a numerical method. Note that all of the results in \cite{Erturk2005747} as a reference are computed under a quite fine mesh $601\times 601$, and hence the propose method has the reasonable efficient performance for this problem especially at high Reynolds numbers.

\begin{figure}[htbp]
\centering
\subfigure[Eddies $\mathrm{BL}2$, $\mathrm{BL}3$]{
\includegraphics[width=4.5cm]{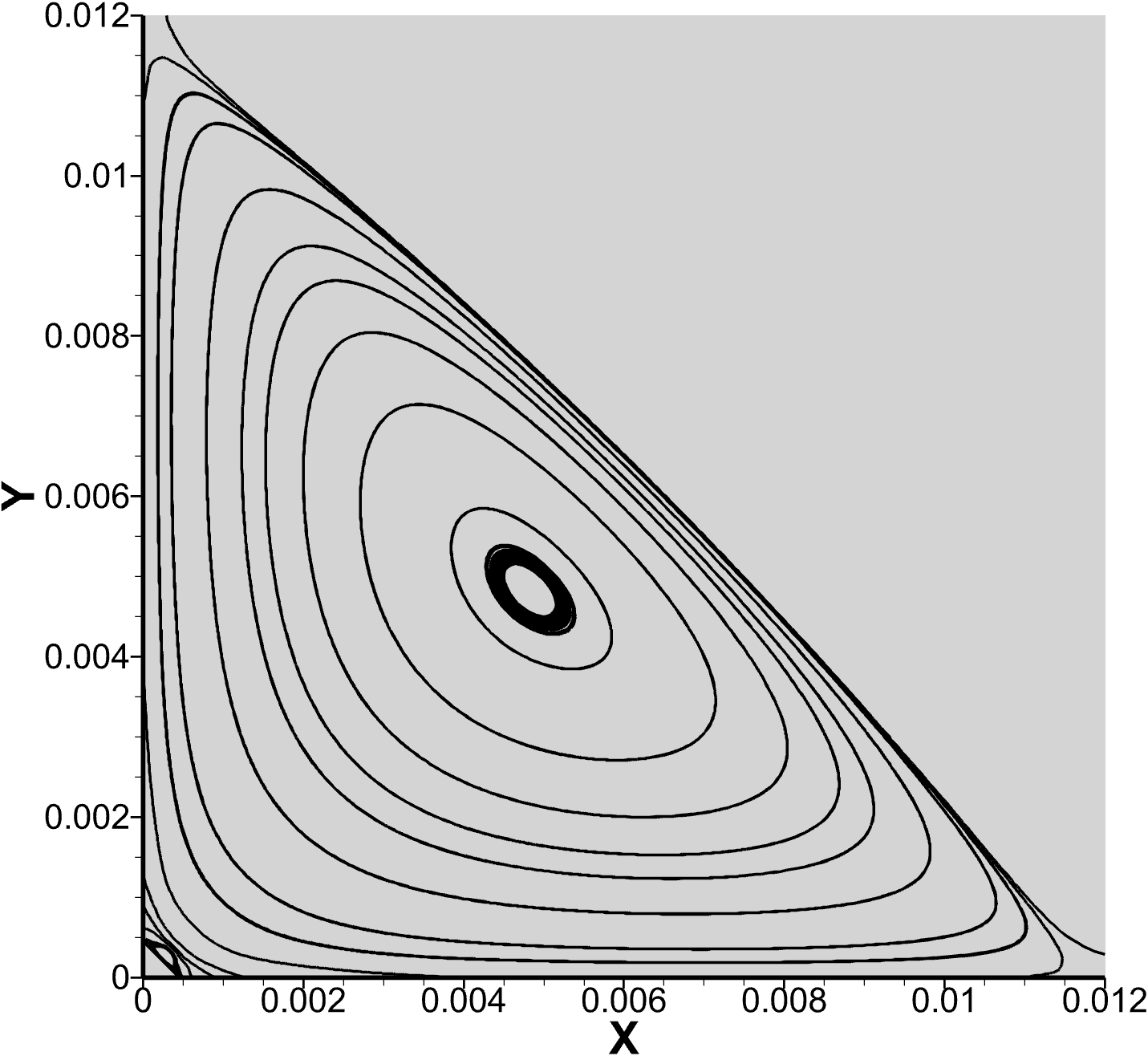}
}
\quad
\subfigure[Eddies $\mathrm{BL}2$, $\mathrm{BL}3$]{
\includegraphics[width=4.5cm]{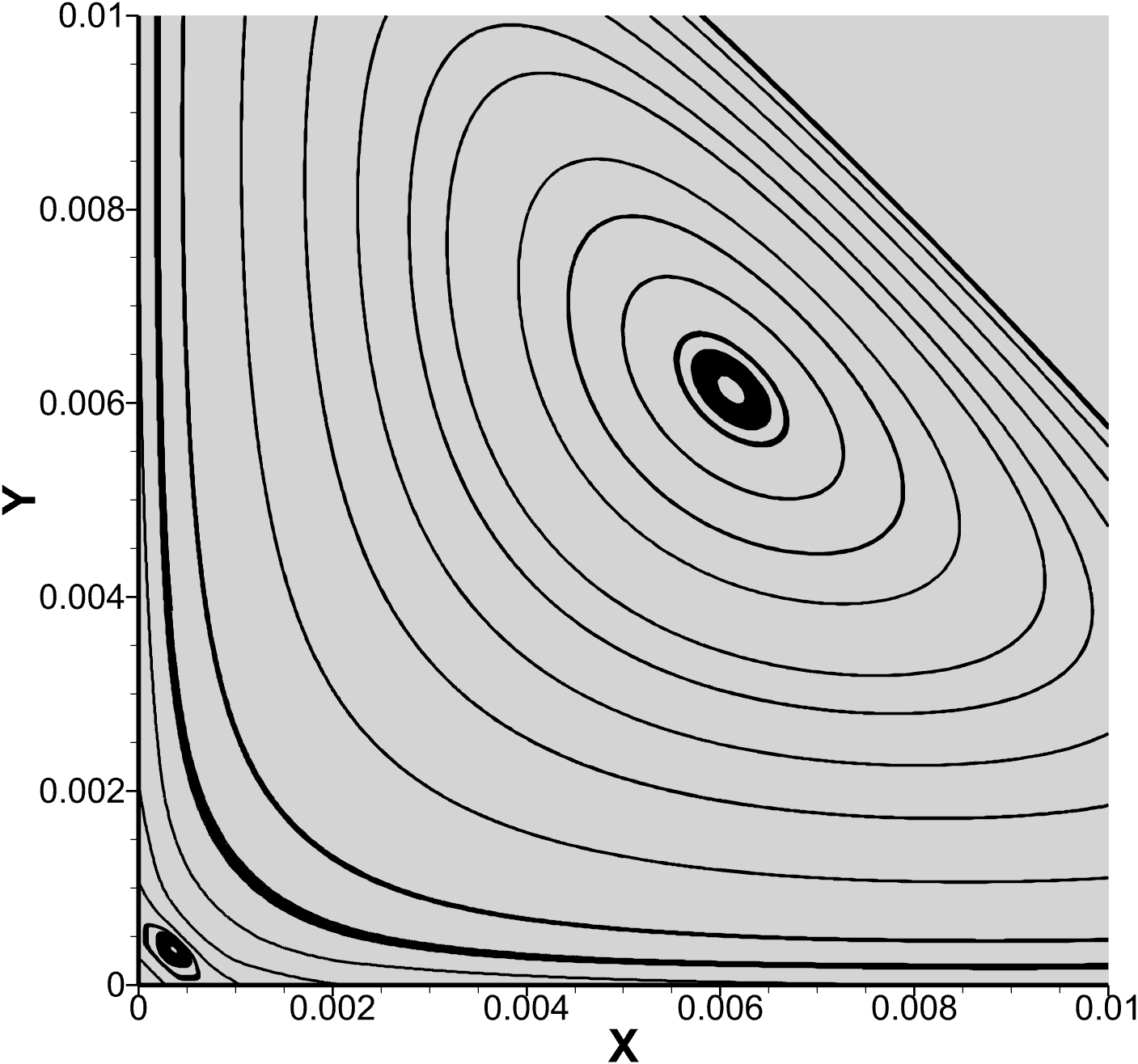}
}
\quad
\subfigure[Eddies $\mathrm{BL}2$, $\mathrm{BL}3$]{
\includegraphics[width=4.5cm]{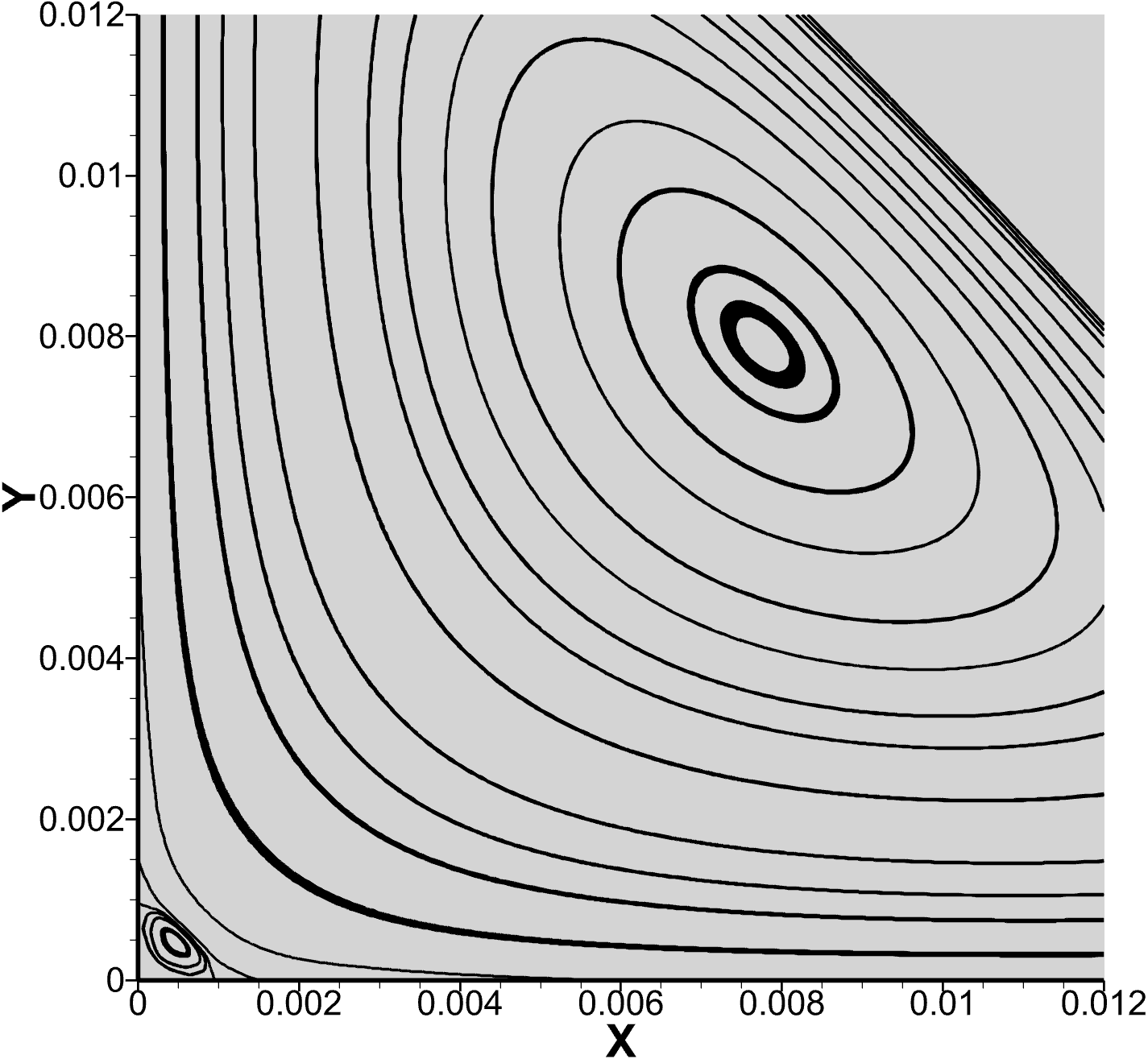}
}

\subfigure[Eddies $\mathrm{BL}2$, $\mathrm{BL}3$]{
\includegraphics[width=4.5cm]{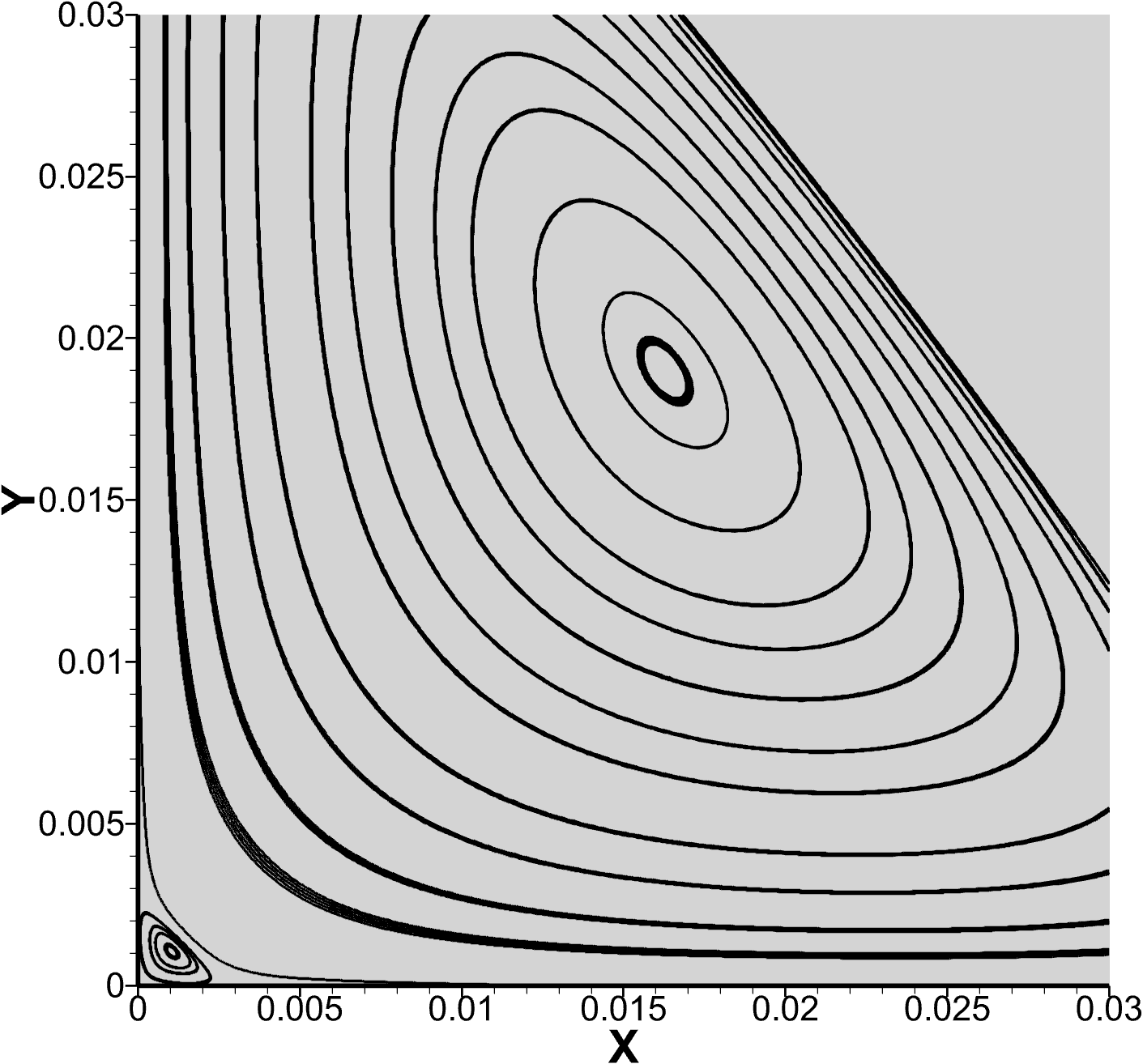}
}
\quad
\subfigure[Eddy $\mathrm{BL}3$]{
\includegraphics[width=4.5cm]{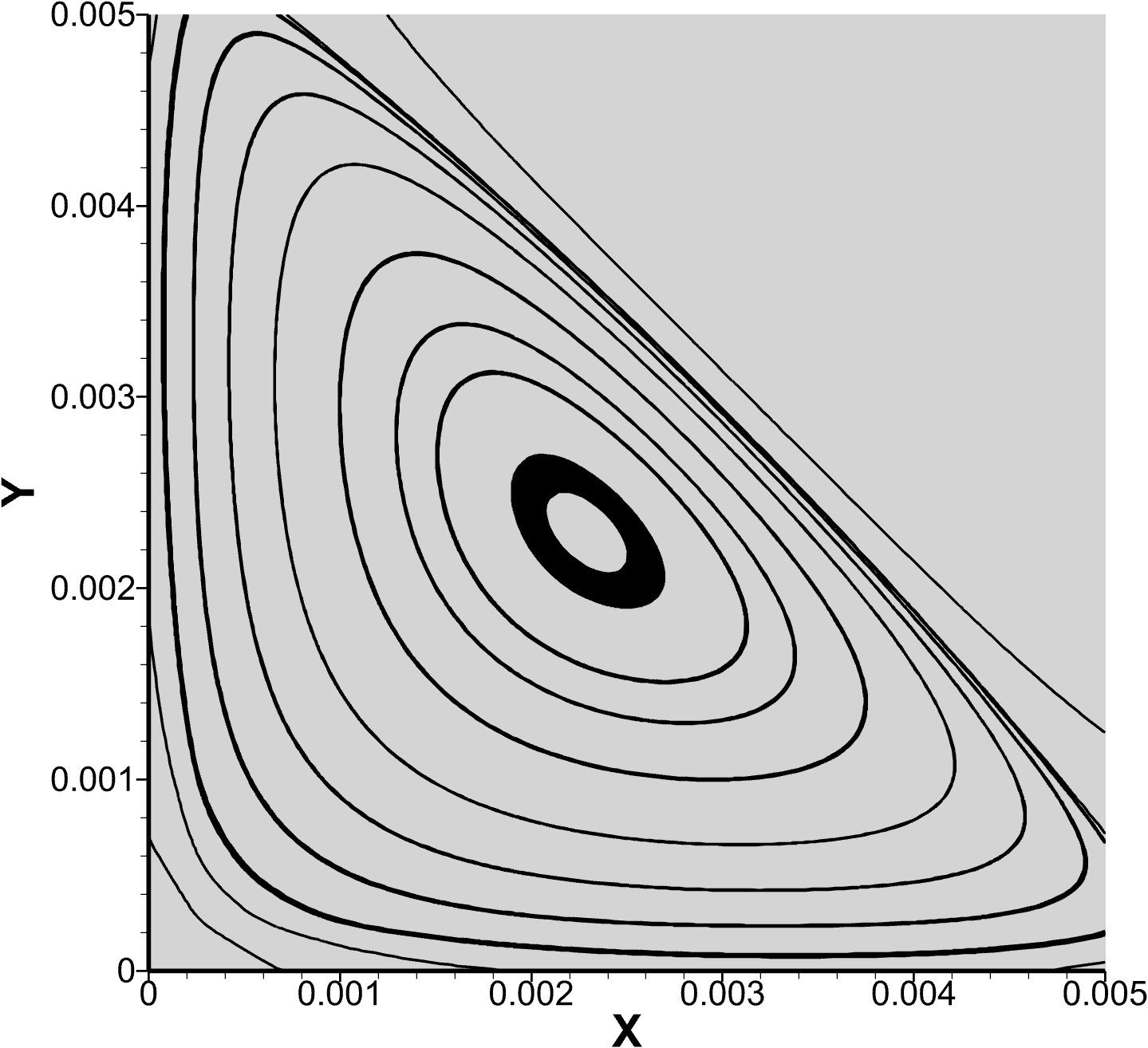}
}
\quad
\subfigure[Eddy $\mathrm{BL}3$]{
\includegraphics[width=4.5cm]{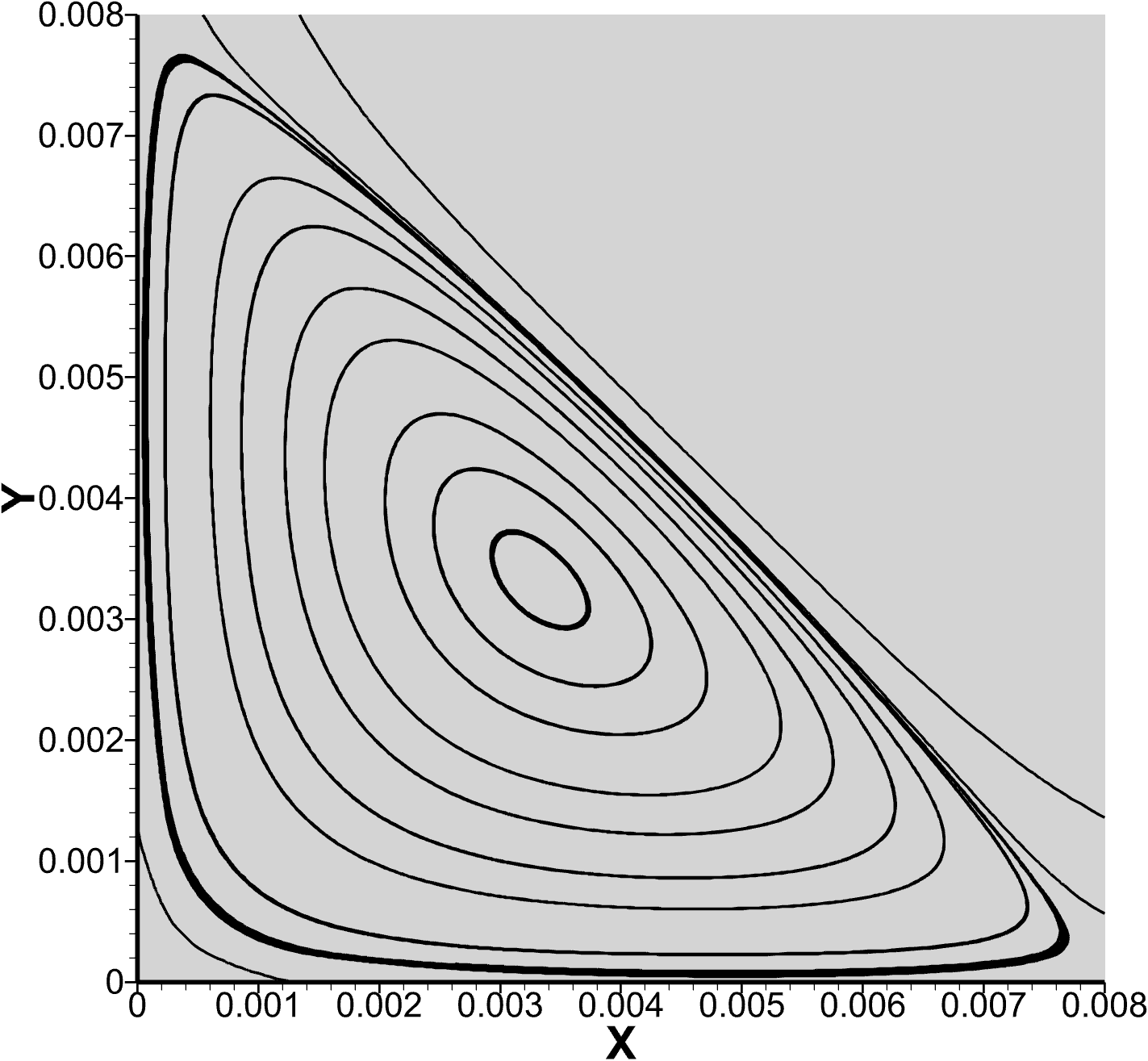}
}
\centering
\caption{\emph{Example 4.} Under the refined mesh, a magnified view of streamline contours at the bottom left corner by $\bm{\mathrm{P}}_3^{\mathrm{bubble}}$-$\mathrm{P}_2^{\mathrm{dc}}$-$\bm{\mathrm{BDM}}_3$ from (a) to (f): $\mathrm{Re}=1000$, $2500$, $5000$, $10000$, $15000$ and $20000$}\label{fig:ex4:mesh80:bl:sl}
\end{figure}

\begin{figure}[htbp]
\centering
\subfigure[Eddies $\mathrm{BR}2$, $\mathrm{BR}3$]{
\includegraphics[width=4.5cm]{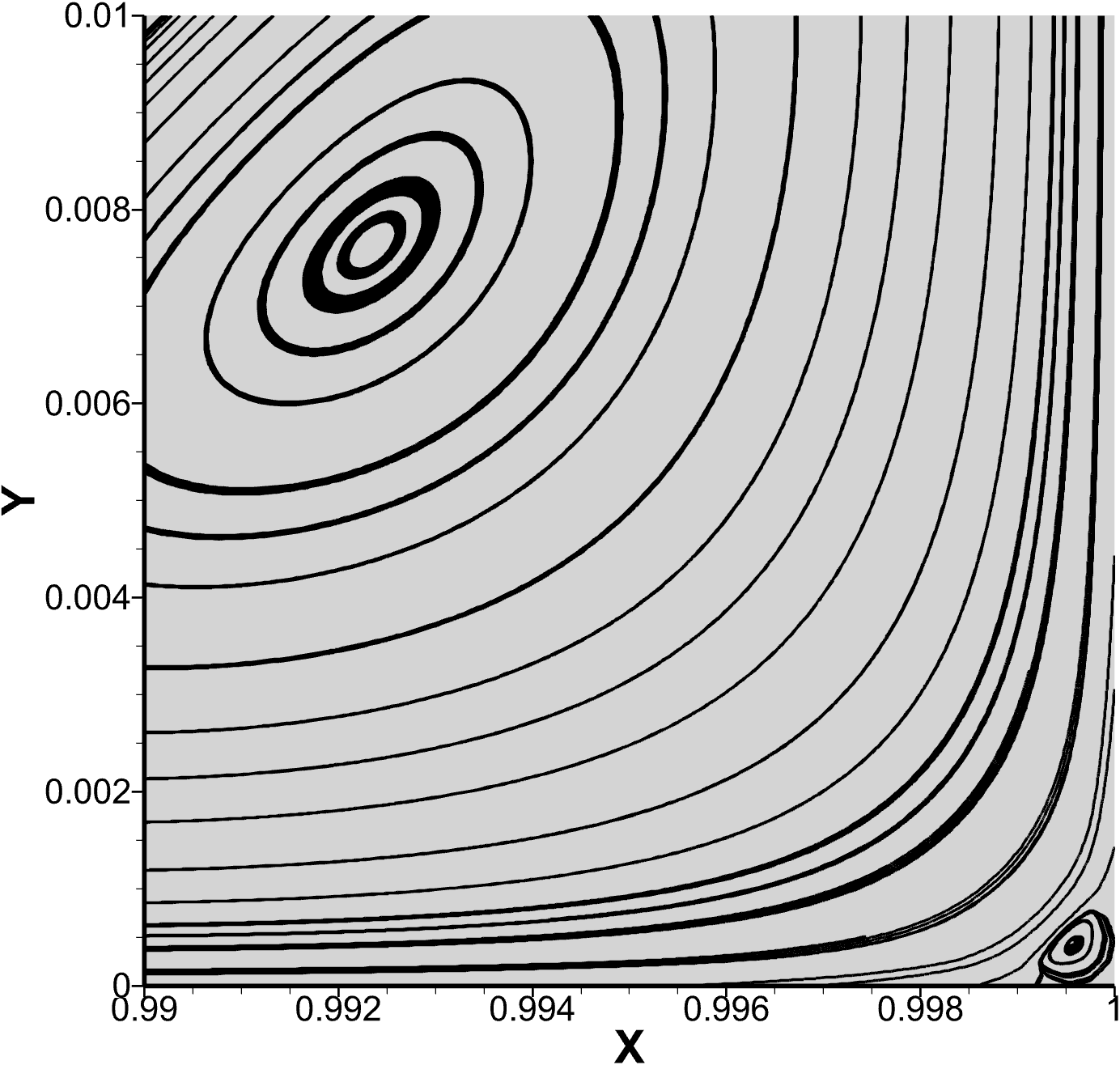}
}
\quad
\subfigure[Eddies $\mathrm{BR}2$, $\mathrm{BR}3$]{
\includegraphics[width=4.5cm]{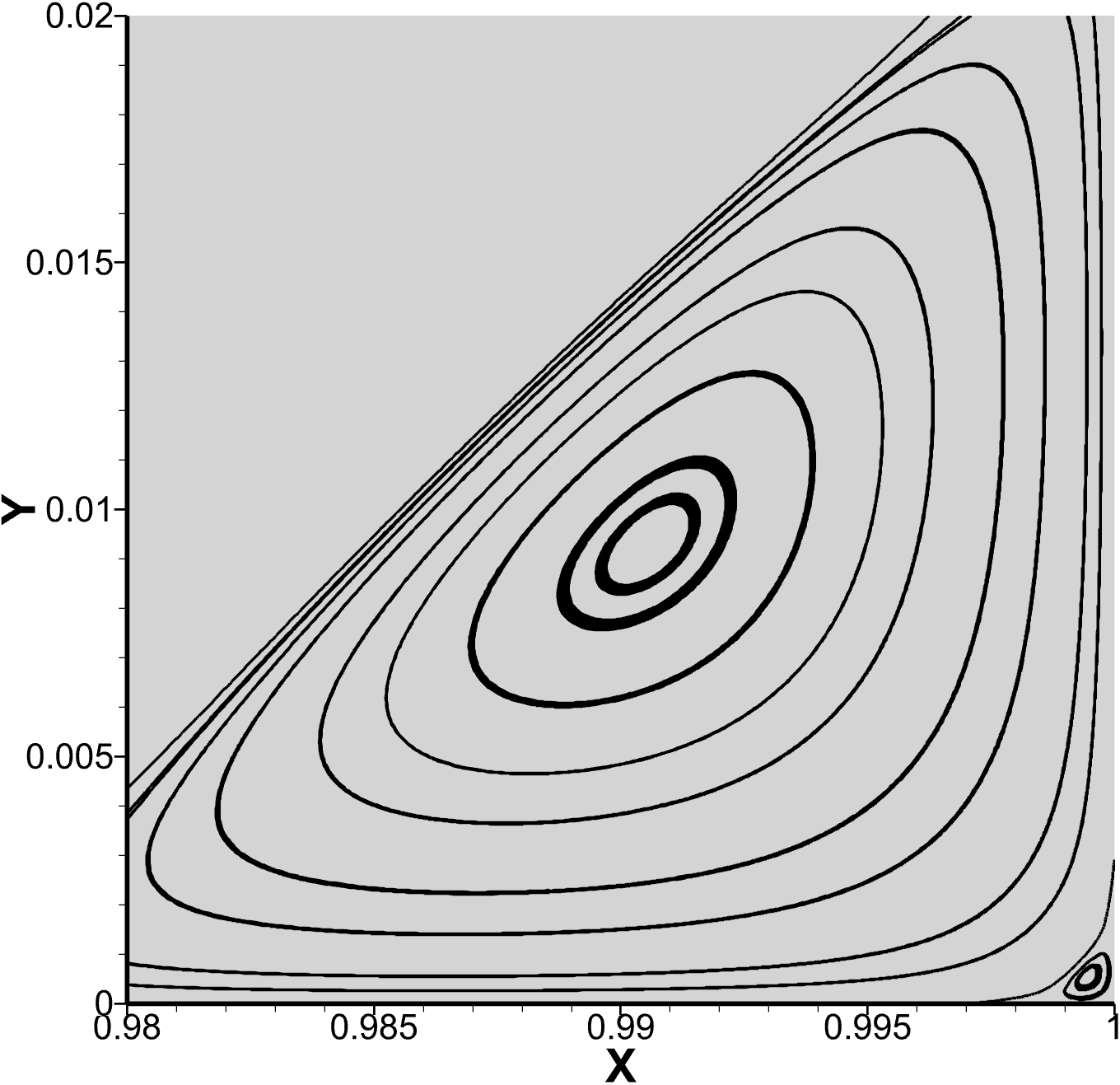}
}
\quad
\subfigure[Eddies $\mathrm{BR}2$, $\mathrm{BR}3$]{
\includegraphics[width=4.5cm]{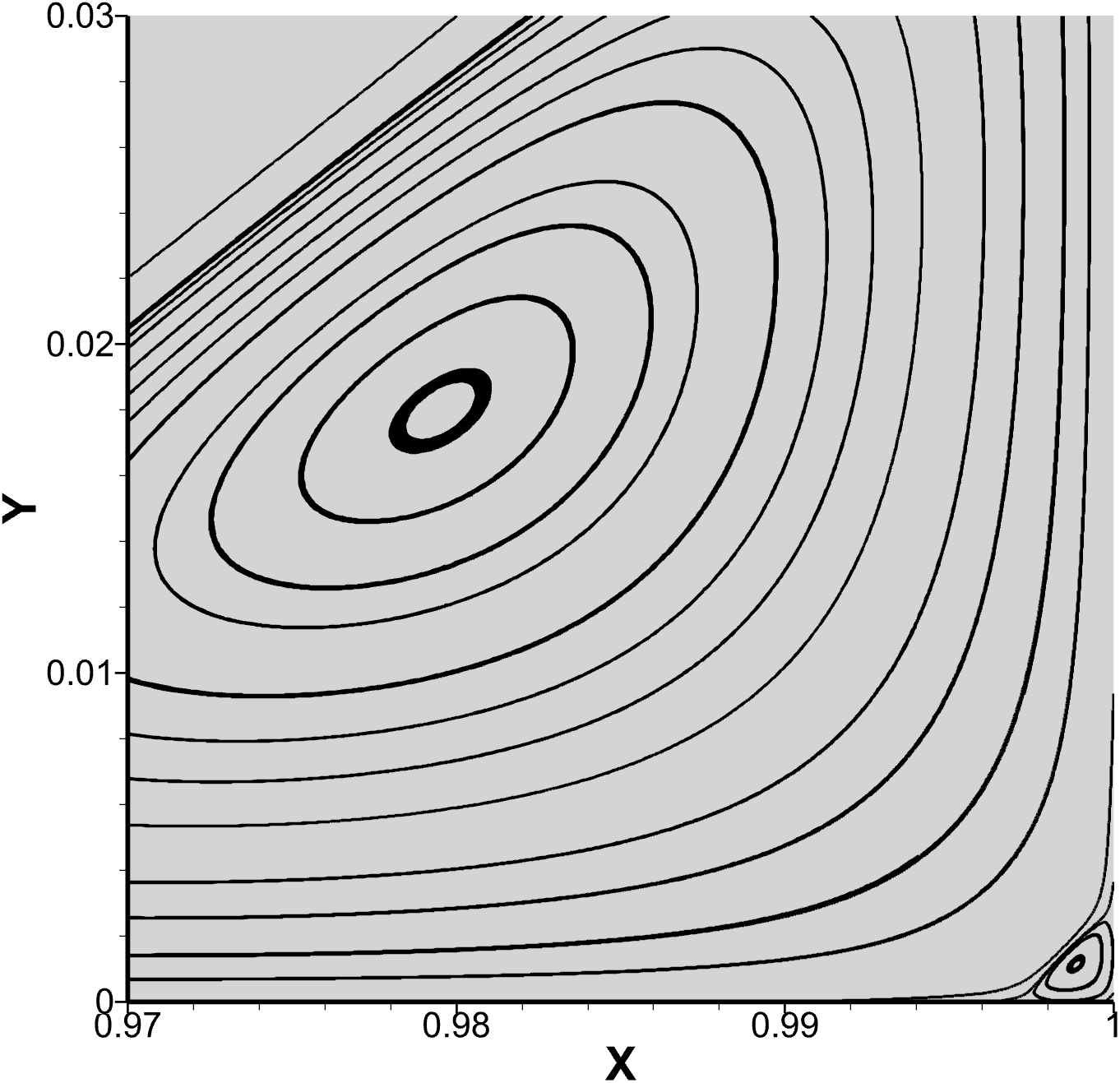}
}

\subfigure[Eddy $\mathrm{BR}3$]{
\includegraphics[width=4.5cm]{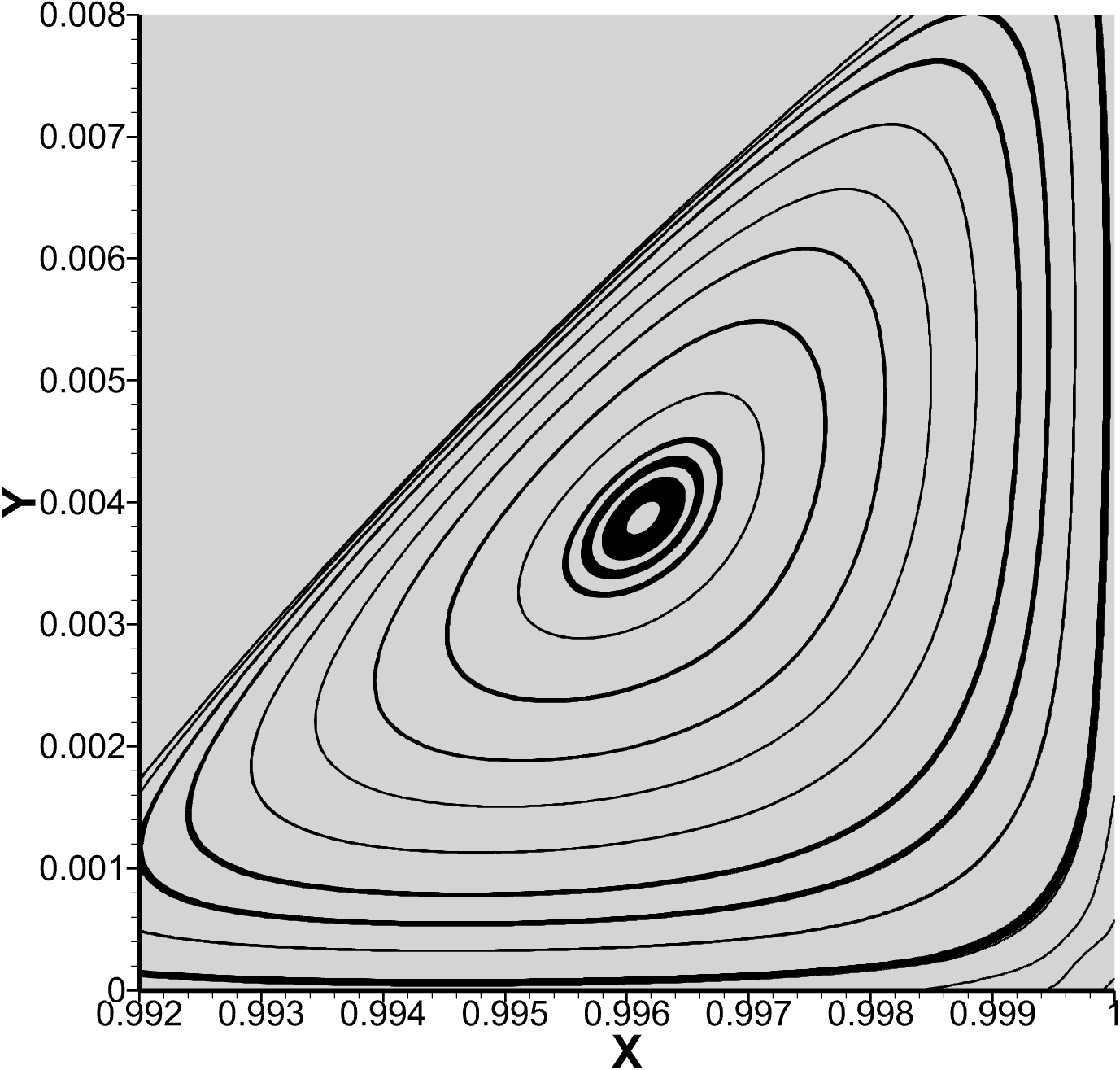}
}
\quad
\subfigure[Eddies $\mathrm{BR}3$, $\mathrm{BR}4$]{
\includegraphics[width=4.5cm]{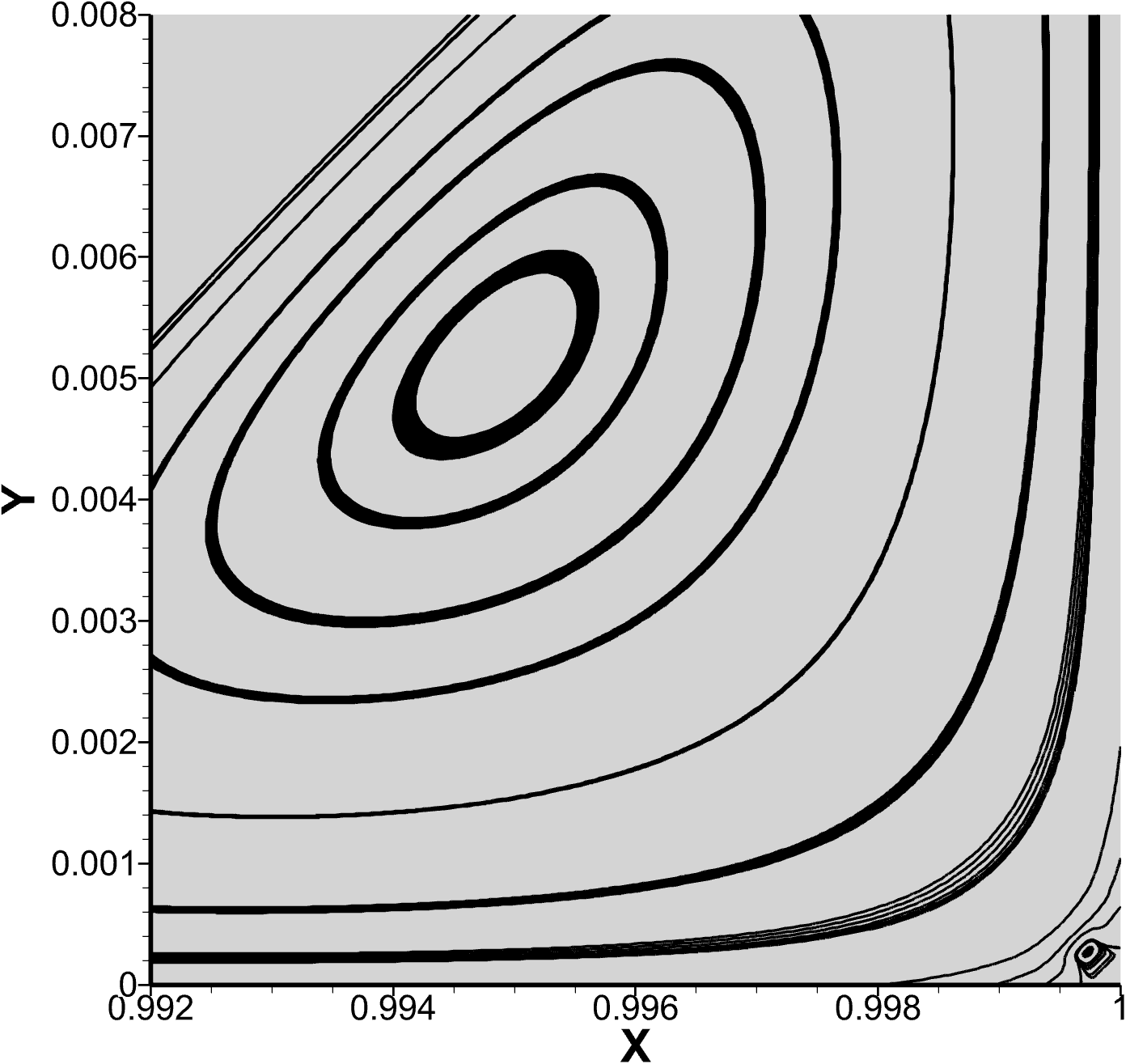}
}
\quad
\subfigure[Eddies $\mathrm{BR}3$, $\mathrm{BR}4$]{
\includegraphics[width=4.5cm]{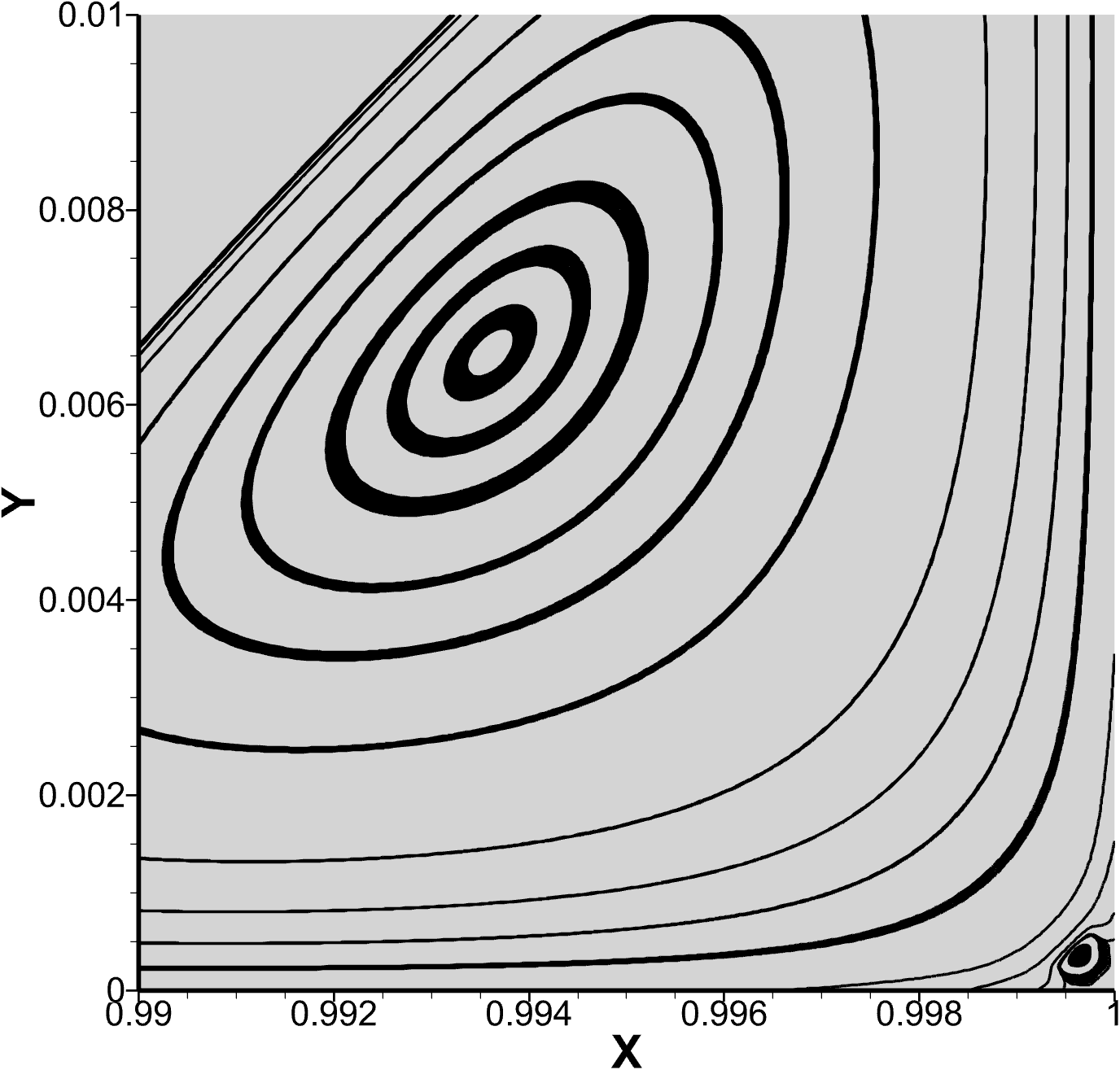}
}
\centering
\caption{\emph{Example 4.} Under the refined mesh, a magnified view of streamline contours at the bottom right corner by $\bm{\mathrm{P}}_3^{\mathrm{bubble}}$-$\mathrm{P}_2^{\mathrm{dc}}$-$\bm{\mathrm{BDM}}_3$ from (a) to (f): $\mathrm{Re}=1000$, $2500$, $5000$, $10000$, $15000$ and $20000$}\label{fig:ex4:mesh80:br:sl}
\end{figure}

Distinctively than many published papers, it can be clearly observed from Figures \ref{fig:ex4:mesh80:bl:sl} and \ref{fig:ex4:mesh80:br:sl} that our computations by $\bm{\mathrm{P}}_3^{\mathrm{bubble}}$-$\mathrm{P}_2^{\mathrm{dc}}$-$\bm{\mathrm{BDM}}_3$ under the refined mesh indicate the appearance of a quaternary vortex at the bottom left corner (BL3) for $\mathrm{Re}=1000$, $2500$ and $5000$; the appearance of a quaternary vortex at the bottom left corner (BR3) for $\mathrm{Re}=1000$ and $2500$; the appearance of a fifth-level vortex at the bottom left corner (BR4) for $\mathrm{Re}=15000$ and $20000$. It means the achievement of our purpose to resolve sub-grid scale vortices with considerable computations.

Then, we shall check some properties of primary vortices at various Reynolds numbers, such as locations of their centers and the values of streamfunction $\phi_h$ at their centers. The streamfunction $\phi_h\in S_h\subset H^1(\Omega)$ is obtained from the discrete velocity solution $\bm{u}_h$ a posteriori by a Poisson equation that

\[
(\mathbf{curl}\,\phi_h, \mathbf{curl}\,\psi_h)=(\bm{u}_h,\mathbf{curl}\,\psi_h),\quad\forall\,\psi_h\in S_h,
\]
where the $H^1$-conforming finite element space $S_h$ is chosen as the fifth-order Lagrange finite element space. Here we also take the corresponding results by a high-order accurate finite difference scheme under a fine mesh $601\times 601$ in \cite{Erturk2005747} as a reference. In Table \ref{tab:ex4:primary:vortex}, for any given Reynolds number $\mathrm{Re}\in\{1000,2500,5000,10000,15000,20000\}$, $(x_1,x_2)$ denotes the location $xy$-coordinates of the primary vortex center and $\phi$ denotes the streamfunction value at the point $(x_1,x_2)$, while $(x_{1,\mathrm{ref}},x_{2,\mathrm{ref}})$ and $\phi_{\mathrm{ref}}$ denote the reference results with the same meanings respectively. Our computations are highly comparable with the reference solutions. Moreover, from not only the qualitative information hidden in Figure \ref{fig:ex4:mesh80:global:sl} but also the quantitative data shown in Table \ref{tab:ex4:primary:vortex}, it follows that the primary vortex center moves towards the cavity's geometric center as the Reynolds number increases, which was also reported by Barragy and Carey \cite{Barragy1997453} and Ghia et al. \cite{Ghia1982387}.

\begin{table}[htbp]
    \centering
    \topcaption{\emph{Example 4}. Under the refined mesh, properties of primary vortices at various Reynolds numbers; streamfunction values and center locations. Comparison between Erturk et al. \cite{Erturk2005747} and $\bm{\mathrm{P}}_3^{\mathrm{bubble}}$-$\mathrm{P}_2^{\mathrm{dc}}$-$\bm{\mathrm{BDM}}_3$}
    \begin{tabular*}{\hsize}{@{}@{\extracolsep{\fill}}lcccccc@{}}
        \toprule
        $\mathrm{Re}$                              &  $1000$     &  $2500$     &  $5000$     & $10000$     & $15000$     & $20000$     \\
        \midrule
        $\phi$                                     & $-0.117697$ & $-0.119444$ & $-0.119360$ & $-0.118355$ & $-0.117364$ & $-0.116484$ \\
        $\phi_{\mathrm{ref}}$ \cite{Erturk2005747} & $-0.118781$ & $-0.121035$ & $-0.121289$ & $-0.120403$ & $-0.119240$ & $-0.118039$ \\
        \hline
        $x_1$                                      & $0.5308$    & $0.5201$    & $0.5154$    & $0.5122$    & $0.5108$    & $0.5100$    \\
        $x_{1,\mathrm{ref}}$ \cite{Erturk2005747}  & $0.5300$    & $0.5200$    & $0.5150$    & $0.5117$    & $0.5100$    & $0.5100$    \\
        \hline
        $x_2$                                      & $0.5653$    & $0.5437$    & $0.5345$    & $0.5297$    & $0.5277$    & $0.5264$    \\
        $x_{2,\mathrm{ref}}$ \cite{Erturk2005747}  & $0.5650$    & $0.5433$    & $0.5350$    & $0.5300$    & $0.5283$    & $0.5267$    \\
        \bottomrule
    \end{tabular*}
    \label{tab:ex4:primary:vortex}
\end{table}

\begin{figure}[htbp]
\centering
\subfigure{
\includegraphics[width=5cm]{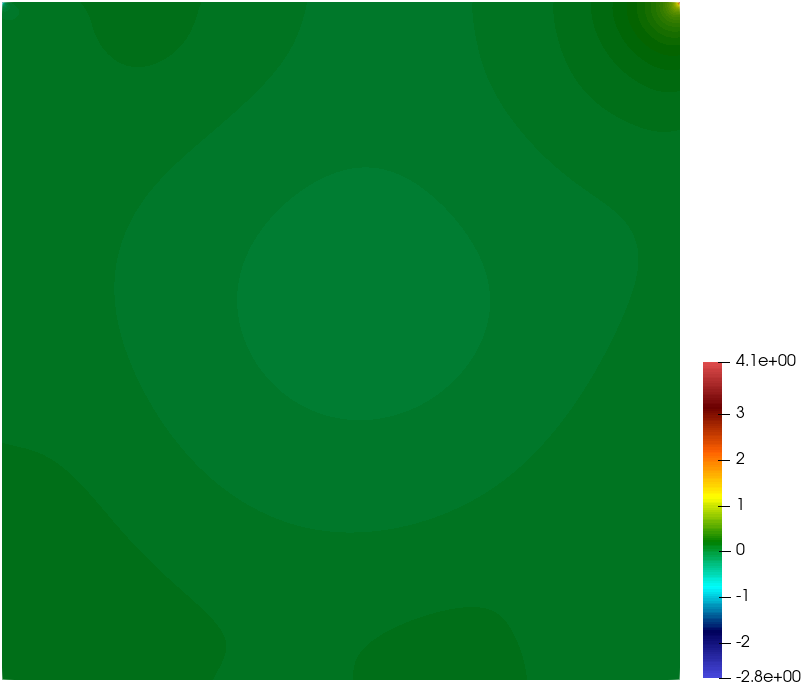}
}
\quad
\subfigure{
\includegraphics[width=5cm]{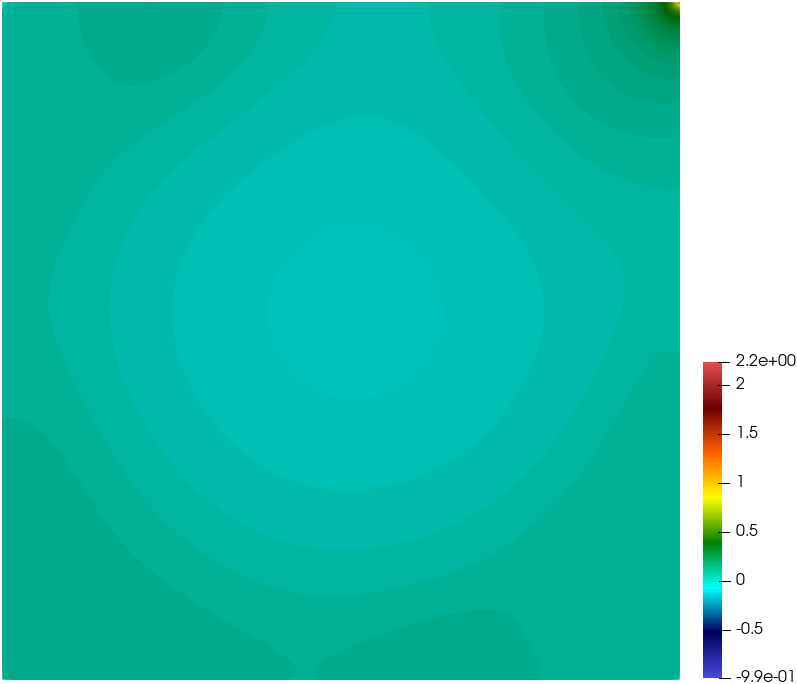}
}
\quad
\subfigure{
\includegraphics[width=5cm]{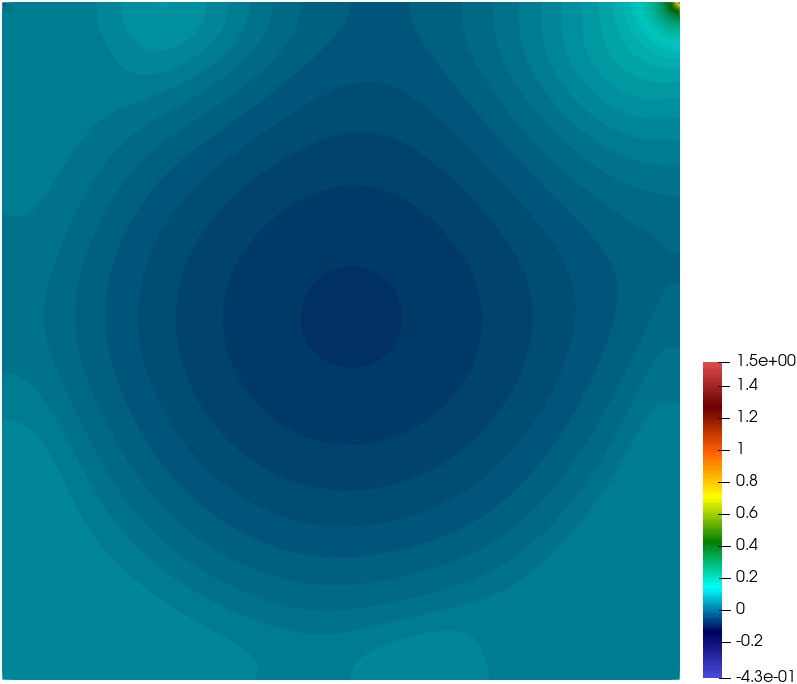}
}

\subfigure{
\includegraphics[width=5cm]{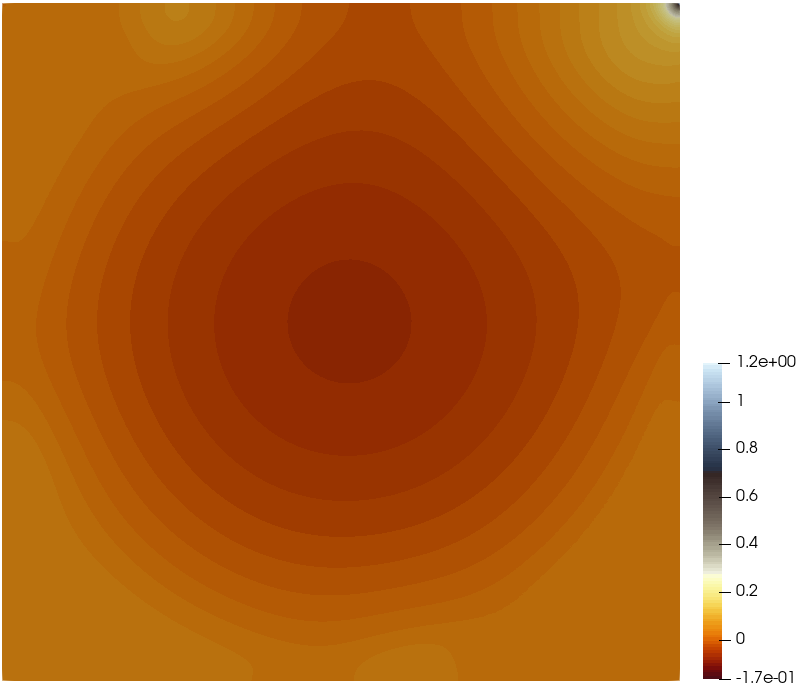}
}
\quad
\subfigure{
\includegraphics[width=5cm]{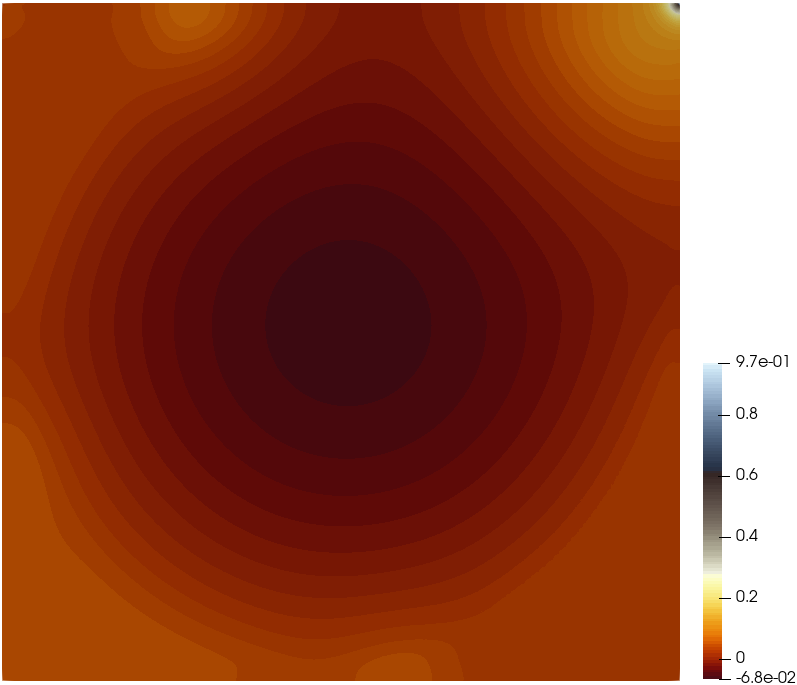}
}
\quad
\subfigure{
\includegraphics[width=5cm]{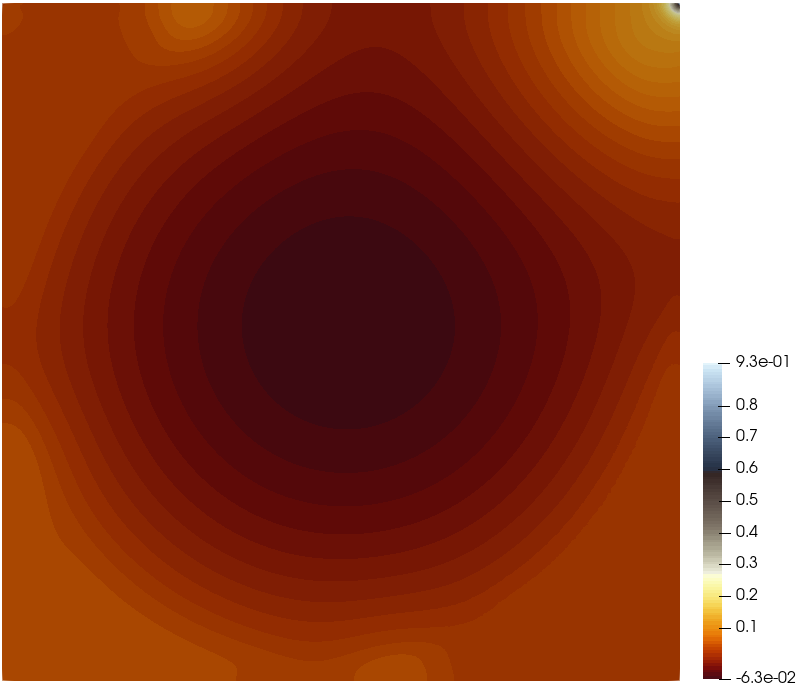}
}
\centering
\caption{\emph{Example 4.} Under the refined mesh, the kinetic pressure distributions by $\bm{\mathrm{P}}_3^{\mathrm{bubble}}$-$\mathrm{P}_2^{\mathrm{dc}}$-$\bm{\mathrm{BDM}}_3$ for $\mathrm{Re}=1000$, $2500$, $5000$, $10000$, $15000$ and $20000$ from top-left to bottom-right}\label{fig:ex4:mesh80:pres}
\end{figure}

Apart from the above many computation results about the velocity, the kinematic pressure distributions at various Reynolds numbers, which can be obtained from \eqref{kin:pres:posteriori}, are displayed in Figure \ref{fig:ex4:mesh80:pres}. By comparison with the authorized results in the literature, in the aspect of contours, the computed kinematic pressure is in good agreement with that of Botella and Peyert \cite{Botella1998421} at $\mathrm{Re}=1000$ and with that of Hachem et al. \cite{Hachem20108643} at $\mathrm{Re}=10000$; in the aspect of value ranges, the computed kinematic pressure almost matches that of Gravemeier et al. \cite{Gravemeier20041323} at $\mathrm{Re}=10000$.

\begin{figure}[htbp]
\centering
\subfigure{
\includegraphics[width=5.5cm]{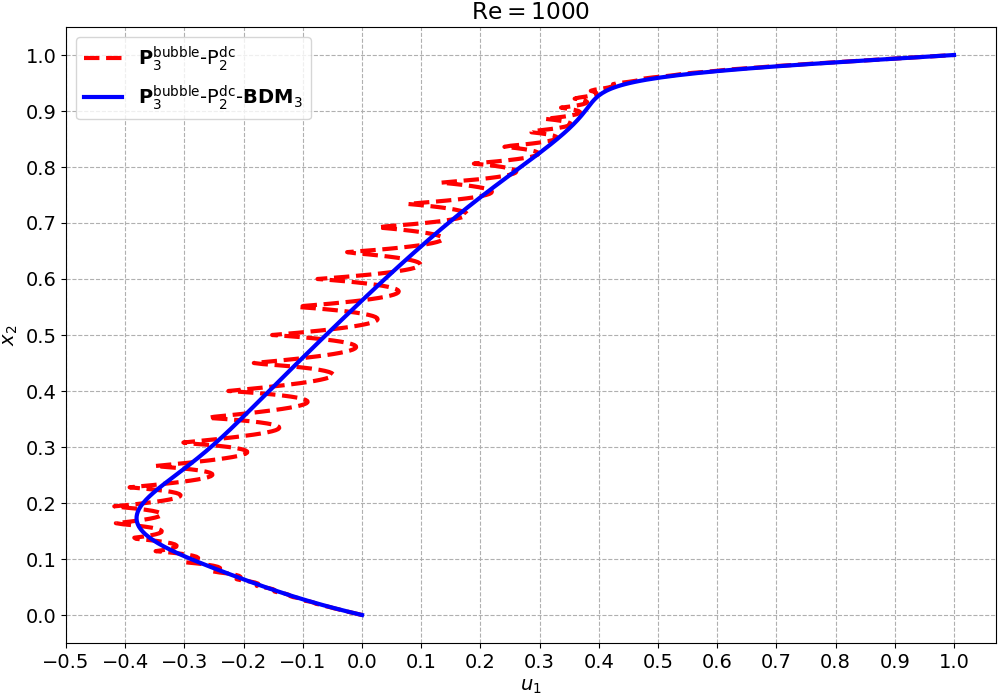}
}
\quad
\subfigure{
\includegraphics[width=5.5cm]{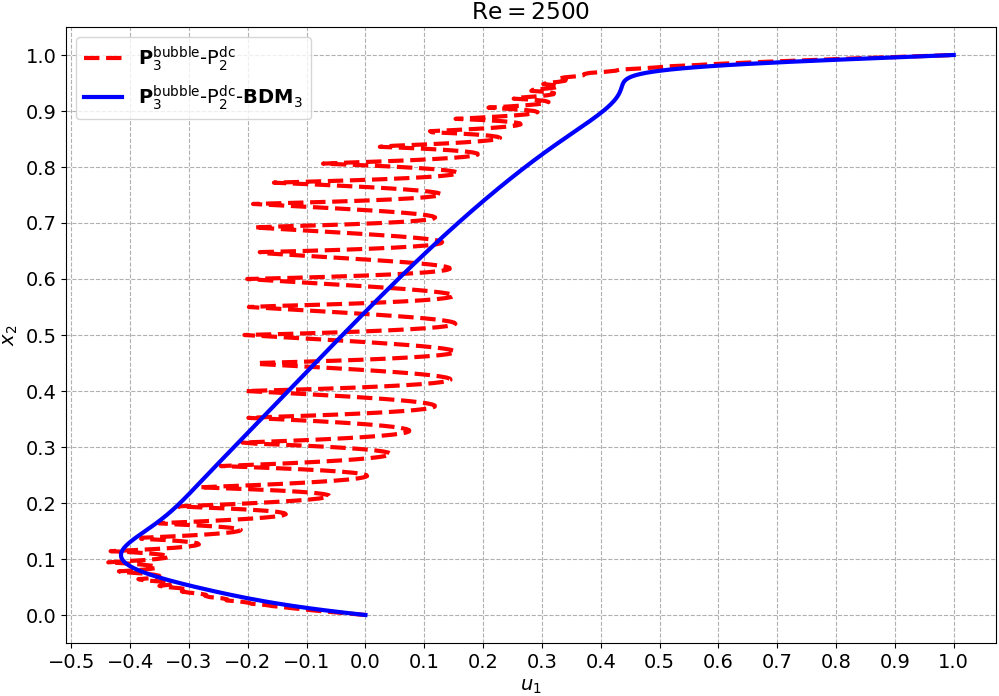}
}
\quad
\subfigure{
\includegraphics[width=5.5cm]{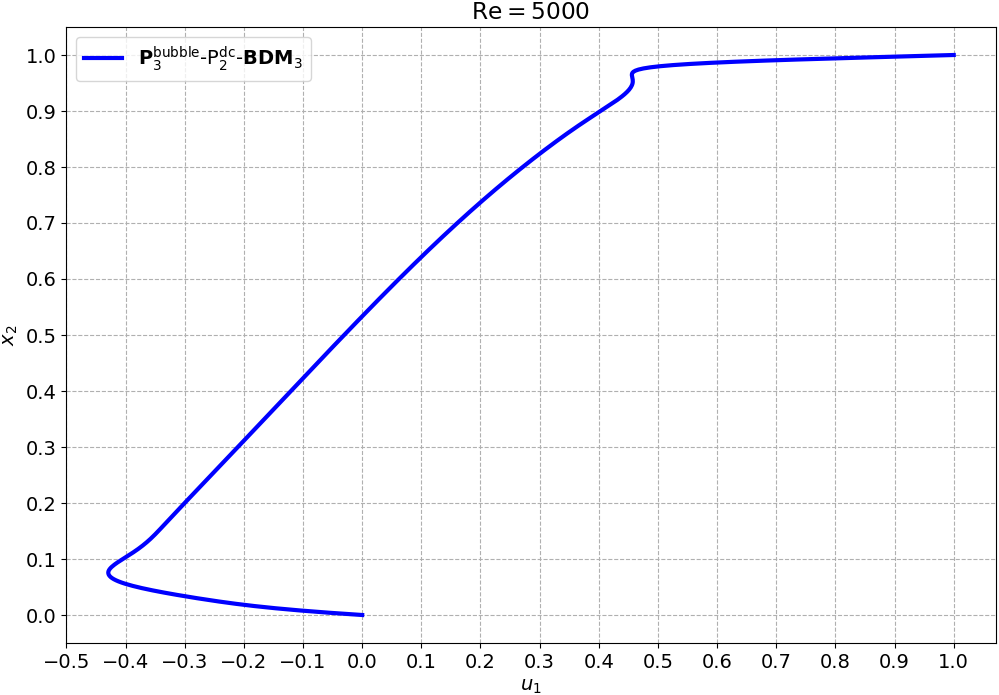}
}

\subfigure{
\includegraphics[width=5.5cm]{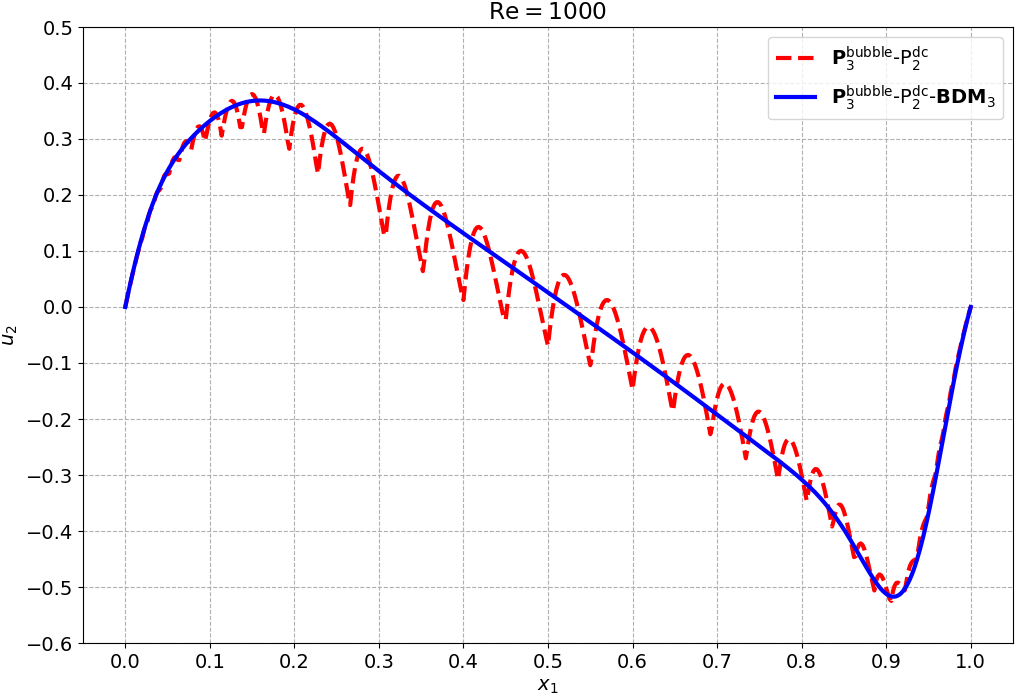}
}
\quad
\subfigure{
\includegraphics[width=5.5cm]{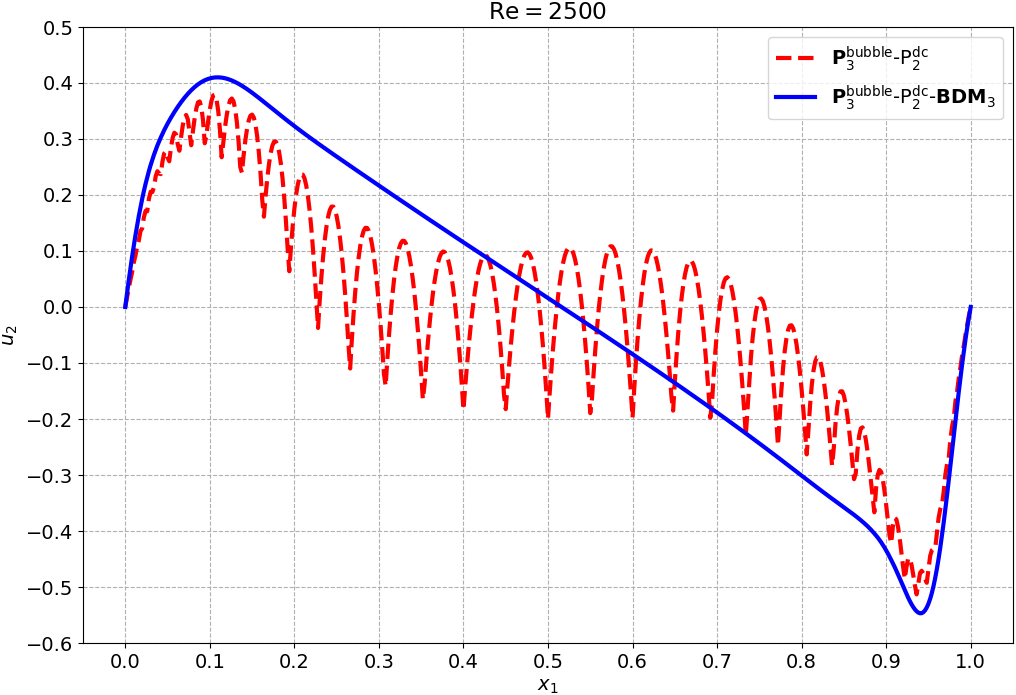}
}
\quad
\subfigure{
\includegraphics[width=5.5cm]{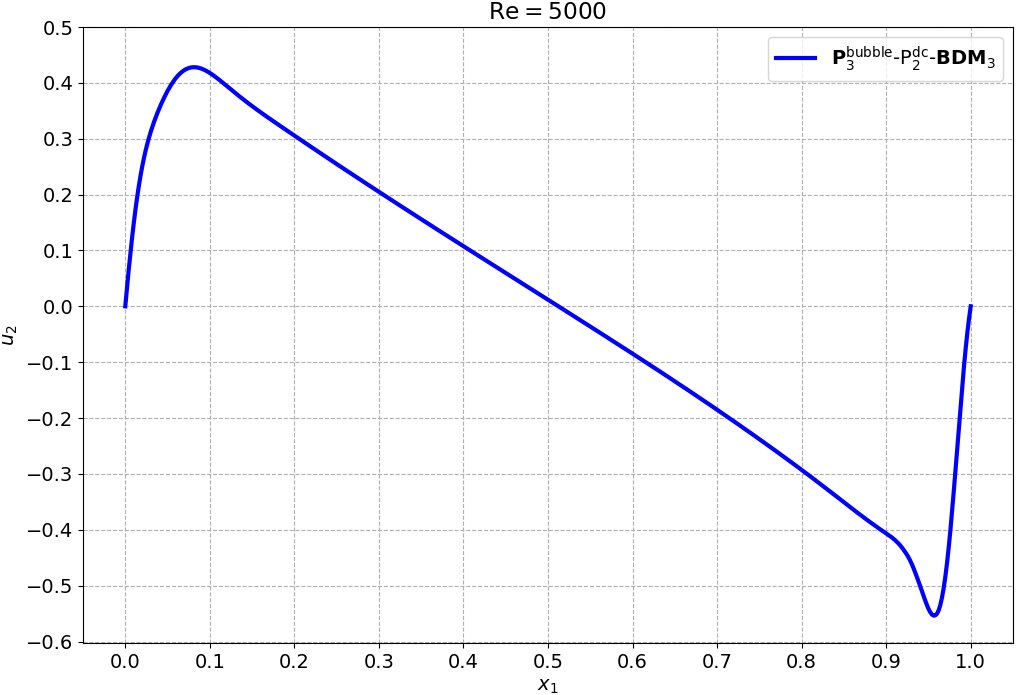}
}
\centering
\caption{\emph{Example 4.} Under the mesh with $M=50$ and $\gamma=2.5$, $u_1$-velocity profiles along the vertical centreline $x_1=0.5$ (upper row) and $u_2$-velocity profiles along the horizontal centreline $x_2=0.5$ (lower row) for the large irrotational body force $\bm{f}=10^5\nabla\psi$. Comparison between $\bm{\mathrm{P}}_3^{\mathrm{bubble}}$-$\mathrm{P}_2^{\mathrm{dc}}$ and $\bm{\mathrm{P}}_3^{\mathrm{bubble}}$-$\mathrm{P}_2^{\mathrm{dc}}$-$\bm{\mathrm{BDM}}_3$ for $\mathrm{Re}=1000$, $2500$ and $5000$ from left to right ($\bm{\mathrm{P}}_3^{\mathrm{bubble}}$-$\mathrm{P}_2^{\mathrm{dc}}$ does not converge at $\mathrm{Re}=5000$)}
\label{fig:ex4:mesh50:P3:compare}
\end{figure}

\begin{figure}[htbp]
\centering
\subfigure{
\includegraphics[width=5cm]{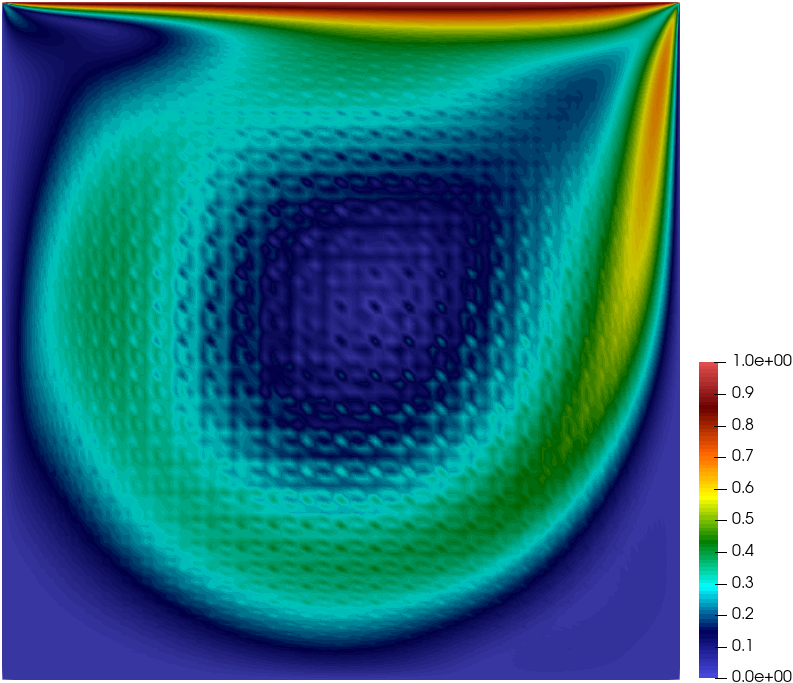}
}
\quad
\subfigure{
\includegraphics[width=5cm]{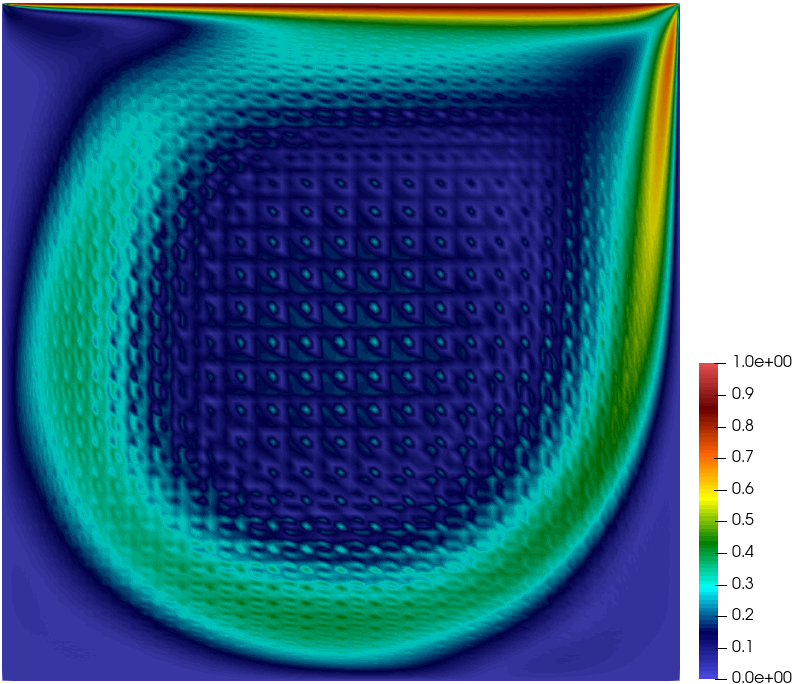}
}
\quad
\subfigure{
\includegraphics[width=5cm]{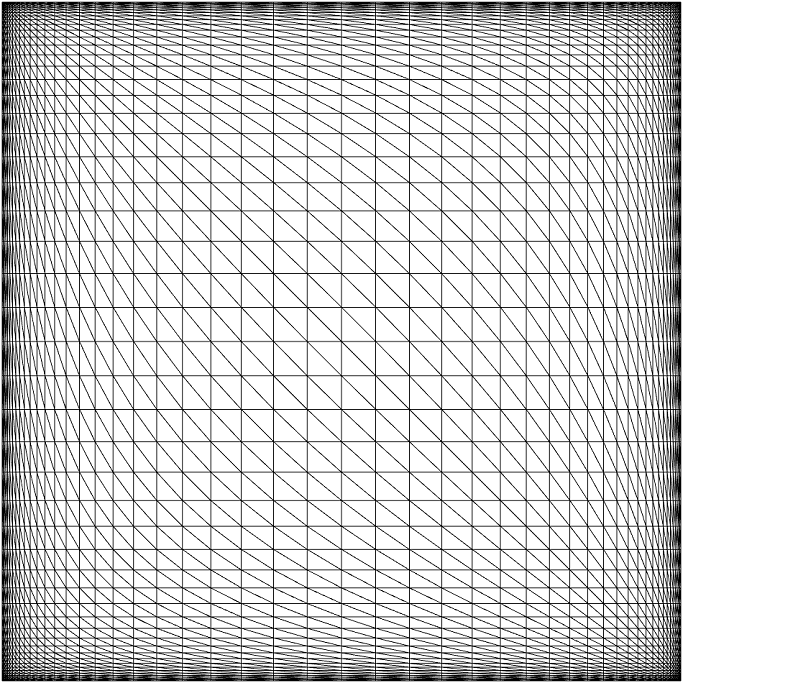}
}

\subfigure{
\includegraphics[width=5cm]{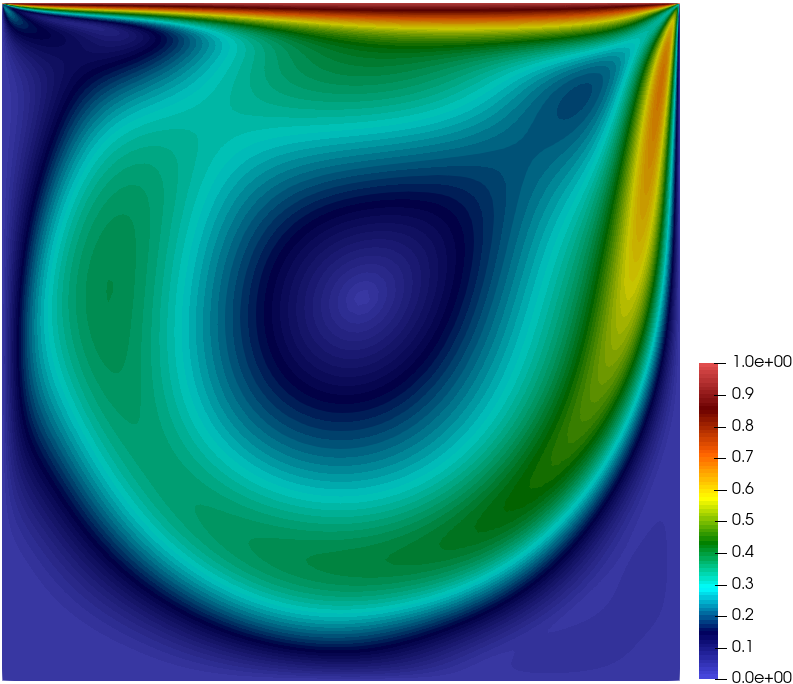}
}
\quad
\subfigure{
\includegraphics[width=5cm]{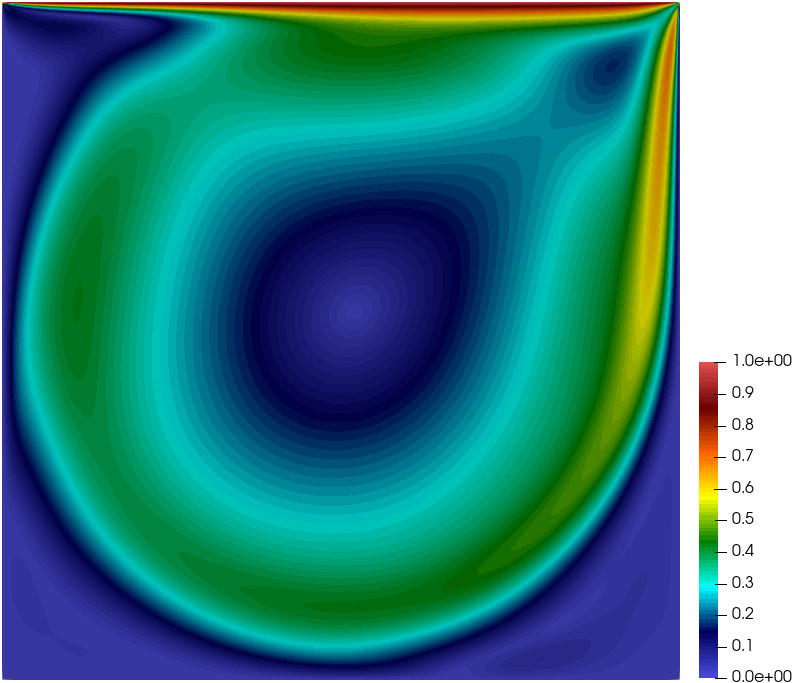}
}
\quad
\subfigure{
\includegraphics[width=5cm]{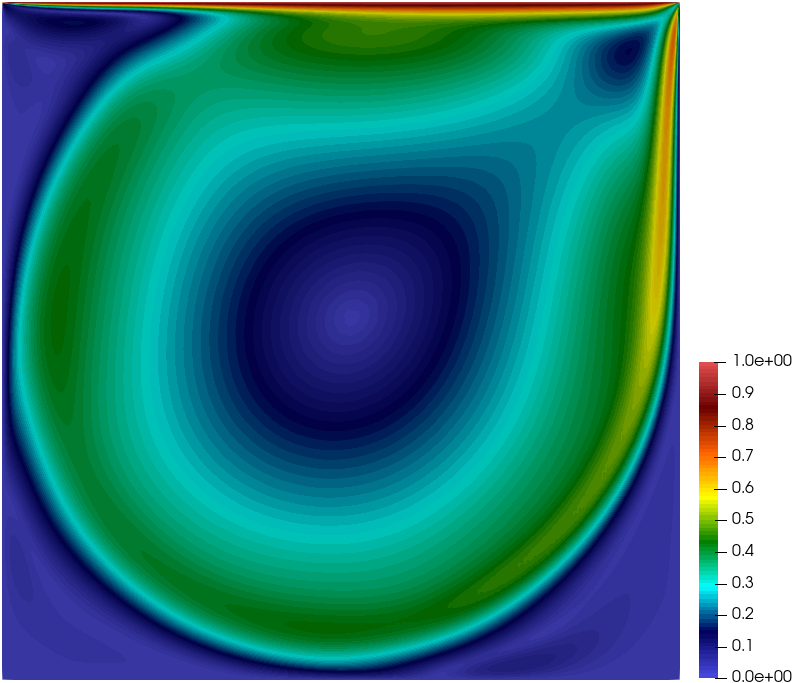}
}
\centering
\caption{\emph{Example 4.} Under the mesh with $M=50$ and $\gamma=2.5$, contours of velocity magnitude for the large irrotational body force $\bm{f}=10^5\nabla\psi$. Comparison between $\bm{\mathrm{P}}_3^{\mathrm{bubble}}$-$\mathrm{P}_2^{\mathrm{dc}}$ (upper row) and $\bm{\mathrm{P}}_3^{\mathrm{bubble}}$-$\mathrm{P}_2^{\mathrm{dc}}$-$\bm{\mathrm{BDM}}_3$ (lower row) for $\mathrm{Re}=1000$, $2500$ and $5000$ from left to right ($\bm{\mathrm{P}}_3^{\mathrm{bubble}}$-$\mathrm{P}_2^{\mathrm{dc}}$ does not converge at $\mathrm{Re}=5000$)}
\label{fig:ex4:mesh50:P3:compare:2}
\end{figure}

Let us end this example by checking the robustness for irrotational body forces. To this end, the same test case will be implemented but with the large right-hand body force

\[
\bm{f}(\bm{x})=10^5\,\nabla\psi,\quad \psi=\frac{1}{3}(x_1^3+x_2^3).
\]
This body force is exactly irrotational, and hence the velocity numerical solutions should not be affected in theory. To verify it, we choose the triangular mesh with $M=50$, $\gamma=2.5$ and the polynomial degree $k=3$, and then compare the performance between $\bm{\mathrm{P}}_3^{\mathrm{bubble}}$-$\mathrm{P}_2^{\mathrm{dc}}$ and $\bm{\mathrm{P}}_3^{\mathrm{bubble}}$-$\mathrm{P}_2^{\mathrm{dc}}$-$\bm{\mathrm{BDM}}_3$ at various Reynolds numbers $\mathrm{Re}=1000$, $2500$ and $5000$ respectively. The striking differences are distinguished in Figure \ref{fig:ex4:mesh50:P3:compare} that both the $u_1$-velocity profiles along a vertical line $x_1=0.5$ and the $u_2$-velocity profiles along a horizontal line $x_2=0.5$ are not affected by the large body force for $\bm{\mathrm{P}}_3^{\mathrm{bubble}}$-$\mathrm{P}_2^{\mathrm{dc}}$-$\bm{\mathrm{BDM}}_3$, while the violent pseudo numerical oscillations occur for $\bm{\mathrm{P}}_3^{\mathrm{bubble}}$-$\mathrm{P}_2^{\mathrm{dc}}$. It is evident that as the Reynolds number increases, the pseudo oscillations by $\bm{\mathrm{P}}_3^{\mathrm{bubble}}$-$\mathrm{P}_2^{\mathrm{dc}}$ are more severe. The same phenomenon is also reported in Figure \ref{fig:ex4:mesh50:P3:compare:2} with respect to the contours of velocity magnitude.

\subsection{Example 5: Laminar flow around a cylinder}

In this validation experiment, we consider the famous Laminar flow past a cylinder \cite{Schafer1996547}, which is also called the Sch\"{a}fer-Turek benchmark. It is popular with the literature about the time-dependent Navier--Stokes equations \cite{Fiordilino2018327, John2004777, Linke2013142, Rong202066}. Our purpose of solving the stationary equations is to investigate the appearance and evolution of the symmetric vortices behind the cylinder as the Reynolds number increases.

The geometry is a 2D channel with a circular obstacle which is positioned (only slightly) off the center of the channel; specifically, the computational domain is set as $\Omega=(0, 2.2)\times (0,0.41)\setminus B_{0.05}(0.2,0.2)$ where $B_r(x_1,x_2)$ denotes a ball with the center $(x_1,x_2)$ and the radius $r$. We choose the body force $\bm{f}=\bm{0}$ and no-slip boundary conditions imposed at all walls. For the inlet and outlet boundary conditions, they are respectively given by
$\bm{u}=(1,0)^\top$ on the inlet and $(p^{\mathrm{kin}}\mathbb{I}-\nu\nabla\bm{u})\bm{n}=\bm{0}$ on the outlet.

\begin{figure}[htbp]
\centering
\includegraphics[width=0.6\textwidth]{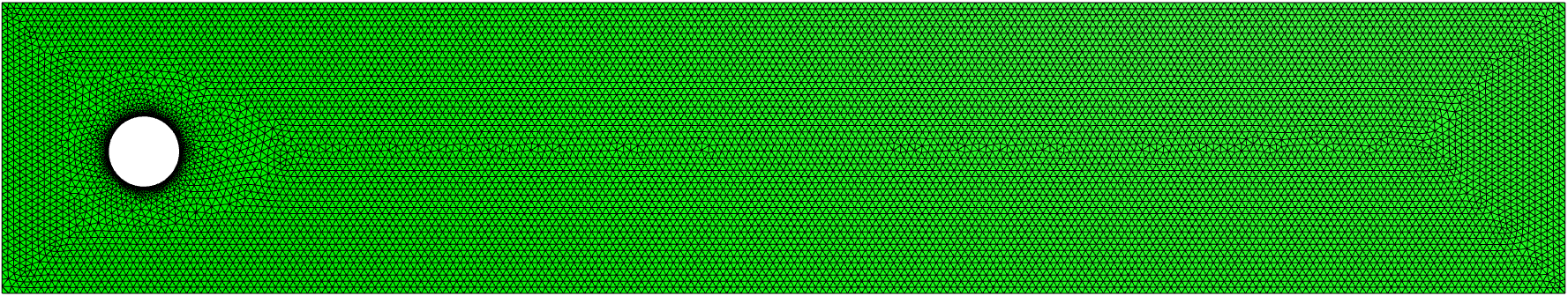}
\caption{\emph{Example 5.} The global triangular mesh plot including the refined grids near the cylinder}\label{fig:ex5:meshplot}
\end{figure}

We consider computations by $\bm{\mathrm{P}}_3^{\mathrm{bubble}}$-$\mathrm{P}_2^{\mathrm{dc}}$-$\bm{\mathrm{BDM}}_3$ performed on a locally refined triangular mesh presented in Figure \ref{fig:ex5:meshplot}. For clarity we denote by $h_l$ the mesh size of the locally refined interface or region and $h_g$ the mesh size elsewhere hereafter. Here we take $h_l=1/500$ and $h_g=1/100$. Recalling Remark \ref{rem:comp:bcs}, with an additional integration on the outlet boundary to \eqref{dis:rot:NS:momentum}, the corresponding Newton iteration should be rewritten as: to find $(\bm{u}_h^{n+1},p_h^{n+1})\in\widetilde{\bm{X}}_h\times Q_h$ such that for all $(\bm{v}_h,q_h)\in\widetilde{\bm{X}}_h\times Q_h$,

\begin{equation}\label{linear:Newton:cylinder}
\begin{split}
&\nu a(\bm{u}_h^{n+1},\bm{v}_h)+b_h(\bm{u}_h^n;\bm{u}_h^{n+1},\bm{v}_h)
+b_h(\bm{u}_h^{n+1};\bm{u}_h^n,\bm{v}_h)+c(\bm{u}_h^n,\bm{u}_h^{n+1},\bm{v}_h)+c(\bm{u}_h^{n+1},\bm{u}_h^n,\bm{v}_h)\\
&+d(\bm{v}_h,p_h^{n+1})+d(\bm{u}_h^{n+1},q_h)
=(\bm{f},\mathcal R(\bm{v}_h))+b_h(\bm{u}_h^n;\bm{u}_h^n,\bm{v}_h)+c(\bm{u}_h^n,\bm{u}_h^n,\bm{v}_h),
\end{split}
\end{equation}
where $c(\bm{w},\bm{z},\bm{v})$, the additional integration on the outlet boundary $\Gamma_N$ (other boundaries are denoted as $\Gamma_D$), is defined as

\[
c(\bm{w},\bm{z},\bm{v}):=\frac{1}{2}\int_{\Gamma_N}(\bm{w}\cdot\bm{z})(\bm{v}\cdot\bm{n})\,\mathrm{d}s,
\]
and the modified finite element space $\widetilde{\bm{X}}_h$ for velocity represents $[\mathcal L_k^1\oplus B_{k+1}]^2\cap\bm{H}_{0,D}^1(\Omega)$ for any integer $k\geqslant 2$. Note that the value of the Reynolds number is figured out by $\mathrm{Re}=\frac{2\times 0.05\times 1}{\nu}=\frac{1}{10\nu}$ here. By adopting the updated Newton iteration \eqref{linear:Newton:cylinder}, the problem is first solved for $\mathrm{Re}=5$, then $\mathrm{Re}=10$, $20$, $40$, $70$, $100$, and then in steps of $25$ until $\mathrm{Re}=200$, with the solution for the previous value of $\mathrm{Re}$ used as initial guess for the next; the Stokes equations are solved to provide the initial guess used at $\mathrm{Re}=5$.

\begin{figure}[htbp]
\centering
\subfigure{
\includegraphics[width=0.45\textwidth]{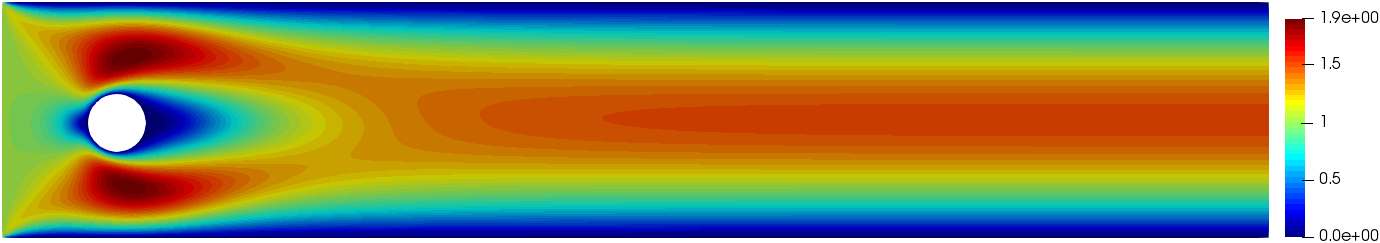}
}
\quad
\subfigure{
\includegraphics[width=0.45\textwidth]{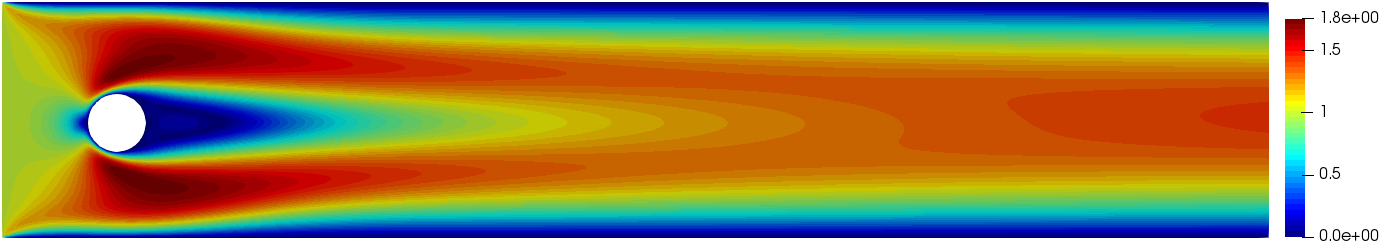}
}

\subfigure{
\includegraphics[width=0.45\textwidth]{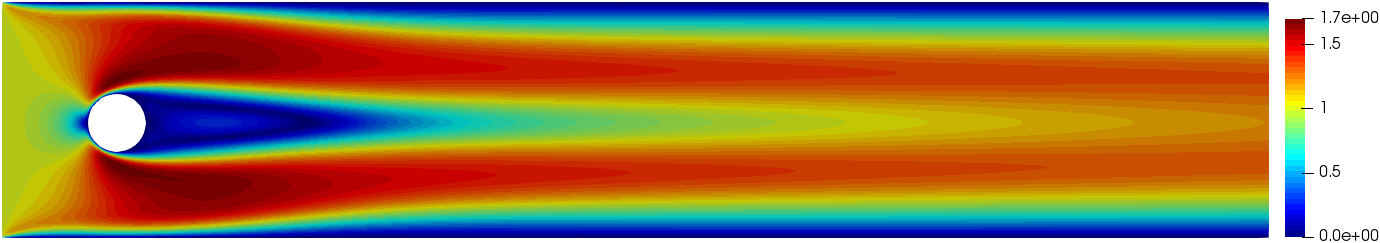}
}
\quad
\subfigure{
\includegraphics[width=0.45\textwidth]{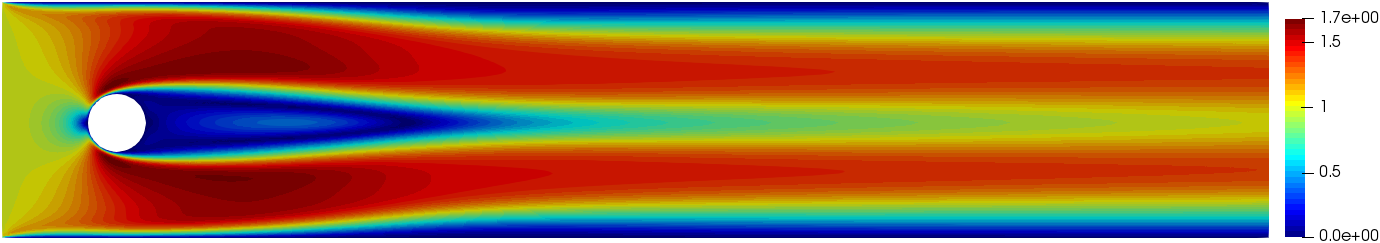}
}
\centering
\caption{\emph{Example 5.} Contours of velocity magnitude for $\mathrm{Re}=10$, $40$, $100$ and $200$ from top-left to bottom-right}\label{fig:ex5:vel:Re:10:200}
\end{figure}

\begin{figure}[htbp]
\centering
\subfigure{
\includegraphics[width=0.45\textwidth]{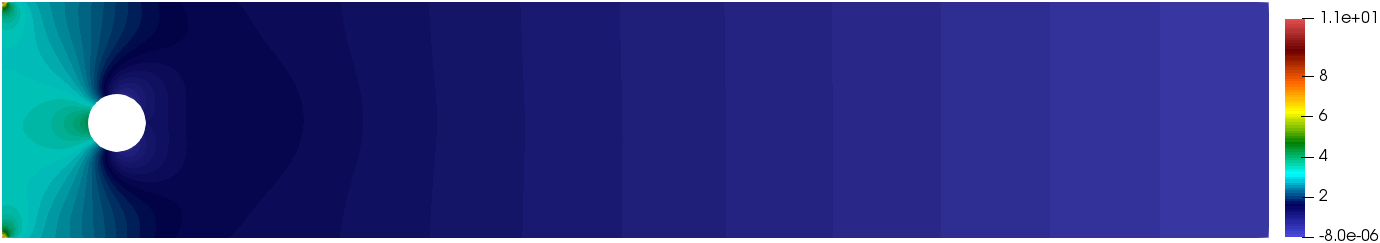}
}
\quad
\subfigure{
\includegraphics[width=0.45\textwidth]{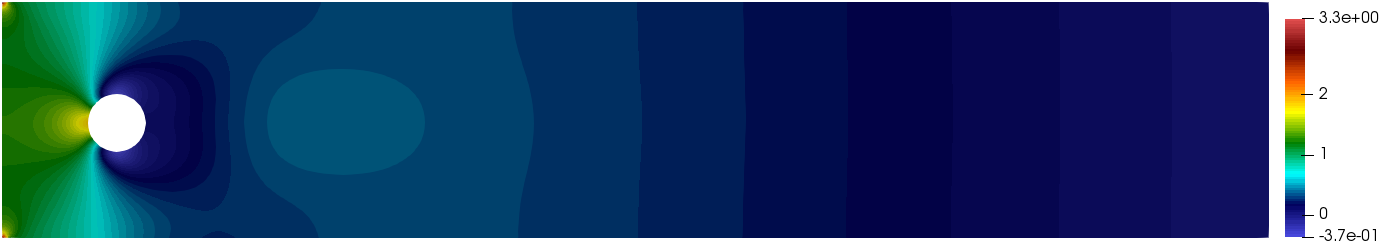}
}

\subfigure{
\includegraphics[width=0.45\textwidth]{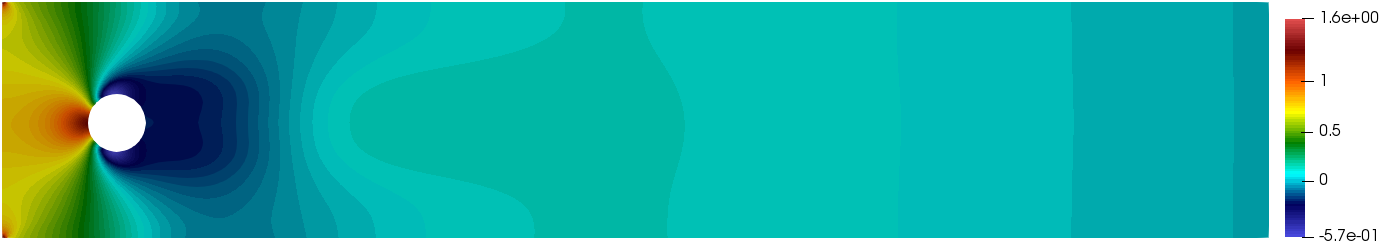}
}
\quad
\subfigure{
\includegraphics[width=0.45\textwidth]{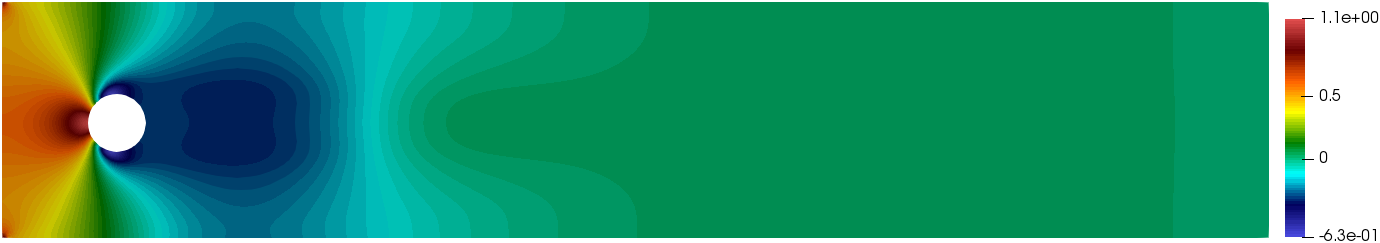}
}
\centering
\caption{\emph{Example 5.} Contours of kinematic pressure for $\mathrm{Re}=10$, $40$, $100$ and $200$ from top-left to bottom-right}\label{fig:ex5:pres:Re:10:200}
\end{figure}

\begin{figure}[htbp]
\centering
\subfigure{
\includegraphics[width=0.45\textwidth]{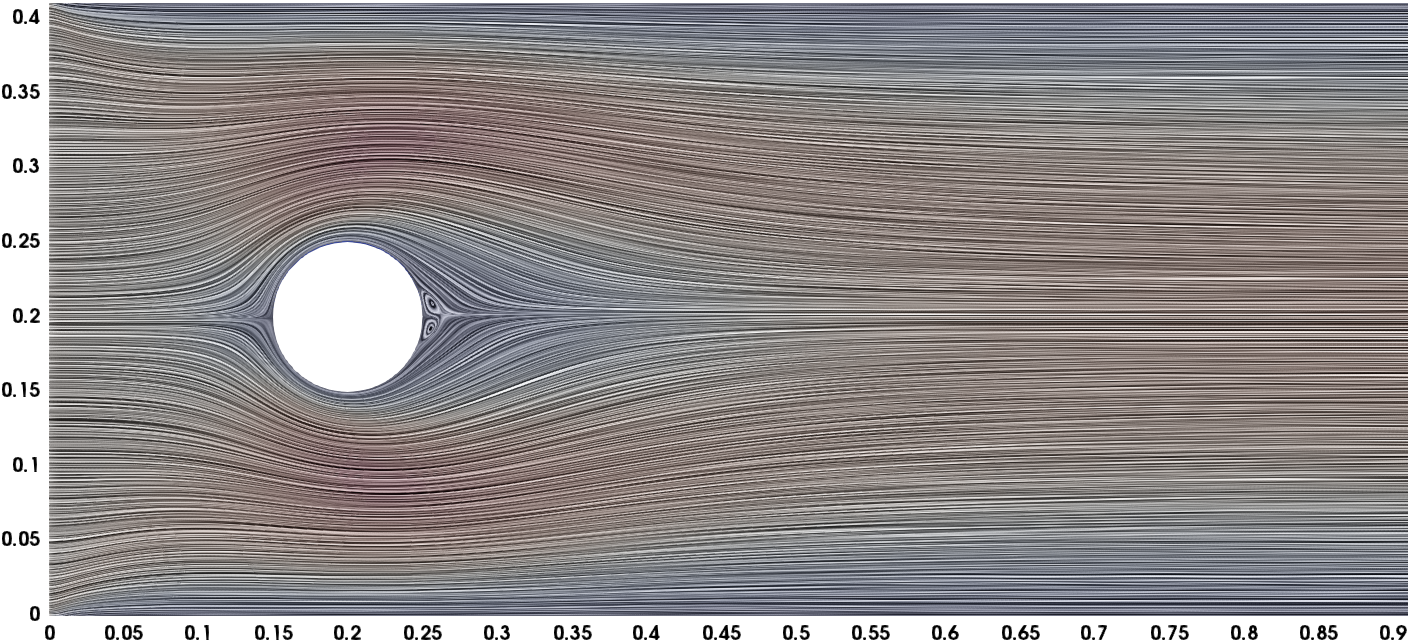}
}
\quad
\subfigure{
\includegraphics[width=0.45\textwidth]{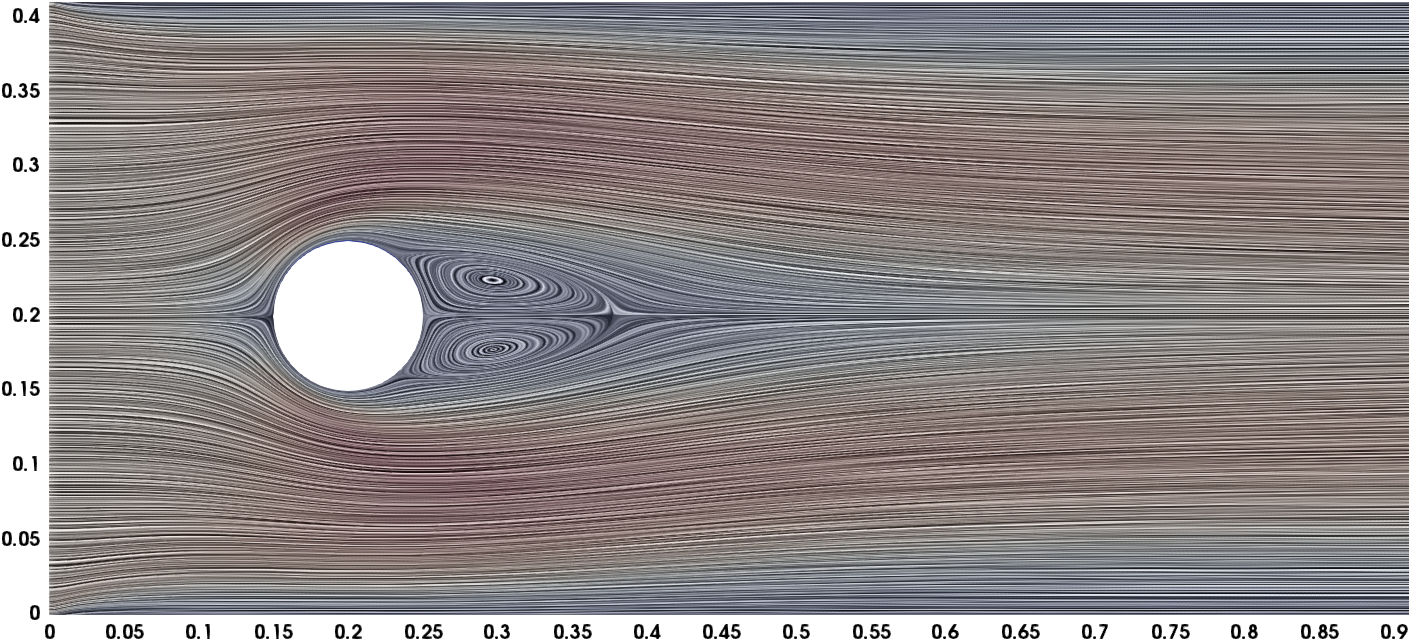}
}

\subfigure{
\includegraphics[width=0.45\textwidth]{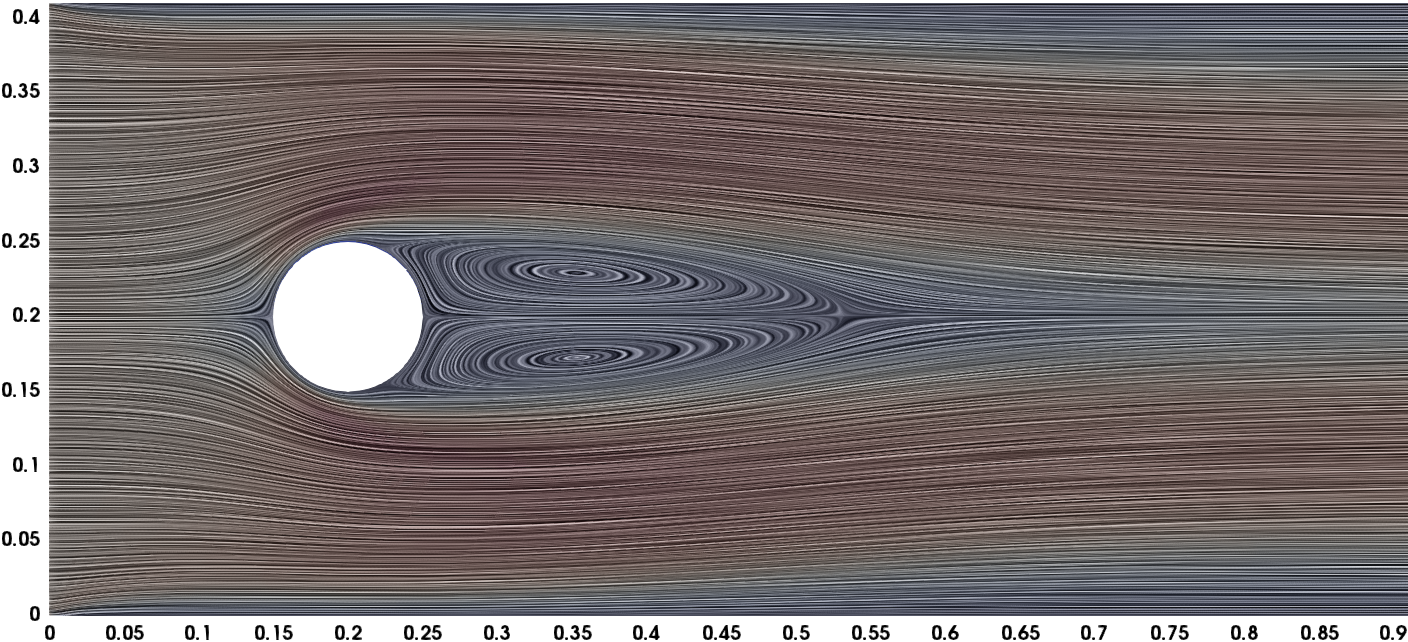}
}
\quad
\subfigure{
\includegraphics[width=0.45\textwidth]{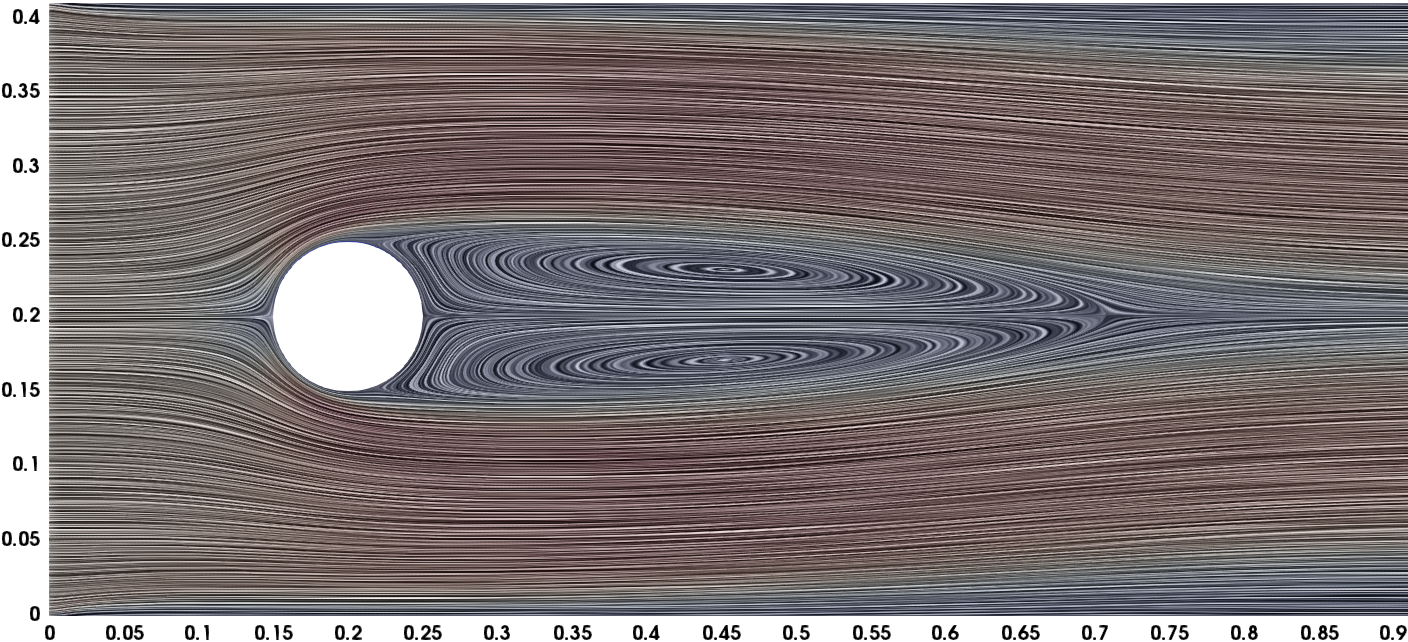}
}
\centering
\caption{\emph{Example 5.} Streamline on colored velocity magnitude distribution in the circulation area for $\mathrm{Re}=10$, $40$, $100$ and $200$ from top-left to bottom-right}\label{fig:ex5:sl:Re:10:200}
\end{figure}

The contours of velocity magnitude and kinematic pressure at various Reynolds numbers are presented in Figures \ref{fig:ex5:vel:Re:10:200} and \ref{fig:ex5:pres:Re:10:200} respectively. The robust symmetry is fully demonstrated in both figures, and the regions of high velocity magnitude extend towards the outlet as the Reynolds number increases, which correspond with the results shown in \cite{Franca2001433}. Moreover, a detailed plot of the streamlines at various Reynolds numbers in the circulation area is provided in Figure \ref{fig:ex5:sl:Re:10:200}. The symmetric vortices make their appearance behind the cylinder at $\mathrm{Re}=10$, which is in good agreement with the reports in Chap. VIII. \cite{VanDyke1975}. As the Reynolds number increases, these two eddies become larger, move towards the outlet, but always keep the symmetry behind the cylinder in the stationary case, which highly match the experimental results shown in \cite{VanDyke2008}.

\subsection{Example 6: Channel flow past a forward-backward facing step}

The benchmark problem is 2D channel flow past a step, which was implemented numerically or experimentally in \cite{Cadou2001825, Fiordilino2018327, Fragos1997495, John2006713, VanDyke2008} through solving the stationary/time-dependent incompressible Navier--Stokes equations. The channel dimensions are $40\times 10$ with a $1\times 1$ step placed five units into the channel from the l.h.s. We choose the body force $\bm{f}=\bm{0}$ and no-slip boundary conditions imposed at all walls. For the inlet and outlet boundary conditions, they are respectively given by
$\bm{u}=(1,0)^\top$ on the inlet and $(p^{\mathrm{kin}}\mathbb{I}-\nu\nabla\bm{u})\bm{n}=\bm{0}$ on the outlet. The solution is interesting and mainly exhibits a smooth velocity distribution with relatively stable eddies in the front of the step and inversely the dynamic eddy formation and detachment occurring behind the step towards the outlet.

\begin{figure}[htbp]
\centering
\includegraphics[width=0.6\textwidth]{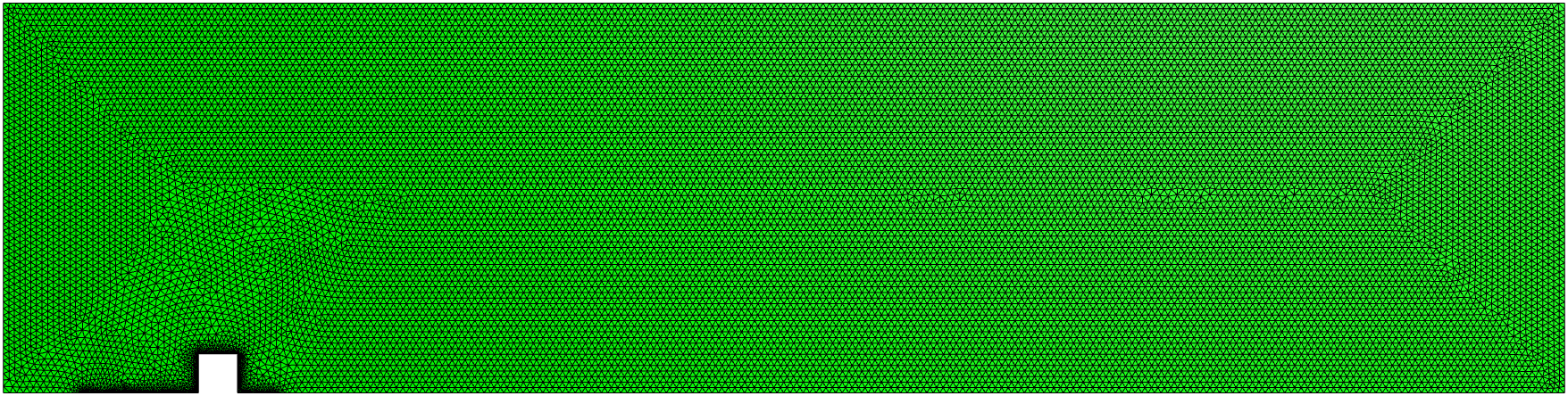}
\caption{\emph{Example 6.} The global triangular mesh plot including the refined grids near the step}\label{fig:ex6:meshplot}
\end{figure}

We consider computations by $\bm{\mathrm{P}}_3^{\mathrm{bubble}}$-$\mathrm{P}_2^{\mathrm{dc}}$-$\bm{\mathrm{BDM}}_3$ performed on a locally refined triangular mesh presented in Figure \ref{fig:ex6:meshplot}. Here we take $h_l=1/30$ and $h_g=1/6$. Note that the Reynolds number is given by $\mathrm{Re}=\frac{1}{\nu}$ here. By adopting the updated Newton iteration \eqref{linear:Newton:cylinder}, the problem is first solved for $\mathrm{Re}=10$, then $\mathrm{Re}=25$, and then in steps of $25$ until $\mathrm{Re}=1000$, with the solution for the previous value of $\mathrm{Re}$ used as initial guess for the next; the Stokes equations are solved to provide the initial guess used at $\mathrm{Re}=10$.

\begin{figure}[htbp]
\centering
\subfigure{
\includegraphics[width=0.45\textwidth]{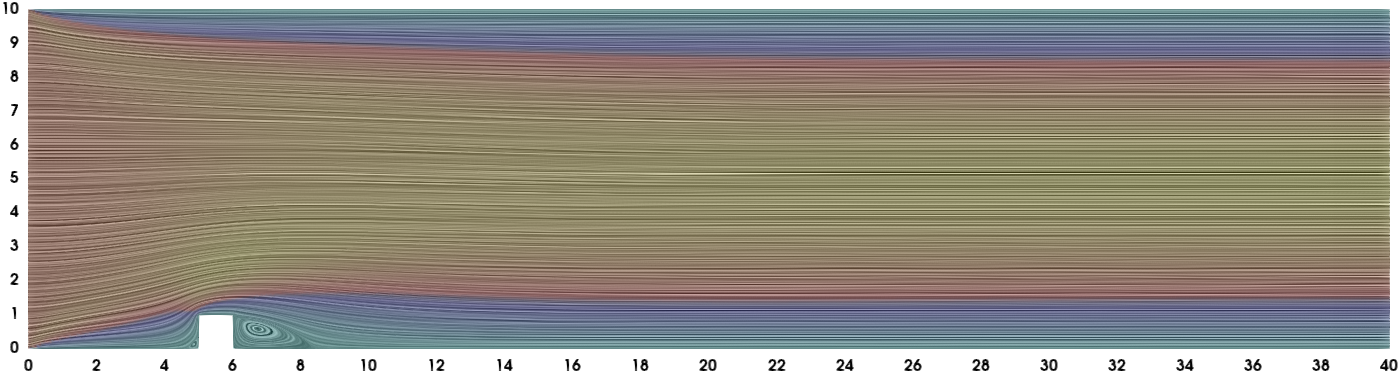}
}
\quad
\subfigure{
\includegraphics[width=0.45\textwidth]{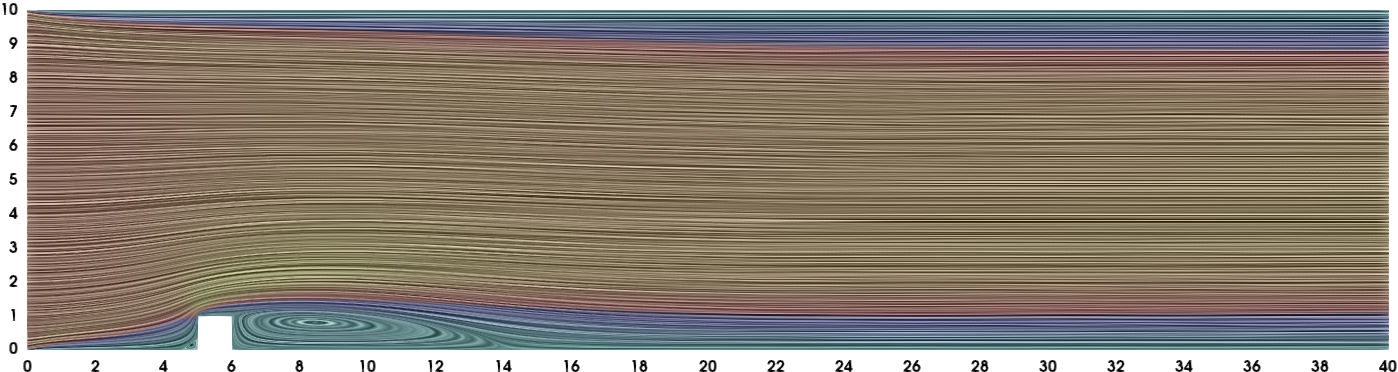}
}

\subfigure{
\includegraphics[width=0.45\textwidth]{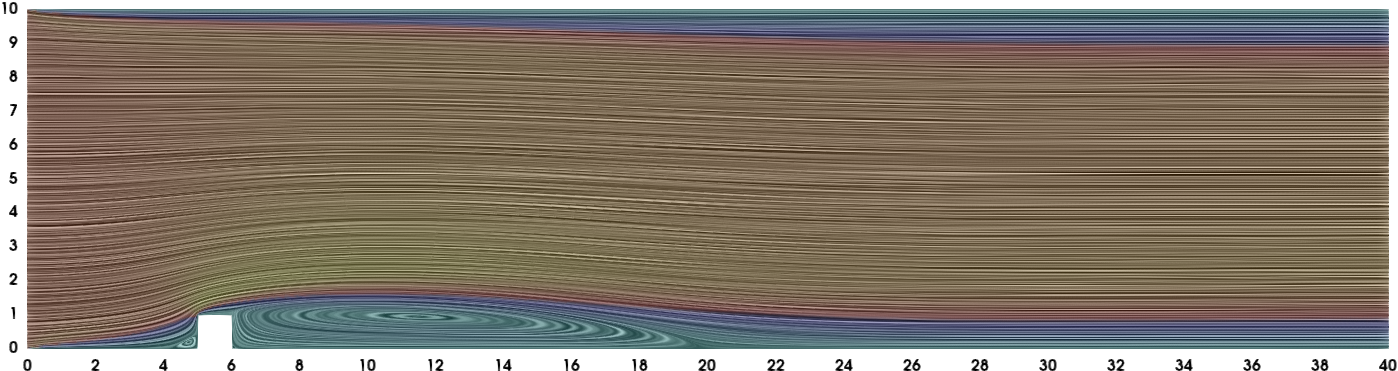}
}
\quad
\subfigure{
\includegraphics[width=0.45\textwidth]{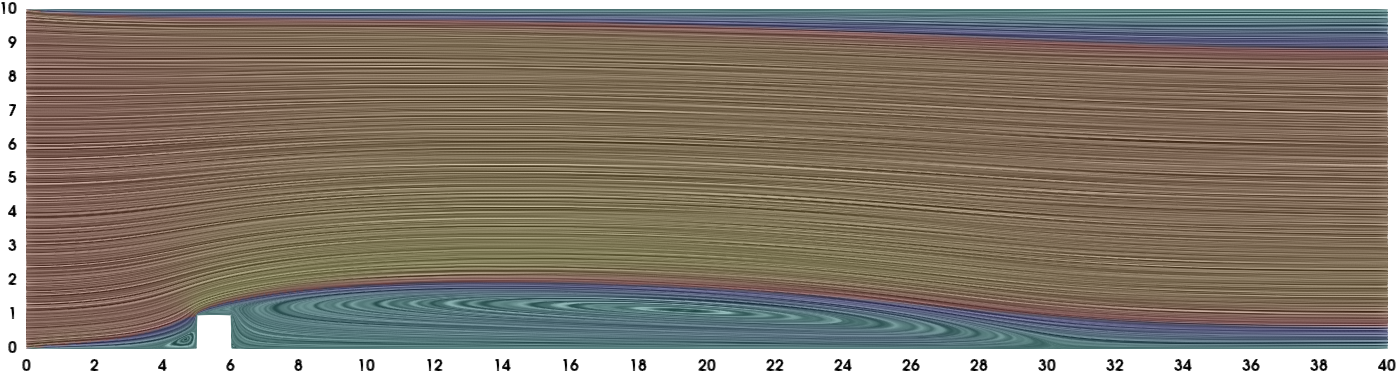}
}
\centering
\caption{\emph{Example 6.} Streamline on colored velocity magnitude distribution for $\mathrm{Re}=10$, $50$, $100$ and $200$ from top-left to bottom-right}\label{fig:ex6:Re:10:200}
\end{figure}

\begin{figure}[htbp]
\centering
\subfigure[Fragos et al. \cite{Fragos1997495}]{
\includegraphics[width=0.47\textwidth]{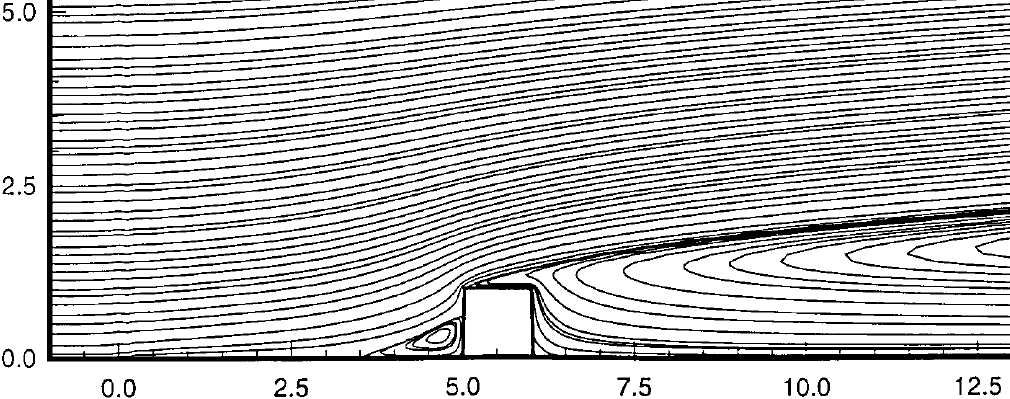}
}
\quad
\subfigure[$\bm{\mathrm{P}}_3^{\mathrm{bubble}}$-$\mathrm{P}_2^{\mathrm{dc}}$-$\bm{\mathrm{BDM}}_3$]{
\includegraphics[width=0.45\textwidth]{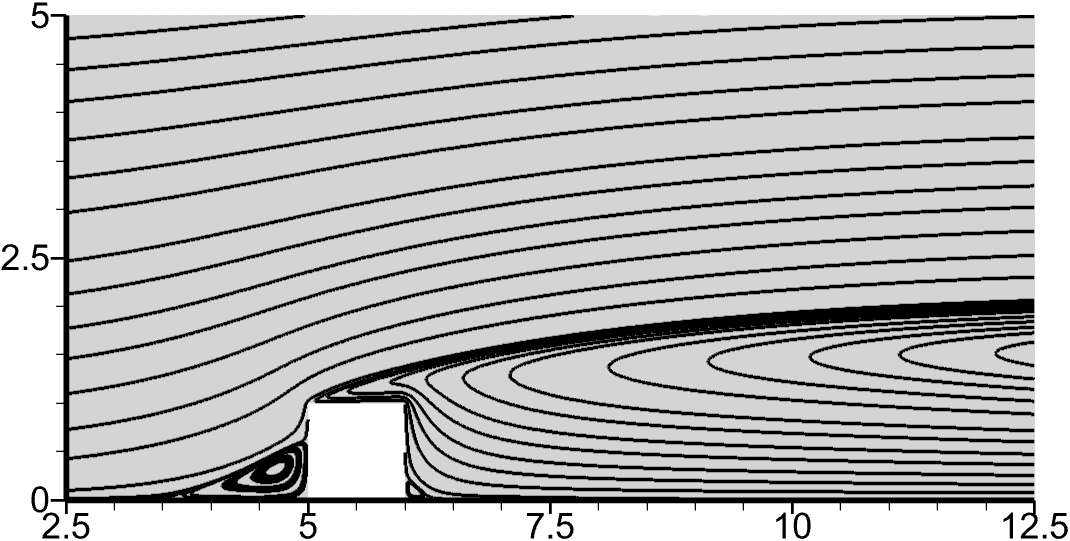}
}
\centering
\caption{\emph{Example 6.} A magnified view of streamline contours near the step for $\mathrm{Re}=275$. Comparison between Fragos et al. \cite{Fragos1997495} (left) and $\bm{\mathrm{P}}_3^{\mathrm{bubble}}$-$\mathrm{P}_2^{\mathrm{dc}}$-$\bm{\mathrm{BDM}}_3$ (right)}\label{fig:ex6:Re275:compare}
\end{figure}

In Figure \ref{fig:ex6:Re:10:200}, it is observed that the streamlines are relatively stable for lower Reynolds numbers ($\mathrm{Re}\leqslant 200$). On the global perspective, there are a growing but stable eddy in the front of the step and a more apparently growing eddy behind the step; the eddy behind the step expands rapidly towards the outlet as the Reynolds number increases. It is worth mentioning that these results match the experimental results shown in \cite{Fragos1997495} involving the sizes, the shapes and the center location $xy$-coordinates of the eddies. In particular, when $\mathrm{Re}\geqslant 275$, in \cite{Fragos1997495} the computational domain is shrunk by about 65\% with the same width but a much shorter length 13.5. We still compare the streamline contours near the step by $\bm{\mathrm{P}}_3^{\mathrm{bubble}}$-$\mathrm{P}_2^{\mathrm{dc}}$-$\bm{\mathrm{BDM}}_3$ for $\mathrm{Re}=275$ with that of Fragos et al. \cite{Fragos1997495}. The comparison is presented in Figure \ref{fig:ex6:Re275:compare}, and one can find two results are virtually identical, which implies the velocity distribution near the step is independent of the channel length at $\mathrm{Re}=275$.

\begin{figure}[htbp]
\centering
\subfigure{
\includegraphics[width=0.45\textwidth]{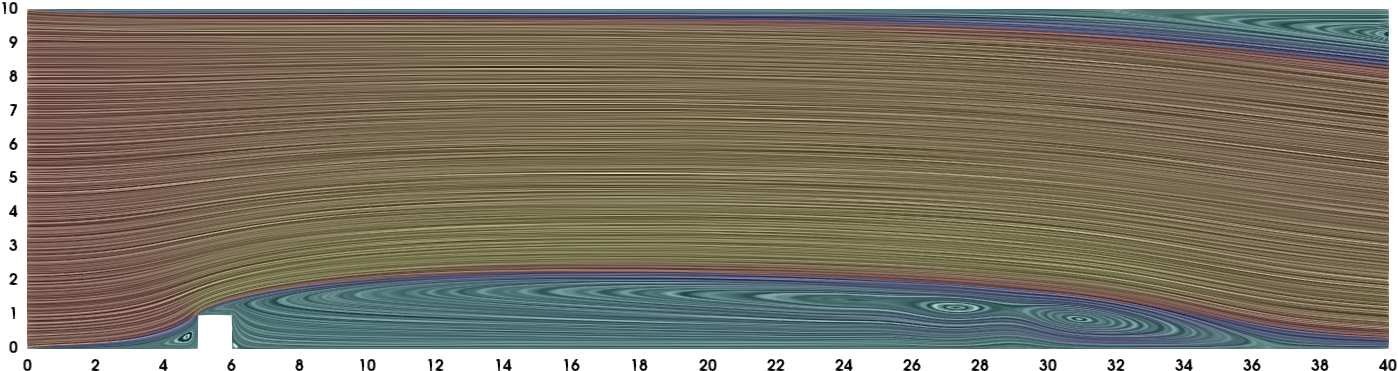}
}
\quad
\subfigure{
\includegraphics[width=0.45\textwidth]{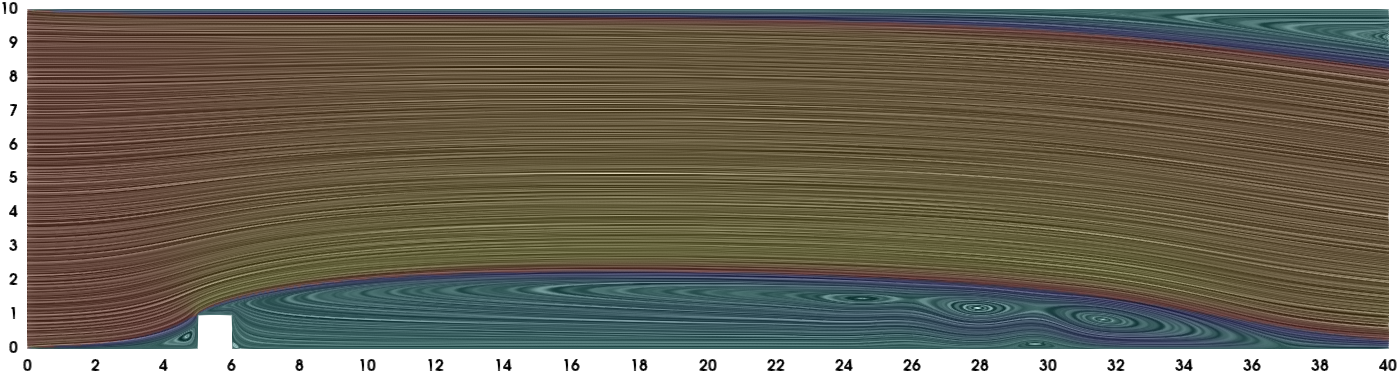}
}

\subfigure{
\includegraphics[width=0.45\textwidth]{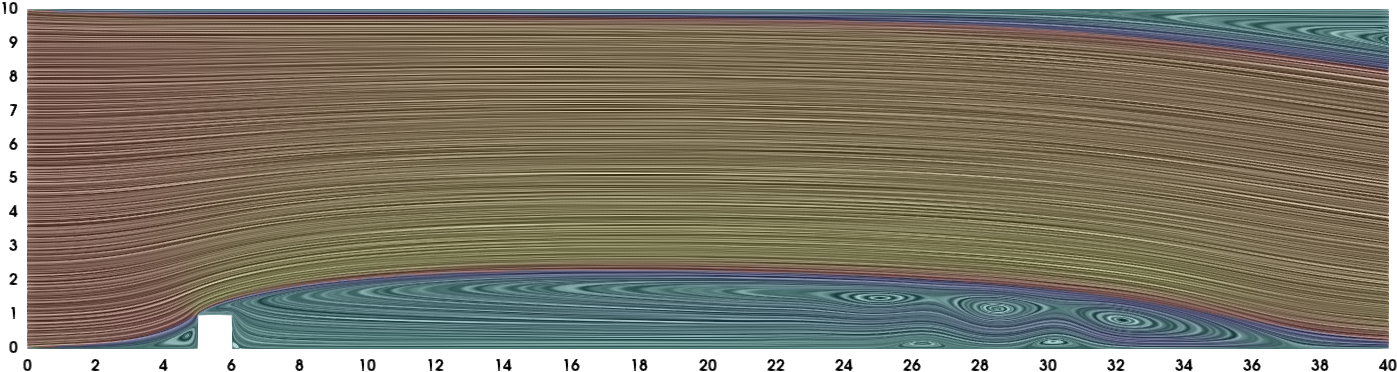}
}
\quad
\subfigure{
\includegraphics[width=0.45\textwidth]{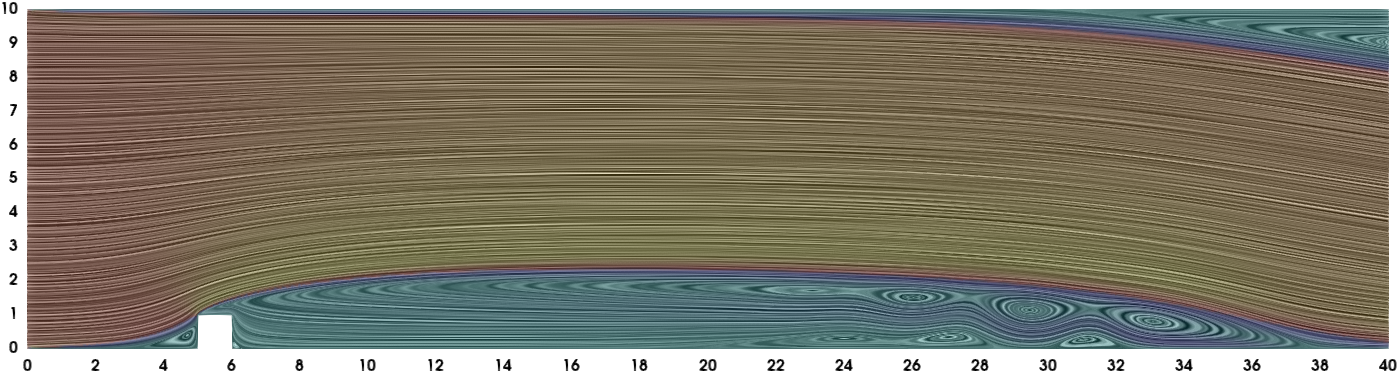}
}
\centering
\caption{\emph{Example 6.} Streamline on colored velocity magnitude distribution for $\mathrm{Re}=400$, $450$, $500$ and $600$ from top-left to bottom-right}\label{fig:ex6:Re:400:600}
\end{figure}

As the Reynolds number further increasing up to 600, a first important difference highlighted by the results shown in Figure \ref{fig:ex6:Re:400:600} that the eddy behind the step starts continuous detachment, which looks like mitosis. In addition, these eddies still move towards the outlet collectively. Another interesting phenomenon is the appearance of an eddy clinging to the bottom at about a quarter of the channel from the r.h.s. at $\mathrm{Re}=450$, and this eddy keeps growing and simultaneously moving towards the outlet. It happens that there is a similar case that the formation of another eddy occurs near the above described eddy at $\mathrm{Re}=500$. But contrarily, its shape and size change rapidly, and the speed of expanding towards the step is visibly faster than towards the outlet. We also note the difference at the top right corner of the channel that an incomplete eddy forms and keeps moving down slowly. According to the numerical results at various Reynolds numbers provided in \cite{Fragos1997495}, we believe that the appearance of this eddy is mainly due to the insufficient length of the channel and the Neumann boundary condition imposed at the outlet.

\begin{figure}[htbp]
\centering
\subfigure{
\includegraphics[width=0.45\textwidth]{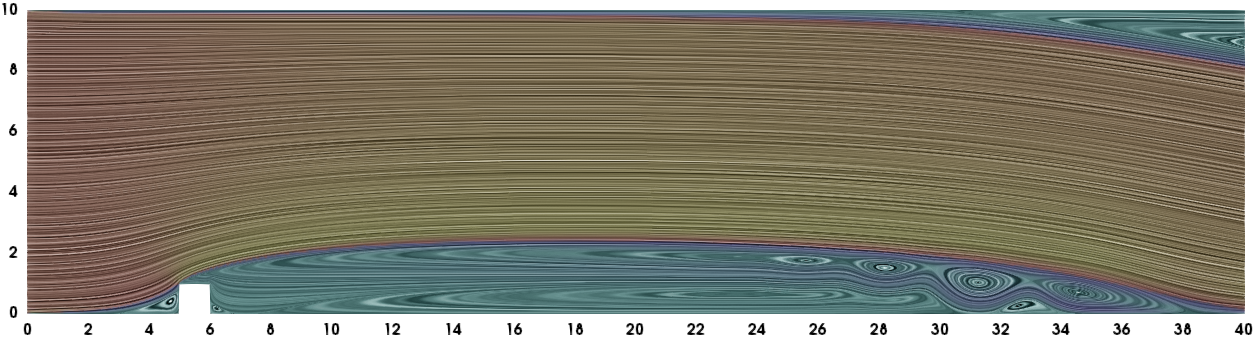}
}
\quad
\subfigure{
\includegraphics[width=0.45\textwidth]{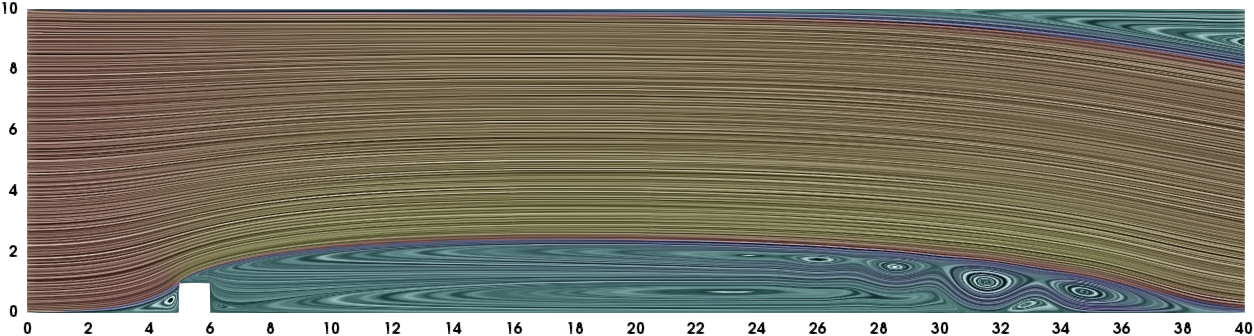}
}

\subfigure{
\includegraphics[width=0.45\textwidth]{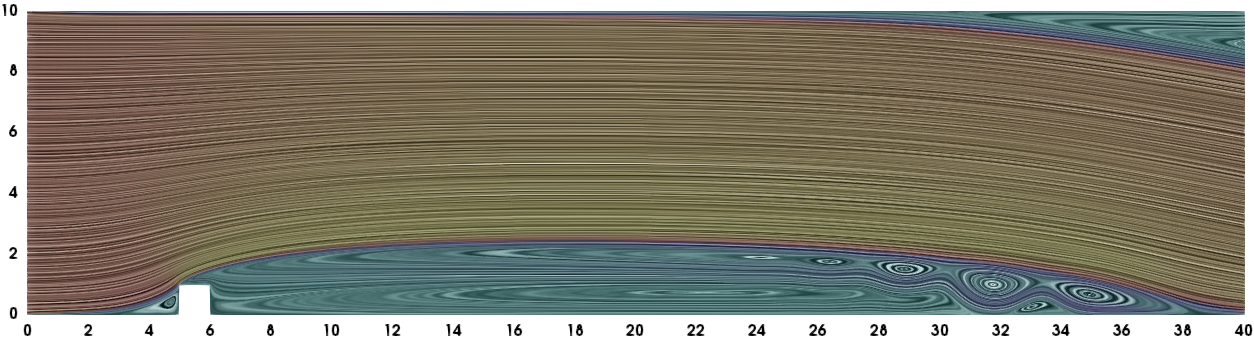}
}
\quad
\subfigure{
\includegraphics[width=0.45\textwidth]{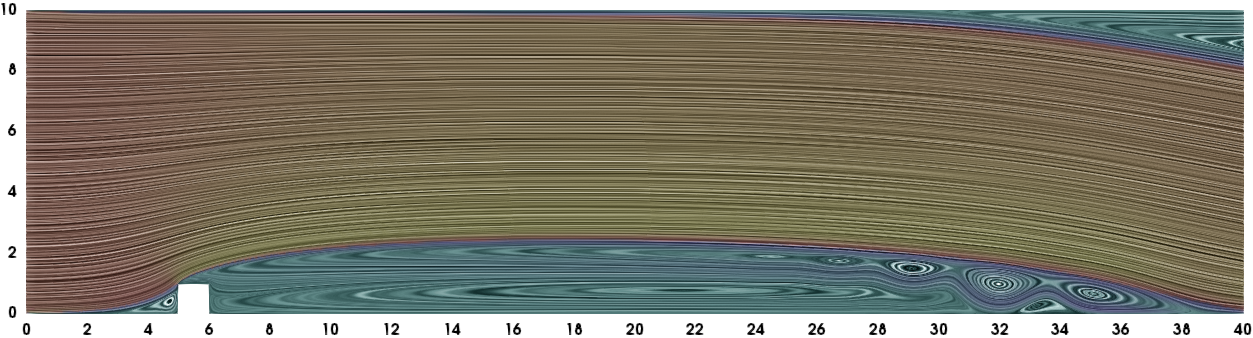}
}
\centering
\caption{\emph{Example 6.} Streamline on colored velocity magnitude distribution for $\mathrm{Re}=850$, $900$, $950$ and $1000$ from top-left to bottom-right}\label{fig:ex6:Re:850:1000}
\end{figure}

\begin{figure}[htbp]
\centering
\subfigure{
\includegraphics[width=4cm]{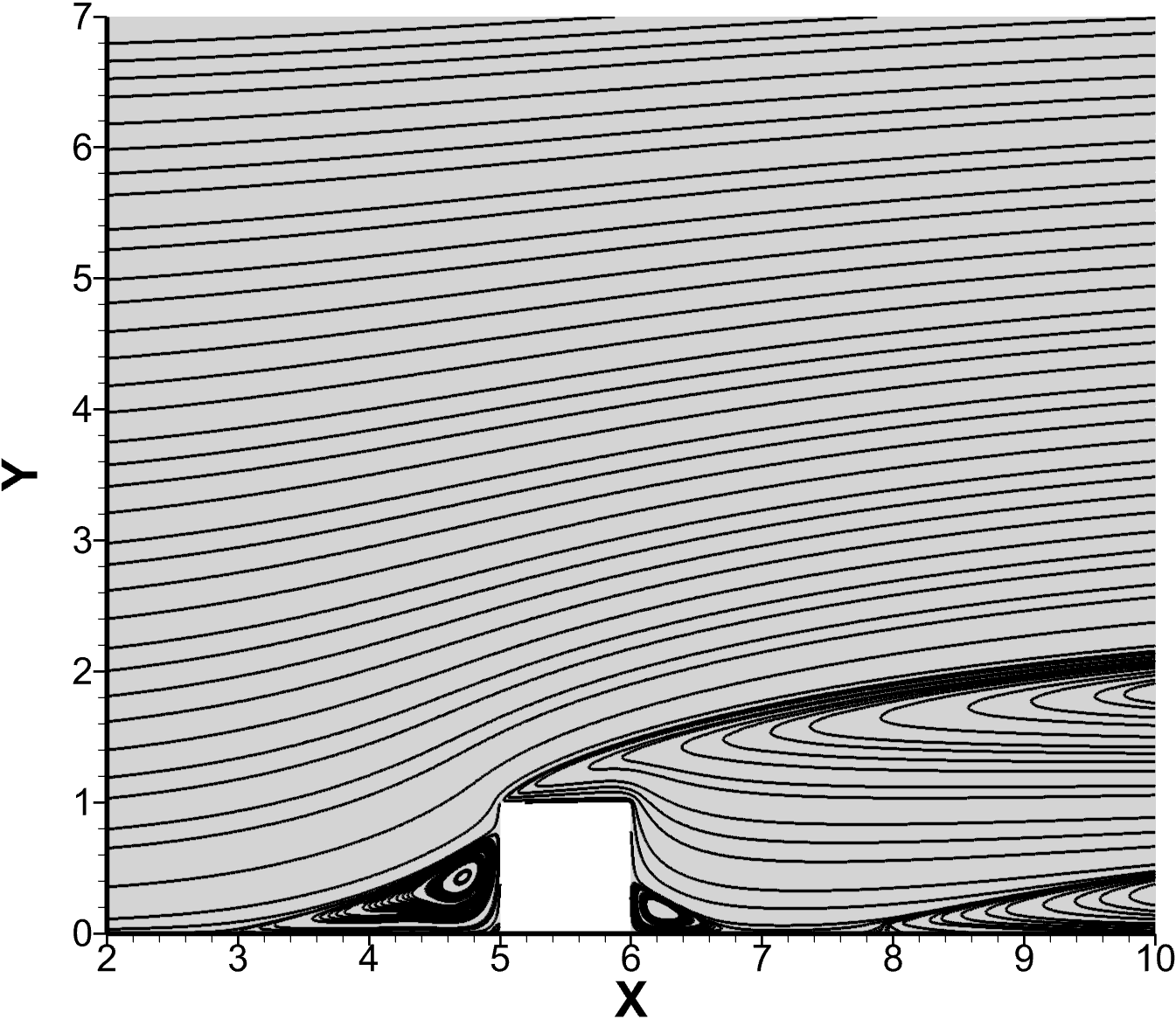}
}
\,
\subfigure{
\includegraphics[width=4cm]{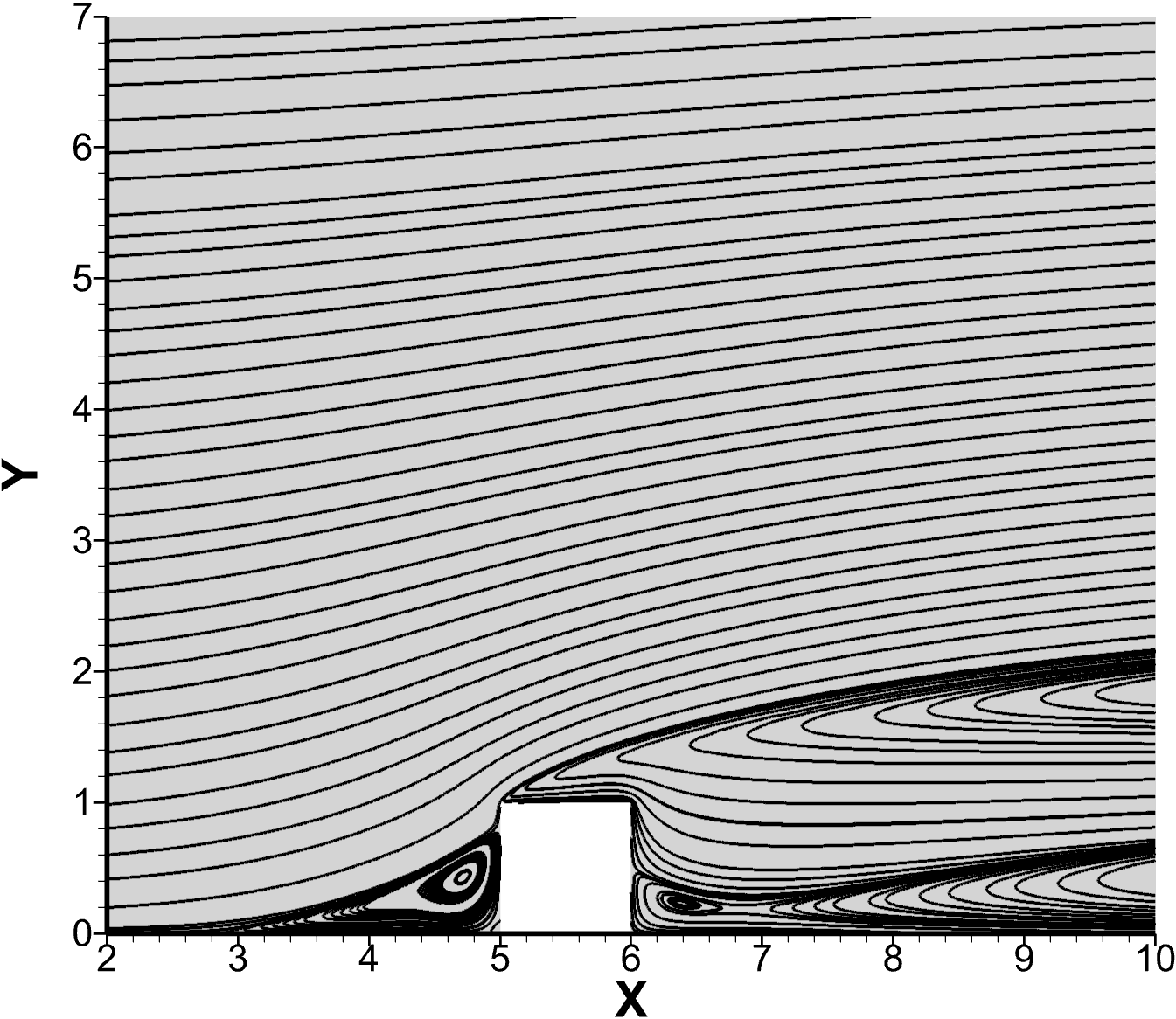}
}
\,
\subfigure{
\includegraphics[width=4cm]{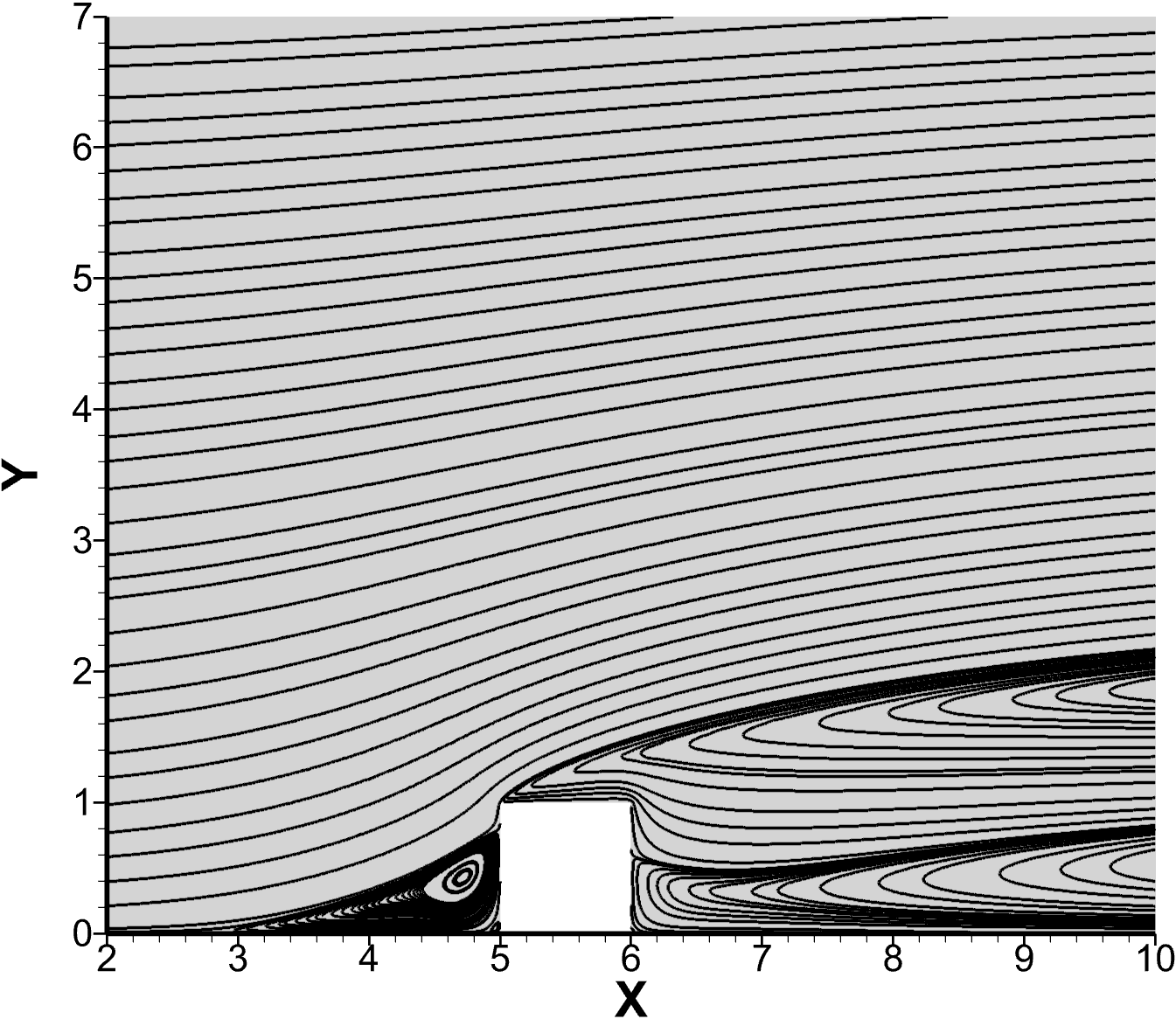}
}
\,
\subfigure{
\includegraphics[width=4cm]{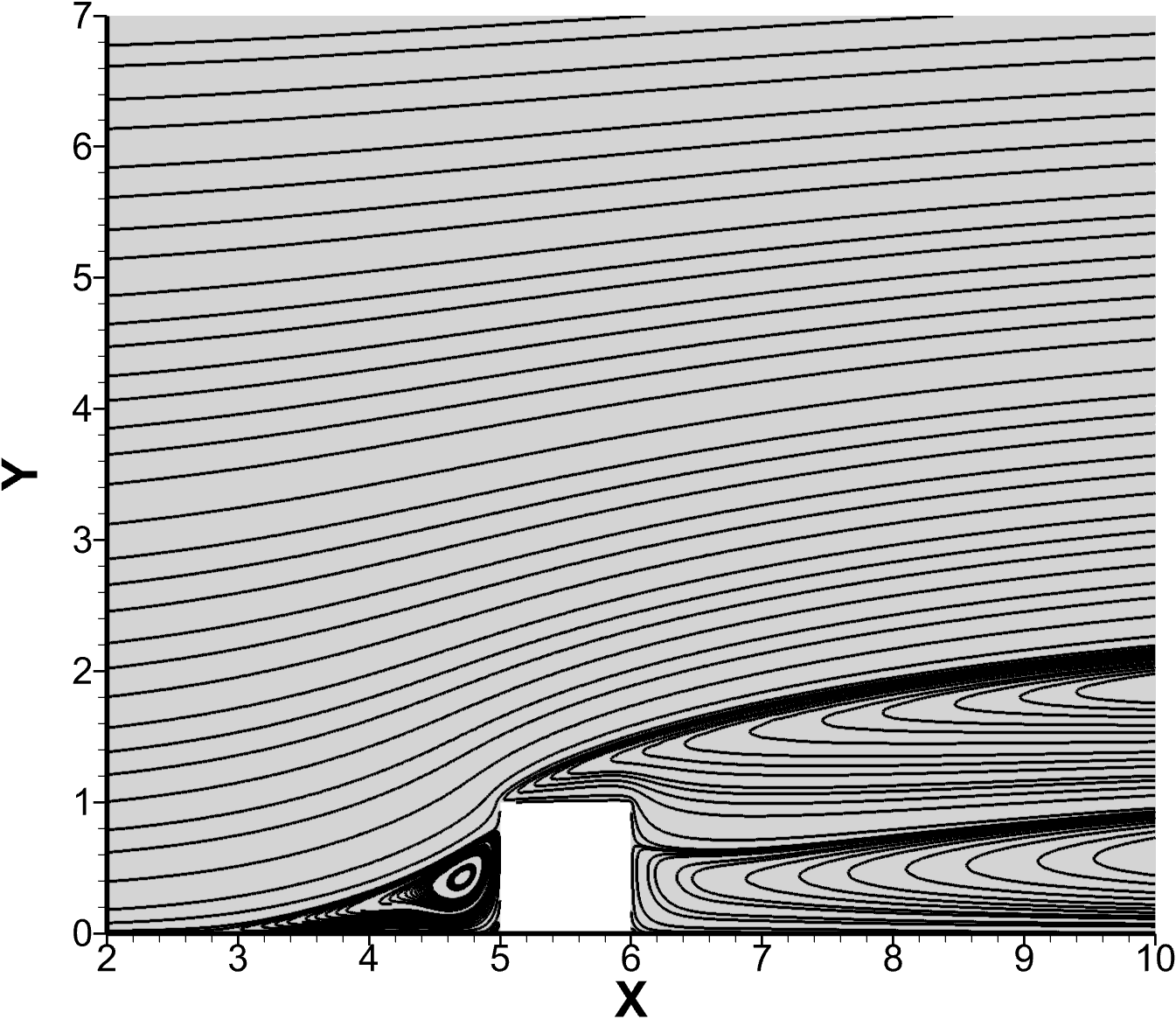}
}
\centering
\caption{\emph{Example 6.} A magnified view of streamline contours near the step for $\mathrm{Re}=850$, $900$, $950$ and $1000$ from left to right}\label{fig:ex6:local:Re:850:1000}
\end{figure}

\begin{figure}[htbp]
\centering
\subfigure{
\includegraphics[width=0.45\textwidth]{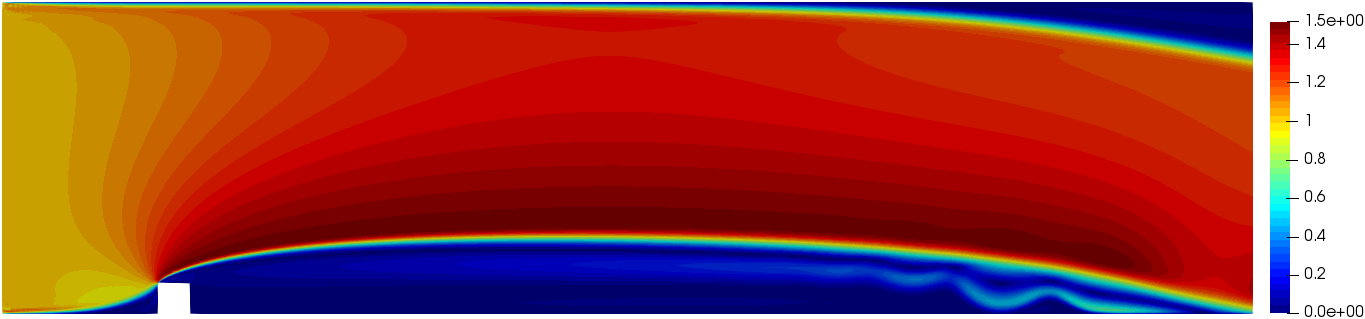}
}
\,
\subfigure{
\includegraphics[width=0.45\textwidth]{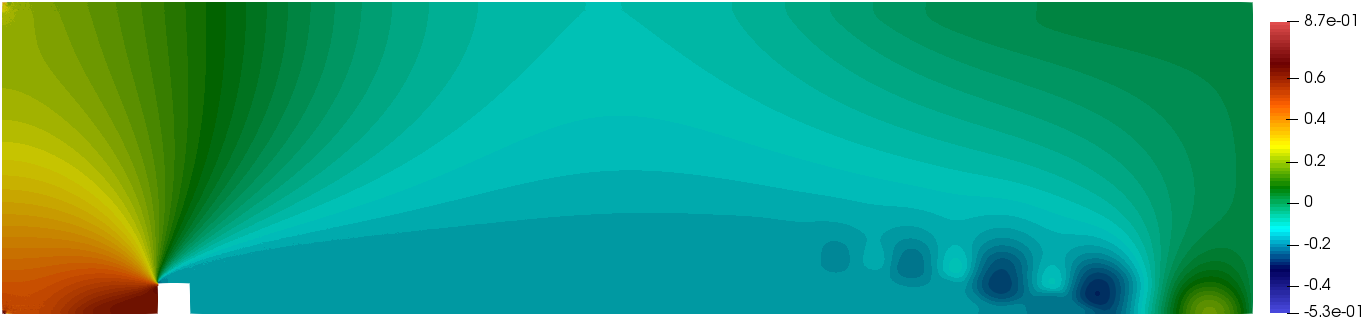}
}
\centering
\caption{\emph{Example 6.} Contours of velocity magnitude (left) and kinetic pressure (right) for $\mathrm{Re}=1000$}\label{fig:ex6:Re1000:vel:pres}
\end{figure}

All of the above descriptions still hold for higher Reynolds numbers ($\mathrm{Re}\geqslant 850$), which are supported in Figure \ref{fig:ex6:Re:850:1000}. In addition, there is an unexpected phenomenon that the eddy, which stems from one debuting at $\mathrm{Re}=500$, coalesces continuously with a smaller eddy clinging to the back of the step, and becomes a large eddy independent of that existing behind the step since the beginning. With a magnified view, a detailed plot of streamline contours near the step for these Reynolds numbers is presented in Figure \ref{fig:ex6:local:Re:850:1000} to observe this evolution with a better perspective. Further the contours of velocity magnitude and kinetic pressure at $\mathrm{Re}=1000$ are also plotted in Figure \ref{fig:ex6:Re1000:vel:pres}. Both of the velocity and the kinetic pressure distributions present dramatic changes near the bottom of the outlet.

\subsection{Example 7: Channel flow past forward-backward facing steps and around square solid objects}

At the end of the last example, it is pointed out that both of the velocity and the kinetic pressure distributions change rapidly near the bottom of the outlet at $\mathrm{Re}=1000$. Hence, we are driven to consider how a geometric change of the channel affects the velocity distribution and/or the kinetic pressure distribution. To this end, the same test case as the last example will be implemented but with a modified channel flow region. In this example, the channel dimensions are $40\times 10$ with two $1\times 1$ steps placed five units into the bottom and the top of channel from the l.h.s. respectively, and two $1\times 1$ solid objects fixed two units over the lower step and two units under the upper step respectively.

\begin{figure}[htbp]
\centering
\includegraphics[width=0.6\textwidth]{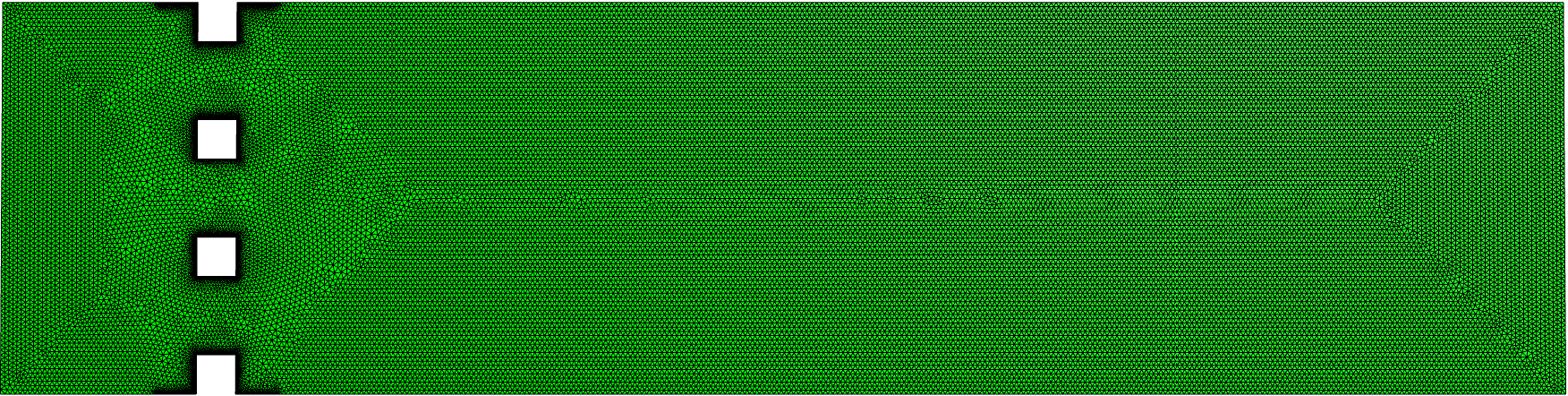}
\caption{\emph{Example 7.} The global triangular mesh plot including the refined grids near the steps and the objects}\label{fig:ex7:meshplot}
\end{figure}

\begin{figure}[htbp]
\centering
\subfigure{
\includegraphics[width=0.45\textwidth]{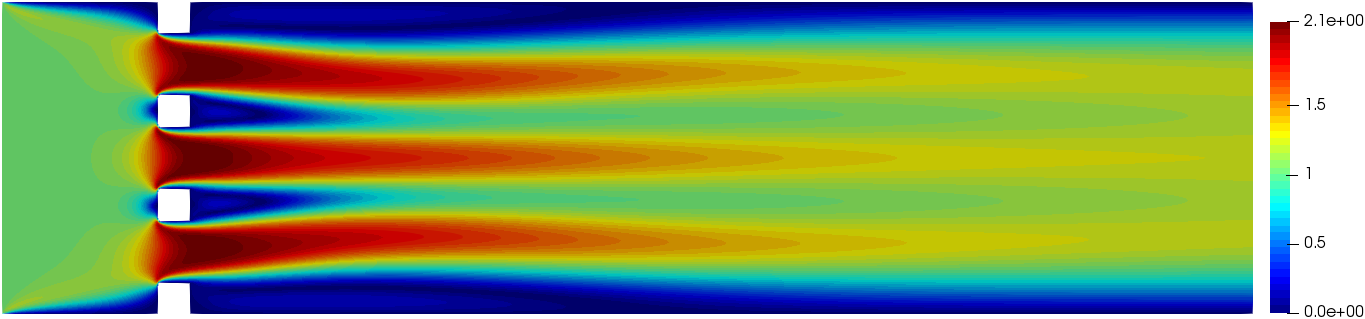}
}
\quad
\subfigure{
\includegraphics[width=0.45\textwidth]{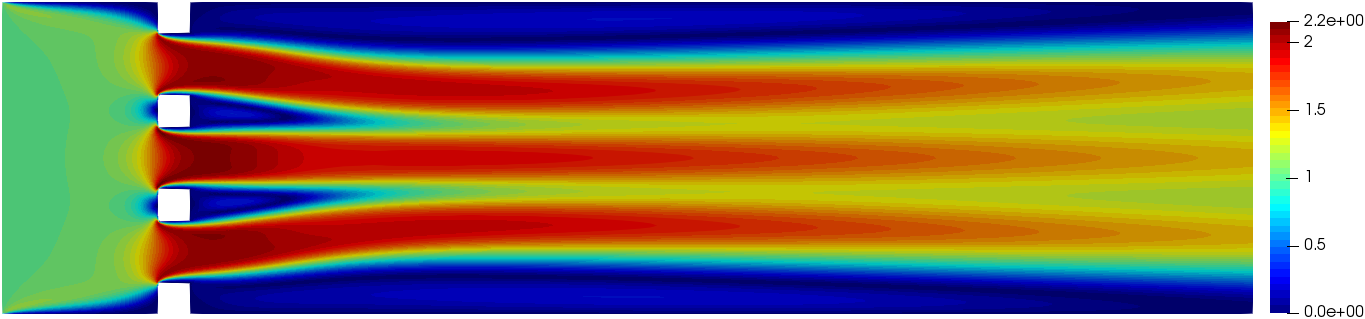}
}

\subfigure{
\includegraphics[width=0.45\textwidth]{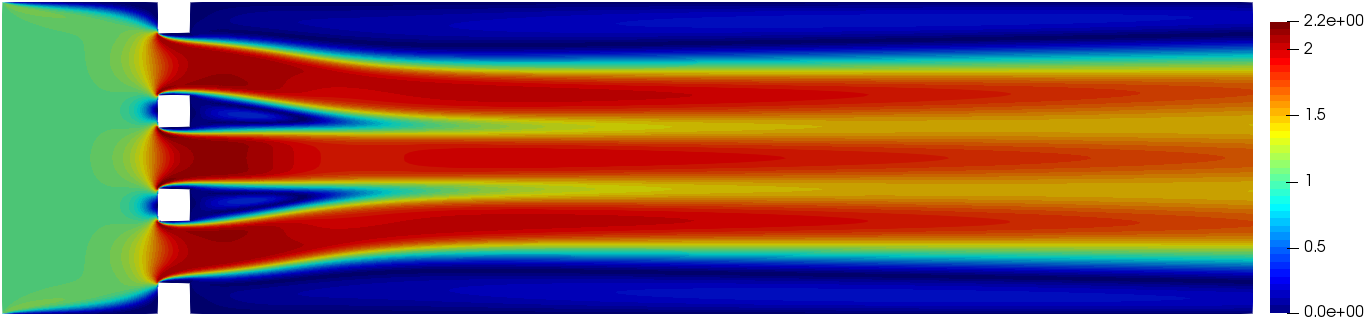}
}
\quad
\subfigure{
\includegraphics[width=0.45\textwidth]{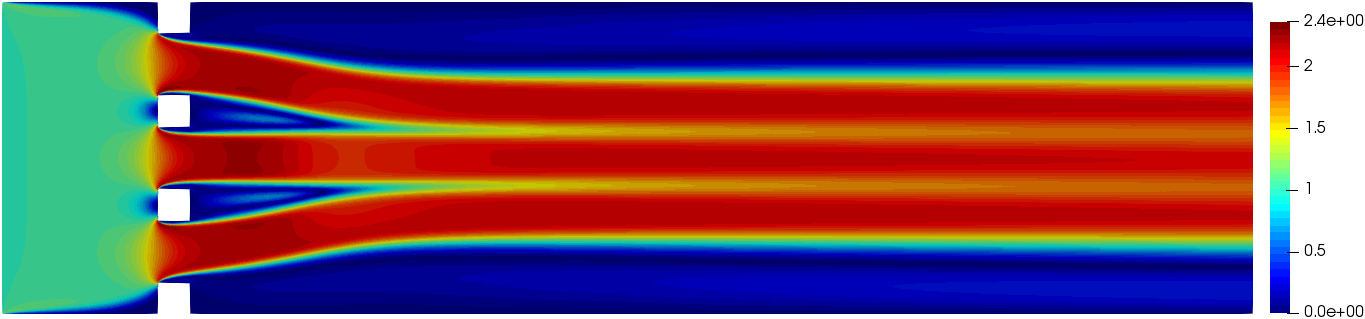}
}
\centering
\caption{\emph{Example 7.} Contours of velocity magnitude for $\mathrm{Re}=50$, $100$, $150$ and $300$ from top-left to bottom-right}\label{fig:ex7:vel:Re:50:300}
\end{figure}

\begin{figure}[htbp]
\centering
\subfigure{
\includegraphics[width=0.45\textwidth]{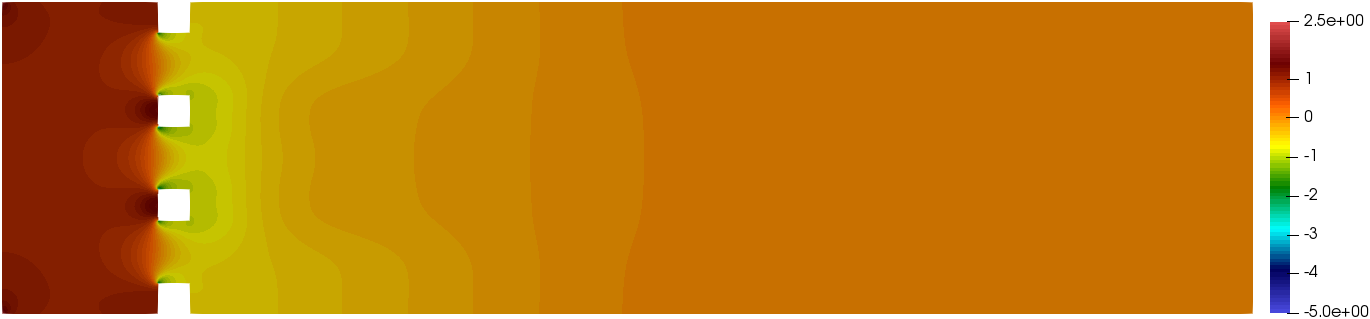}
}
\quad
\subfigure{
\includegraphics[width=0.45\textwidth]{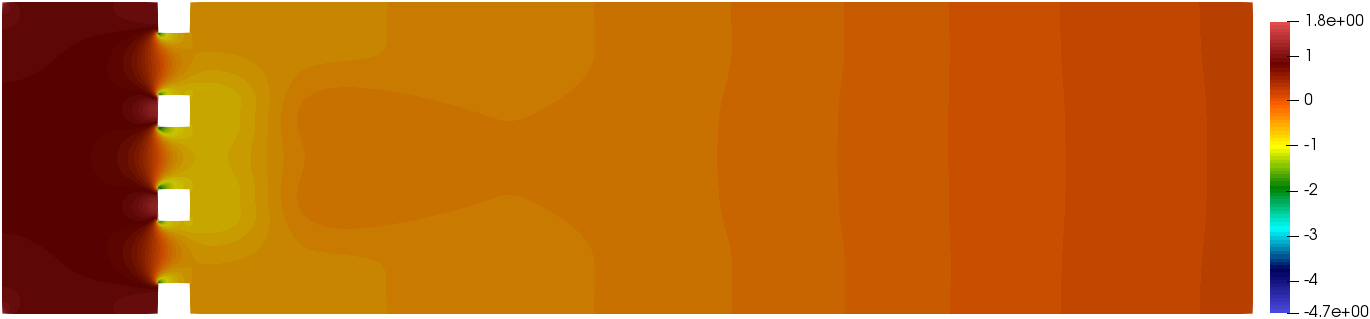}
}

\subfigure{
\includegraphics[width=0.45\textwidth]{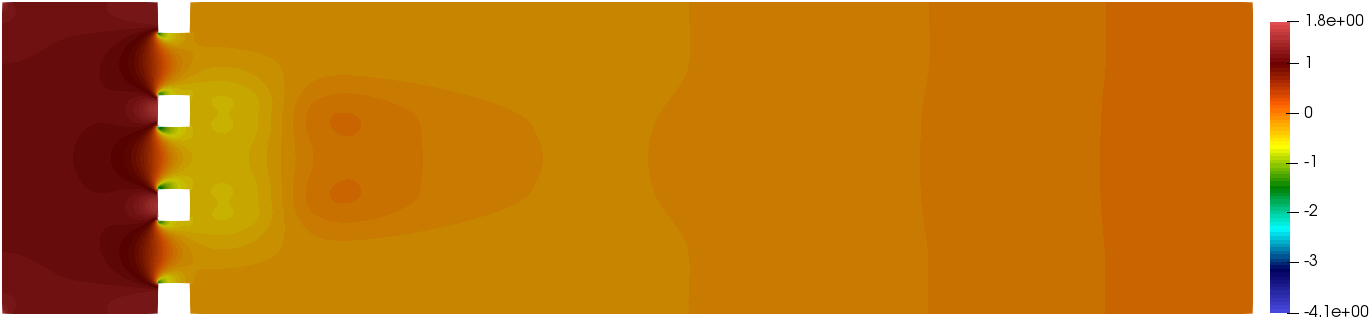}
}
\quad
\subfigure{
\includegraphics[width=0.45\textwidth]{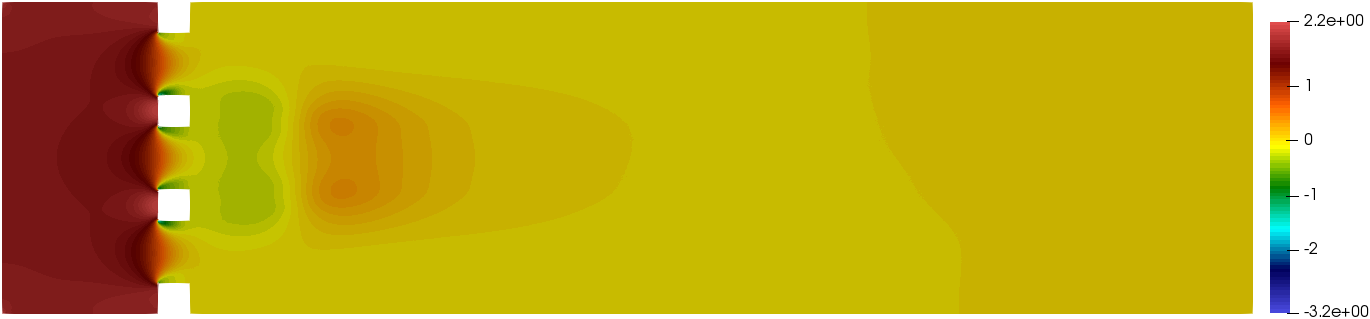}
}
\centering
\caption{\emph{Example 7.} Contours of kinetic pressure for $\mathrm{Re}=50$, $100$, $150$ and $300$ from top-left to bottom-right}\label{fig:ex7:pres:Re:50:300}
\end{figure}

\begin{figure}[htbp]
\centering
\subfigure{
\includegraphics[width=0.45\textwidth]{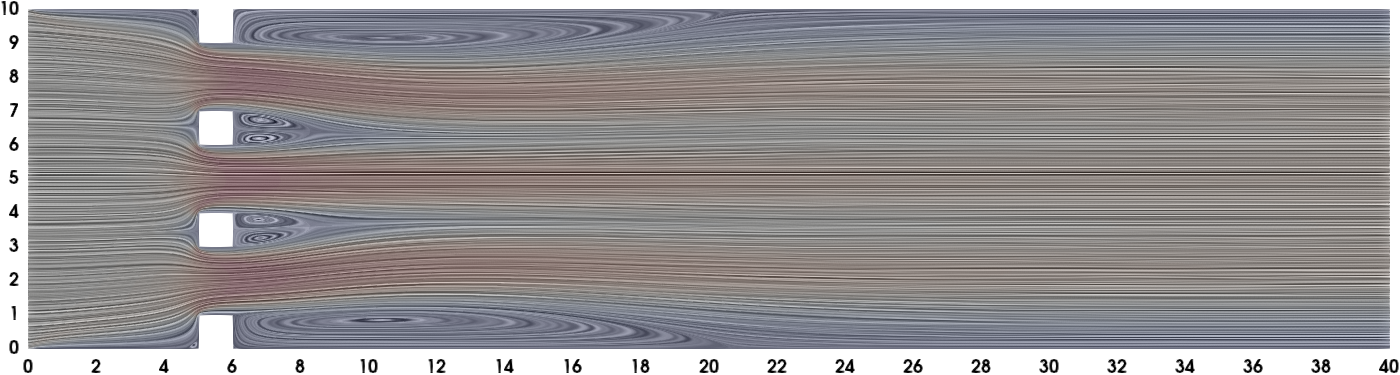}
}
\quad
\subfigure{
\includegraphics[width=0.45\textwidth]{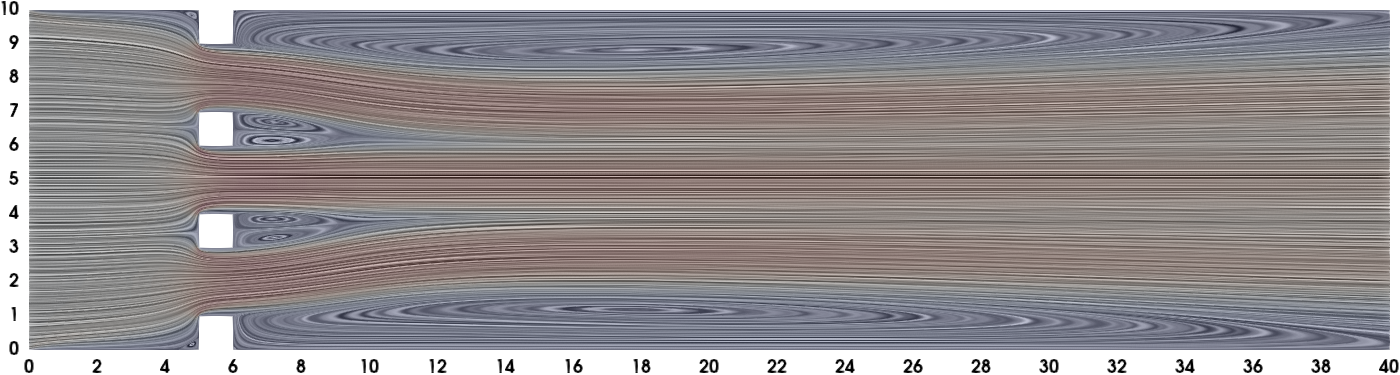}
}

\subfigure{
\includegraphics[width=0.45\textwidth]{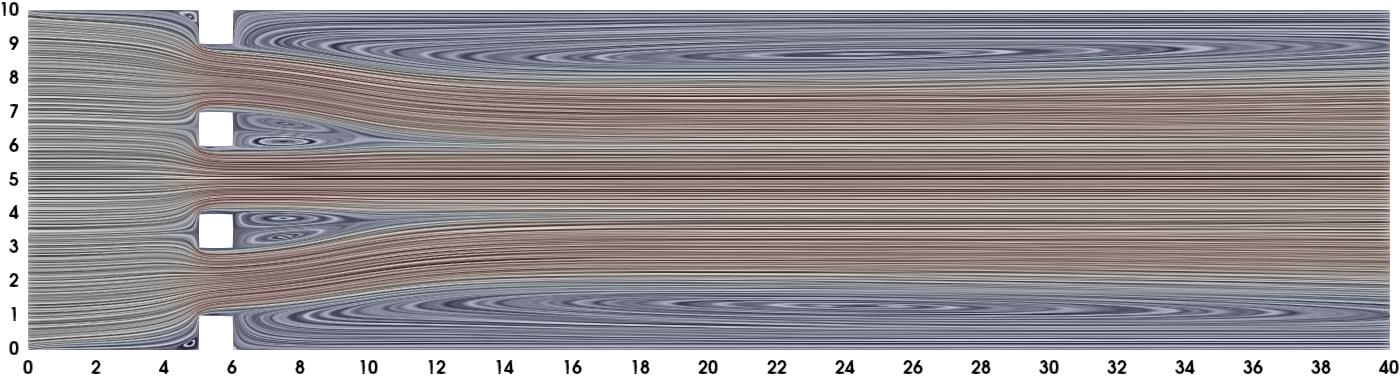}
}
\quad
\subfigure{
\includegraphics[width=0.45\textwidth]{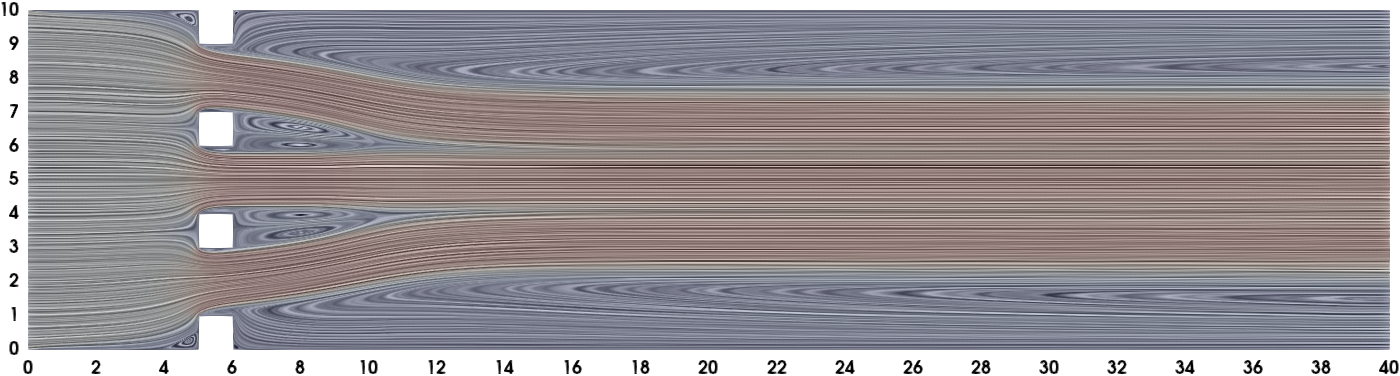}
}
\centering
\caption{\emph{Example 7.} Streamline on colored velocity magnitude distribution for $\mathrm{Re}=50$, $100$, $150$ and $300$ from top-left to bottom-right}\label{fig:ex7:sl:Re:50:300}
\end{figure}

\begin{figure}[htbp]
\centering
\subfigure{
\includegraphics[width=0.45\textwidth]{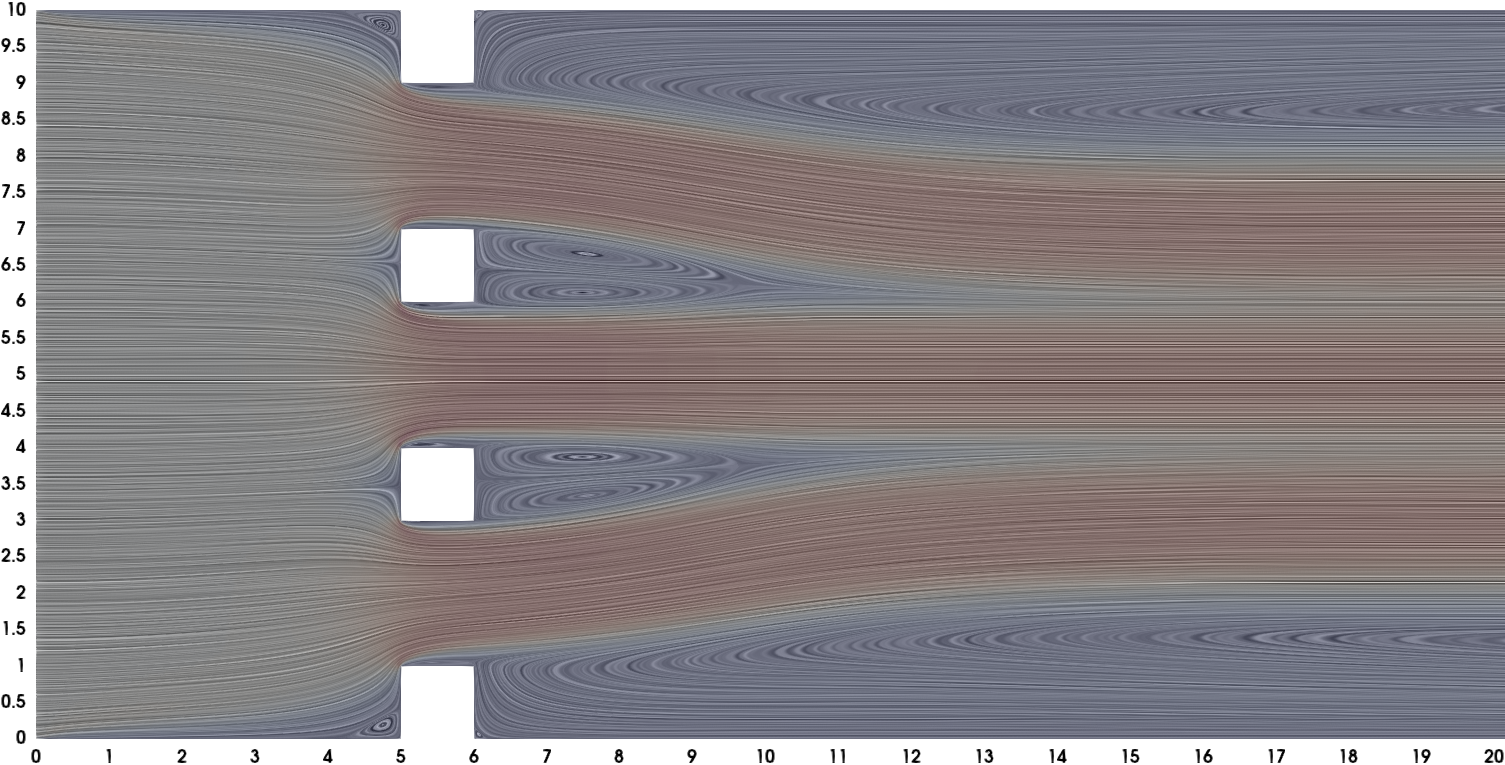}
}
\quad
\subfigure{
\includegraphics[width=0.45\textwidth]{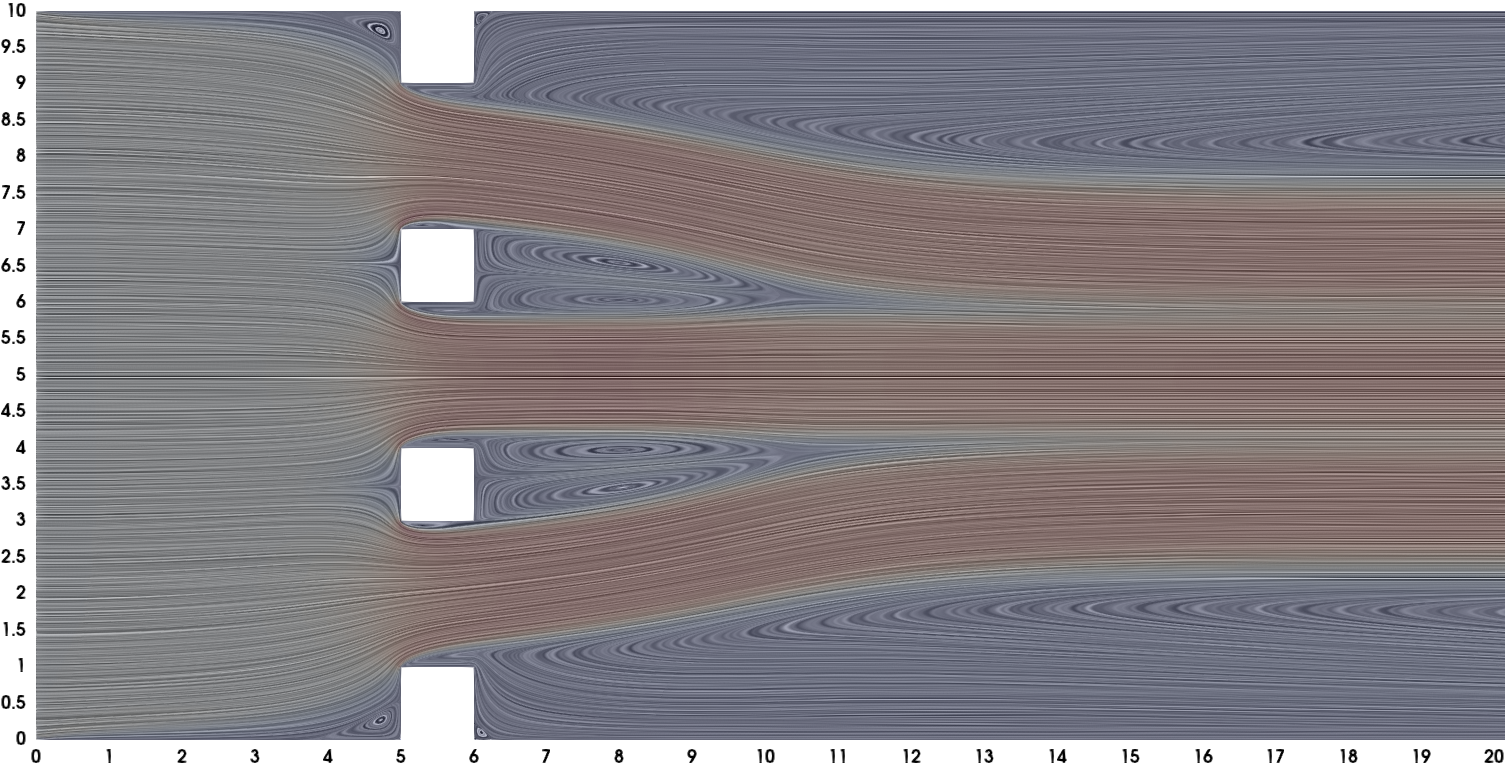}
}
\centering
\caption{\emph{Example 7.} A magnified view of streamline on colored velocity magnitude distribution near the steps and the objects for $\mathrm{Re}=150$ (left) and $\mathrm{Re}=300$ (right)}\label{fig:ex7:sllocal:Re:150:300}
\end{figure}

We consider computations still by $\bm{\mathrm{P}}_3^{\mathrm{bubble}}$-$\mathrm{P}_2^{\mathrm{dc}}$-$\bm{\mathrm{BDM}}_3$ performed on a locally refined triangular mesh presented in Figure \ref{fig:ex7:meshplot}. Here we take $h_l=1/40$ and $h_g=1/8$ to employ a finer mesh. Note that the Reynolds number is still given by $\mathrm{Re}=\frac{1}{\nu}$ here. By adopting the updated Newton iteration \eqref{linear:Newton:cylinder}, the problem is first solved for $\mathrm{Re}=10$, then $\mathrm{Re}=25$, and then in steps of $25$ until $\mathrm{Re}=300$, with the solution for the previous value of $\mathrm{Re}$ used as initial guess for the next; the Stokes equations are solved to provide the initial guess used at $\mathrm{Re}=10$.

The contours of velocity magnitude, and the contours of kinematic pressure, and the streamlines at various Reynolds numbers $\mathrm{Re}=50$, $100$, $150$ and $300$ are presented in Figures \ref{fig:ex7:vel:Re:50:300}, \ref{fig:ex7:pres:Re:50:300} and \ref{fig:ex7:sl:Re:50:300} respectively. Indeed, the added step and two small square solid objects change the distributions dramatically. Specifically, the complete symmetries occur over the whole channel; in the region near the outlet, the velocity magnitude shows the two-level differentiation while the kinematic pressure keeps constant. For the streamlines, the most remarkable feature is that two large eddies behind the steps, whose scale is larger at the same Reynolds number by comparison with that of the last example, compress the ``jet flow'' region in the middle of the channel. Another interesting phenomenon is that the distribution change of the streamlines slows down as the Reynolds number increases. To verify this statement over all the channel, a detailed plot of streamline contours near the steps and the objects for $\mathrm{Re}=150$ and $300$ is presented in Figure \ref{fig:ex7:sllocal:Re:150:300} with a magnified view.

\section{Conclusion}
\label{sec:conc}

Based on the classical mixed method with high-order conforming finite elements, we have developed an efficient pressure-robust method to the rotation form of the stationary incompressible Navier--Stokes equations over shape-regular triangular meshes via constructing a novel skew-symmetric discrete trilinear form containing the velocity reconstruction operator. The proposed method achieves the pressure-independent velocity errors and preserves the same convergence orders of both velocity and pressure as the classical method does. Three numerical examples with exactly analytical solutions have presented to verify the theoretical results and to show the remarkable performance in the aspect of errors compared to the classical method. Furthermore, four practical experiments are implemented to adequately demonstrate the efficiency and robustness of the proposed high-order finite element method.

\section*{References}
\bibliography{mybibfile}

\end{document}